\documentclass[10pt]{book}
\usepackage{remreset}

\makeatletter
  \@removefromreset{section}{chapter}
\makeatother

\usepackage[top=4.2cm, bottom=4.2cm, left=4.0cm, right=4.0cm]{geometry}
\usepackage{amssymb}
\usepackage{latexsym}
\usepackage{parskip}
\usepackage{amsthm}
\usepackage[all]{xy}
\usepackage[T1]{fontenc}
\usepackage{amssymb}
\usepackage{latexsym}
\usepackage{amsthm}
\usepackage[all]{xy}
\usepackage[T1]{fontenc}
\usepackage{chngcntr}

\newtheorem{thm}{Theorem}[section]
\newtheorem{prop}{Proposition}[section]
\newtheorem{lemma}{Lemma} \newtheorem{cor}[thm]{Corollary}

\newtheorem{dfn}{Definition}[section]

\theoremstyle{definition}
 
\newtheorem{exam}{Example}
\newtheorem{remark}{Remark}[section]

\newenvironment{abstract}%
    {\cleardoublepage\thispagestyle{empty}\null\vfill\begin{center}%
     \itshape {\bfseries \upshape \large Abstract}\end{center}}%
    {\vfill\null}

\newenvironment{acknowledgements}%
    {\cleardoublepage\thispagestyle{empty}\null\vfill\begin{center}%
     \itshape {\LARGE Acknowledgements}\end{center}}%
    {\vfill\null}
 \title{{\bf Quenched Asymptotics for the Discrete Fourier Transforms of a Stationary Process}} \author{\vspace{1cm} \large David Barrera\\{} \\ {Department of Mathematical Sciences} \\{  University of Cincinnati}\\ {\small PO Box 210025,
Cincinnati, Oh 45221-0025, USA.}\\{\small  barrerjd@mail.uc.edu}
 }
\date{}

\begin{document}  \maketitle


 
\thispagestyle{empty}
    \null\vspace{\stretch {1}}
        \begin{flushright}
               {\Large \itshape Bernarda y Gonzalo}
        \end{flushright}
\vspace{\stretch{2}}\null

\begin{abstract}
In this dissertation, we show that the Central Limit Theorem and the Invariance Principle for Discrete Fourier Transforms discovered by Peligrad and Wu can be extended to the quenched setting. We show that the random normalization introduced to extend these results is necessary and we discuss its meaning. We also show the validity of the quenched Invariance Principle for fixed frequencies under some conditions of weak dependence. In particular, we show that this result holds in the martingale case. 

The discussion needed for the proofs allows us to show some general facts apparently not noticed before in the theory of convergence in distribution. In particular, we show that in the case of separable metric spaces the set of test functions in the Portmanteau theorem can be reduced to a countable one, which implies that the notion of quenched convergence, given in terms of convergence a.s. of conditional expectations, specializes in the right way in the regular case when the state space is metrizable and second-countable.

We also collect and organize several disperse facts from the existing theory in a consistent manner towards the statistical spectral analysis of the Discrete Fourier Transforms, providing a comprehensive introduction to topics in this theory that apparently have not been systematically addressed in a self-contained way by previous references.
\end{abstract}

\begin{acknowledgements}
\vspace{1cm}
        Thanks to M. Peligrad for introducing me to the problems studied here and guiding my steps through the techniques that made possible to write this monograph. She has been the ideal Ph.D. advisor to me: confident, approachable, generous, and capable of mixing the demands for ``readings and deadlines'' with the permission for ``wondering freedom'' in  amounts that make possible the emergence of creative ideas. 
        
        I would like to extend my gratitude to W. Bryc and Y. Wang, the members of the evaluation committee for my Ph.D. dissertation, for their criticism and suggestions to improve the work here presented and for their academic advise. 
        
        My conversations with C.Dragan were particularly useful and stimulating, and D.Voln\'{y} personally encouraged me to work on the details of the proof of Theorem \ref{nonqueconthe}, pointing out the fact that his proof (with M.Woodroofe) on the non-rotated case was applicable to the discrete Fourier transforms.  
        
        This work was partially supported by the Research Grant no.1512936 from the Division of Mathematical Sciences of the National Science Foundation.
        \end{acknowledgements}

\tableofcontents
\newpage

\section*{Notation}
\pagestyle{plain}
\label{gennot}
\begin{enumerate}
\item {\it The natural numbers.} We will denote by $\mathbb{N}$ the set of natural numbers {starting at zero}. $\mathbb{N}:=\{0,1,2,\dots\}$. We will also use the notation $\mathbb{N}^{*}:=\mathbb{N}\setminus\{0\}$.
\item {\it The space $([0,2\pi),\mathcal{B},\lambda)$.} Throughout this text, $([0,2\pi),\mathcal{B},\lambda)$ will denote, unless otherwise specified, the interval $[0,2\pi)$ seen as probability space with the Borel sigma-algebra $\mathcal{B}$ and the normalized Lebesgue measure $\lambda$. This is, for every $B\in \mathcal{B}$
\begin{equation}
\label{deflamequ}
\lambda(B)=\frac{1}{2\pi}L(B)
\end{equation}
where $L$ is the Lebesgue measure, specified by $L[a,b)=b-a$ for every  real numbers $a<b$.

\item {\it Limits.} Unless otherwise specified, an expression of the form ``$\lim_{n}$'' must be read as ``$\lim_{n\to\infty}$'', and similarly for ``$\limsup_{n}$'' and ``$\liminf_{n}$''.

\item {\it Convergence of series.} Given a sequence $(a_{k})_{k\in\mathbb{Z}}$ of elements in a normed vector space $V$, we say that $\sum_{k\in\mathbb{Z}}a_{k}$ is convergent if $\sum_{k\in\mathbb{N}}a_{k}$ and $\sum_{k\in\mathbb{N}^{*}}a_{-k}$ are convergent (the partial sums have a limit), and we define $\sum_{k\in\mathbb{Z}}a_{k}:=\sum_{k\in\mathbb{N}}a_{k}+\sum_{k\in\mathbb{N}^{*}}a_{-k}.$

\item {\it Measurability.} Given measurable spaces $(\Omega_{1},\mathcal{F})$, $(\Omega_{2},\mathcal{G})$, a function $f:\Omega_{1}\to \Omega_{2}$ is {\it $\mathcal{F}/\mathcal{G}$ measurable} if for every $B\in \mathcal{G}$, $f^{-1}(B)\in\mathcal{F}$. If $(\Omega_{2},\mathcal{G})$ is the space of complex numbers with the Borel sigma algebra $\mathcal{C}$, we will use the term {\it $\mathcal{F}-$measurable function} when referring to an $\mathcal{F}/\mathcal{C}$ measurable function. For a specified $\mathcal{F}$, clear along the discussion, we will speak of a  {\it measurable function} when referring to an $\mathcal{F}-$measurable function.

\item {\it Preimages of sets.} Given  measurable spaces $(\Omega,\mathcal{F})$ and $(\Omega',\mathcal{F}')$, $A\in\mathcal{F}'$ and an $\mathcal{F}/\mathcal{F}'$ measurable function $X$, we will denote by $[X\in A]$ the $\mathcal{F}-$set
$$[X\in A]:=\{\omega\in\Omega:X(\omega)\in A\}.$$

 \item \label{meaequ}{\it Equivalence classes of functions.} If $(\Omega,\mathcal{F},\mu)$ is a measure space and $(\Omega',\mathcal{F}')$ is a measurable space, we say that two $\mathcal{F}/\mathcal{F}'-$measurable functions $X,Y$ are {\it $\mu-$equivalent} if there exists $A\in\mathcal{F}$ with $\mu(\Omega\setminus A)=0$ such that $X(\omega)=Y(\omega)$ for every $\omega\in A$. 
 If $\mu$ is fixed and $X$ is $\mu-$equivalent to $Y$, we call $X$ a {\it version} of $Y$.

\item {\it Two abbreviations}. Here, ``a.s.'' abbreviates ``almost surely'', and  ``a.e''  abbreviates ``almost every'' ({\it not} ``almost everywhere'').

\item {\it $L^{p}$ spaces.} Given a measure space $(\Omega,\mathcal{F},\mu)$, a sigma algebra $\mathcal{F}_{0}\subset \mathcal{F}$, and $0<p<\infty$, $L^{p}_{\mu}(\mathcal{F}_{0})$ denotes the (normed) space of $\mu-$equivalence classes of $p-$integrable functions $X:\Omega\to \mathbb{C}$ that are $\mathcal{F}_{0}-$measurable (with the norm given in the next item). Thus, $X\in L^{p}_{\mu}(\mathcal{F}_{0})$ if and only if (some version of) $X$ is $\mathcal{F}_{0}-$measurable and
$$\int_{\Omega}|X(\omega)|^{p}\, d\mu(\omega)<\infty.$$
If $\mathcal{F}$ is fixed, we will use the notation $L^{p}_{\mu}$ for $L^{p}_{\mu}(\mathcal{F})$. $L^{\infty}_{\mu}(\mathcal{F})$ denotes the (normed) space of $\mu-$equivalence classes of essentially bounded functions: $X\in L^{\infty}_{\mu}(\mathcal{F})$ if there exists $c>0$ such that
$$\mu([|X|>c])=0.$$

\item {\it $L^{p}$ norms.} Given $p>0$ and $X\in L^{p}_{\mu}$,  ``$||X||_{\mu,p}$'' will denote the $L^p-$norm of $X$. This is
\begin{equation}
\label{lpnorequ}
||X||_{\mu,p}:=\left(\int_{\Omega}|X(\omega)|^{p}\,d\mu(\omega)\right)^{{1}/{p}}
\end{equation}
when $p<\infty$, and 
\begin{equation}
\label{linfnorequ}
||X||_{\mu,\infty}:=\inf\{c>0:\mu[|X|>c]=0\}
\end{equation}
when $p=\infty$. 

\item {\it The spaces $l^{p}(\mathbb{Z})$ and $l^{p}(\mathbb{N})$.} If $\Omega=\mathbb{Z}$ or $\Omega=\mathbb{N}$ and $\mu$ is the counting measure ($\mu(\{z\})=1$ for every $z\in \Omega$), we will denote by $l^{p}(\mathbb{Z})$  (resp. $l^{p}(\mathbb{N})$) the space $L^{p}_{\mu}$. Thus  $(a_{k})_{k}$ ($k\in\mathbb{Z}$ or $\mathbb{N}$) belongs to $l^{p}(\mathbb{Z})$ (resp. $l^{p}(\mathbb{N})$) if and only if
\begin{equation}
\label{defainl2equ}
||(a_{k})_{k}||_{\mu,p}^{p}=\sum_{k}|a_{k}|^{p}<\infty.
\end{equation}

\item{\it Random variables and stochastic processes.} A {\it random variable} is a $\mathbb{P}-$equivalence class of measurable functions $X:\Omega\to\mathbb{C}$ defined on some probability space $(\Omega,\mathcal{F},\mathbb{P})$ (note that here random variables are complex valued functions). A {\it stochastic process} is a sequence $(X_{k})_{k}$ of random variables, where $k$ runs over $\mathbb{Z}$ or $\mathbb{N}$.

\item {\it Convergence in distribution.} The convergence in distribution of random elements in a metric space (or of probability measures, or of distribution functions) will be denoted here by ``$\Rightarrow$''. If necessary, we will use the notation ``$\Rightarrow_{n}$'' to indicate that the convergence holds as $n\to\infty$.

\item {\it Characteristic functions.} Given a measurable space $(\Omega,\mathcal{F})$ and $A\in\mathcal{F}$, we will use the notation $I_{A}$ for the {\it characteristic function} of $A$. This is $I_{A}:\Omega\to \{0,1\}$ is given by $I_{A}(\omega)=0$ if $\omega\notin A$ and $I_{A}(\omega)=1$ if $\omega\in A$.

\item {\it Expectation.} If $\mathbb{P}$ is a probability measure other than $\lambda$, we will use the traditional notation ``$E$'' to denote integration with respect to $\mathbb{P}$. Thus for instance $||X||_{_{\mathbb{P},1}}=E[|X|]$ if $X\in L^{1}_{\mathbb{P}}$. 
If we need to specify $\mathbb{P}$, we will use the notation ``$E_{^{_{\mathbb{P}}}}$'' , or some other convenient variation of it, to indicate integration with respect to $\mathbb{P}$.

\item {\it Inner Product in $L^{2}$.} We will also make use of the Hilbert space structure of $L^{2}_{\mu}$, whose inner product 
$\left\langle X,Y\right\rangle_{\mu}:L^{2}_{\mu}\times L^{2}_{\mu}\to [0,\infty)$ is defined by
\begin{equation}
\label{innproL2}
\left\langle X,Y\right\rangle_{\mu}=\int_{\Omega}X(\omega)\overline{Y}(\omega)\,d\mu(\omega),
\end{equation}
where $\overline{Y}(\omega)$ is the conjugate of $Y(\omega)$, and we will say that $X,Y\in L^{2}_{\mu}$ are {\it orthogonal} if $\left\langle X,Y\right\rangle_{\mu}=0$.

\item {\it The one-dimensional torus.} Finally, $\mathbb{T}\subset \mathbb{C}$ denotes the unit circle with the sub-space topology and the Lie-group structure given by multiplication of complex numbers.  

\end{enumerate}

\cleardoublepage
\pagestyle{headings}
\pagenumbering{arabic}

\chapter*{Introduction}
\pagenumbering{arabic}
\addcontentsline{toc}{chapter}{\protect\numberline{} Introduction} \markboth{INTRODUCTION}{INTRODUCTION}

\label{int}

The celebrated {\it Birkhoff's Ergodic Theorem} 
states that if $T:\Omega\to\Omega$ is a measure-preserving transformation on a probability space $(\Omega,\mathcal{F},\mathbb{P})$, $X\in L^{1}_{\mathbb{P}}(\mathcal{F})$, and $X_{k}:=X\circ T^{k}$, then the ergodic averages
$$A_{n}=\frac{1}{n}\sum_{k=0}^{n-1}X_{k}$$
converge $\mathbb{P}-$a.s., as $n\to\infty$, to a function $\hat{X}$ with $\hat{X}\circ T=\hat{X}$ ($\mathbb{P}-$a.s.).

A well-known and easy argument\footnote{See Theorem \ref{ergthedisfoutra} and its proof.} allows one to see that Birkhoff's Ergodic Theorem ``generalizes itself'' in the following way: let
$$S_{n}(\theta):=\sum_{k=0}^{n-1}X_{k}e^{ik\theta}$$
be the {\it $n-$th discrete Fourier Transform} of the process $(X_{k})_{k}$. Then the {\it Fourier averages} 

  
$$A_{n}(\theta):=\frac{1}{n}S_{n}(\theta)$$
converge $\mathbb{P}-$a.s., as $n\to\infty$. This is (a partial statement of) the  {\it pointwise Ergodic Theorem for Discrete Fourier Transforms}, a  version of Birkhoff's Ergodic Theorem whose further analysis has taken  mainly the directions opened by the  following questions: 

{\bf Question 1:} {\itshape Given a version of $X$. Can we choose the (probability one) set of convergence for $A_{n}(\theta)$ {\upshape independent} of $\theta$?}

The answer to this particular question appeared 1941  when Wiener and Wintner found a positive answer\footnote{Though, according to Assani (\cite{ass},p.24), Wiener and Wintner's original proof was flawed, and the first known correct proof of the Wiener-Wintner Theorem is actually due to Furstenberg (\cite{ass}, p.36), who published it in his 1960's monograph \cite{fur60}.} known today as the {\it Wiener-Wintner Theorem}, a result that opened a line of research that would lead to results such at Bourgain's {\it Return Times Theorem} (\cite{bou})  and to the currently very active investigation of convergence theorems for multiple recurrence in Dynamical Systems. These investigations lie at the heart of the connections between probabilistic or ergodic-theoretical techniques and problems in number theory, like Furstenberg's equivalence (and proof) of Szemer\'{e}di's theorem (\cite{fur}). 

{\bf Question 2:} {\it What can be said about the asymptotics of the {\upshape periodogram}}
\begin{equation}
\label{prerequ}
I_{n}(\theta)=\frac{|S_{n}(\theta)|^{2}}{n}\,\,\,\,\,\,\,\,\,\mbox{ \it?}
\end{equation}
The importance of this question came mainly from the research in the direction of the {\it Periodogram Analysis} (or, more widely, {\it Spectral Analysis}) of Time Series, a technique started by Schuster in 1898  (\cite{sch}) that would become the standard tool for the identification of statistically significant frequencies in time series of observations and has been widely applied in the Physical and Social Sciences. Several papers appeared through the 20th and 21st century addressing this and related questions in different cases important for the applications\footnote{For a review of some of them see the introduction to \cite{pelwu} and \cite{Wu} and the references therein.}, some of them departing from the elementary fact that the periodogram is the square of the modulus of $\sqrt{n}A_{n}(\theta)$, and therefore that the investigation of the Periodogram's asymptotics can be seen as a particular instance of the question about the speed of convergence of $A_{n}(\theta)$.


\subsection*{The Central Limit Theorem for Discrete Fourier Transforms}


The 2010 paper \cite{pelwu} by Peligrad and Wu, devoted to the asymptotics of $\sqrt{n}A_{n}(\theta)$, is a remarkable step in this direction of the research on Spectral Analysis. It is shown there that, if we take into account a certain $T-$filtration\footnote{See definitions \ref{defadafil}  and \ref{minadafildef} in this monograph.} associated to the process $(X_{k})_{k}$, the assumptions necessary to prove the Central Limit Theorem for $A_{n}(\theta)$ -this is, that $\sqrt{n}A_{n}(\theta)$ is asymptotically normal- can be basically reduced to the minimal ones plus a certain regularity condition (see Definition \ref{regcondef})\footnote{As the reader will wee, such condition is actually unnecessary in the quenched setting,  because of the ``random centering'' needed for the corresponding results.}. The results on that paper contain many of the precedent ones as special cases. They are also stated in the setting typical for the investigation of quenched limit theorems and gave rise to the main questions addressed in this monograph.  

Without going now into details, it is important to notice that, in contrast to the $\mathbb{P}-$a.s. convergence of $A_{n}(\theta)$, to pass from asymptotic results for $\sqrt{n}A_{n}$ to asymptotic results for (the complex-valued) process $\sqrt{n}A_{n}(\theta)$, a further analysis of the joint distributions of its real and imaginary parts is needed. Upon addressing this problem, one realizes that the  ``frequencies'' (values of $\theta$) associated to the ``square root'' of the {\it point spectrum} of the Koopman (composition) operator induced by the map $T$ (definitions \ref{defkooope} and \ref{poispe}) have the remarkable property of being the ``generic'' set of exceptional frequencies in which the asymptotics of $\sqrt{n}A_{n}(\theta)$ can fail to be (2-dimensional) normal with independent entries. For this and other reasons, the point spectrum of the Koopman operator will play an essential role in the results to be presented here, and the exposition starts with the basic definitions and properties related to it.

\subsection*{The Invariance Principle for Discrete Fourier Transforms}

One can roughly summarize Peligrad and Wu's Central Limit Theorem by saying that, under regularity, the distribution of $\sqrt{n}A_{n}(\theta)$ is (indeed) asymptotically normal with independent real and imaginary parts for {\it $\lambda-$a.e} {\it fixed} frequency $\theta$. Peligrad and Wu's paper addresses also the problem of the invariance principle but, in contrast to the case corresponding to the Central Limit Theorem,  the authors show the weaker statement that the asymptotic distribution of $W_{n}(\theta,t):=\sqrt{n}A_{\left\lfloor nt\right\rfloor}(\theta)$ corresponds to that of a random function of the form $(\theta,\omega)\mapsto f(\theta)(B_{1}(\omega)+iB_{2}(\omega))$ (with {\it random} parameters $\omega$ {\it and} $\theta$) where $B_{1}$ and $B_{2}$ are independent Brownian motions\footnote{Actually $f(\theta)=\sigma(\theta)/\sqrt{2}$, where $\theta\mapsto \sigma^{2}(\theta)$ is the {\it spectral density} of $(X_{k})_{k}$ with respect to the normalized Lebesgue measure. See Section \ref{autcovfunandspeden} for details.}. The underlying probability law is therefore $(\lambda\times\mathbb{P})\,W_{n}^{-1}$: the parameter $\theta$ is considered only ``in average'' in this case.  

The bottom line of the problem when trying to prove the Invariance Principle for fixed frequencies with these methods lies in the lack of a maximal inequality {general enough} as to pass from the martingale approximations for $\sqrt{n}A_{n}(\theta)$ to martigale approximations for $t\mapsto \sqrt{n}A_{\left\lfloor nt\right\rfloor}(\theta)$. On the other side, integrating over $\theta$ allows us to apply Hunt and Young's inequality (Theorem \ref{hunyou} in this monograph), which actually has a role in the fixed frequency approximations, to bypass this problem. To the date, it is not known whether the Central Limit Theorem of Peligrad and Wu for fixed frequencies can be extended to a corresponding invariance principle without additional assumptions.

\subsection*{The Problem of Quenched Convergence}

The results listed before are stated in the context of {\it stationary} sequences. 
There is a certain form of non-stationarity that is very important in the applications and has grown as a topic of intensive research during the last twenty years. In the context of i.i.d. sequences it can be introduced in the following way\footnote{The formal definition is Definition \ref{defquecon} in page \pageref{defquecon}.}: let $(\zeta_{k})_{k\in\mathbb{Z}}$ be an i.i.d. sequence, let $f=f(...,z_{-1},z_{0})$ be a measurable real-valued function defined on the space of complex-valued sequences indexed by the non-positive integers with the product sigma algebra, and consider the (stationary) stochastic process $(X_{k})_{k}$ given by
$$X_{k}:=f(\dots,\zeta_{k-1},\zeta_{k})$$
for all $k\in\mathbb{N}$.

Assume, for the sake of the discussion, that $EX_{0}=0$ and $EX_{0}^{2}=1$, and  that we have proved the Central Limit Theorem for the stationary process $(X_{k})_{k}$, so that
$$\frac{1}{\sqrt{n}}\sum_{k=0}^{n-1}X_{k}\Rightarrow_{n} N(0,1).$$

{\bf Question} (A question on Quenched Convergence.) {\it If we {\upshape fix} a  point $a=(\dots, a_{-1},a_{0})$ in the domain of $f$ and consider the (nonstationary) process $(X_{a,k})_{k}$ given by
$$X_{a,k}=f(\dots,a_{-1},a_{0},\zeta_{1},\dots,\zeta_{k}),$$
does the ``same'' Central Limit Theorem (still) hold for $(X_{a,k})_{k}$?}

This is, {\it can we assert that
$$\frac{1}{\sqrt{n}}\sum_{k=0}^{n-1}X_{a,k}\Rightarrow_{n} N(0,1)\,\,\,\,\mbox{\it ?}$$}

The idea behind the notion of quenched convergence is whether we can give asymptotics for a stochastic process ``started at a point'', or with ``initial conditions''. In this particular example  ``quenched convergence'' means an affirmative answer to the question above for almost every $(\dots, a_{-1},a_{0})$ with respect to the law of $(\dots,\zeta_{-1},\zeta_{0})$.

Note that the law of the modified process $(X_{a,k})_{k}$ is typically singular with respect to the law of the stationary process $(X_{k})_{k}$ (for instance $\mathbb{P}[X_{a,0}=f(a)]=1$ but $\mathbb{P}[X_{0}=f(a)]$ is typically equal to zero), and therefore we cannot give affirmative answers to the question above based on arguments of dominating measures.\footnote{See \cite{bilconpromea}, Theorem 14.2 for an example of this technique.} 

The formal notion of {quenched} convergence, which captures  the  question above and other versions of it, is actually {strictly stronger} than the notion of convergence in distribution\footnote{See Remark \ref{remdefquecon} in page \pageref{remdefquecon}.}: every process converging in the quenched sense converges in ``the annealed'' sense, but the reciprocal is not true, even in the specific setting of the question above. This will be stated ``abstractly'' in Section \ref{necrancen} and proved (in the setting of functions of i.i.d. sequences, and for the corresponding normalized Fourier averages $\sqrt{n}A_{n}(\theta)$) in Chapter \ref{chapronecrancen}. Note also that this not is obvious: consider for instance the ($m+1-$dependent) case in which $f$ has ``finite memory'',  $f=f(z_{-m},\cdots,z_{0})$ for some $m>0$.\footnote{Or, in a more strict language, $f((z_{-k})_{k\in\mathbb{N}})=f((z_{-k}')_{k\in\mathbb{N}})$ 
for any two sequences whose terms coincide for $0\leq k \leq m$.}

\subsection*{A Growing Trend}

The problem of quenched convergence was not intensively studied during the 20th century, though it has been long recognized as an important requirement in the theory of statistical inference for Markov processes\footnote{See for instance the note preceding (1.8) in \cite{bilstainf}. See also Example \ref{staseqasfunmarcha} in page \pageref{staseqasfunmarcha} in this monograph for technical details on the relationship between quenched convergence -as presented here- and convergence with respect to the transition measures induced by the Kernel of a stationary Markov chain.}. 

Results on quenched convergence can be traced to at least 1968 with Billingsley's quenched Invariance Principle for $\phi-$mixing processes (\cite{bilconpromea68}, Theorem 20.4). Other results in this direction appeared sporadically\footnote{According to the remarks in \cite{cunmerpel}, the paper \cite{gorlif} deserves special mention in this respect, since it started the investigation of these results in the sense of Markov operators. See also \cite{barpelpel} for a slightly more detailed account of the results in this direction before 2001.},  but an inflexion point came with the paper \cite{derlin}  published in 2001 by Derrienic and Lin, which was inspired by a question raised by Kipnis and Varadhan in 1987 (\cite{kipvar}, Remark 1.7) and gave rise to a considerable amount of new research\footnote{In order of appearance, some examples are \cite{wuwoo}, \cite{cun},  \cite{merpelpel}, \cite{cunpel}, \cite{cunvol}, \cite{cunmer}, \cite{volwoononque}, \cite{volwoo},  \cite{barpel}, and \cite{barpelpel}.} on the validity of the Central Limit Theorem for functions of  Markov chains when the chain starts at a point. In a more informal way it can be asserted that, nowadays, the word ``quenched'' is becoming a common sound in the conferences and meetings of specialists in Probability.

\subsection*{The Content of this Monograph}

This work  presents the first series of results on quenched limit theorems for the discrete Fourier transforms of a stationary process. 

While the original purpose was to limit the exposition to the minimal amount of material necessary to fully understand the results presented in the series of papers \cite{barqueinv}, \cite{baranexa}, and \cite{barpel},  and therefore to refer the reader to the existing literature for the background theory, I found particularly difficult to navigate between the many references needed to carry on the proofs of the results in question while maintaining at the same time a clear perspective of the mathematical ground in which these arguments rest. For that reason, a chapter on ``Background Theory'', Chapter \ref{chabacthe}, was inserted. While it is my desire that it can serve as a quick introduction for anyone interested in reading the series of papers started by \cite{Wu} in the direction of the speed of convergence for the discrete Fourier transforms, the results presented in this chapter are not original
, and my motivation to present them was to pave the way to a clear exposition in  further sections. I have tried to keep the references to the literature containing the original proofs  even in the cases in which, for pedagogical reasons, I decided to rewrite them. This was not always possible though, and I must advance my apologies to any reader who finds a proof by a different author without a reference, expecting that (s)he believes in the unintentional nature of my omission.

Chapter \ref{chacondis}  covers issues related to convergence in distribution. The following reasons lead me to insert these topics as part of this monograph: first, although perhaps obvious for the expert that knows the real-valued case, the notion of convergence in distribution for {\it complex-valued} {cadlag} functions is not easy to find in the mainstream literature, thus I decided that it was wise to devote a few pages explaining how this notion can be understood via an obvious extension of the Skorohod metric to the complex-valued case, and how the techniques used for real-valued functions {\it indeed} apply to the complex-valued ones. For the same reason, I also considered important to explain why some well-known convergence of types theorems can be carried over to the complex-valued case, and to give them as statements of convergence of random variables instead of distribution functions. The ``transport theorem'' in Section \ref{trathe}, borrowed from an external source, was inserted in order to make the monograph more self-contained.

The ``refinement of the Portmanteau theorem'' (Section \ref{refporthe}) deserves, on the other side, special mention. It came out after many hours confronting a certain question that has some resemblance to the one giving rise to the Wiener-Wintner theorem: when facing the problem of passing from the ``fixed frequency'' to the ``averaged frequency'' limit theorems, which in the annealed case can be trivially solved by integrating with respect to the parameter $\theta$, one has to deal with the fact that, in the quenched case, the (probability one) set of decomposing measures with respect to which the results hold for a fixed $\theta$ may change with $\theta$, and therefore one has to be more careful when performing integrations over the (uncountable) set $[0,2\pi)$ of parameters $\theta$. While this can be done via arguments involving interpretations of Fubini's theorem, I found more illuminating and clear to use the language of conditional expectations in this case, but in order to succeed with this way one is finally lead to ask whether the set of test functions in the Portmanteau theorem can be reduced to a countable one. The answer is ``yes'' in the separable case (and it is what this ``refinement'' deals with), and the consequences for the theory of quenched convergence  pay off, in my opinion, the short digression.

Chapter \ref{quecon} presents the definitions and elementary properties related to the notion of quenched convergence: adapted $T-$filtrations, regular conditional expectations, and the interactions between the product measures and the regular conditional expectation with respect to the  product of two sigma algebras. These notions constitute the elementary ``grammar'' necessary for the results on quenched convergence presented here and for their proofs, and are usually taken for granted along the papers in the literature. For this reason, this is also a chapter aimed to introduce the beginner to these techniques.

This chapter presents also several examples related to the existence of regular conditional expectations. For the most part, they belong to the standard literature, but since some of the constructions along the references are not given in terms of invertible Dynamical Systems, I considered appropriate to spend some energy explaining how the corresponding results are indeed possible if we  restrict ourselves to the invertible case. In particular, we obtain a representation of a stationary stochastic process as a sequence of functions of a stationary Markov Chain preserving the invertibility of the underlying shift operator.  The construction can  be easily adapted to show that {any} stochastic process admits a representation as a function of a (possibly nonstationary) Markov Chain.

The main contributions of this monograph are contained in Chapter \ref{resandcom}. In summary, it is shown there that the limit theorems by Peligrad and Wu admit quenched versions under some ``intuitively obvious'' modification (the ``random centering''). It is shown that this  modification is necessary, and some quenched invariance principles for fixed frequencies are also provided. 

The rest of the monograph is devoted to prove the results in Chapter \ref{resandcom}. Since it is not possible to make comments about this without going into technical details we will just mention two things: 

{\it First}, the reader is invited to note that, in a certain sense, all the quenched results given here for fixed frequencies are just interpretations of corresponding results for (non-rotated) partial sums (including, after a ``creative'' step, the proof of Theorem \ref{quecltfoutra}), and therefore we can consider the investigation of asymptotics for the discrete Fourier transforms $S_{n}(\theta)$ for $\theta$ fixed  ({almost}) as a particular case of the investigation of the asymptotics for partial sums of (complex-valued) random variables. 

{\it Second}. For the reader familiar with these techniques, note that the standard application of maximal inequalities to pass from martingale approximations leading to the Central Limit Theorem to corresponding approximations leading to the Invariance Principle encounters an additional obstacle here: in the nonstationary setting, maximal inequalities are scarce. This is the bottom line behind the necessity of weak-$L^{p}$ spaces along our proofs, and it provides a further reason to call for the (already growing) investigation of maximal inequalities for nonstationary processes.


\part{Background and Results}

\chapter{Background Theory}
\label{chabacthe}

In this chapter we will survey the background theory necessary to justify our forthcoming discussions and to settle a solid ground for them. Specifically, we will be concerned with presenting the objects that motivate the questions leading to the main results in this monograph, leaving aside for later chapters the discussions relative to the methods of our proofs. 

Most of the results presented in this chapter are part of the literature and the reader is referred to the corresponding reference for their proofs. Nonetheless, we will go through detailed discussions whenever the clarity of the arguments would be affected otherwise.

This chapter is organized as follows: in Section \ref{kooopeanditspoispe} we discuss the notion of the {\it Koopman operator} (Definition \ref{defkooope}) associated to a measure preserving transformation, emphasizing the discussion on the structure of its point spectrum. These notions will show up later along the proofs of our main results, particularly in the steps involving asymptotic finite-dimensional distributions. 

Then, in Section \ref{ranelel2}, we will present some results necessary to clarify the construction of the approximating martingales whose asymptotics will be transferred to the processes under consideration. This will require a short review of results from classical Harmonic Analysis and a visit to the problem of measurability for functions defined by limits.

Section \ref{ratcongenergthe} presents a result (Theorem \ref{ergtheposdunsch}) that seems to be implicit in the literature but whose pieces are somehow disperse. This theorem gives rise to a result (Theorem \ref{ergthedisfoutra}) that generalizes the pointwise and $L^{p}$ ergodic theorems to discrete Fourier transforms in a very natural way, justifying the investigation of its rate of convergence via the Central Limit Theorem. To reach this result we have to introduce a technical notion,  the {\it ``extension to the product space''} of a random variable and a measure-preserving transformation (see Definition \ref{extprospa} and the discussion following it), that will be important for some of the steps in the forthcoming proofs of our main results. We also introduce the basics of {\it weak $L^{p}-$spaces}, which will be needed later along the proofs from Chapter \ref{resandcom}.

In Section \ref{adafil} we settle the ground for the forthcoming discussions about ``quenched convergence''. In particular, we will establish (Proposition \ref{relconexpkooope}) the interaction between the Koopman operator and the conditional expectations  with respect to the corresponding elements in the filtration of an adapted process, a fact that will be crucial for our proofs. We will also introduce important notions such as that of a (strictly) ``stationary process'' (Definition \ref{staproL2}), ``left'' and ``right'' sigma-algebras (Definition \ref{taisigalgdef}), and the ``model'' example of linear processes (Example \ref{exalinpro}). We conclude with two ergodic theorems (Theorem \ref{ergtheadafil} and Corollary \ref{dedmerpel}) that will be of utter importance when discussing the quenched asymptotic distributions associated to the normalized Fourier averages, and in particular to understand the role of the point spectrum in the statements of the results to be presented in Chapter \ref{resandcom}.

Finally, in Section \ref{autcovfunandspeden}, we will introduce the notion of the {\it autocovariance function} (Definition \ref{defautcovfun}) and the {\it spectral density} (Definition \ref{defspeden}) of a stationary square-integrable process, whose estimation justify much of the research in the directions explored along this monograph. We will also introduce the notion of {\it regular processes} (Definition \ref{regcondef}), which will be essential for some of our proofs. Our discussion will lead us to the (annealed) limit theorems of Peligrad and Wu (theorems \ref{cltpelwu} and \ref{pelwufclt}), whose extension to the quenched setting is one of the main purposes of this work.

\section{The Koopman Operator and its Point Spectrum}
\label{kooopeanditspoispe}
In this section we present the notion of the {\it Koopman operator} associated to a measure preserving transformation on a probability space, and we introduce the analytic facts about it that will be of use along the proofs of the results present in this monograph.

\subsection{Definitions and General Properties}

Let us start by recalling the notion of a {\it measure preserving transformation}.

\medskip

\begin{dfn}[Measure Preserving Transformation]
\label{defmeapretra}
Given a measure space $(\Omega,\mathcal{F},\mu)$, a {\itshape\bfseries measure preserving transformation} $T:\Omega\to \Omega$ is an $\mathcal{F}/\mathcal{F}$-measurable map such that for every $A\in\mathcal{F}$
\begin{equation}
\label{premeaequ}
\mu(T^{-1}A)=\mu(A).
\end{equation}
\end{dfn} 

We will restrict our attention in this monograph to measure-preserving transformations on probability spaces, but some of the notions presented below can be extended to more general measure spaces.

In particular, measure preserving transformations and their dynamics will be of utter importance to codify the notion of stationary processes used along this work (see Definition \ref{staproL2} below). To settle the ground for the upcoming discussions let us introduce now the notion of the {\it Koopman operator} associated to a measure preserving transformation. 

\bigskip

\begin{dfn}[Koopman Operator]
\label{defkooope}
Let $(\Omega,\mathcal{F},\mathbb{P})$ be a probability space and let $T:\Omega\to\Omega$ be  a measure preserving transformation. Given $p> 0$ we define the {\itshape\bfseries Koopman operator} $T:L^{p}_{\mathbb{P}}\to L^{p}_{\mathbb{P}}$ by 
$$TY:=Y\circ T.$$
\end{dfn} 
\begin{remark}
\label{remunadefkooope}
Note that we are using the same notation for the transformation $T$ and its associated Koopman operator. This should not be a source of confusion in what follows: ``$TU$'' must be interpreted as $U\circ T$ when $U$ is a random variable, and as the image of $U$ under $T$ when $U$ is a subset of $\Omega$. Similarly, $T^{-1}U$ should be understood as $U\circ T^{-1}$ if $U$ is a random variable when $T$ is  invertible and bimeasurable, and as $T^{-1}(U)$, the inverse image of $U$ under $T$, if $U$ is a subset of $\Omega$.
\end{remark}

Since many of our forthcoming proofs depend on spectral properties of the Koopman operator associated to a measure-preserving transformation, we will start by presenting some elementary facts related to the eigenvalues of these operators. Let us start by a formal introduction of these objects.

\medskip

\begin{dfn}[Point Spectrum of $T$]
\label{poispe}
With the notation in Definition \ref{defkooope}, denote by
$$Spec_{p}(T):=\{\alpha\in \mathbb{C}: \,\mbox{\upshape there exists $q> 0$ and  $X \in L^{q}_{\mathbb{P}}\setminus\{0\}$ with\,\,\,} TX=\alpha X \}.$$
$Spec_{p}(T)$ is called the {\bfseries  point spectrum} of $T$, and any element of $Spec_{p}(T)$ is called an {\bfseries eigenvalue} of $T$.  
\end{dfn}

\medskip

\begin{remark}
\label{speinsT}
Note that if $p>0$, $T$ is an isometry in $L^{p}_{\mathbb{P}}$: $(E[|TX|^{p}])^{1/p}=(E[|X|^{p}])^{1/p}$. In particular, $Spec_{p}(T)\subset \mathbb{T}$.
\end{remark}

\medskip

The following proposition shows that  the definition of $Spec_{p}(T)$ can be recast by restricting $T$ to $L^{q}_{\mathbb{P}}$ for a fixed $q>0$. 

\medskip

\begin{prop}[Persistence of $Spec_{p}(T)$]
\label{boueigfun}
In the setting of definitions \ref{defkooope} and \ref{poispe} denote, for every $q>0$ and $\alpha\in \mathbb{T}$
\begin{equation}
\label{defeigspaequ}
V_{\alpha}^{q}:=\{X\in L^{q}_{\mathbb{P}}: TX=\alpha X\}. 
\end{equation}
and let $V_{\alpha}:=\cup_{q>0}V_{\alpha}^{q}$. Then the following statements are equivalent
\begin{enumerate}
\item $\alpha \in Spec_{p}(T)$. 
\item $V_{\alpha}\neq\{0\}$.
\item $V_{\alpha}\cap L^{\infty}_{\mathbb{P}}\neq \{0\}$.
\end{enumerate}
In particular, given $q>0$, $Spec_{p}(T)$ is the set of eigenvalues of the Koopman operator $T:L^{q}_{\mathbb{P}}\to L^{q}_{\mathbb{P}}$. 
\end{prop}

{\bf Proof:} Only $2. \Rightarrow 3.$ requires a proof. 

Indeed, note that if $0\neq Y \in V_{\alpha}$ is given, then from $|TY|=|\alpha Y|=|Y|$ $\mathbb{P}-$a.s. it follows that for all $M\geq 0$
$$TI_{[|Y|\leq M]}=I_{[T|Y|\leq M]}=I_{[|TY|\leq M]}=I_{|Y|\leq M}$$
$\mathbb{P}-$a.s. and therefore, choosing $M$ such that $0\neq I_{[|Y|\leq M]}$,
$$T(YI_{[|Y|\leq M]})=(TY)(TI_{[|Y|\leq M]})=\alpha Y I_{[|Y|\leq M]}.$$
Thus $X:=YI_{[|Y|\leq M]}\in V_{\alpha}$. Since clearly $X\in L^{\infty}_{\mathbb{P}}\setminus \{0\}$ this gives the desired conclusion.\qed

Throughout this monograph, we will be mainly concerned with {\it ergodic transformations}. Ergodic transformations enjoy some special properties and, in some sense, they are the building blocks of any measure preserving transformation (see for instance Theorem 6 in \cite{qua}). The definition is the following.

\medskip

\begin{dfn}[Ergodic Transformation]
\label{ergtra}
Let $(\Omega,\mathcal{F},\mathbb{P})$ be a probability space. A transformation $T:\Omega\to \Omega$ is called {\itshape\bfseries ergodic} if it is measure preserving and every $T-$invariant set in $\mathcal{F}$ is ``trivial''. This is, if for every $A\in\mathcal{F}$: $T^{-1}A=A$ implies that $\mathbb{P}(A)\in\{0,1\}$. Equivalently, if for all $A\in\mathcal{F}$: $TI_{A}=I_{A}$ implies that $EI_{A}\in\{0,1\}$.
\end{dfn}

Ergodicity is a spectral property: among its well-known characterizations the following one will be of special interest for us.

\medskip

\begin{lemma}[Ergodicity as a Spectral Property]
\label{equergonesimeigval}
A measure-preserving transformation $T$ on a probability space is ergodic if and only if $1$ is a simple eigenvalue of $T$: if $X$ satisfies $TX=X$ then $X$ is (a.e. equal to a fixed) constant.
\end{lemma}
{\bf Proof:} See for instance \cite{einandwar}, Proposition 2.14. \qed

\medskip

Our attention along this work will be mainly focused on Koopman operators associated to {ergodic} transformations on a probability space. To give a first consequence of the ergodic hypothesis note the following: according to the first line in the proof of Proposition \ref{boueigfun}, if $Y$ is an eigenvector of $T$ then $|Y|$ is $T-$invariant, and therefore constant if $T$ is ergodic (Lemma \ref{equergonesimeigval}). This gives the following result. 

\bigskip

\begin{prop}[Circularity of Eigenfunctions]
\label{eigfunareuni}
Assume that $T$ is ergodic and $\alpha\in \mathbb{T}$ is given: if $Y$ satisfies $TY=\alpha Y$ then $|Y|$ is constant.
\end{prop}





Even more is true: the following proposition implies that, when $T$ is ergodic, the eigenfunctions of $T$ are unique up to multiplication by a scalar. Note also that in this case $Spec_{p}(T)$ is more than just a subset of $T$.

\medskip

\begin{prop}[Group Structure of $Spec_{p}(T)$]
\label{spesubgro}
With the notation in Proposition \ref{boueigfun}, and assuming $T$ is ergodic, $Spec_{p}(T)$ is a subgroup of $\mathbb{T}$, and every element in $Spec_{p}(T)$ is a simple eigenvalue of $T$.
\end{prop}

{\bf Proof:} The proposition consists of two statements, which we proceed to prove now.

{\it $Spec_{p}(T)$ is a group.} Since clearly $1\in Spec_{p}(T)$ (consider any constant function $X$), it suffices to see that if $\alpha_{1}\in Spec_{p}(T)$ and $\alpha_{2} \in Spec_{p}(T)$, then $\alpha_{1}\overline{\alpha_{2}}\in Spec_{p}(T)$.

Let us prove it: given $\alpha_{1},\alpha_{2}\in Spec_{p}(T)$ and nonzero functions $ X_{1}\in V_{\alpha_{1}}$ and $X_{2}\in V_{\alpha_{2}}$, note that, since $|X_{1}|$ and $|X_{2}|$ are constant non-zero functions (Proposition \ref{eigfunareuni}),  $X_{1}X_{2}$ is (also) nonzero, and that
$$T(X_{1}\overline{X_{2}})=TX_{1}\overline{TX_{2}}=\alpha_{1}\overline{\alpha_{2}}X_{1}X_{2}.$$
In particular $\alpha_{1}\overline{\alpha_{2}}\in Spec_{p}(T)$, as claimed.

{\it The eigenvalues are simple.} If $X,Y$ are (non-zero) eigenfunctions associated to $\alpha\in Spec_{p}(T)\subset \mathbb{T}$, the argument just given shows that $X\overline{Y}$ is an eigenfunction of $T$ associated to $1$. Since $T$ is ergodic, there exists a constant $c\in\mathbb{C}$ such that $X\overline{Y}=c$, $\mathbb{P}-$a.s. It follows (multiply by $Y$) that $|Y|^{2}X=cY$ and therefore, since $0<|Y|$ is  constant, there exists a constant $\beta$($=c/||Y||_{_{\mathbb{P},\infty}}^{2}$) such that $X=\beta Y$.\qed

\subsection{Separability and Cardinality of the Point Spectrum}
\label{remspeiscou}

Let us recall now the following well known definition:

\medskip

\begin{dfn}[Countably generated sigma-algebras]
\label{deffcougen}
A sigma algebra $\mathcal{F}$ is {\it countably generated} if there exists a countable family of sets $\mathbb{A}=\{A_{k}\}_{k\in\mathbb{Z}}\subset \mathcal{F}$ such that $\sigma(\mathbb{A})=\mathcal{F}$. 
\end{dfn}

\medskip

This is the case if, for instance, $\mathcal{F}$ is the Borel sigma algebra of a separable metric space $(S,d)$, or if $\mathcal{F}$ is the sigma algebra generated by a countable family of random elements in a separable space. If $(\Omega,\mathcal{F},\mathbb{P})$ is a probability space and $\mathcal{F}$ is countably generated, $L^{p}_{\mathbb{P}}$ is separable for every $p\in [1,\infty)$ (see \cite{bilpromea}, Theorem 19.2). 

Now, it is a standard exercise to prove that the separability of $L^{2}_{\mathbb{P}}$ (or more generally, of any Hilbert space) is equivalent to the existence of a countable orthonormal basis of $L^{2}_{\mathbb{P}}$: a set $\{Y_{k}\}_{k\in\mathbb{Z}}\subset L^2_{\mathbb{P}}$ of mutually orthogonal elements whose linear span is dense in $L^{2}_{\mathbb{P}}$. In particular, if $\mathcal{F}$ is countably generated, $L^{2}_{\mathbb{P}}$ admits at most countably many mutually orthogonal one-dimensional subspaces: for any family $\{Y_{j}\}_{j\in J}\subset L^{2}_{\mathbb{P}}$ of mutually orthogonal elements with $E[|Y_{j}|^2]=1$, the balls centered at $Y_{j}$ with radius $1$ are mutually disjoint, which restricts the cardinality of $J$ to a countable one if $L^{2}_{\mathbb{P}}$ is separable.

Recall the notation introduced in Proposition \ref{boueigfun} and note that, since $T$ is measure preserving, the spaces $V_{\alpha}$ are mutually orthogonal: given $\alpha_{1}\in V_{\alpha_{1}}$ and $\alpha_{2}\in V_{\alpha_{2}}$, $$E[Y_{1}\overline{Y}_{2}]=E[T[Y_{1}\overline{Y}_{2}]]=\alpha_{1}\overline{\alpha}_{2}E[Y_{1}\overline{Y}_{2}]$$
which implies that either $\alpha_{1}=\alpha_{2}$ or $E[Y_{1}\overline{Y}_{2}]=0$.

From these observations the following follows at once.

\medskip

\begin{prop}[Cardinality of $Spec_{p}(T)$]
\label{prospecou}
Let $(\Omega,\mathcal{F},\mathbb{P})$ be a probability space. If $\mathcal{F}$ is countably generated then for every measure-preserving transformation $T:\Omega\to\Omega$, $Spec_{p}(T)$ is countable. In particular 
\begin{equation}
\label{eigangmeaone}
\lambda(\{\theta\in[0,2\pi):e^{i\theta}\in Spec_{p}(T)\})=0.
\end{equation}
\end{prop}

\section{Random Elements in $L^{2}$}
\label{ranelel2}

In this section, we will introduce the results from Harmonic Analysis that will be used along the monograph. In particular, we will  show how to use Carleson theorem (Theorem \ref{carthe}) to show that a random function $\omega\mapsto f_{\omega}$ in $L^{2}_{\lambda}$ (see Definition \ref{raneledef}) defined on $(\Omega,\mathcal{F},\mathbb{P})$ induces a random function  $([0,2\pi),\mathcal{B},\lambda)\to L^{2}_{\mathbb{P}}$ if the sigma algebra $\mathcal{F}$ is countably generated (Theorem \ref{thedua}), a construction that will be important when justifying that the approximating martingales present along the proofs of the results in Chapter \ref{resandcom} are well defined. 

On doing so, we will stop to discuss the measurability of a function defined by limits in a complete and separable metric space (Section \ref{fundeflimsec}). We will also introduce the notion of discrete Fourier Transforms (Definition \ref{defdisfoutra}) of a stochastic process, a generalization of the notion of partial sums that is at the heart of the results presented in this work. 

\subsection{Functions Defined by Limits}
\label{fundeflimsec}
In this section, we will discuss the issue of the measurability for a map given by pointwise convergence of random functions in a metric space, and we will define the notion of ``limit function'' for an a.s convergent sequence of random elements in a complete and separable metric space in an unambiguous way. The results and definitions introduced here will be used, several times in an implicit way, along the discussions involving functions defined by (a.e.) convergent sequences.

We begin our discussion introducing the following technical notion.
\medskip

\begin{dfn}[Distance to a set, $\epsilon-$Neighborhood]
\label{defdistoaset}
If $(S,d)$ is a metric space with metric $d$, then for any given $x\in S$ and $A\subset S$ we define the {\itshape \bfseries distance from $x$ to $A$} by
\begin{equation}
\label{disxtoa}
d(x,A):=\inf_{a\in A}d(x,a),
\end{equation}
and we define the
 {\itshape\bfseries $\epsilon-$neighborhood of $A$}, $A^{\epsilon}$, as the (open) set 
 \begin{equation}
 \label{defenei}
 A^{\epsilon}:=\{x\in S:d(x,A)<\epsilon\}. 
 \end{equation}
\end{dfn}

Assume that $(S,d)$ is a (nonempty) metric space. In addition assume that $(S,d)$ is  complete and separable\footnote{The assumption of completeness is made to guarantee  that the set of points where a given sequence of functions converges is measurable. The assumption of separability is made to guarantee that $\mathcal{S}\otimes\mathcal{S}$ is the Borel sigma-algebra of $S\times S$ (see Appendix M10 in \cite{bilconpromea}), so that the distance function $d:S\times S\to [0,\infty)$, which is  continuous with respect to the product topology, is  $\mathcal{S}\otimes\mathcal{S}-$measurable, and  for any two given $\mathcal{F}/\mathcal{S}$ measurable functions $f,g$, the function $\omega\mapsto d(f(\omega),g(\omega))$ is $\mathcal{F}-$measurable.}, let $\mathcal{S}$ be the Borel sigma-algebra of $S$, and fix $s\in S$. If $(f_{n})_{n\in\mathbb{N}}$ is a sequence of $\mathcal{F}/\mathcal{S}$ measurable functions,  define $C_{_{(f_{n})_{n}}}$ as the (measurable) set  where $(f_{n})_{n\in\mathbb{N}}$ is a Cauchy sequence. Thus
\begin{equation}
\label{defcauseqfn}
C_{_{(f_{n})_{n}}}:=\bigcap_{m\in \mathbb{N}^{*}}\bigcup_{n\in\mathbb{N}}\bigcap_{k\geq n}[d(f_{n},f_{k})<{1}/{m}],
\end{equation}
 and define the ``limit function'' $\lim_{n} f_{n}$ by
\begin{equation}
\label{deflimfn}
\omega\mapsto \left\{ 
  \begin{array}{c l}
     \lim_{n}f_{n}(\omega)&, \quad \textrm{if $\omega\in C_{_{(f_{n})_{n}}}$ } \\
    s &, \quad \textrm{if $\omega\notin C_{_{(f_{n})_{n}}}$}\\
    
  \end{array} \right.
  \end{equation}

Now remember the well known definition of the $\liminf$ of a family of (measurable) sets $\{A_{n}\}_{n\in\mathbb{N}}\subset {\mathcal{F}}$:
$$\liminf_{n}A_{n}:=\bigcup_{n\in\mathbb{N}}\bigcap_{k\geq n}A_{k},$$
and note that that for every closed set $F\subset S$
\begin{equation}
\label{preimacloset}
(\lim_{n}f_{n})^{-1}(F)=\left\{ 
  \begin{array}{l  l}
    C_{_{(f_{n})_{n}}} \cap \bigcap_{m\in\mathbb{N}^{*}}\liminf_{n}[f_{n}\in F^{1/m}] & , \quad \textrm{if $s\notin F$ } \\
    (\Omega\setminus C_{_{(f_{n})_{n}}})\cup(C_{_{(f_{n})_{n}}} \cap \bigcap_{m\in\mathbb{N}^{*}} \liminf_{n}[f_{n}\in F^{1/m}] ) &, \quad \textrm{if $s\in F.$}\\
    
  \end{array} \right .
  \end{equation} 
The measurablility of these sets, together with the $\pi-\lambda$ theorem (applied to the set of elements $A\in\mathcal{S}$ such that $(\lim_{n}f_{n})^{-1}(A)\in\mathcal{F}$) give at once the following result.

\medskip

\begin{prop}[Measurability of Limit Functions]
\label{measlimfun}
Let $(S,d)$ be a complete and separable (nonempty) metric space with Borel sigma algebra $\mathcal{S}$. Given any sequence $(f_{n})_{n\in\mathbb{N}}$ of $\mathcal{F}/\mathcal{S}$ measurable functions defined on some measurable space $(\Omega,\mathcal{F})$, the function $\lim_{n}f_{n}$ defined by (\ref{deflimfn}) is $\mathcal{F}/\mathcal{S}-$measurable. 
\end{prop}

Finally note that if $\mathbb{P}$ is a probability measure on $(\Omega,\mathcal{F})$ and $\mathbb{P}(\Omega\setminus C_{_{(f_{n})_{n}}})=0$, then the $\mathbb{P}-$equivalence class of $\lim_{n}f_{n}$ is independent of the choice of $s$. 

Let us formalize this in the following definition

\medskip

\begin{dfn}[Functions Defined by Limits]
\label{deflimfun}
In the context of Proposition \ref{measlimfun}, assume that $\mathbb{P}$ is a probability measure on $(\Omega,\mathcal{F})$, and that $\mathbb{P}(\Omega\setminus C_{_{(f_{n})_{n}}})=0$. We define {\itshape \bfseries the} limit function (also denoted by) $\lim_{n}f_{n}$ as the $\mathbb{P}-$equivalence class of functions represented by  $\lim_{n}f_{n}$.
\end{dfn}

Let us finish this section by reminding the formal notion of a random element in a metric space.

\medskip

\begin{dfn}[Random Elements and their Law]
\label{raneledef}
If $(S,d)$ is a metric space with Borel sigma algebra $\mathcal{S}$, a {\itshape \bfseries random element of} $S$ is an $\mathcal{F}/\mathcal{S}$ measurable function $V:\Omega\to S$ from some probability space $(\Omega,\mathcal{F},\mathbb{P})$ to $(S,\mathcal{S})$. If $V$ is a random element on $S$, the {\itshape \bfseries law of $V$} is the probability measure $\mathbb{P}V^{-1}$ on $\mathcal{S}$ defined by
$$\mathbb{P}{V}^{-1}(A)=\mathbb{P}[V\in A]$$
for all $A\in\mathcal{S}$. 
\end{dfn}

During the rest of this section, we will focus our attention on random elements in $L^{2}_{\lambda}$. This is, $\mathcal{F}/\mathcal{S}$ measurable functions $V:\Omega\to L^{2}_{\lambda}$ from a probability space $(\Omega,\mathcal{F},\mathbb{P})$ to the space $(S,d)$ of square-integrable functions in $([0,2\pi),\mathcal{B},\lambda)$ with the $L^{2}_{\lambda}$ norm. 

\subsection{The Fourier Transform of an Integrable Function}
\label{thefoutraofanint}

Let us begin this section by reminding the notion of the {\it Fourier transform} of a function $f\in L^{1}_{\lambda}$, which is the building block for the representation by Fourier series of elements in $L^{2}_{\lambda}$ (or, under an appropriate notion of convergence, of elements in $L^{1}_{\lambda}$).

\medskip

\begin{dfn}[Fourier Transform]
\label{foutra}
Given $f\in L^{1}_{\lambda}$, $\hat{f}:\mathbb{R}\to \mathbb{C}$ will denote  the {\bfseries Fourier transform of $f$}, which is defined by
\begin{equation}
\label{deffoutra}
\hat{f}(x)=\int_{0}^{2\pi}f(\theta)e^{-ix\theta}d\lambda(\theta).
\end{equation}
\end{dfn}

Our first goal is to describe in which sense the Fourier transform of a function allows us to represent it in a convenient way. The first step towards this goal is to define the {\it Fourier partial sums} of a function in $L^{1}_{\lambda}$.

\medskip

\begin{dfn}[Fourier Partial Sums]
\label{parfousum}
For a given $n\in\mathbb{N}^{*}$, the {\bfseries $n-$th Fourier partial sum of a function $f\in L^{1}_{\lambda}$ {\it at a  frequency} $\theta \in [0,2\pi)$} is defined by
\begin{equation}
\label{disfoutraf}
S_{f,n}(\theta):=\sum_{k=1-n}^{n-1}\hat{f}(k)e^{ik\theta}.
\end{equation}
\end{dfn}
In 1966 Lennart Carleson (\cite{Car}) proved the following celebrated result, establishing that the Fourier series representation of a function in $L^{2}_{\lambda}$ is convergent almost surely.\footnote{This result is also true for functions in $L^{p}_{\lambda}$ with $p>1$. See for instance \cite{lac}.}

\medskip

\begin{thm}[Carleson]
\label{carthe}
Let $f\in L^{2}_{\lambda}$ and let $S_{f,n}(\theta)$ be defined by (\ref{disfoutraf}), then
$$f(\theta)=\lim_{n} S_{f,n}(\theta)$$
in the sense of Definition \ref{deflimfun}. This is: there exists a set $I_{f}$ with $\lambda(I_{f})=1$ such that for every $\theta\in I_{f}$, $\lim_{n}S_{f,n}(\theta)=f(\theta)$. 
 \end{thm}

Now, given $f\in L^{2}_{\lambda}$, Parseval's theorem (\cite{gra}, Proposition 3.1.16, (3)) establishes that 
$$\int_{0}^{2\pi}|f(\theta)|^{2}d\lambda(\theta)=\sum_{n\in \mathbb{Z}} |\hat{f}(n)|^{\,2}$$
and, reciprocally, Plancherel's theorem (\cite{gra}, Proposition 3.1.16, (2) and (4)) establishes that for any $(c_{k})_{k\in \mathbb{Z}}\in l^{2}(\mathbb{Z})$, the map ($\lambda-$a.e) given by
\begin{equation}
\label{fouseql2n}
\theta\mapsto \sum_{k\in \mathbb{Z}}c_{k}e^{ik\theta}
\end{equation}
defines a (unique) element $f\in L^{2}_{\lambda}$, with Fourier coefficients $\hat{f}(k)=c_{k}$. These observations can be summarized in the following proposition.

\medskip

\begin{prop}[Representation of $L^{2}_{\lambda}$]
\label{bijcorL2l2}
The correspondence $L^2_{\lambda}\to l^{2}(\mathbb{Z})$ given by $f\mapsto (\hat{f}(k))_{k\in\mathbb{Z}}$ is (well defined and) bijective.
 \end{prop}

\medskip

\begin{remark}
\label{corL2l2isaniso} 
Note that, by Parseval's Theorem, the correspondence given in Proposition \ref{bijcorL2l2} is an isometry of metric spaces.\footnote{This is actually the content of \cite{gra}, Proposition 3.1.16, (4).}
\end{remark}

\subsection{A Duality Theorem}
In virtue of Proposition \ref{bijcorL2l2} and Remark \ref{corL2l2isaniso} we can think of functions in $L^{2}_{\lambda}$ just as elements in $l^{2}(\mathbb{Z})$. In particular, a random function in $L^{2}_{\lambda}$ can be thought of as a measurable map $$\mathbf{Y}:(\Omega,\mathcal{F},\mathbb{P})\to l^{2}(\mathbb{Z})$$
where $(\Omega,\mathcal{F},\mathbb{P})$ is a probability space. Concretely, given a random function $\omega\mapsto f_{\omega}$ of $L^{2}_{\lambda}$, take 
$$\mathbf{Y}(\omega):=(\hat{f}_{\omega}(k))_{k\in\mathbb{Z}},$$
where $\hat{f}_{\omega}$ is the Fourier Transform of $f_{\omega}$ (Definition \ref{foutra}). 

Reciprocally, 
since $l^{2}(\mathbb{Z})$ is separable (see the discussion in Section \ref{remspeiscou}), a random element in $L^{2}_{\lambda}$ is specified by any sequence $(Y_{k})_{k\in\mathbb{Z}}$ of random variables defined on $(\Omega,\mathcal{F},\mathbb{P})$, provided that
\begin{equation}
\label{conranfunl2}
\sum_{k\in\mathbb{Z}}|Y_{k}|^{2}<\infty \,\,\,\,\,\,\,\,\,\mbox{$\mathbb{P}$-a.s.}
\end{equation}
where $\mathbf{Y}(\omega):=0$ if $\omega$ does not belong to the set where (\ref{conranfunl2}) converges.

More can be said in this case: since for $\mathbb{P}-$almost every $\omega$, the series 
$$\sum_{k\in\mathbb{Z}}Y_{k}(\omega)e^{ik\theta}$$
is $\lambda-$a.e convergent, the $\mathcal{B}\otimes\mathcal{F}-$set 
$$A:=\{(\theta,\omega)\in [0,2\pi)\times \Omega: \sum_{k\in\mathbb{Z}}Y_{k}(\omega)e^{ik\theta} \,\,\,\mbox{is convergent}\}$$
satisfies $\lambda\otimes\mathbb{P}(A)=1$, and an application of Fubini's theorem shows that there exists a set $I_{\mathbf{Y}}\subset[0,2\pi)$ with $\lambda(I_{\mathbf{Y}})=1$ satisfying following property: for every $\theta\in I_{\mathbf{Y}}$ there exists $\Omega_{\theta}$ with $\mathbb{P}(\Omega_{\theta})=1$ such that the series
\begin{equation}
\label{ranfunfixfre}
\sum_{k\in \mathbb{Z}}Y_{k}e^{ik\theta}
\end{equation}
converges for all $\omega\in \Omega_{\theta}$. 

If we assume in addition that, for a given $\theta\in I_{\mathbf{Y}}$ (or in a set $I_{\mathbf{Y}}'\subset I_{\mathbf{Y}}$ with $\lambda(I_{\mathbf{Y}}')=1$)
\begin{equation}
\label{consupL2}
E\left[\sup_{n\in \mathbb{N}} \left|\sum_{k= 1-n}^{n-1}Y_{k}(\omega)e^{ik\theta}\right|^{2}\right]<\infty
\end{equation}
then, by Lebesgue's dominated convergence theorem, the function given by
$$\omega\mapsto \sum_{k\in\mathbb{Z}}Y_{k}(\omega)e^{ik\theta}$$
belongs to $L^{2}_{\mathbb{P}}$. In particular we have the following result.

\medskip

\begin{thm}[Duality of Random Elements in $L^{2}$]
\label{thedua}
Let $(\Omega,\mathcal{F},\mathbb{P})$ be a probability space. Assume that $\mathcal{F}$ is countably generated (Definition \ref{deffcougen}), and let $\mathbf{Y}:(\Omega,\mathcal{F},\mathbb{P})\to l^{2}(\mathbb{Z})$ be a random element of $l^{2}(\mathbb{Z})$. If (\ref{consupL2}) holds for $\lambda-$a.e $\theta$,  then the map ($\lambda-$ a.e) defined  by
\begin{equation}
\label{mapduatheequ}
\theta\mapsto \sum_{k\in\mathbb{Z}}Y_{k}e^{ik\theta}
\end{equation}
(where the series is defined in the $\mathbb{P}-$a.s sense) is a random element $([0,2\pi),\mathcal{B},\lambda)\to L^{2}_{\mathbb{P}}$ of $L^{2}_{\mathbb{P}}$. 

If \,$Y:(\Omega,\mathcal{F},\mathbb{P})\to L^{2}_{\lambda}$ is a random element of $L^{2}_{\lambda}$ and we denote by $\hat{f}$ the Fourier transform of $f$ (Definition \ref{foutra}), the same statement holds by taking $\mathbf{Y}(\omega)= (Y_{k}(\omega))_{k\in\mathbb{Z}}:=(\widehat{Y(\omega)}(k))_{k\in\mathbb{Z}}$. 
\end{thm}

{\bf Proof:} First: since $\mathcal{F}$ is countably generated, $L^{2}_{\mathbb{P}}$ is (complete and) separable (see the comments following Definition \ref{deffcougen}). 

Let now $\mathcal{L}^{2}_{\mathbb{P}}$ denote the Borel sigma-algebra of $L^{2}_{\mathbb{P}}$.  Only the $\mathcal{B}/\mathcal{L}^{2}_{\mathbb{P}}$ measurability of (\ref{mapduatheequ}) is left to prove, which will follow if we can prove that the convergence of (\ref{ranfunfixfre}) in the $\mathbb{P}-$a.s sense (which is guaranteed for $\lambda-$a.e $\theta$) together with (\ref{consupL2}) implies the convergence of (\ref{mapduatheequ}) in the $L^{2}_{\mathbb{P}}-$sense for $\lambda-$a.e $\theta$. \footnote{More precisely, note that for given $N\in\mathbb{N}$, the map $f_{N}:[0,2\pi)\to L^{2}_{\mathbb{P}}$ given by
$$f_{N}(\theta):=\sum_{|k|\leq N}Y_{k}e^{ik\theta}$$
is $\mathcal{B}/\mathcal{L}^{2}_{\mathbb{P}}$ measurable (it is indeed continuous), and  that if (\ref{mapduatheequ}) makes sense as a limit in $L^{2}_{\mathbb{P}}$ for $\lambda-$a.e $\theta$, then it is indeed the same as the map $f:=\lim_{N}f_{N}$ (in the sense of Definition \ref{deflimfun}).}

To see this we can argue as follows: by the a.s convergence  of (\ref{ranfunfixfre}) 
 $$\lim_{N}|\sum_{|k|\leq N}Y_{k}e^{ik\theta}-\sum_{k\in \mathbb{Z}}Y_{k}e^{ik\theta}|=0$$
$\mathbb{P}-$a.s for $\lambda-$a.e $\theta$, and since for every $N\in \mathbb{N}$ (and every such $\theta$)
$$|\sum_{|k|\leq N}Y_{k}e^{ik\theta}-\sum_{k\in \mathbb{Z}}Y_{k}e^{ik\theta}|^2\leq 2\sup_{n\in\mathbb{N}}|\sum_{|k|< n}Y_{k}e^{ik\theta}|^{2},$$
the dominated convergence theorem, together with (\ref{consupL2}), imply that
$$ \lim_{N}E|\sum_{|k|\leq N}Y_{k}e^{ik\theta}-\sum_{k\in \mathbb{Z}}Y_{k}e^{ik\theta}|^{2}=0$$
for $\lambda-$a.e $\theta$, as desired.

The last statement follows at once from the previous one combined with Proposition \ref{bijcorL2l2}.\qed
 .
\subsection{Duality via Decay of Second Moments}

In this section we will give a sufficient condition (see (\ref{conmomy}) below) to guarantee the fulfillment of (\ref{consupL2}), and therefore the validity of the conclusion of Theorem \ref{thedua}. We introduce also the notion of the ($n-$th){\it discrete Fourier transform} of a stationary process, whose normalized asymptotic behavior is the main topic of this work. 

\subsection*{A Maximal Inequality, the Discrete Fourier Transforms}
The following result is another classical tool in Harmonic Analysis (we give here a particular version sufficient for our purposes).

\medskip

\begin{thm}[A Maximal Inequality]
\label{hunyou}
There exists a constant $C$ with the following property: for all $f\in L^{2}_{\lambda}$
\begin{equation}
\label{hunyoumaxine}
\int_{0}^{2\pi}\sup_{n\in\mathbb{N}^{*}}|S_{f,n}(\theta)|^{2} d\lambda(\theta)\leq 
C\sum_{k\geq 0}|\hat{f}(k)|^{2}, 
\end{equation}
where $S_{f,n}(\theta)$ is the $n-$th Fourier partial sum of $f$ at $\theta$ (see (\ref{disfoutraf})) and $\hat{f}$ denotes the Fourier transform of $f$ (Definition \ref{foutra}). 
\end{thm}

{\bf Proof:} See \cite{hunyou}.\qed

From now on, we will refer to the inequality (\ref{hunyoumaxine}) as {\it Hunt and Young's maximal inequality}.


Now consider the following condition on a stochastic process $(Y_{k})_{k\in\mathbb{Z}}$ defined on a probability space $(\Omega,\mathcal{F},\mathbb{P})$:
\begin{equation}
\label{conmomy}
\sum_{k\in \mathbb{Z}}||Y_{k}||_{_{\mathbb{P},2}}^{2}<\infty.
\end{equation}

Note that, under this condition, (\ref{conranfunl2}) is satisfied (the function $\sum_{k\in \mathbb{Z}}|Y_{k}|^{2}$ is actually integrable by the monotone convergence theorem), and $(Y_{k})_{k\in\mathbb{Z} }$  is therefore a random element of $l^{2}(\mathbb{Z})$.  Even more, by Theorem \ref{hunyou}, there exists a constant $C$ such that
\begin{equation}
\label{hunyouran}
\int_{0}^{2\pi}\sup_{n} |\sum_{k=1-n}^{n-1}Y_{k}e^{ik\theta}|^{2}d\lambda(\theta)\leq C\sum_{k\in\mathbb{Z}}|Y_{k}|^{2} \,\,\,\,\,\,\,\,\,\,\mbox{$\mathbb{P}$-a.s.}
\end{equation}
More precisely, (\ref{hunyouran}) holds on the set of $\mathbb{P}-$measure one
$$[\sum_{k\in\mathbb{Z}}|Y_{k}|^{2}<\infty].$$
Integrating with respect to $\mathbb{P}$, and using Fubini's theorem we get that, under (\ref{consupL2}),
$$\int_{0}^{2\pi}E[\sup_{n} |\sum_{k=1-n}^{n-1}Y_{k}e^{ik\theta}|^{2}]\,d\lambda(\theta)\leq C\sum_{k\in \mathbb{Z}}E[|Y_{k}|^{2}]<\infty.$$
In particular
\begin{equation}
\label{supparsumL2}
E[\sup_{n}|\sum_{k=1-n}^{n-1}Y_{k}e^{ik\theta}|^{2}]<\infty 
\end{equation}
for $\lambda-$a.e $\theta$. This, combined with Theorem \ref{thedua} gives the following result.

\medskip

\begin{prop}[A Criterion for Duality]
\label{proconsupL2}
If a stochastic  process $(Y_{k})_{k\geq 0}$ defined on $(\Omega,\mathcal{F},\mathbb{P})$ satisfies (\ref{conmomy}) and $\mathcal{F}$ is countably generated,  the function 
\begin{equation}
\label{proconsupL2equ}
\theta\mapsto \sum_{k\in \mathbb{Z}}Y_{k}\,e^{ik\theta}
\end{equation}
defines a random element $([0,2\pi),\mathcal{B},\lambda)\to L^{2}_{\mathbb{P}}$, and there exists $I'\subset[0,2\pi)$ with $\lambda(I')=1$ such that for every $\theta\in I'$,  (\ref{supparsumL2}) is verified and (\ref{proconsupL2equ}) converges $\mathbb{P}-$a.s .
\end{prop}

\medskip

\begin{remark}
\label{nonnecfcougen}
The assumption on $\mathcal{F}$ (being countably generated) is made only to prove the $\mathcal{B}/\mathcal{L}^{2}_{\mathbb{P}}$ measurability  of the map (\ref{proconsupL2equ}) (see the proof of Theorem \ref{thedua}): the existence of the set $I'$  holds regardless of the nature of $\mathcal{F}$.\footnote{Note that for every $N\in\mathbb{N}$ the map $$(\theta,\omega)\mapsto \sum_{|k|\leq N}Y_{k}(\omega)\,e^{ik\theta}$$ is $\mathcal{B}\otimes\mathcal{F}-$measurable. So is $(\theta,\omega)\mapsto f_{n}(\theta,\omega):=\max_{0\leq k\leq n}|\sum_{j=1-k}^{k-1}Y_{k}e^{ik\theta}|^{2}$, and therefore so is  $(\theta,\omega)\mapsto \lim_{n}f_{n}(\theta,\omega)$ (in the sense of (\ref{deflimfn})). 

It is then clear that the last map is  $\mathcal{B}\otimes\mathcal{F}-$measurable and,   under (\ref{conmomy}), it coincides $\lambda\times\mathbb{P}-$a.s with  $\sup_{n\geq 0}|\sum_{j=1-n}^{n-1}Y_{k}(\omega)e^{ik\theta}|^{2}$ by (\ref{hunyouran}). The argument for the existence of $I'$ goes through just as explained.}
\end{remark}


Before continuing with our discussion, let us stop here to introduce the notion of {\it discrete Fourier Transforms} of a stochastic process.

\medskip

\begin{dfn}[Discrete Fourier Transforms]
\label{defdisfoutra}
Given a stochastic  process $(Y_{k})_{k}$ defined on a probability space $(\Omega,\mathcal{F},\mathbb{P})$ and $n\in \mathbb{N}^{*}$, we will define the {\it {\bfseries $n$-th discrete Fourier transform of $(Y_{k})_{k}$ at the frequency $\theta\in [0,2\pi)$}}, $S_{n}(\theta,\cdot):\Omega\to \mathbb{C}$, by
\begin{equation}
\label{nfoutradef}
S_{n}(\theta,\omega):=\sum_{k=0}^{n-1}Y_{k}({\omega})e^{ik\theta}.
\end{equation}
If $(Y_{k})_{k}$ is not clear from the context, we will use the notation $S_{n}((Y_{k})_{k},\theta,\omega)$ for $S_{n}(\theta,\omega)$. If  $\theta$ is fixed, we will  denote by $S_{n}(\theta)$ the random variable $S_{n}(\theta,\cdot)$. If $\theta=0$, we denote by $S_{n}$ the random variable $S_{n}(0,\cdot)$. 
\end{dfn}

\medskip

\begin{remark}[A note on the definition of $S_{n}(\theta)$]
As the reader can see, if $(Y_{k})_{k\in\mathbb{Z}}$ is a process indexed by $\mathbb{Z}$, we are not including the elements $Y_{k}$ for $k<0$ in our notion of discrete Fourier Transforms. 

A plausible alternative may be to sum over the set of indexes $\{1-n,\cdots,n-1\}$ but, while (\ref{nfoutradef}) may certainly be an example of a ``bad definition'' in the framework of a more general theory, we stick to it  here mainly due to the facts that, first, all of our forthcoming discussions will be made under the additional hypothesis that $(X_{k})_{k}$ is strictly stationary (Definition \ref{staproL2}), which allows us to generate the process $(X_{k})_{k\in \mathbb{Z}}$ knowing only the initial function $X_{0}$ and, second, our results will be concerned with the asymptotics related to ``$S_{n}(\theta)-E_{0}S_{n}(\theta)$'' (see Section \ref{adafil} for the corresponding notation), a normalization that would annihilate the summands with negative index in the ``extended'' definition of the discrete Fourier Transforms.

Finally, our theory is concerned with the asymptotics of processes ``with initial time'', an assumption implicit along the tradition of the study of central limit theorems, and even more important here given the heuristics of the notion of convergence that we will deal with (see Section \ref{heuint}). 

There are other practical reasons to keep this definition (for instance: we would encounter problems passing form the ``randomly centered'' to the ``non-centered'' case along the discussion in Section \ref{necrancen} if we adopted the extended definition of $S_{n}(\theta)$) but after all this choice is, to a certain extent, just a matter of taste, and it is a good exercise for the reader to verify which of the proofs concerning $S_{n}(\theta)$ can be carried through with the suggested, more symmetric definition of the discrete Fourier transforms.
\end{remark}

\section{Dunford-Schwartz Operators and the Ergodic Theorem}
\label{ratcongenergthe}

In this section we present the ergodic theorem for positive Dunford-Schwartz operators and its consequent ergodic theorem for discrete Fourier transforms, a result that has interest in itself and justifies the investigation of the validity of the central limit theorem for the normalized averages of the discrete Fourier transforms of a stationary process. We also make a digression towards weak $L^{p}-$spaces, whose weak norms provide a framework that will be of use along the proofs of forthcoming results.

\subsection{The Ergodic Theorem for Positive Dunford-Schwartz Operators}
\label{dunschope}
To begin with, let us recall the definition of a {\it Dunford-Schwartz operator}. 
\medskip

\begin{dfn}[Dunford-Schwartz Operators]
\label{defdunsch}
Let $(\Omega,\mathcal{F},\mu)$ be a measure space. A {\itshape\bfseries Dunford-Schwartz operator} $T:L^{1}_{\mu}\to L^{1}_{\mu}$ is a linear operator with the following property: for every $p\geq 1$ and every $X\in L^{p}_{\mu}\cap L^{1}_{\mu}$
\begin{equation}
\label{equdunsch}
||TX||_{p}\leq ||X||_{p}\, \, .
\end{equation}
\end{dfn}

\begin{remark}
\label{remondunschdef}
It is possible to see (\cite{eisfarhaanag}, Theorem 8.23) that if (\ref{equdunsch}) holds for $p=1$ and $p=\infty$ then $T$ is Dunford-Schwartz. It is also clear that when $\mu(\Omega)<\infty$ (for instance if $\mu$ is a probability measure), this definition is equivalent to the condition that $T$ is a contraction in $L^{p}_{\mu}$ for every $p\geq 1$ (use the well known continuous embedding $L^{p}_{\mu}\subset L^{1}_{\mu}$, valid when $\mu(\Omega)$ is finite). 
\end{remark}

We will also need to make use of the notion of {\it positive operator}.

\medskip
\begin{dfn}[Positive Operators]
\label{defposope}
Let $(\Omega,\mathcal{F},\mu)$ be a measure space and let $T:L^{1}_{\mu}\to L^{1}_{\mu}$ be a bounded linear operator. $T$ is called {\itshape\bfseries positive} if for any $X\in L^{1}_{\mu}$, $T|X|$ is nonnegative.
\end{dfn}
The following theorem arises from a combination of Theorems 8.24, 11.4 and 11.6 in \cite{eisfarhaanag}, together with the fact that the operator $T$ involved in the hypotheses is continuous in the respective $L^{p}$ space.

\medskip 

\begin{thm}[The Mean and Pointwise Ergodic Theorem for Positive Dunford-Schwartz Operators]
\label{ergtheposdunsch} 
Let $(\Omega,\mathcal{F},\mu)$ be a {finite measure} space ($\mu(\Omega)<\infty$) and let $T:L^{1}_{\mu}\to L^{1}_{\mu}$ be a positive Dunford-Schwartz operator (Definitions \ref{defdunsch} and \ref{defposope}). Then for every $p\geq 1$ and every $X\in L^{p}_{\mu}$ there exists $P_{T}X$ with the following properties
\begin{enumerate}
\item $P_{T}X$ is $T-$invariant: $TP_{T}X=P_{T}X$.
\item The Cesaro-averages $(X+\cdots +T^{n-1}X)/n$ converge to $P_{T}X$ $\mu-$a.s and in $L^{p}_{\mu}$:
\begin{equation}
\label{equconcesave}
\lim_{n}\frac{1}{n}\sum_{k=0}^{n-1}T^{k}X=P_{T}X\mbox{\,\,\,\,\,\,\,{\it $\mu-$a.s and in $L^{p}_{\mu}$}.}
\end{equation}
\end{enumerate} 

\end{thm}

{\bf Proof:} Denote by $A_{n}X$ ($n\geq 1$) the corresponding averages in the conclusion of Theorem \ref{ergtheposdunsch}. This is
\begin{equation}
\label{defanequ}
A_{n}X:=\frac{1}{n}\sum_{k=0}^{n-1}T^{k}X.
\end{equation}
For the existence of the limit in statement 2. see to the proofs of Theorems 8.24, 11.4 and 11.6 in \cite{eisfarhaanag}. Denote this limit by $P_{T}X$.

To see that the limit satisfies 1. note that, since $A_{n}X$ converges in $L^{p}_{\mu}$
$$TP_{T}X=\lim_{n}TA_{n}X$$
(here ``$\lim$'' denotes limit in $L^{p}_{\mu}$) and that, since $T^{n}X/n\to 0$ as $n\to\infty$ $\mu-$a.s ($T^{n}X/n=(n+1)A_{n+1}/n-A_{n}$): 
$$\lim_{n}(A_{n}X-TA_{n}X)=\lim_{n}\frac{1}{n}(X-T^{n}X)=0,\mbox{\,\,\,\,\,\,\,{\it $\mu-$a.s.}}\qed$$

\begin{remark}[$P_{T}$ as a projection, a case of orthogonality.]
\label{remptasaproj}
Given a Banach space $B$ with norm $||\cdot||_{B}$, a {\itshape\bfseries projection on $B$} is a continuous linear operator $P:B\to B$ with the property that $P^2=P$. If $P$ is a projection and $V_{P}:=PB$, we say that $P$ {\it projects $B$ onto $V_{P}$}. 

Notice that Theorem \ref{ergtheposdunsch} can be stated in the following way: {\it let $T$ be a nonnegative Dunford-Schwartz operator and, for $p\geq 1$, let $V_{T,p}\subset L^{p}_{\mathbb{P}}$ be the (closed) subspace of $T-$invariant functions ($Y\in V_{T,p}$ if and only if $TY=Y$, $\mathbb{P}-$a.s.). Then the function $P_{T}:L^{p}_{\mathbb{P}}\to V_{T,p}$ given by
$$P_{T}Y=\lim_{n}A_{n}Y$$
is well defined both in the $\mathbb{P}-$a.s. and $L^{p}_{\mathbb{P}}$-senses.}

It is easy to see that, indeed, $P_{T}L^{p}_{\mathbb{P}}=V_{T,p}$, and $P_{T}$ is clearly linear. Since $P_{T}$ is a contraction in $L^{p}_{\mathbb{P}}$ ($||\lim_{n}A_{n}Y||_{p}=\lim_{n}||A_{n}Y||_{p}\leq ||Y||_{p}$), and since $P_{T}^{2}=P_{T}$ (by the $T-$invariance of $P_{T}X$ for every $X$), $P_{T}$ is a projection.

Assume now that $T$ preserves the mean: for every $p\geq 1$ and every $Y\in L^{p}_{\mathbb{P}}$
\begin{equation}
\label{tpresmeanequ}
E[TY]=E[Y],
\end{equation}
and assume in addition that, either $T$ is multiplicative
\begin{equation}
\label{tmulequ}
T[XY]=TXTY
\end{equation}
(this is the case for instance when $T$ is a Koopman operator) or that
\begin{equation}
\label{psemulpro}
T[XTY]=TX\,TY
\end{equation}
(for instance if $T$ is a conditional expectation), whenever the expressions involved make sense. Then we can show that, actually, $P_{T}$ is {\itshape \bfseries orthogonal}. This is, that if $p\in [1,\infty)$ is given, then for every $X\in L^{p}_{\mathbb{P}}$ and $Y\in L^{p/(p-1)}_{\mathbb{P}}$ ($Y\in L^{\infty}_{\mathbb{P}}$ if $p=1$), $E[(X-P_{T}X)P_{T}Y]=0$. 

To do so we notice the following: first, since $P_{T}Z=\lim_{n}A_{n}Z$ in the $L^{1}_{\mathbb{P}}$ sense,
\begin{equation}
\label{premeapt}
E[P_{T}Z]=E[\lim_{n}A_{n}Z]=\lim_{n}E[A_{n}Z]=E[Z]
\end{equation}
for every $Z\in L^{1}_{\mathbb{P}}$. Then, since $P_{T}Y$ is $T-$invariant,
\begin{equation}
\label{linwrtptequ}
T^{n}(XP_{T}Y)=(T^{n}X)P_{T}Y 
\end{equation}
where $n$ is any natural number\footnote{This is obvious under (\ref{tmulequ}),  and to prove it under (\ref{psemulpro}) we proceed by induction: the case $n=0$ is trivial, an assuming that (\ref{linwrtptequ}) holds for a value of $n$: $$T^{n+1}(XP_{T}Y)=T(T^{n}(XP_{T}Y))=T((T^{n}X)P_{T}Y)=T((T^{n}X)T(P_{T}Y))=(T^{n+1}X)(TP_{T}Y)=$$
$$(T^{n+1}X)P_{T}Y.$$}, and therefore $P_{T}XP_{T}Y=P_{T}(XP_{T}Y)$. 

All together,
this gives that for every $p\geq 1$, $X\in L^{p}_{\mathbb{P}}$ and $Y\in L^{p/(p-1)}_{\mathbb{P}}$:
$$E[(X-P_{T}X)(\overline{P_{T}Y})]=E[XP_{T}\overline{Y}]-E[(P_{T}X)(P_{T}\overline{Y})]=E[XP_{T}\overline{Y}]-E[P_{T}(XP_{T}\overline{Y})]=$$
$$E[XP_{T}\overline{Y}]-E[XP_{T}\overline{Y}]=0$$
as claimed.
 
\end{remark}
\subsection{The Ergodic Theorem for Discrete Fourier Transforms}
\label{ergthefoutra}
The results provided in section \ref{dunschope} allow us to generalize the mean and pointwise ergodic theorems to the case of rotated partial sums (discrete Fourier transforms). In particular, this justifies an interpretation of the main results of this monograph as theorems about the ``speed of convergence'' for the normalized averages of the discrete Fourier transforms. 

Let $(\Omega,\mathcal{F},\mathbb{P})$ be a probability space, let $T:\Omega\to \Omega$ be an invertible, bimeasurable measure-preserving transformation, and let $\theta\in [0,2\pi)$ be given. 

Consider the transformation $\tilde{T}_{\theta}:[0,2\pi)\times\Omega\to [0,2\pi)\times\Omega$ specified by
\begin{equation}
\label{exttprospathe}
\tilde{T}_{\theta}(u,\omega)=((u+\theta)mod(2\pi),T\omega).
\end{equation}

Notice that $\tilde{T}_{\theta}$ is simply the product map between the rotation $u\mapsto (u+\theta)mod (2\pi)$ and $T$. This transformation is clearly measure preserving and invertible. 
\bigskip

\begin{dfn}[Extension to the Product Space]
\label{extprospa} 
Let $p\geq 1$ and  $Y\in L^{p}_{\mathbb{P}}$ be given,  we will denote by $\tilde{Y}$ the extension of 
$Y$ to $[0,2\pi)\times\Omega$ given by the following formula: 
$$\tilde{Y}(u,\omega)=e^{iu}Y(\omega).$$
\end{dfn}

It is clear that $\tilde{Y}\in L^{p}_{\lambda\times\mathbb{P}}$ and that the $L^{p}_{\lambda\times\mathbb{P}}$ norm of this extension is the same as the $L^{p}_{\mathbb{P}}$ norm of $Y$. Note also that if $\tilde{T}_{\theta}$ is given by (\ref{exttprospathe}), then for all $k\in\mathbb{Z}$, 
\begin{equation}
\label{relkooext}
\tilde{T}_{\theta}^{k}\tilde{Y}=\widetilde{e^{ik\theta}T^{k}Y}
\end{equation}
Note that $\tilde{T}_{\theta}$, seen as an operator in $L^{p}_{\lambda\times\mathbb{P}}$ for $p\in[1,+\infty)$ (namely $Z\in L^{p}_{\lambda\times\mathbb{P}} \mapsto Z\circ \tilde{T}_{\theta}$: the Koopman operator associated to $\tilde{T}_{\theta}$), is a positive contraction for every $p$. It follows from Theorem \ref{ergtheposdunsch} that there exists a $\tilde{T}_{\theta}-$invariant function $\tilde{P}_{\theta}\tilde{Y}$ such that
$$\frac{1}{n}\sum_{k=0}^{n-1}\tilde{T}_{\theta}^{k}\tilde{Y}(u,\omega)=\frac{e^{iu}}{n}S_{n}((T^{k}Y)_{k\in\mathbb{N}},\theta)(\omega)\to_{n} \tilde{P}_{\theta}\tilde{Y}(u,\omega)$$
$\lambda\times\mathbb{P}-$a.s and in $L^{p}_{\lambda\times\mathbb{P}}$. Fixing $u_{0}$ such that $\sum_{k=0}^{n-1}\tilde{T}_{\theta}^{k}\tilde{Y}(u_{0},\cdot)/n$ converges $\mathbb{P}-$a.s we see that, if $$P_{\theta}Y(\omega):=e^{-iu_{0}}\tilde{P}_{\theta}\tilde{Y}(u_{0},\omega)$$
then, necessarily
$$\tilde{P}_{\theta}\tilde{Y}(u,\omega)=e^{iu}P_{\theta}Y(\omega), \,\,\,\,\,\,\,\mbox{$\lambda\times\mathbb{P}-$a.s.}$$
In particular, $\widetilde{P_{\theta}Y}=\tilde{P}_{\theta}\tilde{Y}$.

Finally note that, since $\widetilde{{P}_{\theta}{Y}}$ is $\tilde{T}_{\theta}-$invariant
$$TP_{\theta}Y=e^{-iu}e^{-i\theta}\tilde{T}_{\theta}\widetilde{{P}_{\theta}({Y})}=e^{-iu}e^{-i\theta}\tilde{T}_{\theta}(\tilde{P}_{\theta}\tilde{Y})=e^{-i\theta}e^{-iu}\tilde{P}_{\theta}\tilde{Y}=e^{-i\theta}P_{\theta}Y.$$
This proves the following result.

\bigskip

\begin{thm}[The Ergodic Theorem for Discrete Fourier Transforms]
\label{ergthedisfoutra}
Let $(\Omega,\mathcal{F},\mathbb{P})$ be a probability space, $T:\Omega\to \Omega$ a measure-preserving transformation, $p\geq 1$, $Y\in L^p_{\mathbb{P}}$, and denote (also) by $T$ the Koopman operator associated to $T$ (Definition \ref{defkooope}) and by $S_{n}(Y,\theta)$ the $n-$th discrete Fourier transform of the process $(T^{k}Y)_{k\in\mathbb{N}}$ (Definition \ref{defdisfoutra}). Then for every $\theta\in [0,2\pi)$ there exists a  function $P_{\theta}Y\in L^{p}_{\mathbb{P}}$ with the following properties
\begin{enumerate}
\item $TP_{\theta}Y=e^{-i\theta}P_{\theta}Y$, $\mathbb{P}-$a.s.
\item $S_{n}(Y,\theta)/n\to_{n} P_{\theta}Y$, $\mathbb{P}-$a.s. and in $L^{p}_{\mathbb{P}}$.
\end{enumerate}
If $T$ is ergodic $|P_{\theta}Y|$ is constant and $P_{\theta}Y$ is unique up to a scalar multiple. This is: if $Y_{1}, Y_{2}\in L^{p}_{\mathbb{P}}$, then  there exists $c\in \mathbb{C}$ such that $P_{\theta}Y_{1}=cP_{\theta}Y_{2}$.
\end{thm}

{\bf Proof:} By the preceding discussion, only the last statement requires a proof, but this follows at once from Proposition \ref{eigfunareuni} and Colollary \ref{spesubgro}.\qed

\medskip

\begin{remark}[$P_{\theta}$ as an orthogonal projection]
\label{ptkooopecas}
It is easy to see that, in general, 
$$E_{_{\mathbb{P}}}[X\overline{Y}]=E_{_{\lambda\times\mathbb{P}}}[\tilde{X}\overline{\tilde{Y}}]. $$ 
By Remark \ref{remptasaproj}, $\tilde{P}_{\theta}$ is an orthogonal projection onto the subspace $V_{\tilde{T}_{\theta}}\subset L_{\lambda\times\mathbb{P}}^{p}$ of functions that are invariant with respect to $\tilde{T}_{\theta}$. It follows that for every $X\in L^{p}_{\mathbb{P}}$, $Y\in L^{p/(p-1)}_{\mathbb{P}}$
$$E_{_{\mathbb{P}}}[(X-P_{\theta}X)\overline{P_{\theta}Y}]=
E_{_{\lambda\times\mathbb{P}}}[(\tilde{X}-\tilde{P}_{\theta}\tilde{X})\overline{\tilde{P}_{\theta}\tilde{Y}}]=0.$$
This is: for fixed $p\geq 1$, $P_{\theta}$ is the orthogonal projection onto $V_{\theta,p}$, where
$$V_{\theta,p}:=\{Y\in L_{\mathbb{P}}^{p}: TY=e^{-i\theta}Y\}.$$
In particular, taking $\theta=0$, we get the classical statement of the mean and pointwise ergodic theorems for stationary sequences. 
\end{remark}

We remark also the following corollary.
\bigskip

\begin{cor}
\label{corcontozer}
With the notation of Theorem \ref{ergthedisfoutra}, if $e^{-i\theta}\notin Spec_{p}(T)$ (equivalently, if $e^{i\theta}\notin Spec_{p}(T)$), then
$$\frac{1}{n}S_{n}(Y,\theta)\to 0 \,\,\,\,\,\,\,\,\mbox{$\mathbb{P}-$a.s and in $L^{p}_{\mathbb{P}}$}\,.$$
\end{cor}

{\bf Proof:} This is a trivial consequence of the definition of $Spec_{p}(T)$ and the statement {\it 1.} in Theorem \ref{ergthedisfoutra}.\qed

What is the speed of convergence of the averages in Theorem \ref{ergthedisfoutra} and Corollary \ref{corcontozer}? By considering the case $\theta=0$ (the ``classical'' case) we see that this question does not admit an answer valid for any given $Y\in L^{p}_{\mathbb{P}}$ but, as shown by Peligrad and Wu in \cite{pelwu}, the central limit theorem (CLT) holds {\it for $\lambda$-a.e frequency $\theta\in [0,2\pi)$} under the additional (standard) assumptions $Y\in L^{2}_{\mathbb{P}}$, $EY=0$, and a certain regularity condition (see (\ref{reg})) to be discussed later \footnote{Peligrad and Wu's paper was preceded by Wu's \cite{Wu}, in which the same result is proved under the additional assumption $$\sum_{k>0}\frac{||E[X_{k}|\mathcal{F}_{0}]||_{_{2}}}{k}<\infty, $$
where for every $k\in\mathbb{Z}$, $X_{k}=f(\xi_{k})$ with   $(\xi_{k})_{k\in\mathbb{Z}}$ a stationary Markov Chain and $\mathcal{F}_{0}=\sigma(\xi_{k})_{k\leq 0}$ (see Section \ref{examarpro}  for a discussion related to this setting). 
}, this is Theorem \ref{cltpelwu} in this monograph. Peligrad and Wu show also  that under these hypotheses the functional CLT (FCLT) also holds for {\it averaged} frequencies (\cite{pelwu}, Theorem 2.1). This is Theorem \ref{pelwufclt} in page \pageref{pelwufclt}.

Our main goal in this work is to adjust these CLTs to the quenched setting, and to give more precise information about the nature of the asymptotic distribution for a given frequency $\theta\in[0,2\pi)$. We will also see that these ``adjustments'' are not simply an extension of the results just mentioned: not every centered process in $L^{2}_{\mathbb{P}}$ for which the CLT or the FCLT is valid admits quenched asymptotic distributions: the discrete Fourier transforms have to be (randomly) centered to remain orthogonal to the subspace of functions that are measurable with respect to the initial sigma-field (see Chapter \ref{quecon} for the terminology). For a precise description of these results see Chapter \ref{resandcom}.
 
\subsection{Dunford-Schwartz Operators and the Weak $L^{p}-$spaces}
\label{dunschweal2}

The purpose of this section is to prove that, for Dunford-Schwartz operators, a variation of Hunt and Young's maximal inequality (Theorem \ref{hunyou}) is available if
we refine the norm in $L^{p}_\mathbb{\mu}$ to the {\it weak} (induced) norm in $L^{p,\infty}$ for $p>1$. This result will be of importance to prove approximations involving the maxima of (normalized) partial sums via techniques akin to those involving Doob's maximal inequality. 
We start this section by recalling the notion of the {\it weak $L^{p}-$spaces, $L^{p,\infty}$}.

\medskip

\begin{dfn}[Weak $L^{p}$ spaces]
\label{weaLpdef}
Let $(\Omega,\mathcal{F},\mu)$ be a measure space, and given a measurable function $Y:\Omega\to\mathbb{C}$ and $0<p<\infty$, let $[Y]_{\mu,p}$ be given by
\begin{equation}
\label{deflpweanor}
[Y]_{\mu,p}:=\sup_{\alpha>0}(\alpha^{p}\mu([|Y|>\alpha])).
\end{equation}
We define the {\bfseries weak $L^p-$space associated to $\mu$}, $L^{p,\infty}_{\mu}$, as the topological space set-theoretically given by
\begin{equation}
\label{weaLpdef}
L^{p,\infty}_{\mu}:=\{Y: [Y]_{\mu,p}<\infty\},
\end{equation}
and whose topology is induced by the quasi-norm $[\,\cdot\,]_{\mu,p}$. If $p=\infty$ we define the $L^{\infty}$ weak space by $L^{\infty,\infty}_{\mu}:=L^{\infty}_{\mu}$. 
\end{dfn}

\medskip 

\begin{remark}
\label{remtoplpwea}
See section 1.1 in \cite{gra} for more details about $[\,\cdot\,]_{\mu,p}$. We point out in particular that, as stated in Exercise 1.1.12 in that book, the space $L^{p,\infty}_{\mu}$ is metrizable for every $p>0$ (and normable for $p>1$, a fact that we are just about to use).
\end{remark}

Markov's classical inequality shows that if $p>0$ is given, then for all  $Y\in L^{p}_{\mu}$, $[Y]_{p,\mu}\leq ||Y||_{p,\mu}$, so that $L^{p}_{\mu}$ is contained in the {\it weak $L^{p}$-space}, $L^{p,\infty}_{\mu}$. The inclusion  $L^{p}_{\mu}\subset L^{p,\infty}_{\mu}$ is continuous.

Even more (see for instance \cite{gra}, p.13, Exercise 1.1.12): if $p>1$ there exists a norm, $|||\,\,\,|||_{p,\mu}$ on $L^{p,\infty}_{\mu}$ with respect to which $L^{p,\infty}_{\mu}$ is a Banach space,  satisfying 
$$[\,\,\,]_{p,\mu}\leq|||\,\,\,|||_{p,\mu} \leq \frac{p}{p-1}[\,\,\,]_{p,\mu}.$$
In particular, if $p>1$ and $Y$ is any measurable function
\begin{equation}
\label{equlplpw}
(1-\frac{1}{p})|||Y|||_{p,\mu}\leq [Y]_{p,\mu} \leq ||Y||_{p,\mu}, 
\end{equation}

(with the convention $||Y||_{p,\mu}=\infty$ if $Y\notin L^{p}_{\mu}$, and analogously if $Y\notin L^{p,\infty}_{\mu}$).

Let $T:L^{1}_{\mu}\to L^{1}_{\mu}$ be a positive Dunford-Schwartz operator (Definitions  \ref{defdunsch} and \ref{defposope})   and define, for every $Y\in L^{1}_{\mu}$, 
\begin{equation}
\label{defmaxavedunsch}
M_{T}Y:=\sup_{n\in \mathbb{N}}\frac{1}{n}|\sum_{j=0}^{n-1} T^{n}Y|,
\end{equation} 
then (see \cite{kre}, Lemma 6.1, p.51) for every $\alpha>0$ the following Markov-type inequality holds
\begin{equation}
\label{maxmarine}
\mu([M_{T}|Y|>\alpha]) \leq\frac{1}{\alpha}E[|Y|\,I_{[M_{T}|Y|>\alpha]}]\leq \frac{1}{\alpha}E|Y|.
\end{equation}

Therefore, for $Y\in L^{p,\infty}_{\mu}$
$$(1-\frac{1}{p})|||(M_{T}|Y|^p)^{1/p}|||_{p,\mu}\leq  [(M_{T}|Y|^{p})^{1/p}]_{p,\mu}:=(\sup_{\alpha>0}\alpha^p \mu[(M_{T}|Y|^p)^{\frac{1}{p}}>\alpha])^{1/p} \leq ||Y||_{p,\mu},$$
where for the last inequality we used (\ref{maxmarine}).

We  summarize this discussion in the following proposition: 

\medskip

\begin{prop}
\label{prodominfnorpnor}
Let $(\Omega,\mathcal{F},\mu)$ be a measure space, let $T:L^{1}_{\mu}\to L^{1}_{\mu}$ be a positive Dunford-Schwartz operator (Definitions \ref{defdunsch} and \ref{defposope}), and define $M_{T}Y$ as in (\ref{defmaxavedunsch}). Then for every $p>1$  and every $Y\in L^{p,\infty}_{\mu}$
\begin{equation}
\label{equl2wl2}
|||(M_{T}|Y|^p)^{1/p}|||_{p,\mu}\leq \frac{p}{p-1}||Y||_{p,\mu}.
\end{equation}
\end{prop}

\section{$T-$Filtrations and Adapted Processes}
\label{adafil}

In this section we discuss the notions of  {\it $T$-filtrations} and {\it adapted processes}. We shall also briefly discuss, in a heuristic language, how this notion codifies the idea of ``initial conditions'' for a given stationary process, an idea that will be formalized in a precise way and used in Chapter \ref{quecon}.  
The setting of adapted filtrations will be a fundamental part of the assumptions present along the main results of this work. 

\subsection{Definitions and Examples}

Let us begin our discussion by giving the definition of a stationary  process and an ergodic process, a family of processes for which the main results of this monograph are devoted. 

\medskip

\begin{dfn}[Stationary Processes, Ergodic Processes]
\label{staproL2}
A stochastic process $(X_{n})_{n\in\mathbb{Z}}$ defined on a probability space $(\Omega,\mathcal{F},\mathbb{P})$  is called {\bfseries\itshape stationary} if there exists an invertible, bimeasurable, measure-preserving transformation $T:\Omega\to\Omega$ such that for all $k\in\mathbb{Z}$, $X_{k}=T^{k}X_{0}$. The process is called {\itshape\bfseries ergodic} if $T$ is ergodic. 
If $X_{0}\in L^{p}_{\mathbb{P}}$ (for some $p>0$), we say that $(X_{k})_{k\in\mathbb{Z}}$ is a {\bfseries $p-$integrable process}.
\end{dfn}

Let $(X_{k})_{k\in\mathbb{Z}}=(T^{k}X_{0})_{k\in\mathbb{Z}}$ be a stationary stochastic process defined on $(\Omega,\mathcal{F},\mathbb{P})$ (Definition \ref{staproL2}), let $\mathcal{F}_{0}\subset\mathcal{F}$ be a sub-sigma algebra of $\mathcal{F}$, and consider the following properties: 
\begin{enumerate}
\label{conadafil}
\item[{\itshape\bfseries F1}.] $T^{-1}$ is $\mathcal{F}_{0}$ measurable: $\mathcal{F}_{0}\subset T^{-1}\mathcal{F}_{0}$, where for all $k\in\mathbb{Z}$ $$T^{-k}\mathcal{F}_{0}:=\{A\in\mathcal{F}: T^{k}A\in\mathcal{F}_{0}\}=\{T^{-k}B:B\in\mathcal{F}_{0}\}.$$
(here $T^{k}$ denotes the $k-$fold composition of $T$).

\item[{\itshape\bfseries F2}.] $X_{0}$ is $\mathcal{F}_{0}-$measurable: $\sigma(X_{0})\subset\mathcal{F}_{0}$.
\end{enumerate}

Notice that if {\it F\,{1}.} is satisfied, the sequence $(\mathcal{F}_{k})_{k\in\mathbb{Z}}$ of sub-sigma algebras of $\mathcal{F}$ defined by
\begin{equation}
\label{deffk}
\mathcal{F}_{k}:=T^{-k}\mathcal{F}_{0}
\end{equation} 
is nondecreasing ($\mathcal{F}_{k}\subset\mathcal{F}_{k+1}$ for all $k\in\mathbb{Z}$), and that $X_{k}$ is $\mathcal{F}_{k}-$measurable if {\it F2.} holds. 

\medskip

\begin{dfn}[$T-$filtrations, adapted processes]
\label{defadafil}
Let $(\Omega,\mathcal{F},\mathbb{P})$ be a probability space and let $T:\Omega\to\Omega$ be an invertible, bimeasurable, measure-preserving transformation.

\begin{enumerate}
\item A {\bfseries $T-$filtration} is a filtration of the form (\ref{deffk}), where $\mathcal{F}_{0}$ satisfies {\it F1}.
\item The process $(T^{k}X_{0})_{k\in\mathbb{Z}}$ is {\bfseries adapted to $(\mathcal{F}_{k})_{k\in\mathbb{Z}}$} if {\it F2.} holds. 
\end{enumerate}

\end{dfn}

\medskip

Clearly, the ``trivial'' filtrations specified by $\mathcal{F}_{0}=\{\emptyset,\Omega\}$ and by $\mathcal{F}_{0}=\mathcal{F}$ are $T-$filtrations (any process defined in $(\Omega,\mathcal{F},\mathbb{P})$ is adapted to the second one), thus the existence of these objects poses no serious questions. 

Now note that the family of $T-$filtrations admits a partial order in the following way: given two $T-$filtrations $(\mathcal{F}_{k})_{k\in\mathbb{Z}}$ and $(\mathcal{G}_{k})_{k\in\mathbb{Z}}$, $(\mathcal{F}_{k})_{k\in\mathbb{Z}}\leq (\mathcal{G}_{k})_{k\in\mathbb{Z}}$ if $\mathcal{F}_{0}\subset\mathcal{G}_{0}$. With regards to this order, any stationary process $\mathbf{X}=(X_{k})_{k\in\mathbb{Z}}$ admits a {\it minimal} (and unique) adapted filtration in the obvious way: 

\medskip

\begin{dfn}[Minimal Adapted Filtration] 
\label{minadafildef}
Given a stationary process $(X_{k})_{k\in\mathbb{Z}}=(T^{k}X_{0})_{k\in\mathbb{Z}}$ (Definition \ref{staproL2}), define $\mathcal{M}_{0}$ by
\begin{equation}
\label{exiadafil}
\mathcal{M}_{0}:=\cap_{\alpha}\mathcal{G}_{\alpha}
\end{equation}
where the intersection runs over the sub-sigma algebras $\mathcal{G}_{\alpha}\subset\mathcal{F}$ for which $\sigma(X_{0})\subset \mathcal{G}_{\alpha}\subset T^{-1}\mathcal{G}_{\alpha}$. Then the filtration $(\mathcal{M}_{k})_{k\in\mathbb{Z}}:=(T^{-k}\mathcal{M}_{0})_{k\in\mathbb{Z}}$ is the {\bfseries minimal adapted filtration} associated to $(X_{k})_{k\in\mathbb{Z}}$: it is the smallest $T-$filtration for which $(X_{k})_{k\in\mathbb{Z}}$ is adapted (Definition \ref{defadafil}).
\end{dfn}
\medskip

To verify that $(\mathcal{M}_{k})_{k\in\mathbb{Z}}$ is indeed a $T-$filtration notice the following: it is clear that if $\{\mathcal{G}_{\alpha}\}_{\alpha}$ is the family described in Definition \ref{minadafildef} then 
$$\mathcal{M}_{0}\subset \cap_{\alpha}T^{-1}\mathcal{G}_{\alpha},$$ 
and note that for any given   $A\in \cap_{\alpha}T^{-1}\mathcal{G}_{\alpha}$, if $A=T^{-1}A_{\alpha}$ ($A_{\alpha}\in\mathcal{G}_{\alpha}$), then $A_{\alpha}=TA$, which proves ($A_{\alpha}$ does not depend on $\alpha$) that 
$$\cap_{\alpha}T^{-1}\mathcal{G}_{\alpha}\subset T^{-1}\cap_{\alpha}\mathcal{G}_{\alpha}=:T^{-1}\mathcal{M}_{0}.$$
The minimality of $(\mathcal{M}_{k})_{k\in\mathbb{Z}}$ among the adapted filtrations and the uniqueness of $\mathcal{M}_{0}$ are clear from the definition. 

Let us give a further definition, which we will need in subsequent sections.

\medskip

\begin{dfn}[Left and Right sigma-algebras]
\label{taisigalgdef}
If $(\mathcal{F}_{k})_{k\in\mathbb{Z}}$ is a $T-$filtration (Definition \ref{defadafil}), we define the {\bfseries left} and {\bfseries right} sigma algebras $\mathcal{F}_{-\infty}$, $\mathcal{F}_{+\infty}$ by
\begin{equation}
\label{deffmininf}
\mathcal{F}_{-\infty}:=\cap_{k\in\mathbb{Z}}\mathcal{F}_{k} \mbox{\,\,\,\,\,\,\,\,\,\,\it and \,\,\,\,\,\,\,\,\,\,} \mathcal{F}_{+\infty}:=\sigma(\cup_{k\in\mathbb{Z}}\mathcal{F}_{k}).
\end{equation}
\end{dfn}

To illustrate the notion of $T-$filtrations and adapted processes, it is convenient to look at an example, being perhaps the simplest nontrivial one that of {\it linear processes}. 

\medskip

\begin{exam}[Bernoulli Shifts and $T-$filtrations. Linear Processes]
\label{exalinpro}

Consider the space $\Omega=\mathbb{C}^{\mathbb{Z}}$ and, for every $j\in\mathbb{Z}$, let $x_{j}:\Omega\to\mathbb{C}$ be the projection on the $j-$th coordinate: for every $\omega=(\omega_{k})_{k\in\mathbb{Z}}\in\mathbb{C}^{\mathbb{Z}}$
\begin{equation}
\label{prokcoo}
x_{j}(\omega)=\omega_{j}.
\end{equation}
Let $\mathcal{F}$ be the sigma-algebra generated by the finite dimensional cylinders in $\Omega$. This is, by sets of the form 
\begin{equation}
\label{deffindimcyl}
H_{n,k,A}=\{\omega\in\Omega:(x_{n}(\omega),\dots,x_{n+k}(\omega))\in A\}
\end{equation}
where $(n,k)\in \mathbb{Z}\times \mathbb{N}$ and $A$ is a Borel set in $\mathbb{C}^{k+1}$.

Given a sequence $(\xi_{k})_{k\in\mathbb{N}}$ of random variables defined on a probability space $(\Omega',\mathcal{F}',\mathbb{P}')$, consider the map $\xi:\Omega'\to\Omega$ given by $\xi(\omega')=(\xi_{k}(\omega'))_{k\in\mathbb{Z}}$. By the $\pi-\lambda$ theorem, there exists a unique probability measure $\mathbb{P}$ in $(\Omega,\mathcal{F})$ such that, for every set $H_{n,k,A}$ as in (\ref{deffindimcyl}),
\begin{equation}
\mathbb{P}H_{n,k,A}=\mathbb{P}'\xi^{-1}H_{n,k,A}.
\end{equation}

If the sequence $(\xi_{k})_{k\in\mathbb{Z}}$ is stationary, in the sense that $\mathbb{P}'\xi^{-1}H_{n,k,A}$ is independent of $n$ for every fixed $k$ and $A$ (for instance if $(\xi_{k})_{k\in\mathbb{Z}}$ is i.i.d.), then the left shift $T:\Omega\to \Omega$, specified by $x_{k}(T\omega)=x_{k+1}(\omega)$ is an invertible bimeasurable, measure preserving transformation. In any case, if we define for every $k\in\mathbb{Z}$
\begin{equation}
\label{tfillinpro}
\mathcal{F}_{k}:=\sigma(x_{j})_{j\leq k}.
\end{equation}
then it is clear that  $(\mathcal{F}_{k})_{k\in\mathbb{Z}}$ is a $T-$filtration.

Note that, in this setting, the sequence of coordinate functions $(x_{k})_{k\in\mathbb{Z}}$ is a copy (in distribution) of $(\xi_{k})_{k\in\mathbb{Z}}$, thus we can replace ``$\xi_{k}$'' by ``$x_{k}$'' when referring to issues about distribution. 

Assume, in addition to stationarity, that $x_{0}\in L^{2}_{\mathbb{P}}$ (therefore $x_{k}\in L^{2}_{\mathbb{P}}$ for every $k\in\mathbb{Z}$), that $Ex_{0}=0$ and that the $x_{k}$'s are orthogonal: 
$$E[x_{k}\overline{x}_{l}]=\delta_{k,l}E[|x_{0}|^2]$$
where $\delta_{k,l}$ is the Kronecker $\delta-$function ($\delta_{k,l}\in\{0,1\}$, and $\delta_{k,l}=0$ if and only it $k\neq l$).

In this setting, given any sequence $(a_{k})_{k\in\mathbb{Z}}\in l^{2}({\mathbb{Z}})$ and any $k\in \mathbb{Z}$, the function
\begin{equation}
\label{deflinpro}
X_{k}(\omega)=\sum_{j\in\mathbb{Z}}a_{j}x_{k-j}(\omega)
\end{equation}
is well defined in the $L^{2}_{\mathbb{P}}$ sense: for any $k\in\mathbb{Z}$ and $N\in\mathbb{N}$, if $J\subset \mathbb{Z}$ is a finite set with $[-N,N]\cap \mathbb{Z}\subset J$ and $J':=J\setminus [-N,N]$
$$E[|\sum_{j\in J'}a_{j}x_{k-j}|^{2}]\leq E[|x_{0}|^2]\sum_{|j|> N}|a_{j}|^2.$$
which guarantees the convergence (well definition) of $X_{k}$  because $\sum_{j\in\mathbb{Z}}|a_{j}|^{2}<\infty$ ($L^{2}_{\mathbb{P}}$ is complete). Note also that for every $k\in\mathbb{Z}$, $X_{k}=T^{k}X_{0}$, and that $(X_{k})_{k\in\mathbb{Z}}$ is $(\mathcal{F}_{k})_{k\in\mathbb{Z}}-$adapted provided that $a_{j}=0$ for every $j<0$. We will explore a particular case of this example when proving Theorem \ref{nonqueconthe} in Chapter \ref{resandcom}.
\end{exam}
 
\medskip
 
\begin{remark}
\label{remrepprospecden}
Note that, by replacing $a_{k}$ by $a_{k}'=E[|x_{0}|^2]a_{k}$ and $x_{k}$ by $x_{k}/|| x_{k}||_{_{\mathbb{P},2}}$, we can assume without loss of generality that $E[|x_{k}|^2]=1$. 
\end{remark}

\subsection{Heuristic Interpretation} 
\label{heuint}

In a heuristic language, taking $\mathcal{F}_{0}$ as the ``information available to an observer'',  we require that $T$ preserves the information in $\mathcal{F}_{0}$ in order to obtain a $T-$filtration: any set of the form $TA$ for $A\in\mathcal{F}_{0}$ {\it still} belongs to $\mathcal{F}_{0}$. Pulling this heuristic further, we can thus think of $\mathcal{F}_{0}$ as the {\it deterministic part} of the dynamical system $T:\Omega\to\Omega$: an observer capable of knowing all the elements in $\mathcal{F}_{0}$ can follow their evolution under $T$ in a deterministic way. 

It is also usual to interpret $\mathcal{F}_{0}$ as ``the information from the past'', an interpretation that is particularly meaningful in the case of linear processes or, more generally, in the setting of functions of stationary Markov Chains (see Example \ref{examarproexa}). Example \ref{exalinpro} allows us to see how this naturally makes sense: $\mathcal{F}_{0}$, in this case, is the sigma algebra generated by all the coordinates  of the process up to the time $k=0$. 

Now,  a $T-$filtration is adapted to $(X_{k})_{k\in\mathbb{Z}}=(T^{k}X_{0})_{k\in\mathbb{Z}}$ if the information provided by $X_{0}$ is deterministic: the observer 
is able to know the outcome of the process at the time $k=0$. In this setting we can think of an adapted process as a process ``with given initial conditions'': the outcome of $X_{0}$ is known at the moment of running the process.

{\it How does this knowledge affect the asymptotics related to $(X_{k})_{k\in\mathbb{Z}}$?} This is, in broad terms, the question addressed by the notion of {\it quenched convergence}, to be discussed in Chapter \ref{quecon}. In short, and following the traditional interpretation of ``conditioning'', we will codify the ``influence'' of this knowledge by means of the conditional expectation with respect to $\mathcal{F}_{0}$.

\subsection{Interactions with the Koopman Operator}

To begin with this section let us  prove the following result.

\medskip

\begin{prop}
\label{relconexpkooope}
Let $(\mathcal{F}_{k})_{k\in\mathbb{Z}}$ be a $T-$filtration (Definition \ref{defadafil}), let $T:L^{1}_{\mathbb{P}}\to L^{1}_{\mathbb{P}}$ be the corresponding Koopman operator (Definition \ref{defkooope}) and for every $k\in\mathbb{Z}$, denote by $E_{k}$ the conditional expectation with respect to $\mathcal{F}_{k}$: for every $Y\in L^{1}_{\mathbb{P}}$ and $k\in\mathbb{Z}$:
\begin{equation}
\label{defek}
E_{k}Y:=E[Y|\mathcal{F}_{k}].
\end{equation}
Then for every $k,r\in\mathbb{Z}$
\begin{equation}
\label{relconexpkooopeequ}
T^{r}E_{k}=E_{k+r}T^{r}
\end{equation}
as operators in $L^{1}_{\mathbb{P}}$.
\end{prop}

\medskip

{\bf Proof:} Given $Y\in L^{1}_{\mathbb{P}}$, since clearly $T^{r}E_{k}Y$ is $\mathcal{F}_{r+k}-$measurable, we need to prove that for all $A\in\mathcal{F}_{k+r}$,
$$E[(T^{r}E_{k}Y)I_{A}]=E[(T^{r}Y)I_{A}].$$
To do so, let $A'=T^{r+k}A$. Notice that $A'\in\mathcal{F}_{0}$, and therefore $T^{-k}A'\in\mathcal{F}_{k}$. Using this and the fact that $T$ is measure preserving we see that
$$E[(T^{r}E_{k}Y)I_{A}]=E[T^{r}[(E_{k}Y)I_{T^{-k}A'}]]=E[E_{k}[YI_{T^{-k}A'}]]=E[YI_{T^{-k}A'}]$$
$$=E[(T^{r}Y)I_{T^{-(r+k)}A'}]]=E[(T^{r}Y)I_{A}]$$
as desired.\qed

\medskip

Now notice the following: assume that, for a given $k\in\mathbb{N}$, $(E_{0}T)^{k}=E_{0}T^{k}$, then: 
$$(E_{0}T)^{k+1}=E_{0}T(E_{0}T)^{k}=E_{0}TE_{0}T^{k}=E_{0}E_{1}T^{k+1}=E_{0}T^{k+1}$$
which shows, by induction on $k$, that for every $k\in\mathbb{N}^{*}$, 
\begin{equation}
\label{eque0tk}
(E_{0}T)^{k}=E_{0}T^{k}
\end{equation}
The operator $E_{0}T$ satisfies the following ergodic theorem.

\medskip

\begin{thm}[An Ergodic Theorem for Adapted $T-$filtrations]
\label{ergtheadafil}
In the context of Theorem \ref{ergthedisfoutra} and Proposition \ref{relconexpkooope}, given $p\geq 1$ and $Y\in L^{p}_{\mathbb{P}}$:
\begin{equation}
\label{limaveeot}
\lim_{n}\frac{1}{n}\sum_{k=0}^{n-1}E_{0}T^{k}Ye^{ik\theta}=E_{0}P_{\theta}Y, \mbox{\,\,\,\,\,\,\,\it $\mathbb{P}-$a.s and in $L^{p}_{\mathbb{P}}$.}
\end{equation}
\end{thm}

{\bf Proof:} Let $p$ and $Y$ be as in the given hypothesis. The convergence $\mathbb{P}-$a.s and in $L_\mathbb{P}^{p}$ follows via the following argument, similar to the one given for the proof of Theorem \ref{ergthedisfoutra}. The details are left to the reader.

{\it Convergence.} With the notation introduced in Definition \ref{extprospa} and the discussion following it, and defining 
$$\tilde{E}_{0}:=E[\,\cdot\,|\mathcal{B}\otimes\mathcal{F}_{0}]$$
(where the conditional expectation is with respect to $\lambda\times\mathbb{P}$), we can observe that (\ref{eque0tk}) holds with $\tilde{E}_{0}$ in place of $E_{0}$ and $\tilde{T}_{\theta}$ in place of $T$. An application of Theorem \ref{ergtheposdunsch} to the positive Dunford-Schwartz operator $\tilde{E}_{0}\tilde{T}_{\theta}$ allows one to see that
$$\frac{1}{n}\sum_{k=0}^{n-1}\tilde{E_{0}}\tilde{T}^{k}_{\theta}\tilde{Y}=\frac{1}{n}\sum_{k=0}^{n-1}\tilde{E_{0}}\widetilde{e^{ik\theta}T^{k}Y}\to_{n} {P_{\tilde{E}_{0}\tilde{T}_{\theta}}\tilde{Y}},$$
$\lambda\times\mathbb{P}-$a.s. and in $L^{p}_{\lambda\times\mathbb{P}}\,$. This implies that there exists a function\footnote{Actually given by  $\omega\mapsto e^{-iu}P_{\tilde{E}_{0}\tilde{T}_{\theta}}\tilde{Y}(u,\omega)$,
but this is not the representation that we are looking for.} 
 $P_{E_{0},T,\theta}Y:\Omega\to \mathbb{C}$ such that  
 $$\frac{1}{n}\sum_{k=0}^{n-1}E_{0}T^{k}Ye^{ik\theta}\to_{n} P_{E_{0},T,\theta}Y$$
$\mathbb{P}-$a.s. and in $L^{p}_{\mathbb{P}}$.   

{\it Limit function.} To identify $P_{E_{0},T,\theta}Y$ we use the continuity of $E_{0}$ as a linear operator in $L^{p}_{\mathbb{P}}\,$: since, according to Theorem \ref{ergthedisfoutra}, $S_{n}(Y,\theta)\to P_{\theta}Y$ in $L^{p}_{\mathbb{P}}$ as $n\to\infty$, $E_{0}S_{n}(Y,\theta)$ converges in $L^{p}_{\mathbb{P}}$ as $n\to\infty$, and
$$\lim_{n}E_{0}[S_{n}(Y,\theta)]=E_{0}[\lim_{n}S_{n}(Y,\theta)]=E_{0}P_{\theta}Y,$$
as claimed.\qed

\medskip

\begin{cor}
\label{dedmerpel}
In the context of Theorem \ref{ergtheadafil}. If $T$ is ergodic,
\begin{equation}
\label{dedmerpelequ}
\lim_{n}\frac{1}{n}\sum_{k=0}^{n-1}E_{0}T^{k}Y=E[Y], \mbox{\,\,\,\,\,\,\,\it $\mathbb{P}-$a.s and in $L^{p}_{\mathbb{P}}$.}
\end{equation}
\end{cor}

{\bf Proof:} Immediate from Remark \ref{ptkooopecas} and the fact that, in this case, $P_{0}Y=EY$, $\mathbb{P}$-a.s.\qed

\section{The Autocovariance Function and the Spectral Density}
\label{autcovfunandspeden}

In this section we discuss the notions of the {\it autocovariance function} and the {\it spectral density} of a stationary process. 

The autocovariance function is of relevance both in the theoretical and applied aspects of the theory of stochastic processes because it encodes the covariance structure of a given process (allowing inferences about, for instance, rates of decay), and its estimation is part of the study carried out here. 

To give a method to estimate the values of the autocovariance function we will introduce the closely related notion of the spectral density,  which can be computed studying the asymptotic behavior of the normalized $L^{2}-$norms of the discrete Fourier transforms (see Theorem \ref{spedenasalim} below).

\subsection{The Autocovariance Function}

To make the discussion clear let us start by recalling the following definition.

\medskip

\begin{dfn}[Nonnegative Definite Function]
\label{nonnegdef}
A function $\gamma:\mathbb{Z}\to \mathbb{C}$ is {\itshape\bfseries nonnegative definite} if for all vectors $\mathbf{c}=(c_{1},\dots,c_{n})\in \mathbb{C}^{n}$
\begin{equation}
\label{nonnegdefdef}
\sum_{i,j=1}^{n}c_{i}\gamma(i-j)\overline{c}_{j}\geq 0.
\end{equation}
\end{dfn}

These functions happen to be an important object in the study of the spectral properties of stationary sequences. The essential connection with this topic is {\it Herglotz's Theorem}:

\medskip

\begin{thm}[Herglotz's Theorem]
\label{herthe}
A function $\gamma:\mathbb{Z}\to \mathbb{C}$ is nonnegative definite (Definition \ref{nonnegdef}) if and only if  there exists a nondecreasing, right continuous bounded function $F:[-\pi,\pi)\to [0,+\infty)$ such that $F(0)=0$ and
$$\gamma(n)=\int_{-\pi}^{\pi}e^{in\theta}\,d\mu_{F}(\theta)$$
where $\mu_{F}$ denotes the measure induced by $F$: $\mu_{F}((a,b]):=F(b)-F(a)$ for all $[a,b)\subset [-\pi,\pi)$.
\end{thm}

{\bf Proof:} This is Theorem 4.3.1 in \cite{broanddav}.\qed

\medskip

\begin{remark}
\label{remherthe}
By the periodicity of the functions $\theta\mapsto e^{in\theta}$, the conclusion of this theorem remains valid if we substitute $[-\pi,\pi)$ by any interval of length $2\pi$. If we consider for instance the interval $[0,2\pi)$ with the Borel sigma-algebra $\mathcal{B}$ and the measure $$\mu'_{F}(A)=\mu_{F}(A\cap[0,\pi))+\mu_{F}((A-2\pi)\cap[-\pi,0))$$
 for all $A\in\mathcal{B}$ (where $A-2\pi:=\{a-2\pi:a\in A\}$), then the statement of Theorem \ref{herthe} remains valid replacing $[-\pi,\pi)$ by $[0,2\pi)$ and $\mu_{F}$ by $\mu_{F}'$.
\end{remark}
The connection of this theorem with the theory of stationary stochastic processes is made via the notion of the {\it autocovariance function}.

\medskip

\begin{dfn}[Autocovariance Function]
\label{defautcovfun}
Given a stationary square-integrable process $(X_{n})_{n\in \mathbb{Z}}$ (Definition \ref{staproL2}), we define the {\itshape\bfseries autocovariance function} $\gamma:\mathbb{Z}\to \mathbb{C}$ by
\begin{equation}
\label{defcovfun}
\gamma(n):=E[(X_{0}-EX_{0})(\overline{X}_{n}-\overline{EX_{n}})]=E[(X_{0}-EX_{0})(\overline{X}_{n}-\overline{EX_{0}})].
\end{equation}
\end{dfn}

\medskip

\begin{remark}
\label{autcovgivallcov}
Note that $\gamma(\cdot)$ encodes all the covariances of the process $(X_{n})_{n\in\mathbb{Z}}$: given integers $j, k$
$$E[(X_{j}-EX_{j})(\overline{X}_{k}-\overline{EX_{0}})]=E[T^{j}[(X_{0}-EX_{0})(\overline{X}_{k-j}-\overline{EX_{0}})]]=$$
$$E[(X_{0}-EX_{0})(\overline{X}_{k-j}-\overline{EX_{0}})]=\gamma(k-j).$$
\end{remark}

\medskip

Note also that the autocovariance function is {hermitian} ($\gamma(n)=\overline{\gamma(-n)}$) and  {nonnegative definite} (Definition \ref{nonnegdefdef}): given $\mathbf{c}=(c_{1},\dots,c_{n})\in \mathbb{C}^{n}$, and denoting by $\mathbb{X}_{n}=(X_{1}-EX_{1},\dots,X_{n}-EX_{n})$
$$\sum_{i,j}c_{i}\gamma(i-j)\overline{c}_{j}=E[\sum_{i,j}c_{i}(X_{i}-EX_{i})(\overline{X}_{j}-\overline{EX_{j}})\overline{c}_{j}]=E[|\mathbf{c}\cdot\mathbb{X}_{n}|^{2}]\geq 0$$
Herglotz's theorem implies therefore the following.

\medskip

\begin{prop}[Existence of the Spectral Measure]
\label{exispemea}
Given a stationary square-integrable process $\mathbf{X}=(X_{k})_{k\in\mathbb{Z}}$ (Definition \ref{staproL2}) there exists a (finite) measure $m_{_\mathbf{X}}$  on $([0,2\pi),\mathcal{B})$ such that the autocovariance function $\gamma:\mathbb{Z}\to \mathbb{C}$ of $(X_{k})_{k\in\mathbb{Z}}$ (Definition \ref{defautcovfun}) is given by
\begin{equation}
\label{forautcovfuncom}
n\mapsto \gamma(n)=\int_{0}^{2\pi}e^{in\theta}dm_{_\mathbf{X}}(\theta).
\end{equation}
\end{prop}
{\bf Proof:} Use Herglotz's theorem (Theorem \ref{herthe}) together with Remark \ref{remherthe}.\qed






\subsection{The Spectral Density}
\label{speden}

Our goal in this section is to connect the notion of the {\it spectral density} of a stationary process with the asymptotic theory of discrete Fourier transforms. To begin with, let us start by recalling the {\it F\'{e}jer-Lebesgue Theorem.}

\medskip

\begin{thm}[F\'{e}jer-Lebesgue]
\label{fejlebthe}
Let $f\in L^{1}_{\lambda}$ be given and denote by $\hat{f}$ the Fourier transform of $f$ (Definition \ref{foutra}). Then the sequence of functions $(C_{n}f)_{n\geq 0}$ defined by
\begin{equation}
\label{defsign}
C_{n}f(\theta):=\frac{1}{n}\sum_{j=0}^{n-1}\sum_{k=-j}^{j}\hat{f}(k)e^{ik\theta}
\end{equation}
converges to $f$ $\lambda-$a.s. 
\end{thm}

{\bf Proof:} See the proof of Theorem 3.3.3 in \cite{gra}. 
\qed

\medskip

\begin{remark}[$L^{p}_{\lambda}$ convergence in Theorem \ref{fejlebthe}]
\label{remconlp}
According to the referred proof in \cite{gra}, the convergence in (\ref{defsign}) holds in the $L^{p}_{\lambda}$ sense if $p>1$ and $f\in L^{p}_{\lambda}$: in such case there exists a constant $C_{p}$ such that 
$$||\sup_{n\in\mathbb{N}^{*}}|C_{n}f|||_{\lambda,p}\leq C_{p}||f||_{\lambda,p}.$$
In particular, $(C_{n}f-f)_{n\in\mathbb{N}^{*}}$ is  dominated in $L^{p}_{\lambda}$ (by $2\sup_{n\in\mathbb{N}^{*}}|C_{n}f|$), and the dominated convergence theorem implies that $||C_{n}f-f||_{\lambda,p}\to 0$ as $n\to\infty$.
\end{remark}

We  saw in Proposition \ref{exispemea} that for a stationary square-integrable process $\mathbf{X}=(X_{k})_{k\in\mathbb{Z}}$ there exists a measure $m_{_\mathbf{X}}$ on $([0,2\pi),\mathcal{B})$ such that
\begin{equation}
\label{equspemea}
E[(X_{0}-EX_{0})(\overline{X}_{k}-\overline{EX_{0}})]=\int_{0}^{2\pi} e^{ik\theta}d\,m_{_\mathbf{X}}(\theta)
\end{equation}
for all $k\in\mathbb{Z}$.

Note that, {\it if $F$ is absolutely continuous with respect to $\lambda$} and 
\begin{equation}
\label{defspeden}
f(\theta):=\frac{d m_{_\mathbf{X}}}{d \lambda} (\theta)
\end{equation}
is the Radon-Nikodym derivative of $m_{_\mathbf{X}}$ with respect to $\lambda$, then it follows from (\ref{equspemea}) that, if we denote by $\hat{f}$ the Fourier transform of $f$, then
 $$\hat{f}(-k)=\int_{0}^{2\pi}f(\theta) e^{ik\theta}d\lambda(\theta)=E[(X_{0}-EX_{0})(\overline{X}_{k}-\overline{EX_{0}})]=\gamma(k)$$
 This is: {\it the autocovariance function of $(X_{k})_{k\in\mathbb{Z}}$ is given by the sequence of the negative Fourier coefficients of $f$}. This justifies the following definition.

\medskip

\begin{dfn}[Spectral Density]
\label{defspeden}
 We say that a stationary square-integrable process $(X_{k})_{k\in\mathbb{Z}}$ (Definition \ref{staproL2}) admits a {\itshape\bfseries spectral density} if there exists a nonnegative function $f\in L^{1}_{\lambda}([0,2\pi))$ such that for every $k\in \mathbb{Z}$
\begin{equation}
\label{spedenequ}
\hat{f}(-k)=\gamma(k),
\end{equation}
where $\hat{f}$ denotes the Fourier transform of $f$ (Definition \ref{foutra}) and  $\gamma$ is the autocovariance function of $(X_{k})_{k\in\mathbb{Z}}$ (Definition \ref{defautcovfun}).

\end{dfn}

\medskip

\begin{remark}
\label{autcovgivbyfoucoe}
It is an immediate consequence of Theorem \ref{fejlebthe} that if $(X_{k})_{k\in\mathbb{Z}}$ admits a spectral density $f$ then it is unique (up to a set of $\lambda-$measure zero). Note also that if the process $(X_{k})_{k\in\mathbb{Z}}$ is real-valued and admits a spectral density $f$, then $\hat{f}(k)=\gamma(k)$ for all $k\in\mathbb{Z}$ ($\gamma$ is hermitian and real valued, i.e., even).
\end{remark}





\subsection{Regular Processes}
Let us make now a short digression that will allow us to relate the notion of $T-$filtrations to the existence of the spectral density.

\medskip 

\begin{dfn}[Regularity of an Adapted Process]
\label{regcondef}
Let $(\mathcal{F}_{k})_{k\in\mathbb{Z}}$ be a $T-$filtration and let $(X_{k})_{k\in\mathbb{Z}}=(T^{k}X_{0})_{k\in\mathbb{Z}}$ be a $(\mathcal{F}_{k})_{k\in\mathbb{Z}}-$adapted stationary square-integrable process (Definition \ref{defadafil}). The process is called {\bfseries regular} (with respect to $(\mathcal{F}_{k})_{k\in\mathbb{Z}}$) if $E[X_{0}|\mathcal{F}_{-k}]$ converges to $0$ in $L^{2}_{\mathbb{P}}$. This is, if
\begin{equation}
\label{reg}
\lim_{k\to\infty} E[|E[X_{0}|\mathcal{F}_{-k}]|^2]=0.
\end{equation}

\end{dfn}

\begin{remark}
\label{remequreg}
Recall the notation introduced in Definition \ref{taisigalgdef} and Proposition \ref{relconexpkooope}. Since for every $p\geq 1$, $k\in\mathbb{Z}$ and $Y\in L^{p}_{\mathbb{P}}$ the process $(E_{-j+k}Y)_{j\geq 0}$ is a reverse martingale in $L^{p}_{\mathbb{P}}$, the reverse martingale convergence theorem (see Theorem 5.6.1 and Exercise 5.6.1 in \cite{dur}) and the continuity of the Koopman operator $T$ imply that the following equalities hold  both $\mathbb{P}-$a.s and in $L^{p}_{\mathbb{P}}$:
\begin{equation}
\label{equinvprotai}
T^{k}E_{-\infty}Y=T^{k}\lim_{j\to\infty}E_{-j+k}Y=\lim_{j\to\infty}T^{k}E_{-j+k}Y=\lim_{j\to\infty}E_{-j}T^{k}Y=E_{-\infty}T^{k}Y.
\end{equation} 
 In particular, (\ref{reg}) is equivalent to the following condition: {\it for every $k\in\mathbb{Z}$}
 \begin{equation}
 \label{regconequ}
 E[X_{k}|\mathcal{F}_{-\infty}]=0 \,\,\,\,\,\,\mbox{($\mathbb{P}-$a.s.).}
 \end{equation} 
 \end{remark}

Now notice the following: if $(\mathcal{M}_{k})_{k\in\mathbb{Z}}$ is the minimal adapted filtration associated to $(X_{k})_{k\in\mathbb{Z}}$ (Definition \ref{minadafildef}) and $\mathcal{M}_{-\infty}$ is its left sigma algebra of this filtration (Definition \ref{taisigalgdef})  then, since (\ref{regconequ}) is equivalent to the condition $E[X_{k}I_{A}]=0$ for every $A\in\mathcal{F}_{-\infty}$ and $\mathcal{M}_{-{\infty}}\subset \mathcal{F}_{-\infty}$, we have that if $(X_{k})_{k\in\mathbb{Z}}$ is regular with respect to $(\mathcal{F}_{k})_{k\in\mathbb{Z}}$, then for every $k\in\mathbb{Z}$
\begin{equation}
\label{regdesminfilequ}
E[X_{k}|\mathcal{M}_{-\infty}]=0.
\end{equation}

In virtue of Remark \ref{remequreg} this gives the following result.

\medskip

\begin{prop}[Regularity and Minimal Adapted Fitrations]
\label{regandminfil}
If a process $(X_{k})_{k\in\mathbb{Z}}$ is regular with respect to some (adapted) filtration $(\mathcal{F}_{k})_{k\in\mathbb{Z}}$ (Definition \ref{regcondef}), then it is regular with respect to its minimal adapted filtration $(\mathcal{M}_{k})_{k\in\mathbb{Z}}$ (Definition \ref{minadafildef}).
 \end{prop}
 
Since the minimal adapted filtration is unique, this shows that the notion of regularity of an adapted process can be made ``unambiguous'' if we declare a process ``regular'' if it is regular with respect to its minimal adapted filtration.

\subsection{On the Existence of the Spectral Density}

What stationary processes $\mathbf{X}=(X_{k})_{k\in\mathbb{Z}}$ admit a spectral density? First, as stated by Theorem 31.28 in \cite{bilpromea}, the existence of $f$ is equivalent to the {\it absolute continuity} (in the sense of real calculus) of the distribution function $F_{_\mathbf{X}}$ of  $m_{_\mathbf{X}}$ ($F_{_\mathbf{X}}(t):=m_{_\mathbf{X}}((-\infty,t])$): $f$ exists if and only if for every $\epsilon>0$ there exists $\delta>0$ such that if $\{[a_{k},b_{k}]\}_{k=1}^{n}$ is any collection of disjoint intervals contained in $[0,2\pi)$ for which $\sum_{k=1}^{n}(b_{k}-a_{k})<\delta$, then  $\sum_{k=1}^{n}(F_{_\mathbf{X}}(b_{k})-F_{_\mathbf{X}}(a_{k}))<\epsilon $.

An interesting question is how to characterize the existence of the spectral density in terms of rates of decay of the values of the autocovariance function $\gamma$. 

To be more precise, note first that by the proof of Proposition \ref{exispemea}, $m_{_\mathbf{X}}$ is characterized by the equation
 (\ref{forautcovfuncom}). 
 
Now, as proved in \cite{broanddav}, Corollary 4.3.1 (together with Remark \ref{remherthe} above), {\it every} function $\gamma:\mathbb{Z}\to\mathbb{C}$ that can be represented in the form
\begin{equation}
\label{equfoutramea}
\gamma(k)=\int_{0}^{2\pi}e^{ik\theta }d{\mu}(\theta)
\end{equation}
 for some finite measure $\mu$ on $([0,2\pi),\mathcal{B})$, is the autocovariance function of some stationary square-integrable process. Note again that $\gamma(-k)$ is (by definition) the $k-$th Fourier coefficient of the measure $\mu$.

Thus, since any finite measure on $([0,2\pi),\mathcal{B})$ is determined by the sequence of its Fourier coefficients (\cite{bhaway}, Proposition 6.3), the problem of the absolute continuity of $m_{_\mathbf{X}}$ (for any $\mathbf{X}$) is equivalent to the following question:

{\bf Question: }{\it Let $\mu$ be a finite measure on $[0,2\pi)$ and let $\gamma:\mathbb{Z}\to\mathbb{C}$ be given by (\ref{equfoutramea}). What conditions on the sequence $(\gamma(k))_{k\in\mathbb{Z}}$ are necessary and/or sufficient to guarantee that $\mu$ is absolutely continuous with respect to $\lambda$?}

Any answer to this question has an immediate translation to a criterion about the existence of the spectral density of a stationary process $(X_{k})_{k\in\mathbb{Z}}$ in terms of the sequence of its covariances $(\gamma(k))_{k\in\mathbb{Z}}=(E[(X_{0}-EX_{0})(\overline{X}_{k}-\overline{EX_{0}})])_{k\in\mathbb{Z}}$. 

The following criterion is just a reformulation of one of the equivalences of Theorem 1 in \cite{lifpel}.

\medskip

\begin{thm}[Absolute Continuity via Fourier Coefficients]
\label{absconviafoucoe}
Let $\mu$ be a finite measure in $[0,2\pi)$ and define $\gamma:\mathbb{Z}\to\mathbb{C}$ by (\ref{equfoutramea}). Then the following are equivalent
\begin{enumerate}
\item $\mu$ is absolutely continuous with respect to $\lambda$.
\item There exists a sequence of complex numbers $(a_{k})_{k\in\mathbb{Z}}\in l^{2}_{\mathbb{Z}}$ such that for all $n\in\mathbb{Z}$
\begin{equation}
\label{necsufconabsconspemea}
\gamma(n)=\sum_{j\in\mathbb{Z}}a_{j}\overline{a}_{j+n}.
\end{equation}
\end{enumerate}
\end{thm}

\medskip

\begin{remark}
\label{remabsconviafoucoe}
The convolution $a*b$ between two sequences $a=(a_{k})_{k\in\mathbb{Z}}$ and $b=(b_{k})_{k\in\mathbb{Z}}$ in $l^{2}(\mathbb{Z})$ is equal to the sequence ${a}*{b}=((a*b)(k))_{k\in\mathbb{Z}}\in l^{1}(\mathbb{Z})$ given by 
$$({a}*{b})(k)=\sum_{k\in\mathbb{Z}}a_{j}b_{j-k}.$$
Using the fact that every function $f\in L^{1}_{\lambda}$ is the product of two functions in $L^{2}_{\lambda}$ (consider a branch of the square root), it is possible to show that the correspondence $l^{2}(\mathbb{Z})\times l^{2}(\mathbb{Z})\to L^{1}_{\lambda}$ given by $(a,b)\mapsto f_{(a,b)}$ where 
$$f_{(a,b)}(\theta):=\lim_{n}\frac{1}{n}\sum_{k=0}^{n-1}\sum_{j=-k}^{k}(a*b)(j)e^{ik\theta}$$
is surjective (the limit is $\lambda-$a.e well defined by Theorem \ref{fejlebthe}). In this language, Theorem \ref{absconviafoucoe} can be understood as the statement that every absolutely continuous measure in $\mathbb{T}$ corresponds to a sequence 
of the form $a*\overline{a}$ for some $a\in l^{2}(\mathbb{Z})$, where $\overline{a}=(\overline{a}_{k})_{k\in\mathbb{Z}}$ is the conjugate sequence of $a$. The corresponding density is actually given by $f_{(a,\overline{a})}$, which is the same as the function
$$\theta\mapsto |\sum_{k\in\mathbb{Z}} a_{k}e^{ik\theta}|^{2}.$$ 

\end{remark}

In the context of Theorem \ref{absconviafoucoe}, it is possible to give other sufficient conditions implying the absolute continuity of $\mu$ with respect to $\lambda$. Assume for instance that $(\gamma(k))_{k\in\mathbb{Z}}\in l^{2}(\mathbb{Z})$. Then by Proposition \ref{bijcorL2l2}  the function
$$f(\theta)=\sum_{k\in\mathbb{Z}}\gamma(k)e^{ik\theta}$$
is the spectral density of $(X_{k})_{k\in\mathbb{Z}}$. In particular, the condition
\begin{equation}
\label{secmoncon}
\sum_{k\geq 0}|E[(X_{0}-EX_{0})(\overline{X}_{k}-\overline{EX_{0}})]|^{2}<\infty
\end{equation}
implies the existence of the spectral density. For other sufficient conditions see for instance \cite{cunmerpel} and the references therein.



Our next result is the following:

\medskip

\begin{prop}[Spectral Density via Regularity]
\label{regcon}
Every regular process (Definition \ref{regcondef}) admits a spectral density.
\end{prop}

{\bf Proof:}\footnote{This argument follows the proof of Theorem 3 in \cite{lifpel}.} With the notation (\ref{defek})  define, for every $k\in\mathbb{Z}$ and $Y\in L^{1}_{\mathbb{P}}$,  
\begin{equation}
\label{defpkequ}
\mathcal{P}_{k}Y:=(E_{k}-E_{k-1})Y.
\end{equation}
Remember now that $L^{2}_{\mathbb{P}}(\mathcal{F}_{k})\subset L^{2}_{\mathbb{P}}$ denotes the subspace of functions that are measurable with repsect to $\mathcal{F}_{k}$, and denote also by
\begin{equation}
\label{deflfkminlfkminone}
V_{k}:=L^{2}_{\mathbb{P}}(\mathcal{F}_{k})\ominus L^{2}_{\mathbb{P}}(\mathcal{F}_{k-1})
\end{equation}
the orthogonal complement of $L^{2}_{\mathbb{P}}(\mathcal{F}_{k-1})$ in $L^{2}_{\mathbb{P}}(\mathcal{F}_{k})$. Then, by defining $\mathcal{F}_{-\infty}$ as in (\ref{deffmininf}), we see that
\begin{equation}
\label{decl2f0}
L^{2}_{\mathbb{P}}(\mathcal{F}_{k})\ominus L^{2}_{\mathbb{P}}(\mathcal{F}_{-\infty})=\bigoplus_{j\leq k}V_{j}
\end{equation}
and that $\mathcal{P}_{k}$, restricted to $L^{2}_{\mathbb{P}}$, is just the orthogonal projection on the space $V_{k}$. In particular, since under (\ref{reg}), $X_{k}\in L^{2}_{\mathbb{P}}(\mathcal{F}_{k})\ominus L^{2}_{\mathbb{P}}(\mathcal{F}_{-\infty})$ (see Remark \ref{remequreg}), we have that for every $k\in\mathbb{Z}$
$$X_{k}=\sum_{j\in \mathbb{Z}}\mathcal{P}_{-j}X_{k}=\sum_{j\geq -k  }\mathcal{P}_{-j}X_{k}$$
and therefore, by orthogonality and Proposition \ref{relconexpkooope},
\begin{equation}
\label{pitequ}
E[|X_{0}|^2]=\sum_{k\geq 0}E[|\mathcal{P}_{-k}X_{0}|^2]=\sum_{k\geq 0}E[|\mathcal{P}_{0}X_{k}|^{2}].
\end{equation}
It follows from Proposition \ref{proconsupL2} that the function $f:[0,2\pi)\to[0,\infty)$ specified by
$$f(\theta)=E[|\sum_{k\geq 0}\mathcal{P}_{0}X_{k}e^{ik\theta}|^2] $$
is well defined. More precisely: for $\lambda-$a.e $\theta$ the integrand converges $\mathbb{P}-$a.s and the integral (with respect to $\mathbb{P}$) makes sense. 

We claim that $f$ is the spectral density of $(X_{k})_{k\in\mathbb{Z}}$.


Fix $k\in\mathbb{Z}$ and begin by noticing that, by orthogonality and Proposition \ref{relconexpkooope},
$$E[X_{0}\overline{X_{-k}}]=E[(\sum_{j\in\mathbb{Z}}\mathcal{P}_{-j}X_{0})(\overline{\sum_{l\in \mathbb{Z}}\mathcal{P}_{-l}X_{-k}})]=
\sum_{j\geq 0}E[(\mathcal{P}_{-j}X_{0})(\overline{\mathcal{P}_{-j}X_{-k}})]=$$
\begin{equation}
\label{expfoucoe}
\sum_{j\geq 0}E[(\mathcal{P}_{0}X_{j})(\overline{\mathcal{P}_{0}X_{j-k}})].
\end{equation}
Our goal is thus to prove that $f$ is integrable and that for every $k\in\mathbb{Z}$, $\hat{f}(k)$ is equal to the last term in (\ref{expfoucoe}).

To begin with, define $\Omega_{1}$ as the set of probability one
\begin{equation}
\label{defomeoneexispeden}
\Omega_{1}:=\{\omega\in\Omega:\sum_{k\geq 0}|\mathcal{P}_{0}X_{k}(\omega)|^{2}<\infty\}.
\end{equation}

\begin{enumerate}
\item {\it The function $f$ is integrable.} 
By (\ref{pitequ}) and Carleson's Theorem (Theorem \ref{carthe})
the function 
$$\theta\mapsto \sum_{k\geq 0}\mathcal{P}_{0}X_{k}(\omega)e^{ik\theta}$$ 
is well defined (the series converges $\lambda-$a.s) for every $\omega\in\Omega_{1}$. It follows from the dominated convergence theorem and (\ref{hunyouran}) ({see also the line following that equation}) 
that for every $\omega\in\Omega_{1}$
$$\int_{0}^{2\pi}|\sum_{k\geq 0}\mathcal{P}_{0}X_{k}(\omega)e^{ik\theta}|^{2}d\lambda(\theta)=\sum_{k\geq 0}|\mathcal{P}_{0}X_{k}(\omega)|^{2}<\infty.$$
Now, by Tonelli's theorem, Proposition \ref{relconexpkooope} and the monotone convergence theorem
$$\int_{0}^{2\pi}f(\theta)d\lambda(\theta)=E[\int_{0}^{2\pi}|\sum_{k\geq 0}\mathcal{P}_{0}X_{k}e^{ik\theta}|^2d\lambda(\theta)]=\sum_{k\geq 0}E[|\mathcal{P}_{0}X_{k}|^{2}]=$$
$$\sum_{k\geq 0}E[|\mathcal{P}_{-k}X_{0}|^{2}]=E[|X_{0}|^{2}],$$
which shows that $f\in L^{1}_{\lambda}$. Note also that this proves the required equality $\hat{f}(0)=E[|X_{0}|^2]$.

\item {\it The Fourier coefficient $\hat{f}(k)$ is given by the last term of  (\ref{expfoucoe})}. Fix $k\in\mathbb{Z}$ and note first that, by (\ref{hunyouran}) and the dominated convergence theorem, the following holds: for every $\omega\in\Omega_{1}$,
$$\int_{0}^{2\pi}e^{-ik\theta}|\sum_{j\geq 0}\mathcal{P}_{0}X_{j}(\omega)e^{ij\theta}|^{2}d\lambda(\theta)=\lim_{N}\int_{0}^{2\pi}e^{-ik\theta}|\sum_{j=0}^{N} \mathcal{P}_{0}X_{j}(\omega)e^{ij\theta}|^{2}d\lambda(\theta)=$$
$$\lim_{N}\int_{0}^{2\pi}\sum_{j=0}^{N}\sum_{l=0}^{N} (\mathcal{P}_{0}X_{j}(\omega))(\overline{\mathcal{P}_{0}X_{l}(\omega)})e^{i(j-l-k)\theta}d\lambda(\theta)=\lim_{N}\sum_{j=0}^{N}(\mathcal{P}_{0}X_{j})(\overline{\mathcal{P}_{0}X_{j-k}})=$$
$$\sum_{j\geq 0}(\mathcal{P}_{0}X_{j}(\omega))(\overline{\mathcal{P}_{0}X_{j-k}(\omega)}).$$
Now notice the following: by the Cauchy-Schwartz inequality and Proposition \ref{relconexpkooope},
$$\sum_{j\geq 0}E[|(\mathcal{P}_{0}X_{j})(\overline{\mathcal{P}_{0}X_{j-k}})|]\leq \sum_{j\geq 0}(E[|\mathcal{P}_{0}X_{j}|^{2}])^{1/2}(E[|\mathcal{P}_{0}X_{j-k}|^{2}])^{1/2}\leq$$
$$(\sum_{j\geq 0}E[|\mathcal{P}_{0}X_{j}|^{2}])^{1/2}(\sum_{l\geq 0}E[|\mathcal{P}_{0}X_{l-k}|^{2}])^{1/2}=(\sum_{j\geq 0}E[|\mathcal{P}_{0}X_{j}|^{2}])^{1/2}(\sum_{l\geq -k}E[|\mathcal{P}_{0}X_{l}|^{2}])^{1/2}=$$
$$(\sum_{j\geq 0}E[|\mathcal{P}_{-j}X_{0}|^{2}])^{1/2}(\sum_{l\geq -k}E[|\mathcal{P}_{-l}X_{0}|^{2}])^{1/2}= \sum_{j\geq 0}E[|\mathcal{P}_{0}X_{j}|^{2}]=E[|X_{0}|^{2}],$$
 
and therefore, by the dominated convergence theorem
\begin{equation}
\label{casi}
E[\sum_{j\geq 0}(\mathcal{P}_{0}X_{j}(\omega))(\overline{\mathcal{P}_{0}X_{j-k}(\omega)})]=\sum_{j\geq 0}E[(\mathcal{P}_{0}X_{j}(\omega))(\overline{\mathcal{P}_{0}X_{j-k}(\omega)})].
\end{equation}

This information allows us to finish the proof: an application of Fubini's theorem and (\ref{casi}) gives
$$\hat{f}(k):=\int_{0}^{2\pi}f(\theta)e^{-ik\theta}d\lambda(\theta)=\int_{0}^{2\pi}E[|\sum_{j\geq 0}\mathcal{P}_{0}X_{j}e^{ij\theta}|^2]e^{-ik\theta}d\lambda(\theta)=$$
$$ \int_{0}^{2\pi}E[(\sum_{j\geq 0}\mathcal{P}_{0}X_{j}e^{ij\theta})(\sum_{l\geq 0}\overline{\mathcal{P}_{0}{X}_{j}}e^{-ij\theta})]e^{-ik\theta}d\lambda(\theta)=$$
$$E[\int_{0}^{2\pi}\sum_{j,l\geq 0}(\mathcal{P}_{0}X_{j})(\overline{\mathcal{P}_{0}X_{l}})e^{i(j-l-k)\theta}d\lambda(\theta)]=$$
$$E[\sum_{k\geq 0}(\mathcal{P}_{0}X_{j})\overline{(\mathcal{P}_{0}X_{j-k})}]=\sum_{k\geq 0}E[(\mathcal{P}_{0}X_{j})(\overline{\mathcal{P}_{0}X_{j-k}})].$$
which is the desired expression.\qed
\end{enumerate}




\medskip

\begin{remark}
\label{remintmar}
Note that the function
\begin{equation}
\label{defd0}
D_{0}(\theta)=\sum_{k\geq 0}\mathcal{P}_{0}X_{k}e^{ik\theta}
\end{equation}
is an adapted martingale difference with respect to $(\mathcal{F}_{k})_{k\geq 0}$: $D_{0}(\theta)$ is $\mathcal{F}_{0}-$measurable and, for every $k\geq 1$, $E_{k-1}T^{k}D_{0}(\theta)=T^{k}E_{-1}D_{0}(\theta)=0$. In subsequent proofs we will show that if for every $n\geq0$ we define
\begin{equation}
\label{defMn}
M_{n}(\theta):=\sum_{k=0}^{n-1}T^{k}D_{0}(\theta)e^{ik\theta},
\end{equation}
then the quenched asymptotics of $(M_{n}(\theta))_{n\geq 0}$ can be transported to corresponding results for $(S_{n}(\theta)-E_{0}S_{n}(\theta))_{n\geq 0}$: the study of quenched limit theorems for adapted martingales will therefore play an essential role in some of the proofs of the forthcoming results.
\end{remark}  

\subsection{Estimating the Spectral Density} 

According to Theorem \ref{ergthedisfoutra} in Section \ref{ergthefoutra}, for a stationary square-integrable process $(X_{k})_{k\in\mathbb{Z}}=(T^{k}X_{0})_{k\in\mathbb{Z}}$ (Definition \ref{staproL2}), the  averaged discrete Fourier transforms $S_{n}(\theta)/n$ (see Definition \ref{defdisfoutra}) converge almost surely and in $L^{2}_{\mathbb{P}}$  to a function $P_{\theta}X_{0}$ with the property that
$$TP_{\theta}X_{0}=e^{-i\theta}P_{\theta}X_{0}.$$
In particular, as we stated in Corollary \ref{corcontozer}, $P_{\theta}X_{0}=0$ when $e^{-i\theta}\notin Spec_{p}(T)$ (Definition \ref{poispe}).  Even more is true: when $\mathcal{F}$ is countably generated, Proposition \ref{prospecou} implies that the set 
\begin{equation}
\label{setofexclimequzer}
\{\theta\in[0,2\pi): P_{\theta}X_{0}\neq 0\}
\end{equation}
has $\lambda-$measure zero.

We will go a bit further now by showing that, if $(X_{k})_{k\in\mathbb{Z}}$ is a centered process and admits a spectral density then it is given by the asymptotic variance of the properly normalized discrete Fourier transforms. The result is the following:

\medskip

\begin{thm}[Spectral Density as an Asymptotic Variance]
\label{spedenasalim}
Let $(X_{k})_{k\in \mathbb{Z}}$ be a stationary process (Definition \ref{staproL2}) with $E[X_{0}]=0$,  and assume that $(X_{k})_{k\in\mathbb{Z}}$ admits a spectral density $\sigma^{2}:[0,2\pi)\to [0,+\infty)$ (Definition \ref{defspeden}). Then 
\begin{equation}
\label{forspeden}
\sigma^{2}(\theta)=\lim_{n}\frac{1}{n}E[|S_{n}(\theta)|^{2}]
\end{equation}
in the sense of Definition \ref{deflimfun} (this is, (\ref{forspeden}) holds for $\lambda-$a.e $\theta$), where $S_{n}(\theta)$ is the $n-$th discrete Fourier transform of $(X_{k})_{k\in\mathbb{Z}}$ at $\theta$ (Definition \ref{defdisfoutra}).

\end{thm}

{\bf Proof:} Let us start with the following observation: given a sequence of complex numbers $(c_{k})_{k\in\mathbb{Z}}$, it is not hard to see that
\begin{equation}
\label{equcesmea}
\sum_{j=0}^{n-1}\sum_{k=0}^{n-1}c_{j-k}=\sum_{j=-(n-1)}^{n-1}(n-|j|)c_{j}=\sum_{j=0}^{n-1}\sum_{k=-j}^{j}c_{k}.
\end{equation}

Now assume that $(X_{k})_{k\in\mathbb{Z}}=(T^{k}X_{0})_{k\in\mathbb{Z}}$ is a stationary process with $E[X_{0}]=0$. 
Denoting by $\gamma:\mathbb{Z}\to\mathbb{C}$ the autocovariance function of $(X_{k})_{k\in\mathbb{Z}}$ ($\gamma(k)=E[X_{0}\overline{X_{k}}]$) and taking $$c_{j}:=E[X_{0}\overline{X}_{-j}]e^{ij\theta}=\gamma(-j)e^{ij\theta}$$ it follows from (\ref{equcesmea}) that

$$E|\frac{1}{\sqrt{n}}S_{n}(\theta)|^{2}=\sum_{j=0}^{n-1}\sum_{k=0}^{n-1}E[X_{j}\overline{X}_{k}]e^{ij\theta}e^{-ik\theta}=\frac{1}{n}\sum_{j=0}^{n-1}\sum_{k=0}^{n-1}{E[X_{0}\overline{X}_{-(j-k)}]}e^{i(j-k)\theta}=$$
$$
=\frac{1}{n} \sum_{j=0}^{n-1}\sum_{k=-j}^{j}E[X_{0}\overline{X}_{-k}]e^{ik\theta}. $$

and the conclusion follows at once from Definition \ref{defspeden} and Theorem \ref{fejlebthe}.\qed

\subsection*{Conclusive Remarks}
\label{conrem}
As we have suggested along the previous sections, this monograph is devoted to the (quenched) asymptotic behavior of the normalized averages 
\begin{equation}
\label{defavedisfoutra}
A_{n}(\theta,\omega)=\frac{1}{n}\sum_{k=0}^{n-1}X_{k}(\omega)e^{ik\theta}=\frac{1}{n}S_{n}(\theta)
\end{equation}
of the discrete Fourier transforms of a stationary, square-integrable centered process $(X_{k})_{k\in\mathbb{Z}}=(T^{k}X_{0})_{k\in\mathbb{Z}}$ (Definition \ref{staproL2}) defined on a probability space $(\Omega,\mathcal{F},\mathbb{P})$. We have seen (Corollary \ref{corcontozer}, Proposition \ref{prospecou}, and Theorem \ref{spedenasalim}) that when $\mathcal{F}$ is countably generated (Definition \ref{deffcougen}) and $(X_{k})_{k\in\mathbb{Z}}$ admits a spectral density $\sigma^{2}:[0,2\pi)\to [0,+\infty)$ (Definition \ref{defspeden}),  there exists a set $I\subset [0,2\pi)$ of $\lambda-$measure one such that for every $\theta\in I$, $$A_{n}(\theta,\cdot)\to_{n\to \infty} 0, \mbox{\,\,\,\,\,\,\,$\mathbb{P}-$\it a.s. and in $L^{2}_{\mathbb{P}}$.}$$ 
and 
$$E[|\sqrt{n}A_{n}(\theta,\cdot)|^2]\to_{n} \sigma^{2}(\theta)
.$$  

\medskip

The next obvious task in this direction is to explore the validity of the central limit theorem for the normalized averages $(\sqrt{n}A_{n}(\theta,\cdot))_{n\geq 0}$. This is, we would  like to have an answer to (for instance) the following questions: if  $N(0,\Sigma)$ denotes a centered normal ($2-$dimensional) random variable with covariance matrix $\Sigma$, 

{\bf Question 1:} {\it Can we choose the set $I$ with the property that for every $\theta\in I$}
\begin{equation}
\label{equconnorave}
\sqrt{n}A_{n}(\theta,\cdot)\Rightarrow_{n\to \infty} N(0,\Sigma(\theta))?
\end{equation}
{\it and if so...}

{\bf Question 2:} {\it Can we actually prove the functional form of this convergence?} This is, {\it can we prove that $I$ can be chosen such that the  sequence of random functions $(B_{n}(\theta,\cdot))_{n\geq 0}$ defined on $[0,+\infty)$ by
\begin{equation}
\label{defbn}
B_{n}(\theta,\omega)(t)=\frac{S_{\left\lfloor nt \right\rfloor}(\theta)}{\sqrt{n}}(\omega)
\end{equation}
converge weakly (in a sense to be specified later) to a $2-$dimensional Brownian motion $B(\theta)$?}

Peligrad and Wu, in \cite{pelwu}, answered the first question positively for real-valued processes under the regularity condition (\ref{reg}). Their result is the following:

\medskip

\begin{thm}[CLT for Discrete Fourier Transforms]
\label{cltpelwu}
If $(\mathcal{F}_{k})_{k\in\mathbb{Z}}$ is a $T-$filtration and  $(X_{k})_{k\in\mathbb{Z}}=(T^{k}X_{0})_{k\in\mathbb{Z}}$ is a stationary real-valued, $(\mathcal{F}_{k})_{k\in\mathbb{Z}}-$adapted, and regular process (Definitions \ref{staproL2}, \ref{defadafil} and \ref{regcondef}) with spectral density $\theta\mapsto \sigma^{2}(\theta)$ (Definition \ref{defspeden}), there exists $I\subset [0,2\pi)$ with $\lambda(I)=1$ such that for every $\theta\in I$, (\ref{equconnorave}) holds, where
\begin{equation}
\label{covmatasynor}
\Sigma(\theta)=\left[\begin{array}{cc}
\sigma^{2}(\theta)/2 & 0\\
0& \sigma^{2}(\theta)/2
\end{array}\right].
 \end{equation}
\end{thm}

{\upshape Equivalently}, for every $\theta\in I$, 
$$\sqrt{n}A_{n}(\theta,\cdot)\Rightarrow_{n} (\sigma^{2}(\theta)/2)^{1/2}(N_{1}+iN_{2}),$$
where $N_{1},N_{2}$ are standard (real valued) independent normal random variables.

With regards to the second question, Peligrad and Wu were able to provide an ``intermediate'' answer by considering the random functions $B(\theta,\cdot)$ as random elements defined on the probability space $$([0,2\pi)\times\Omega,\mathcal{B}\otimes\mathcal{F},\lambda\times\mathbb{P})$$
To present their result, recall first that a two dimensional or {\it complex} standard Brownian motion is a random function of the form $B_{1}+iB_{2}$ where $B_{1},B_{2}$ are independent standard Brownian motions. Peligrad and Wu's invariance principle is the following.

\medskip

\begin{thm}[FCLT (averaged frequency) for Discrete Fourier Transforms]
\label{pelwufclt}
Consider $D([0,1],\mathbb{C})$, the space of cadlag  complex-valued functions $[0,1)\to\mathbb{C}$ endowed with the Skorohod topology (see Section \ref{concomvalcadfun}) and, in the context of Theorem \ref{cltpelwu}, let $B=B_{1}+iB_{2}$ be a  two-dimensional standard Brownian motion defined on a probability space $(\Omega',\mathcal{F}',\mathbb{P}')$. 
Then the random functions $$V_{n}:([0,2\pi)\times\Omega,\mathcal{B}\otimes\mathcal{F},\lambda\times\mathbb{P})\to D([0,1],\mathbb{C})$$
specified by (\ref{defbn}) converge weakly (i.e, in distribution) to the random function $([0,2\pi)\times \Omega',\mathcal{B}\otimes\mathcal{F}',\lambda\times\mathbb{P}')\to D([0,1],\mathbb{C})$ defined by
$$B(\theta,\omega'):=\frac{\sigma(\theta)}{\sqrt{2}}B(\omega').$$
\end{thm}

The main purpose of this work is to extend these theorems to similar results about {\it quenched} convergence, a notion that will be explained and explored in Chapter \ref{quecon}. It is important to mention that the results that we will obtain can be considered as generalized versions -under the respective hypotheses-  of Theorems \ref{cltpelwu} and \ref{pelwufclt}, in the sense that quenched convergence is, as we shall see, a notion of convergence stronger than convergence in distribution (we will also be able to deduce the convergence in Theorem \ref{cltpelwu} from our quenched central limit theorem for Fourier transforms, Theorem \ref{quecltfoutra}). We will also see that our quenched results are strictly stronger than the ones given here: while it is not difficult to show in a general sense that quenched convergence is strictly stronger than convergence in distribution (see Example \ref{exaqueconvscondis} in page \pageref{exaqueconvscondis}), for the processes under consideration a more careful analysis will be needed (see Theorem \ref{nonqueconthe} and the comments following it). 

\chapter{Convergence in Distribution}
\label{chacondis}

In this chapter we discuss several topics related to the notion of convergence in distribution of random elements in a metric space $S$, specialized to the cases that we will address in subsequent chapters. This will be important in order to understand the methods involved the proofs of the results presented in Chapter \ref{resandcom}.

First we will prove, in Section \ref{refporthe}, a further equivalence to the Portmanteau theorem (Theorem \ref{refporthesta}), valid in the case in which the underlying metric space $S$ is separable. 
  The important point of this reduction is that the family of functions to be tested for verifying convergence in distribution can be reduced in this case to a countable one, a fact that we will use to prove that our notion of quenched convergence (to be given in Chapter \ref{quecon}) specializes ``in the right way'' to the case under our consideration (this is, to the setting of regular conditional expectations). 

Then, in Section \ref{concomvalcadfun},  we will discuss briefly the notion of convergence in $D[[0,\infty),\mathbb{C}]$, in order to set the ground for the proofs of the forthcoming results involving the convergence in distribution of complex valued {\it cadlag} functions. The discussion will give rise to Theorem \ref{theconindoinfty}, which will be the key to proceed when addressing the proofs of the invariance principles in Chapter \ref{resandcom}.

We will then present, in Section \ref{contyp}, some results about convergence of types, which will be used on our discussions relating the possible ``quenched'' and ``annealed'' limits of a stochastic process, a discussion that will be essential to prove that the annealed limit theorems of Fourier transforms inspiring our results cannot themselves be extended to quenched ones: a random normalization is necessary (see Theorem \ref{nonqueconthe}, and the discussion preceding its statement, in page \pageref{nonqueconthe}).

In the section devoted to ``random elements and product spaces'' (Section \ref{raneleprospa}) we will address the relationship within the convergence in distribution of a sequence of random elements $(Z_{k})_{k}$ depending on two (random) parameters $(\theta,\omega)$ for a.e. fixed $\theta$ and the convergence in distribution of this sequence on the product space of the domain of the parameters. This discussion will serve later to clarify the hierarchy between the invariance principles under our consideration.

Finally, in Section \ref{trathe}, we will present a result (Theorem \ref{limlimlem}) used along our arguments in order to transport the asymptotic distributions of the processes under our consideration from suitable martingale approximations. 


\section{A Refinement of the Portmanteau Theorem}
\label{refporthe}

Throughout this section $(S,\mathcal{T})$ will denote a topological space with topology $\mathcal{T}$. If $S$ is a metric space, we will use the notation $(S,d)$, where $d:S\times S\to [0,+\infty)$ is the corresponding metric.  

To begin with, remember the notion of a {\it Urysohn function}.

\medskip

\begin{dfn}[Urysohn Function]
\label{uryfun} Given two closed, disjoint sets $F_{0},F_{1}$ in a perfectly normal
 topological space (for instance, any metric space) $(S,\mathcal{T})$, a function 
\begin{equation}
\label{equuryfun}
U(F_{0},F_{1}):S\to [0,1]
\end{equation}
is called a {\bfseries Urysohn function} if it is continuous, $U^{-1}\{0\}=F_{0}$ and $U^{-1}\{1\}=F_{1}$.
\end{dfn}

\medskip

\begin{remark}
\label{rempernorspa}
The existence of Urysohn functions for every two disjoint closed sets is the axiom characterizing perfectly normal spaces, and it is a standard fact from general topology that metrizable spaces are perfectly normal. We also recall the following: {\it for a perfectly normal space every closed set $F$ is a $G_{\delta}-$set}: there exists a countable family $\{G_{k}\}_{k\in\mathbb{N}}$ of open sets such that $F=\cap_{k\in\mathbb{N}}G_{k}$.
\end{remark}

The following definition is introduced for technical purposes.

\medskip

\begin{dfn}[Co-base]
\label{defcobas}
Given a topological space $(S,\mathcal{T})$, let us call a collection $\{F_{j}\}_{j\in J}$ of closed subsets of $S$ a {\bfseries co-base} if $\{S\setminus F_{j}\}_{j\in J}$ is a base of $\mathcal{T}$.
\end{dfn}

Note that if $(S,d)$ is separable it admits a co-base that is also a 
$\pi-$system (consider the finite intersections on any co-base). In the sequel, $\mathbf{C}^{b}(S)$ {\bfseries denotes the space of continuous and bounded functions $f:S\to \mathbb{R}$}. If needed, we will consider it also as a metric space via the {\bfseries \itshape uniform norm}
\begin{equation}
\label{uninor}
||f||_{_{\infty}}:=\sup_{s\in S}|f(s)|,
\end{equation} 
for every $f\in \mathbf{C}^{b}(S)$.

\medskip

\begin{thm}[A Refinement of the Portmanteau Theorem]
\label{refporthesta}
Let $S$ be a separable metric space,  let $\{F_{n}\}_{n\in\mathbb{N}}$ be a co-base of $S$ which is also a $\pi-$system, and let $X_{n}$, $X$ ($n\in \mathbb{N}$) be random elements of $S$ (Definition \ref{raneledef})\footnote{Note that the $X_{n}$'s are not necessarily defined on the same probability space.}. Then the following two statements are equivalent
\begin{enumerate}
\item For every $f\in \mathbf{C}^{b}(S)$, 
$$\lim_{n} Ef(X_{n})=Ef(X).$$
\item For every $k\in \mathbb{N}$, every rational $\epsilon>0$, and some Urysohn function   $U_{k,\epsilon}=U(S\setminus F_{k}^{\epsilon},F_{k})$
$$\lim_{n}EU_{k,\epsilon}(X_{n})=EU_{k,\epsilon}(X),$$
where $F_{k}^{\epsilon}$ is given according to Definition \ref{defdistoaset}.
\end{enumerate}

\end{thm}

\medskip

As stated before, the importance of this theorem for our purposes resides in the fact that it allows us to reduce the family of test functions in Portmanteau's Theorem to a countable one, a fact that will be exploited in the proof of Proposition \ref{prounicon} in Chapter \ref{quecon}. 

{\bf Proof of Theorem \ref{refporthesta}:} Denote by $P_{n}$ the law of $X_{n}$ and by $P$ the law of $X$ (see Definition \ref{raneledef}). Since 1. clearly implies 2. it suffices to see, by the Portmanteau Theorem (\cite{bilconpromea}, Theorem 2.1), that if 2. is true then for any given closed set $F$
$$\limsup_{n}{P}_{n}F\leq PF.$$
If for some $k$, $F=F_{k}$, this is a consequence of the inequality
$$I_{F}\leq U_{k,\epsilon}\leq I_{F^{\epsilon}},$$
the hypothesis in 2. and the continuity from above of finite measures. 

If $F$ is an arbitrary closed set, say $F=\cap_{j\in J}F_{j}$ for some $J\subset \mathbb{N}$, and if we define for all $k\in \mathbb{N}$, $J_{k}:=J\cap [0,k]$ and $A_{k}:=\cap_{j\in J_{k}} F_{j}$ then, since $A_{k}\in\{F_{n}\}_{n}$,
$$\limsup_{n} P_{n}F\leq \limsup_{n} P_{n} A_{k} \leq PA_{k} $$ 
for all $k$. By letting $k\to \infty$ we get the desired conclusion.\qed

\medskip

\begin{remark}
\label{remporthepernorspa}
We remark that the Portmanteau theorem can be extended to the context of random elements in abstract perfectly normal spaces (with their Borel sigma algebra) if one interprets ``convergence in distribution'' as the fulfillment of the hypothesis 1. of Theorem \ref{refporthesta}. This can be seen by following the arguments in \cite{bilconpromea} and using the fact that every closed set is a $G_{\delta}$ set (Remark \ref{rempernorspa}). In this context, Theorem \ref{refporthesta} corresponds to the second-countable case.
\end{remark}

\medskip

\section{Convergence of Complex-valued {Cadlag} Functions}
\label{concomvalcadfun}

This monograph contains results about convergence in distribution (under several measures) of random elements of $D[[0,\infty),\mathbb{C}]$: the space of functions $f:[0,\infty)\to \mathbb{C}$ that are continuous from the right and have left-hand limits at every point ({\itshape\bfseries cadlag} functions in $\mathbb{C}$). This space is an algebra with the operation of multiplication and addition given by the usual pointwise operations between complex functions, and it is a ($\mathbb{C}$ or $\mathbb{R}$-)vector space with the usual operation of multiplication by constants regarded as constant functions. 

To clarify the notions behind our results about convergence in $D[[0,\infty),\mathbb{C}]$, let us start in the following way: first, {\bfseries denote by  $(D[[0,\infty)],d)$ the space of {\itshape real-valued} cadlag functions} endowed with the Skorohod distance $d$ defined in  \cite{bilconpromea}, (16.4), which we proceed to describe now for the sake of completeness. 

\subsection*{Definition of the Skorohod Distance ({\itshape real-valued} case)}
Fix $m\in \mathbb{N}^{*}$, and consider the following definitions
\begin{enumerate}
\label{defskodis}
\item  First, define the family 
\begin{equation}
\Lambda_{m}:=\{\varphi:[0,m]\to [0,m]: \mbox{$\varphi$ is surjective, nondecreasing, and $||\varphi||_{m}<\infty$}\},
\end{equation}
where 
$$||\varphi||_{m}=\sup_{0\leq s<t\leq m}\left|\log\frac{\phi(t)-\phi(s)}{t-s} \right|.$$
Note in particular that for every $\varphi\in\Lambda_{m}$, $\varphi(0)=0$, $\varphi(m)=m$, and $\varphi$ is continuous.
\item Now consider the {\it Skorohod distance} $d_{m}$ in the space $D([0,m])$ of  real-valued {\it cadlag} functions with domain $[0,m]$: for every $w_{1},w_{2}\in D([0,m])$
\begin{equation}
\label{skodisdm}
d_{m}(w_{1},w_{2})=\inf_{\phi\in{\Lambda_{m}}}\{||\phi||_{m}\vee ||w_{1}-w_{2}\circ\phi||\}
\end{equation}
where $||\cdot||$ denotes the corresponding uniform norm in $D([0,m])$:
$$||w||=\sup_{0\leq t\leq m}|w(t)|.$$
\item Finally denote by $r_{m}$ the {\it restriction operator} $r_{m}: D([0,\infty))\to D([0,m])$
given by 
\begin{equation}
\label{resope}
(r_{m}w)(t)=w(t),
\end{equation}
define  $g_{m}:[0,\infty)\to [0,1]$ by
\begin{equation}
\label{defgm}
g_{m}(t)=\left\{
\begin{array}{ll}
1 & ,0\leq t\leq m-1 \\
m-t &, m-1< t\leq m \\
0 &, m<t
\end{array}\right.
\end{equation}
and denote, for every $w\in D([0,\infty))$
\begin{equation}
\label{defsmores}
w^{k}:=r_{k}(g_{k}w).
\end{equation} 
\end{enumerate}
With these notations we define $d$ as follows: given $w_{1},w_{2}\in D([0,\infty))$
\begin{equation}
\label{defskodis}
d(w_{1},w_{2}):=\sum_{k\geq 1}2^{-k}(1\wedge d_{k}(w_{1}^{k},w_{2}^{k})).
\end{equation}
  
\subsection{The topology of $D[[0,\infty),\mathbb{C}]$}

Let $S=D[[0,\infty),\mathbb{C}]$. The bijection $D[[0,\infty),\mathbb{C}]\to D[0,\infty)\times D[0,\infty)$ given by 
$$w=Re(w)+iIm(w)\mapsto (Re(w), Im(w))$$
allows us to regard $S$ as a topological space whose topology is the topology generated by the product Skorohod topology of $(D[[0,\infty)],d)$, this is, by the product of the topologies induced by  the metric (\ref{defskodis}). 
This topology is metrizable: it is induced by   the {\it product Skorohod metric} denoted (also) by
$d:S\times S \to S$ and given by
\begin{equation}
\label{defskodiscom}
d(w_{1},w_{2}):=((d(Re(w_{1}),Re(w_{2})))^2+(d(Im(w_{1}),Im(w_{2})))^2)^{1/2}
\end{equation}
where ``$d$'', at the right-hand side, is given by (\ref{defskodis}).

\medskip

\begin{dfn}[The space $\mbox{$D[[0,\infty),\mathbb{C}]$}$]
\label{skospacomvalcadfun}
The Skorohod distance in $S=D[[0,\infty),\mathbb{C}]$ is the distance $d$ defined by (\ref{defskodiscom}). $(S,d)$ is the space of {\bfseries cadlag complex-valued functions on $[0,\infty)$}.
\end{dfn}

\medskip

\begin{remark}[A Criterion for Measurability]
\label{remcrimead0infc}
Let $\mathcal{D}_{\infty,\mathbb{C}}$ be the Borel sigma algebra on $S=D[[0,\infty),\mathbb{C}]$, let $(\Omega,\mathcal{F})$ be a measurable space, and let $X:\Omega\to S$ be a given function.

By our definition  of the topology of  $S$, to prove that $X$ is $\mathcal{F}/\mathcal{D}_{\infty,\mathbb{C}}$ measurable it suffices to see the measurability of the real and imaginary parts of $X$. This observation, combined with the argument in \cite{bilconpromea}, p.84, and with Theorem 16.6 in that book shows that $X$ is $\mathcal{F}/\mathcal{D}_{\infty,\mathbb{C}}-$measurable if and only there exists a dense set $T\subset [0,\infty)$ such that for every $t\in T$, $\omega\mapsto X(\omega)(t)$ is $\mathcal{F}-$measurable. 

\end{remark}

Finally let us point out that, since separability and completeness ascend to the product space (with the product metric), and $(D[[0,\infty)],d)$ is separable and complete (\cite{bilconpromea}, Theorem 16.3), we have the following proposition.

\medskip

\begin{prop}
\label{comdoinfty}
The space of cadlag complex-valued functions on $[0,\infty)$ (Definition \ref{skospacomvalcadfun}) is separable and complete.
 \end{prop}

\subsection{Convergence on $D[[0,\infty),\mathbb{C}]$}
\label{condoinfc}
The space of cadlag complex-valued function admits, as any other metric space, a notion of convergence in distribution. To prove that an actual sequence of random elements in this space converges in distribution we will use the theoretical framework explained in \cite{bilconpromea} for convergence of real-valued cadlag functions, whose arguments can be transported to the case of complex-valued functions without major difficulties. More precisely, we will prove convergence in $D[[0,\infty),\mathbb{C}]$ via the following facts:

\begin{enumerate}
\item {\it Generic Idea.} Let us start by recalling the generic idea: remember that a sequence $(P_{n})_{n\in \mathbb{N}}$ of probability measures  on a metric space  is {\it tight} if for any $\epsilon>0$ there exists a compact set $K$ such that for every $n\in \mathbb{N}$, $P_{n}(K)>1-\epsilon$, and that when the space is separable and complete, tightness is equivalent to the {\it relative compactness} of $(P_{n})_{n\in\mathbb{N}}$ (see \cite{bilconpromea}, Theorems 5.1 and 5.2): $(P_{n})_{n\in\mathbb{N}}$ is tight if and only if for every (strictly increasing) sequence $(n_{k})_{k\in\mathbb{N}}$ of natural numbers there exists a subsequence $(n_{k'})_{k'}$ and a probability measure $P'$ with $P_{n_{k'}}\Rightarrow P'$. It follows that a tight sequence is convergent if $P'$ is independent of the given (sub)sequences. Since tightness is a necessary condition for convergence of measures, a way of addressing proofs of convergence in distribution is to give criteria for tightness and criteria to identify asymptotic distributions so that, in practice, one proves that a sequence of probability measures is convergent by proving that it is tight and that there exists a unique subsequential distribution via these criteria. 
\item {\it Criteria for Tightness.} Now, if $W:(\Omega,\mathcal{F},\mathbb{P})\to D[[0,\infty),\mathbb{C}]$ is a random element of $D[[0,\infty),\mathbb{C}]$, then the inequalities
$$\mathbb{P}[W\notin K_{1}\times K_{2}]\leq\mathbb{P}[Re(W)\notin K_{1}]+\mathbb{P}[Im(W)\notin K_{2}]\leq 2\mathbb{P}[W\notin K_{1}\times K_{2}]$$
show that, given a sequence $(W_{n})_{n\in\mathbb{N}}$ of random elements in $D[[0,\infty),\mathbb{C}]$, $(Re(W_{n}))_{n}$ and $(Im(W_{n}))_{n}$ are tight if and only if $(W_{n})_{n}$ is tight. Of course, this argument shows (the well known fact) that a sequence of random elements in the product of two metric spaces is tight if and only if the component  sequences are tight. 

The important observation is that we can prove tightness in $D[[0,\infty),\mathbb{C}]$ by applying criteria for tightness in $D[[0,\infty)]$ to the real and imaginary parts of any given random sequence of cadlag complex-valued functions. 

\item {\it Asymptotic Distributions.} By an adaptation of the arguments in \cite{bilconpromea68} and \cite{bilconpromea}, it is possible to show that the finite dimensional distributions are a {\it separating class} in $D[[0,\infty),\mathbb{C}]$: if for every $t_{1}\leq\dots\leq t_{n}$ we denote by $\pi_{t_{1}\dots t_{k}}:D[[0,\infty)]\to\mathbb{C}^{k}$ the projection 
\begin{equation}
\label{proope}
\pi_{t_{1}\dots t_{k}}(w):=(w(t_{1}),\dots,w(t_{k})),
\end{equation} 
then two measures ${P}_{1}$ and ${P}_{2}$ in $D[[0,\infty),\mathbb{C}]$ coincide if and only if there exists a dense subset $T\subset [0,\infty)$ such that for every $0\leq t_{1}\leq\dots\leq t_{n}$ in $T$ the  $n$-th dimensional distributions ${P}_{j}\pi_{t_{1}\dots t_{k}}^{-1}$ ($j=1,2$) on $\mathbb{C}^{n}=\mathbb{R}^{2n}$ are the same. 

\item\label{resdom}{\it Restriction of the domains.} Let $P$ be a probability measure in $D[[0,\infty),\mathbb{C}]$ and, for $m> 0$, consider the space $D([0,m],\mathbb{C})$ of cadlag complex-valued functions on $[0,m]$ with the Skorohod distance defined by identifying $D([0,m],\mathbb{C})=D[[0,m]]\times D[[0,m]]$ and extending $d_{m}$ (see (\ref{skodisdm})) to the product space as above. If $m$ is such that 
\begin{equation}
\label{conconpoi}
{P}\{w:\lim_{t\to m^{-}}w(t)\neq w(m)\}=0
\end{equation}
and $r_{m}: D[[0,\infty),\mathbb{C}]\to D[[0,m],\mathbb{C}]$ is the restriction operator defined above ($(r_{m}w)(t)=w(t)$), then the hypothesis ${P}_{n}\Rightarrow {P}$ in $D[[0,\infty),\mathbb{C}]$ implies that $\mathbb{P}_{n}r_{m}^{-1}\Rightarrow \mathbb{P}r_{m}^{-1}$ in $D[[0,m],\mathbb{C}]$, and the following ``converse'' holds: if $(m_{k})_{k}$ is a sequence increasing to infinity such that (\ref{conconpoi}) holds for all $m=m_{k}$, then the hypothesis $${P}_{n}r_{m_{k}}^{-1}\Rightarrow {P}r_{m_{k}}^{-1} \mbox{\,\,\,\,\,\,\, \it for every $k\in\mathbb{N}$ }$$
(on $D[[0,m_{k}],\mathbb{C}]$)  implies that ${P}_{n}\Rightarrow{P}$. This observation allows us to prove convergence in $D[[0,\infty),\mathbb{C}]$ by restricting our attention to $D[[0,m],\mathbb{C}]$.

\end{enumerate}
Let us finish this section by giving a more concrete criterion for convergence in $D[[0,\infty),\mathbb{C}]$. The proof will be just briefly sketched using the facts recalled here and referring to additional arguments from \cite{bilconpromea68} and \cite{bilconpromea}.

\medskip

\begin{thm}[Criterion of Convergence]
\label{theconindoinfty}
Let $P_{n},P$  be probability measures in $D[[0,\infty),\mathbb{C}\,]$ ($n\in\mathbb{N}$) and consider, for every $t>0$, the set
$$J_{t}=\{w\in D[[0,\infty),\mathbb{C}\,]: \lim_{s\to t-}w(s)\neq w(t)\}.$$
Then the set $A_{P}$ of nonnegative numbers $t$ such that $P{J_{t}}=0$ has a countable complement (in $[0,\infty)$), and if $(P_{n})_{n\in\mathbb{N}}$ is tight, ${P}_{n}\Rightarrow P$ if and only if for every $t_{1}\leq \cdots\leq t_{k}$ in $A_{P}$
 \begin{equation}
 \label{confindimdis}
 P_{n}\pi_{t_{1}\cdots t_{k}}^{-1}\Rightarrow_{n} P\pi_{t_{1}\cdots t_{k}}^{-1}.
 \end{equation}
  \end{thm}

\medskip

{\bf Proof (sketch):} To see that $A_{P}$ has a countable complement show, using the argument in \cite{bilconpromea68}, p.124 that for every given $m\geq 0$, $[0,m]\setminus T_{P}$ is countable.

To prove the second statement assume that $(P_{n})_{n\in\mathbb{N}}$ is tight, and start by considering a strictly increasing sequence $(s_{k})_{k\in\mathbb{N}}$ of elements in $A_{P}$ with $\lim_{k}s_{k}=\infty$.  

Given $m>0$,  denote by $r_{m}$ the restriction operator given by (\ref{resope}). By 4. above, $P_{n}\Rightarrow P$ if and only if given $m\in\{s_{k}\}_{k\in \mathbb{N}}$, $P_{n}r_{m}^{-1}\Rightarrow Pr_{m}^{-1}$ as $n\to \infty$.

Now, by an adaptation of Theorem 15.1 in \cite{bilconpromea68}: $P_{n}r_{m}^{-1}\Rightarrow Pr_{m}^{-1}$ if and only if  $(P_{n}r_{m}^{-1})_{n\in\mathbb{N}}$ is tight and for every $0\leq t_{1}\leq\cdots\leq t_{k}\leq m$
\begin{equation}
\label{confindimdism}
P_{n}(\pi_{t_{1}\dots t_{k}}^{m}r_{m})^{-1}\Rightarrow P(\pi_{t_{1}\dots t_{k}}^{m}r_{m})^{-1},
\end{equation}
where $\pi_{t_{1}\dots t_{k}}^{m}:D([0,m],\mathbb{C})\to\mathbb{C}^{k}$ denotes the projection specified by (\ref{proope}) (the superindex ``$m$'' is introduced to indicate the domain).

The tightness of $(P_{n}r_{m}^{-1})_{n\in\mathbb{N}}$ follows by an adaptation of the argument at the beginning of the proof of  Theorem 16.7  in \cite{bilconpromea} (in short: $r_{m}$ is continuous because $m\in A_{P}$, and since $(P_{n})_{n\in\mathbb{N}}$ is tight, $(P_{n}r_{m}^{-1})_{n\in\mathbb{N}}$ is tight by the mapping theorem), and therefore it suffices to prove, by Theorem 15.1 in \cite{bilconpromea68} again, that 
the convergence (\ref{confindimdis}) for every $0\leq t_{1}\leq\cdots\leq t_{k}$ is equivalent to the convergence (\ref{confindimdism}) for every $m\in A_{P}$ and every $0\leq t_{1}\leq\cdots\leq t_{k}\leq m$, but this is  just a consequence of the equality 
$$\pi_{t_{1}\dots t_{k}}^{m}r_{m}=\pi_{t_{1}\dots t_{k}}$$
(where $m\geq t_{k}$) and the fact that $\lim_{n}s_{n}=\infty$.\qed 


\section{Convergence of Types}
\label{contyp}

In this section we present some facts about Convergence of Types in a form that is convenient for the proofs given along this monograph. Let us start by recalling the notion of a {\it non-degenerate} distribution function.

\medskip

\begin{dfn}[Nondegenerate Distribution Function]
\label{nondegdisfun}
A distribution function $F$ is {\itshape\bfseries non-degenerate} if it is not the indicator function of some interval $[a,+\infty)$. This is, if it does not correspond to a constant random variable.
\end{dfn}

Our arguments in this section will be mostly based on the following {\it Convergence of Types} theorem (\cite{bilpromea}, Theorem 14.2). In accordance with the notation introduced at the beginning, if $F_{n}$ and $F$ are (probability) distribution functions, ``$F_{n}\Rightarrow_{n} F$'' will denote pointwise convergence at the continuity points of $F$ (convergence of distribution functions).

\medskip

\begin{lemma}[Convergence of Types]
\label{bilcontyp}
Let $F_{n}$, $F$ and $G$ be distribution functions,  and let $a_{n},u_{n},b_{n}, v_{n}$ be constants with $a_{n}>0,\,u_{n}>0$. If $F$, $G$ are non-degenerate, $F_{n}(a_{n}x+b_{n})\Rightarrow F(x)$, and $F_{n}(u_{n}x+v_{n})\Rightarrow G(x)$ then there exist $a=\lim_{n}a_{n}/u_{n}$, $b=\lim_{n}(b_{n}-v_{n})/u_{n}$ and $G(x)=F(ax+b)$.  
\end{lemma}

Note that, necessarily, $a>0$ (as otherwise $G$ would be constant).

We will translate this statement to a statement about convergence of stochastic processes (with a restricted choice of $u_{n},v_{n}$ see Proposition \ref{quecon} below), which we will be able to extend to the complex valued case. 

\subsection{Preliminary Facts}
\label{prefac}

To begin with, we remind the following elementary facts, here capital letters ($U_{n}$, $V_{n}$, etc) denote (real-valued) random variables and ``$\to_{P}$'' denotes convergence in probability. For more details see for instance Section 25 in \cite{bilpromea}.
\begin{enumerate}

\item If $a$ is constant then $U_{n}\Rightarrow a$ if and only if $U_{n}\to_{P} a$. 

\item If $U_{n}\Rightarrow W$ and $V_{n}\to_{P} 0$  then $U_{n}+V_{n}\Rightarrow W$.

\item If $(a_{n})_{n}$ is a sequence of constant functions then $a_{n}\Rightarrow A$ if and only if $a=\lim_{n}a_{n}$ exists (and therefore $A=a$ a.s.).

\end{enumerate}

These facts will be used along the proof of the forthcoming results in this section, which will be useful when addressing the issues of the relationship between ``annealed'' convergence and ``quenched'' convergence of a sequence of random variables.

\subsection{Convergence of Types Results}

In this section all the processes under consideration will be assumed real-valued, unless otherwise specified. We will simply write ``$V_{n}\Rightarrow V$'' to indicate that the stochastic process $(V_{n})_{n\in\mathbb{N}}$ converges in distribution to $V$ as $n\to\infty$.

\medskip

\begin{lemma}
\label{lemcondeg}
If  $Y_{k}\Rightarrow Y$ and $\{c_{k}\}_{k}\subset \mathbb{R}$ are such that $Y_{k}+c_{k}\Rightarrow 0$, then $Y=-\lim_{k}c_{k}$. In particular, $Y$ is a constant function. 
\end{lemma}

{\bf Proof:} Note that $c_{k}=-Y_{k}+(Y_{k}+c_{k})\Rightarrow -Y$ because $Y_{k}+c_{k}\Rightarrow 0$. Now use  3. in Section \ref{prefac}. \qed

\medskip

\begin{cor}
\label{degcas}
If $X$ is not constant, $X_{n}\Rightarrow X$, and $a_{n}$, $b_{n}$ are such that $a_{n}X_{n}+b_{n}\Rightarrow 0$, then $a_{n}\to 0$ and $b_{n}\to 0$.
\end{cor}

{\bf Proof:} If $0<a:=\limsup_{n}a_{n}\leq \infty$ and $a_{n_{k}}\to_{k\to \infty} a$ with $a_{n_{k}}>0$, then applying Lemma \ref{lemcondeg} with $Y_{k}=X_{n_{k}}$ and $c_{k}=b_{n_{k}}/a_{n_{k}}$ we conclude that $X$ is constant. This proves that, necessarily, $\limsup_{n}a_{n}\leq 0$. A similar argument shows that $\liminf_{n}a_{n}\geq 0$, and therefore $\lim_{n}a_{n}=0$.

The fact that $b_{n}\to 0$ follows from here applying Lemma \ref{lemcondeg} again, because $a_{n}X_{n}\Rightarrow 0$ . \qed  

\medskip

These results give rise to the following proposition

\medskip

\begin{prop}
\label{proquecon}
If $X$ is not constant, $X_{n}\Rightarrow X$ and $a_{n}>0$, $b_{n}$ are such that $a_{n}X_{n}+b_{n}\Rightarrow Y$, then there exists $a=\lim_{n}a_{n}$, $b=\lim_{n}b_{n}$ and, therefore, $Y=aX+b$ (in distribution).
\end{prop}

{\bf Proof:} If $Y$ is constant then, from $a_{n}X_{n}+b_{n}-Y\Rightarrow 0$ (see 1. in Section \ref{prefac}) it follows, via Corollary \ref{degcas}, that $\lim_{n}a_{n}=0$ and $\lim_{n}b_{n}=Y$.

If $Y$ is not constant we apply Lemma \ref{bilcontyp} with $F_{n}$, $F$, and $G$ the distribution functions of $X_{n}$, $X$ and $Y$ respectively, and with $u_{n}=1$, $v_{n}=0$.\qed

\medskip
 
\begin{remark} Taking $X_{n}=1$ (the constant function), $a_{n}=n$, and $b_{n}=-n$, we see that the given restriction on $X$ (to be non constant) is necessary. The asymptotically degenerate case is nonetheless covered by the following proposition (note the additional restriction on the coefficient of $X_{k}$).

\end{remark}

\medskip

\begin{prop}
\label{degcascon}
If $X$ is constant, $X_{k}\Rightarrow X$, and $X_{k}+c_{k}\Rightarrow Z$ then $c=\lim_{k}c_{k}$ exists and therefore $Z=X+c$ (in distribution). 
\end{prop}

{\bf Proof:} Use $X+c_{k}=(X-X_{k})+X_{k}+c_{k}\Rightarrow Z$ by 2. in Section \ref{prefac}. The conclusion follows from the item 3. there.

These results can be transported to the case of complex-valued random elements. More concretely.
\medskip
  \begin{prop}[Convergence of Types for complex-valued Random Variables]
\label{procontypcom} 
  Proposition \ref{proquecon} and Proposition \ref{degcascon} remain valid if the processes  involved are complex-valued, provided that the constants $(a_{n})_{n}$ in Proposition \ref{proquecon} are still real and positive (all the other constants can be assumed complex).
 \end{prop} 
 
  {\bf Proof:} To see this for Proposition \ref{proquecon} notice that if $X_{n}, X$ are complex valued, $X_{n}\Rightarrow X$, and $\mathbf{u}\in \mathbb{C}=\mathbb{R}^{2}$ is any vector then, by the mapping theorem
  $$\mathbf{u}\cdot(a_{n}X_{n}+b_{n})=a_{n}(\mathbf{u}\cdot X_{n})+\mathbf{u}\cdot b_{n}\Rightarrow \mathbf{u}\cdot Y$$
  so that, by the real valued case just proved, there exists $a=\lim_{n}a_{n}$ and $b_{\mathbf{u}}=\lim_{n}\mathbf{u}\cdot b_{n}$. Since $\mathbf{u}$ is arbitrary, there actually exists $b=\lim_{n}b_{n}$. 
  
The second conclusion ($Y=aX+b$ in distribution) follows at once from the Cramer-Wold theorem. The argument for Proposition \ref{degcascon} is similar.\qed

\section{Random Elements and Product Spaces}
\label{raneleprospa}

Since we will be concerned with random cadlag functions seen as random elements depending on two random parameters $(\theta,\omega)$ or on a single parameter $\omega$ for $\theta$ fixed (see Chapter \ref{resandcom}), it is convenient to give now the following   proposition.

\medskip

\begin{prop}[Convergence for fixed parameters and on the Product Space]
\label{relconfixfreavefre}
Let $(S,d)$ be a metric space and let $(\Theta,\mathcal{B},\lambda)$ be a probability space.\footnote{This is just {\it some} probability space but, as the notation suggests, we will use only the case $\Theta=[0,2\pi)$ with the Borel sigma algebra and the normalized Lebesgue measure.} Assume that for every $\theta\in\Theta$ and every $n\in\mathbb{N}\cup\{\infty\}$, $V_{n}(\theta)$ is a random element (Definition \ref{raneledef}) in $S$ defined on a probability space $(\Omega_{n},\mathcal{G}_{n},\mathbb{P}_{n})$, and that the function $V_{n}$ given  by $(\theta,\omega)\mapsto V_{n}(\theta)(\omega)$ is measurable with respect to the product sigma-algebra $\mathcal{B}\otimes\mathcal{G}_{n}$. If for $\lambda-$a.e $\theta$  $V_{n}(\theta)\Rightarrow V_{\infty}(\theta)$ as $n\to\infty$, then the random elements $V_{n}:(\Theta\times\Omega_{n},\mathcal{B}\otimes\mathcal{G}_{n},\lambda\times\mathbb{P}_{n})\to S$ converge in distribution to the random element $V_{\infty}:(\Theta\times\Omega_{\infty},\mathcal{B}\otimes\mathcal{G}_{\infty},\lambda\times\mathbb{P}_{\infty})\to S$ as $n\to \infty$. 
\end{prop}

{\bf Proof:} Given any bounded and continuous function $f:S\to\mathbb{R}$, consider the function
$$g_{n}(\theta):=E(f(V_{n}(\theta))-Ef(V(\theta))$$
where (we emphasize again that) ``$E$'' denotes integration with respect to the respective probability measures ($\mathbb{P}_{n}$ and $\mathbb{P}_{\infty}$). By the hypotheses on $V_{n}(\theta)$ and Fubini's theorem, $g_{n}$ is $\mathcal{B}-$measurable, and since $V_{n}(\theta)\Rightarrow V_{\infty}(\theta)$ for $\lambda-$a.e $\theta$, $g_{n}(\theta)\to 0$ as $n\to\infty$, $\lambda-$a.s. It follows from the dominated convergence theorem that
$$\int_{\Theta}g_{n}(\theta)\, d\lambda(\theta)\to_{n} 0 $$
as $n\to\infty$. This is (Fubini's Theorem again), that 
$$\int_{\Theta\times\Omega}f\circ V_{n} \, d(\lambda\times\mathbb{P}_{n})\to_{n} \int_{\Theta\times\Omega}f\circ V_{\infty} \, d(\lambda\times\mathbb{P}_{\infty}),$$
which gives the desired conclusion.\qed

\medskip

The following example shows that the converse of Proposition \ref{relconfixfreavefre} does not hold.

\medskip

\begin{exam}
\label{proconnoascon}
Consider the probability space $([0,2\pi),\mathcal{B},\lambda)$ and let, for every $n\geq 0$, $f_{n}:[0,2\pi)\to [0,\infty)$ be a sequence of ($\mathcal{B}-$measurable) functions with the property that $f_{n}\to 0$ in $L^{1}_{\lambda}$ and for every $\theta\in [0,2\pi)$, $(f_{n}(\theta))_{n\in\mathbb{N}}$ is not convergent. For instance take $f_{0}=I_{[0,2\pi)}$, $f_{1}=I_{[0,\pi)}$, $f_{2}=I_{[\pi,2\pi)}$, $f_{3}=I_{[0,\pi/2)}$, $f_{4}=I_{[\pi/2,\pi)}$, and so on. 

Given any probability space $(\Omega,\mathcal{F},\mathbb{P})$ let $X$ be the constant function $X(\omega)=1$ and consider, for every $n\in\mathbb{N}$, the random variable
$V_{n}:([0,2\pi)\times\Omega,\mathcal{B}\otimes\mathcal{F},\lambda\times\mathbb{P})\to\mathbb{R}$ given by $V_{n}(\theta,\omega)=f_{n}(\theta)X(\omega)=f_{n}(\theta)$. Since $f_{n}\to 0$ in $L^{1}_{\lambda}$, $V_{n}\Rightarrow 0$, but note that since the law of $X$ is the Dirac measure $\delta_{1}$ ($\delta_{1}\{1\}=1$), the sequence of random variables $(V_{n}(\theta,\cdot))_{n\in\mathbb{N}}$ defined on $(\Omega,\mathcal{F},\mathbb{P})$ does not converge in distribution for any $\theta\in[0,2\pi)$ (the law of $V_{n}(\theta,\cdot)$ is $\delta_{f_{n}(\theta)}$). 
\end{exam}

\medskip

As this discussion shows, given a sequence $(V_{n})_{n\in\mathbb{N}}$ as in the statement of Proposition \ref{relconfixfreavefre},  the convergence in distribution of $V_{n}(\theta,\cdot)$ for $\lambda-$a.e $\theta$ is in general a notion {\it stronger} than that of the convergence in distribution of $V_{n}$. We will return to this discussion in Chapter \ref{resandcom}.
  
\section{A Transport Theorem}
\label{trathe}

The last result to be recalled in this chapter, Theorem \ref{limlimlem}, is an improvement due to Dehling, Durieu and Voln\'{y}, of Theorem 3.1 in \cite{bilconpromea} for the case in which the target (state) space is a complete and separable metric space.

\medskip

\begin{thm}[Transport Theorem]
\label{limlimlem}
{\it Let $(S,{d})$ be a complete and separable metric space. Assume that for all natural numbers $r,n$, $X_{r,n}$ and $X_{n}$ are random elements of $S$ defined on the same probability space $(\Omega,\mathcal{F},\mathbb{P})$, and that $X_{r,n}\Rightarrow_{n} Z_{r}$. Then the hypothesis
\begin{equation}
\label{liminf}
\lim_{r}\limsup_{n}\mathbb{P}[d(X_{r,n},X_{n})\geq \epsilon]=0 \mbox{\,\,\,\,\,\,\,\it for all $\epsilon>0$,}
\end{equation}
implies the existence of a random element $X$ of $S$ such that $Z_{r}\Rightarrow_{r} X$  and  $X_{n}\Rightarrow_{n}X$.}
\end{thm}

{\bf Proof:} This is Theorem 2 in \cite{deduvo}.\qed
\medskip

\begin{cor}
\label{corlimlimlem}
In the context of Theorem \ref{limlimlem} denote, for any given $q>0$, $$||Z||_{\mathbb{P},q}:=\left(\int_{\Omega}|Z|^{q}d\mathbb{P}(\omega)\right)^{1/q}.$$ 
If for some $q> 0$
$$\lim_{r}\limsup_{n}||d(X_{r,n},X_{n})||_{\mathbb{P},q}=0$$
and if for all (but finitely many) $r\in\mathbb{N}$ there exists a random element $Z_{r}$ with $X_{r,n}\Rightarrow_{n} Z_{r}$, then there exists a random element $X$ such that $X_{n}\Rightarrow_{n} X$ and $Z_{r}\Rightarrow_{r} X$.
\end{cor}

{\bf Proof:} Apply Markov's inequality to verify the hypothesis of Theorem \ref{limlimlem}. \qed

We will use these results to obtain the asymptotic distributions of the processes under our consideration from suitable martingale approximations.

\chapter[Quenched Convergence and Regular Conditional Expectations]{Quenched Convergence and Regular Conditional Expectations}
\chaptermark{Quenched Convergence...}
\label{quecon}

In this chapter we introduce the notions of {\it quenched convergence with respect to a sigma algebra} and {\it regular conditional expectation}. These notions will settle the formal ground for our discussions on asymptotic limit theorems ``started at a point''. 

Results on quenched convergence -in particular those related to the asymptotics of averages for dependent structures- are the object of intensive research at the moment of writing this monograph (see for instance \cite{baranexa}, \cite{barpel}, \cite{cunmer}, \cite{cunvol},  \cite{volwoo} and the references therein), but they have been in the literature for at least about forty years (see for instance Theorem 20.4 in \cite{bilconpromea68}). These results belong to the category of limit theorems for {nonstationary} processes: in loose terms, they refer to convergence in distribution of a process with respect to a family of random measures that ``integrate'' to a stationary distribution for the process in question.\footnote{See Definition \ref{regconexpdef} and Proposition \ref{proredregcondis} for the precise meaning of this statement.}

The presentation is organized as follows: in Section \ref{quecondefgenrem} we introduce the notions of {\it quenched convergence} of a stochastic process defined on a probability space $(\Omega,\mathcal{F},\mathbb{P})$ with respect to a sub sigma-algebra $\mathcal{F}_{0}$ of $\mathcal{F}$, and the notion of {\it regularity} for conditional expectations, which will be an important assumption along our forthcoming proofs, and we give some basic properties associated to this notions.

Then, in Section \ref{exareg}, we provide some examples of probability spaces and initial sigma algebras that admit a regular conditional expectation. To be more precise, we show  (Example \ref{examarproexa}) that this is the case for the setting of functions of stationary Markov Chains, which encompasses a broad family of the processes present in the applications.

We move then to quickly discuss, in Section \ref{regandquecon},  the relationship between regularity and quenched convergence. We prove there (Proposition \ref{proredregcondis}) that, in the case of a separable state space, quenched convergence with respect to a sigma algebra $\mathcal{F}_{0}$ admitting a regular conditional expectation is the same as convergence in distribution with respect to any family of probability measures decomposing $E[\,\cdot\,|\mathcal{F}_{0}]$, an assumption that is apparently implicit in the literature but whose proof is not present among the visible references.

Finally, we show that the notion of regular conditional expectation behaves well with respect to the product of probability spaces: the product of two sigma-algebras admitting  regular conditional expectations satisfies itself this property, and a decomposition of the expectation with respect to this product sigma-algebra is  given by the product of any two decompositions of the factor algebras (Proposition \ref{proregpro}).

\section{Definitions and General Remarks}
\label{quecondefgenrem}

In this section we introduce the notions of {\it quenched convergence} of a stochastic process and {\it regular conditional expectation} with respect to a sigma algebra.

\subsection{Quenched Convergence}

Let $(Y_{n})_{n\geq 1}$ be a measurable sequence on some metric space $(S,d)$. This is, for every $n$, 
\begin{equation}
\label{defyn}
Y_{n}:(\Omega,\mathcal{F})\to (S,\mathcal{S})
\end{equation}
is an $\mathcal{F}/\mathcal{S}$ measurable function where $(\Omega,\mathcal{F})$ is a (fixed) measure space and $\mathcal{S}$ is the Borel sigma algebra of $S$. Let $\mathbb{P}$ be a given probability measure on $(\Omega,\mathcal{F})$, so that $(\Omega,\mathcal{F},\mathbb{P})$ is a probability space, and denote by ``$\Rightarrow_{\mathbb{P}}$'' the convergence in distribution with respect to $\mathbb{P}$.

The Portmanteau theorem (\cite{bilconpromea}, Theorem 2.2) states, among other equivalences, that if $Y:(\Omega',\mathcal{F}',\mathbb{P}')\to (S,\mathcal{S})$ is a random element of $S$ (Definition \ref{raneledef}) then $Y_{n}\Rightarrow_{\mathbb{P}} Y$ if and only if for every bounded and continuous function $f:S\to \mathbb{R}$
\begin{equation}
\label{usupro}
\int_{\Omega}f\circ Y_{n}(\omega)\, d\mathbb{P}(\omega)\to_{n\to \infty}, \int_{\Omega'}f\circ Y (z) d\mathbb{P}'(z)
\end{equation}

or, in the usual probabilistic notation, if $\lim_{n\to \infty}Ef(Y_{n})=Ef(Y)$, where $E$ is the expectation (Lebesgue integral) with respect to the corresponding probability measures.


\medskip


A stronger kind of convergence, {\it quenched} convergence, can be defined in the following way: 

\medskip

\begin{dfn}[Quenched Convergence]
\label{defquecon}
Let $(\Omega,\mathcal{F},\mathbb{P})$ be a probability space, let $(Y_{n})_{n\in\mathbb{Z}}$ be a sequence of random elements on a metric space $(S,d)$ as in (\ref{defyn}), and let $Y$ be a random element of $(S,d)$ defined on some probability space $(\Omega',\mathcal{F}',\mathbb{P}')$. Fix a sub-sigma algebra $\mathcal{F}_{0}\subset \mathcal{F}$, and denote by $E_{0}$ the conditional expectation with respect to $\mathcal{F}_{0}$. We say that {\bfseries $Y_{n}$ converges to $Y$ in the  quenched sense with respect to $\mathcal{F}_{0}$} if for every bounded and continuous function $f:{S}\to\mathbb{R}$
$$E_{0}[f(Y_{n})]\to_{n} Ef(Y),\,\,\,\,\mbox{\it $\mathbb{P}-$a.s.}$$

\end{dfn}

\medskip

\begin{remark}
\label{remdefquecon}
{As indicated in Section \ref{heuint}, $\mathcal{F}_{0}$ will represent in the practice, in a heuristic language, the ``initial information'' about (or ``the  past'' of) the process $(Y_{n})_{n\in\mathbb{Z}}$.} In most of our discussions it will be clear from the context what $\mathcal{F}_{0}$ is, thus we will simply speak of {\it quenched convergence} when addressing quenched convergence with respect to $\mathcal{F}_{0}$. 
 \end{remark}
 
 Note also the following: since the convergence in Definition \ref{defquecon} is  pointwise convergence of uniformly bounded functions (to a constant value), the dominated convergence theorem guarantees that for every continuous and bounded function $f$
 $$\lim_{n}E[E_{0}f(Y_{n})]=\lim_{n}E[f(Y_{n})]=Ef(Y),$$ 
 thus, certainly, {\it quenched convergence implies convergence in distribution}. 

\medskip

\begin{exam}[Quenched Convergence vs Convergence in Distribution]
\label{exaqueconvscondis} An example showing that the notion of quenched convergence is strictly stronger than convergence in distribution can be constructed by starting from any sequence $(Y_{n})_{n}$ of  $\mathcal{F}_{0}-$measurable functions and noticing that quenched convergence of $Y_{n}$ to $Y$ in this case is the same as  
 $$f(Y_{n})\to Ef(Y),\mbox{\,\,\,\,\,\,\,\, $\mathbb{P}-$\it a.s}$$
 for all continuous and bounded functions $f$, which is not possible if, for instance, $(Y_{n})_{n}$ takes the values $1$ and $0$ infinitely often  $\mathbb{P}-$a.s
. Thus it suffices to consider a sequence $(Y_{n})_{n}$ of random variables that converges in distribution but gives $\mathbb{P}$-a.s a sequence with infinitely many  $0$'s and $1$'s, and then to define $\mathcal{F}_{0}:=\sigma(Y_{n})_{n}$: take for instance the functions $f_{n}$ in Example \ref{proconnoascon} or, for an even  simpler example, consider unit interval with Lebesgue measure as the underlying probability space and, for every $k\in \mathbb{N}$, define $Y_{2k}:=I_{[0,1/2]}$ and $Y_{2k+1}:=I_{(1/2,1]}$.
 \end{exam}

Now recall the following property of conditional expectations (for a proof see for instance Theorem 34.2 (v) in \cite{bilpromea}):
\medskip

\begin{thm}[Dominated Convergence Theorem for Conditional Expectation]
\label{domcontheconexp}
Let $(\Omega,\mathcal{F},\mathbb{P})$ be a probability space, and let $X_{n},X\in L^{1}_{\mathbb{P}}$ be real-valued random variables defined on $(\Omega,\mathcal{F},\mathbb{P})$. If $X_{n}\to X$ {$\mathbb{P}$-}a.s. and there exists $Y\in L^{1}_{\mathbb{P}}$ with $|X_{n}|\leq Y$ ($\mathbb{P}-$a.s) for all $n$, then $E[X_{n}|\mathcal{F}_{0}]\to E[X|\mathcal{F}_{0}]$, $\mathbb{P}-$a.s.
\end{thm}

 Applying this lemma to $X_{n}=f(Y_{n})$ and $X=f(Y)$ we get, in the context of Definition \ref{quecon}, the following property.

\medskip
 
 \begin{prop}[Quenched Convergence with respect to sub sigma-algebras]
 \label{inhquecon}
 If $Y_{n}$ converges to $Y$ in the quenched sense with respect to $\mathcal{F}_{0}$ and $\mathcal{G}_{0}\subset\mathcal{F}_{0}$, then $Y_{n}$ converges to $Y$ in the quenched sense with respect to $\mathcal{G}_{0}$.
 \end{prop}

Note that if $Y_{n}$ converges to $Y$ in the quenched sense, the convergence in distribution of $Y_{n}$ to $Y$ is a consequence of Proposition \ref{inhquecon} by considering $\mathcal{G}_{0}=\{\emptyset,\Omega\}$. Though Example \ref{exaqueconvscondis} shows how the notions of quenched convergence and convergence in distribution differ in general, we will address the problem of non-quenched convergence later,  in the specific context of our quenched results. Concretely, we will see that the processes for which the CLT is known to happen within our discussion do not admit a quenched version without a ``random centering'', corresponding to the usual normalization of the mean transported to the setting of conditional expectation.

\subsection{Regular Conditional Expectations}

We begin this section introducing the notion of {\it regular conditional expectation}, which will allow us to interpret the notion of quenched convergence as a notion of convergence in distribution with respect to a family of measures.

\medskip

\begin{dfn}[Regular Conditional Expectation]
\label{regconexpdef}
Let $(\Omega,\mathcal{F},\mathbb{P})$ be a probability space, let $\mathcal{F}_{0}\subset \mathcal{F}$ be a sub-sigma algebra of $\mathcal{F}$, and denote by $E_{0}$ the conditional expectation with respect to $\mathcal{F}_{0}$.
We say that $E_{0}$ is {\bfseries regular} if there exists a family of probability measures  $\{\mathbb{P}_{\omega}\}_{\omega\in \Omega}$  such that for every integrable $X:(\Omega,\mathcal{F},\mathbb{P})\to \mathbb{R}$, the function defined by
\begin{equation}
\label{defe0regequ}
\omega\mapsto \int_{\Omega}X(z)d\mathbb{P}_{\omega}(z)
\end{equation}
if the integral makes sense\footnote{This will happen over an $\mathcal{F}_{0}-$set of $\mathbb{P}-$measure one, see Remark \ref{remregconexpdef} below.}, and zero otherwise, defines an $\mathcal{F}_{0}-$measurable version of $E_{0}X$. In this case we call $\{\mathbb{P}_{\omega}\}_{\omega\in \Omega}$ a {\bfseries decomposition of $E_{0}$}. 
\end{dfn}

\medskip

\begin{remark}
\label{remregconexpdef}
Note that, in Definition \ref{regconexpdef}, $X$ is an actual $\mathbb{P}-$integrable function, {\it not} a $\mathbb{P}-$equivalence class of functions. 

Note also that if, in the context of Definition \ref{regconexpdef},  $X$ is a bounded function, then the integral in (\ref{defe0regequ}) is well defined for every $\omega\in \Omega$.

Let now $X\in L^{1}_{\mathbb{P}}$ be given, and fix a version (also denoted by) $X$ of this function. If $(X_{n})_{n\in\mathbb{N}}$ is a family of nonnegative simple functions with $X_{n}(\omega)$ increasing to $|X(\omega)|$ for all $\omega\in \Omega$ (see for instance p.254 in \cite{bilpromea}), then the monotone convergence theorem gives that, for every $\omega\in\Omega$
$$\int_{\Omega}X_{n}(z)\,d\mathbb{P}_{\omega}(z)\to_{n}\int_{\Omega}|X(z)|\,d\mathbb{P}_{\omega}(z),$$
where the right hand side is regarded as $\infty$ if $X\notin L^{1}_{\mathbb{P}_{\omega}}$.

For every $n\in\mathbb{N}$, denote by $\tilde{X}_{n}$ the function specified by (\ref{defe0regequ}) with $X_{n}$ in place of $X$. By the definition of $\{\mathbb{P}_{\omega}\}_{\omega\in\Omega}$, $\tilde{X}_{n}$ is an $\mathcal{F}_{0}-$measurable version of $E_{0}X_{n}$, and we have just seen that 
\begin{equation}
\label{coninttoint}
\tilde{X}_{n}(\omega)\to_{n} \int_{\Omega}|X(z)|\,d\mathbb{P}_{\omega}(z)
\end{equation}
for {\it every} $\omega\in \Omega$.

Now, by Theorem \ref{domcontheconexp} and the fact that $\tilde{X}_{n}$ is a version of $E_{0}X_{n}$, 
$$\tilde{X}_{n}\to_{n}E_{0}|X|,$$ 
$\mathbb{P}-$a.s. This, together with (\ref{coninttoint}), implies that
\begin{equation}
\label{intalmthe}
E_{0}[|X|](\omega)=\int_{\Omega}|X(z)|\,d\mathbb{P}_{\omega}(z),
\end{equation}
$\mathbb{P}-$a.s. In particular,
\begin{equation}
\label{intfinas}
\int_{\Omega}|X(z)|\,d\mathbb{P}_{\omega}(z)<\infty
\end{equation}
for $\mathbb{P}-$a.e. $\omega$, provided that the set $\Omega_{X}$ of $\omega$'s where (\ref{intfinas}) holds is $\mathcal{F}-$measurable. But it turns out that $\Omega_{X}$ is {\it indeed} $\mathcal{F}_{0}-$measurable: simply note that, following with the notation along this remark
$$\Omega_{X}:=\bigcup_{n\in\mathbb{N}}\bigcap_{k\in\mathbb{N}}[\tilde{X}_{k}\leq n].$$
Even more is true: denoting by $Y^{+}=YI_{[Y\geq 0]}$ and $Y^{-}=-YI_{[Y<0]}$  the nonnegative and negative parts of $Y$, it is easy to see that if $(X_{n}^{+})_{n\in\mathbb{N}}$ and $(X_{n}^{-})_{n\in\mathbb{N}}$ are sequences of simple functions increasing (respectively) to $X^{+}$ and $X^{-}$ then,  following the definitions explained above, the function in Definition \ref{regconexpdef} is the same as the function
$$\tilde{X}:=\lim_{n}((\widetilde{X^{+}_{n}}-\widetilde{X^{-}_{n}})I_{\Omega_{X}}),$$
which is clearly $\mathcal{F}_{0}-$ measurable and is easily seen to satisfy 
$$E[\tilde{X}I_{A}]=E[XI_{A}]$$
for every $A\in\mathcal{F}_{0}$, being therefore a version of $E_{0}X$. It follows by linearity that the following proposition holds:

\end{remark}

\medskip
  
\begin{prop}
\label{proredregcondis}
In the context of Definition \ref{regconexpdef}, $\{\mathbb{P}_{\omega}\}_{\omega\in\Omega}$ is a decomposition of $E_{0}$ if and only if for every $A\in\mathcal{F}$, $\omega\mapsto \mathbb{P}_{\omega}(A)$ defines an $\mathcal{F}_{0}-$measurable version of $\mathbb{P}[A|\mathcal{F}_{0}]$. 
\end{prop}

Note that, necessarily, the set where the integral makes sense has $\mathbb{P}-$measure one, because
$$E|X|=EE_{0}|X|=\int_{\Omega}\int_{\Omega}|X(z)|d\mathbb{P}_{\omega}(z)d\mathbb{P}(\omega)$$

In other words, (\ref{defe0regequ}) defines an $\mathcal{F}_{0}-$measurable function and  $E_{0}X(\omega)=E^{\omega}X$, $\mathbb{P}-$a.s., where $E^{\omega}$ denotes integration with respect to $\mathbb{P}_{\omega}$. If $\mathcal{F}$ is countably generated and $E_{0}$ is regular then the following uniqueness (up to $\mathbb{P}-$negligible sets) result holds.

\medskip

\begin{prop}
\label{unidec}
In the context of Definition \ref{regconexpdef}, if $\mathcal{F}$ is countably generated and $E_{0}$ is regular, and given any two decompositions $\{\mathbb{P}_{1,\omega}\}_{\omega\in\Omega}$ and $\{\mathbb{P}_{2,\omega}\}_{\omega\in\Omega}$ of $E_{0}$, there exists a set $\Omega_{0}\subset\Omega$ with $\mathbb{P}\Omega_{0}=1$ such that for every $\omega\in\Omega_{0}$, $\mathbb{P}_{1,\omega}=\mathbb{P}_{2,\omega}$.
\end{prop}

{\bf Proof:} Denote by $E_{1,0}$ and $E_{2,0}$ the versions of $E_{0}$ given, respectively, by (integration with respect to) $\{\mathbb{P}_{1,\omega}\}_{\omega\in\Omega}$ and $\{\mathbb{P}_{2,\omega}\}_{\omega\in\Omega}$.

Now, given $A\in\mathcal{F}$, consider the function
\begin{equation}
\label{difdecA}
U_{A}(\omega):=\mathbb{P}_{1,\omega}(A)-\mathbb{P}_{2,\omega}(A)=:E_{1,0}I_{A}(\omega)-E_{2,0}I_{A}(\omega).
\end{equation}
Note that $U_{A}$ is $\mathcal{F}_{0}-$measurable and therefore so is the set $[U_{A}\geq 0]$. In particular
$$\int_{\Omega}U_{A}(\omega)I_{[U_{A}\geq 0]}(\omega)d\mathbb{P}(\omega)=$$
$$\int_{\Omega}E_{1,0}I_{A}(\omega)I_{[U_{A}\geq 0]}(\omega)d\mathbb{P}(\omega)-\int_{\Omega}E_{2,0}I_{A}(\omega)I_{[U_{A}\geq 0]}(\omega)d\mathbb{P}(\omega)=$$
$$\int_{\Omega}E_{0}[I_{A}I_{[U_{A}\geq 0]}](\omega)d\mathbb{P}(\omega)-\int_{\Omega}E_{0}[I_{A}I_{[U_{A}\geq 0]}](\omega)d\mathbb{P}(\omega)=0,$$
and by a similar argument using the set $[U_{A}<0]$ we conclude that there exists $\Omega_{A}$ with $\mathbb{P}\Omega_{A}=1$ such that $\mathbb{P}_{1,\omega}(A)=\mathbb{P}_{2,\omega}(A)$ for every $\omega\in\Omega_{A}$.

Let  $\{A_{k}\}_{k\in\mathbb{N}}$ be a countable $\pi-$system generating $\mathcal{F}$ and let $\Omega_{0}:=\cap_{k\in\mathbb{N}}\Omega_{A_{k}}$. Clearly, $\mathbb{P}\Omega_{0}=1$. 

By the $\pi-\lambda$ theorem (applied to the set of $A\in\mathcal{F}$ such that $\mathbb{P}_{1,\omega}(A)=\mathbb{P}_{2,\omega}(A)$ for every $\Omega\in\Omega_{0}$), $\mathbb{P}_{1,\omega}(A)=\mathbb{P}_{2,\omega}(A)$ for every $A\in\mathcal{F}$ and every $\omega\in\Omega_{0}$. This is, $\mathbb{P}_{1,\omega}=\mathbb{P}_{2,\omega}$ for every $\omega\in\Omega_{0}$
.\qed

\medskip

\begin{remark}
\label{remproone}
For future reference, we will point out the following: in the context of Definition \ref{regconexpdef}, and given a decomposition $\{\mathbb{P}_{\omega}\}_{\omega\in\Omega}$ of $E_{0}$, a set $A$ satisfies $\mathbb{P}A=1$ if and only if $\mathbb{P}_{\omega}A=1$ for $\mathbb{P}-$a.e $\omega$. This is a simple consequence of the equality $$\mathbb{P}A=\int_{\Omega}\mathbb{P}_{\omega}A\,d\mathbb{P}(\omega).$$ 
\end{remark}

\subsection{Regularity and $T-$Filtrations}
Assume that $(\mathcal{F}_{k})_{k\in{\mathbb{Z}}}$ is a given $T-$filtration (Definition \ref{defadafil}), that $E_{0}:=E[\cdot|\mathcal{F}_{0}]$ is regular, and that $\{\mathbb{P}_{\omega}\}_{\omega\in\Omega}$ is a given decomposition of $E_{0}$. How do we relate the conditional expectations $E[\,\cdot\,|\mathcal{F}_{k}]$ (which depend on $\mathbb{P}$) with the conditional expectations induced by $\mathbb{P}_{\omega}$?  The following answer is sufficient for our purposes:

\medskip

\begin{lemma} 
\label{perconexpadacas}
Let $(\mathcal{F}_{k})_{k\in\mathbb{Z}}$ be a $T-$filtration (Definition \ref{defadafil}) and for every $k\in\mathbb{Z}$, denote by $E_{k}$ the conditional expectation with respect to $\mathcal{F}_{k}$ and $\mathbb{P}$. Assume that $\mathcal{F}_{0}$ is countably generated (Definition \ref{deffcougen}), that $E_{0}$ is regular, and that $\{\mathbb{P}_{\omega}\}_{\omega\in\Omega}$ is a decomposition of $E_{0}$ (Definition \ref{regconexpdef}). Denoting further by $E^{\omega}_{k}$ the conditional expectation with respect to $\mathcal{F}_{k}$ and $\mathbb{P}_{\omega}$, the following property holds: for every $\mathbb{P}-$integrable $Y$, every $k\in \mathbb{Z}$, and every $\mathcal{F}_{k}-$measurable version of $E_{k}Y$, there exists $\Omega_{Y}$ with $\mathbb{P}\Omega_{Y}=1$ such that
\begin{equation}
\label{herconexp}
E^{\omega}_{k}Y=E_{k}Y
\end{equation}
$\mathbb{P}_{\omega}-$a.s. for every $\omega\in \Omega_{Y}$. 
\end{lemma}

\medskip

\begin{remark}
\label{remver}
Note that if $Z$ is any version of $E_{k}Y$, (\ref{herconexp}) and Remark \ref{remproone} imply that $E^{\omega}_{k}Y=Z$, $\mathbb{P}_{\omega}-$a.s for $\mathbb{P}-$a.e. $\omega$ (over a set of probability one depending on $Z$). 

\end{remark}
\medskip

{\bf Proof of Lemma \ref{perconexpadacas}:} Fix a version of $Y\in L^{1}_{\mathbb{P}}$. We will prove that for any ($\mathcal{F}_{k}-$measurable) version of $E_{k}Y$, there exists a set $\Omega_{Y}\subset \Omega$ with $\mathbb{P}\Omega_{Y}=1$  such that the following holds: for every $\omega \in \Omega_{Y}$ and every $A\in \mathcal{F}_{k}$
\begin{equation}
\label{equconexp}
\int_{A} Y(z)d\mathbb{P}_{\omega}(z)=\int_{A} E_{k}Y(z)d\mathbb{P}_{\omega}(z),
\end{equation}
this clearly implies the first conclusion.

Fix a ($\mathcal{F}_{k}-$measurable) version of $E_{k}Y$ and notice that for $A$ fixed, a set $\Omega_{Y,A}$ of probability one such that (\ref{equconexp}) holds for all $\omega\in\Omega_{Y,A}$ exists by the property defining the family $\{\mathbb{P}_{\omega}\}_{\omega\in \Omega}$ and because
$$E_{0}[YI_{A}]=E_{0}[(E_{k}Y)I_{A}],$$
$\mathbb{P}-$a.s. Without loss of generality $\Omega_{Y,A}\subset \{\omega\in \Omega: |Y|+|E_{k}Y|\in L^{1}_{\mathbb{P}_{\omega}}\}$ (the last set has $\mathbb{P}-$measure one because $E|Z|=EE_{0}|Z|$ for every $Z\in L^{1}_{\mathbb{P}}$).

Now proceed as follows: let $\{A_{n}\}_{n\in\mathbb{N}}\subset \mathcal{F}_{k}$ be a countable family generating $\mathcal{F}_{k}$  which is also a $\pi-$system and includes $\Omega$ (such a family exists because $\mathcal{F}_{0}$ is assumed countably generated), let $\Omega_{Y}:=\cap_{n\geq 1}\Omega_{Y,A_{n}}$, and let $\mathcal{G}_{k}\subset\mathcal{F}_{k}$ be the family of sets $A \in \mathcal{F}_{k}$ such that ($\ref{equconexp}$) holds for all $\omega \in \Omega_{Y}$. It is easy to see that $\mathcal{G}_{k}$ is a $\lambda-$system and therefore, since it includes $\{A_{n}\}_{n\in\mathbb{N}}$, the $\pi-\lambda$ theorem implies that $\mathcal{G}_{k}=\mathcal{F}_{k}$. Note that $\mathbb{P}\Omega_{0,Y}=1$, and that for all $\omega\in \Omega_{Y}$, (\ref{equconexp}) holds for all $A\in \mathcal{F}_{k}$. 


This gives the proof of the first conclusion. 
The second conclusion (the one about martingales) follows easily from this, together with the fact that  $E|X|^{p}=EE_{0}|X|^{p}$ and therefore $E|X|^{p}<\infty$ if and only if $E^{\omega}|X|^{p}<\infty$ for $\mathbb{P}-$a.e. $\omega$.\qed

\medskip

\begin{cor}
\label{cormarundpome}
In the context of Lemma \ref{perconexpadacas} and denoting further by $E^{\omega}$ the integration with respect to $\mathbb{P}_{\omega}$, if $p\geq 1$ and $D_{0}\in L^{p}_{\mathbb{P}}(\mathcal{F}_{0})$ is such that $E_{-1}D_{0}=0$, there exists a set $\Omega_{0}\subset \Omega$ with $\mathbb{P}\,\Omega_{0}=1$ such that for every $k\geq 1$ and every $\omega\in \Omega_{0}$, $E^{\omega}|T^{k}D_{0}|^{p}<\infty$ and $E_{k-1}^{\omega}T^{k}D_{0}=0$, $\mathbb{P}_{\omega}-$a.s.
\end{cor}

It follows easily that { if $(T^{k}D_{0})_{k\in\mathbb{N}}$ is a $(\mathcal{F}_{k})_{k\in\mathbb{N}}$ adapted (stationary) sequence of martingale differences in $L^{p}_{\mathbb{P}}$, then for $\mathbb{P}-$almost every $\omega$ (over a set depending on fixed versions of $(T^{k}D_{0})_{k\in\mathbb{N}}$), $(T^{k}D_{0})_{k\in\mathbb{N}}$ is a $(\mathcal{F}_{k})_{k\in \mathbb{N}}$ adapted (not necessarily stationary) sequence of martingale differences in $L^{p}_{\mathbb{P}_{\omega}}$}.

{\bf Proof of Corollary \ref{cormarundpome}:} Let $D_{k}:=T^{k}D_{0}$ ($k\in\mathbb{Z}$), and first note that $E_{k-1}D_{k}=0$, $\mathbb{P}-$a.s. for every $k\in\mathbb{Z}$. 

Now let $\Omega_{0,1}$ be a set of probability one such that if $\omega\in \Omega_{0,1}$, $E_{k-1}D_{k}=E^{\omega}_{k-1}D_{k}$ $\mathbb{P}_{\omega}-$a.s. for all $k\geq 1$ (Lemma \ref{perconexpadacas} and Remark \ref{remver}),  let $\Omega_{0,2}$ be a set of probability one with the property that for all $\omega\in \Omega_{0,2}$ and all $k\geq 1$, $E_{k-1}D_{k}=0$ $\mathbb{P}_{\omega}-$a.s. (Remark \ref{remproone}), and let $\Omega_{0,3}$ be a set of probability one such that for all $\omega\in \Omega_{0,3}$ and all $k\geq 0$, $E^{\omega}|D_{k}|^{p}<\infty$ (such a set exists because $\infty>E|D_{k}|^{p}=E[E_{0}|D_{k}|^{p}]$). The set $\Omega_{0}=\cap_{j=1}^{3}\Omega_{0,j}$ satisfies the given conclusion.\qed

\section{Examples of Regularity}
\label{exareg}

In this section we illustrate the notion of regularity by constructing regular conditional expectations in specific settings. The setting in Section \ref{examarpro} is of particular interest due to its generality and its importance along the applications.

Let us start by illustrating the trivial cases:

\medskip

\begin{exam}[Trivial Examples of Regularity]
\label{triexareg}
In the context of Definition \ref{regconexpdef}, if $\mathcal{F}$ includes the singletons $\{\omega\}$ ($\omega\in \Omega$), $\mathcal{F}_{0}=\mathcal{F}$, and for a given $\omega\in\Omega$, $\delta_{\omega}$ denotes the Dirac probability measure at $\omega$ ($\delta_{\omega}\{\omega\}=1$), then $\{\delta_{\omega}\}_{\omega\in\Omega}$ is a decomposition of $E_{0}$. If $\mathcal{F}_{0}=\{\emptyset,\Omega\}$ is the trivial sigma-algebra, then we get a decomposition $\{\mathbb{P}_{\omega}\}_{\omega\in\Omega}$ of $E_{0}$ by taking $\mathbb{P}_{\omega}=\mathbb{P}$  for every $\omega\in \Omega$.
\end{exam}

\medskip

\subsection{Functions of i.i.d. Sequences}
\label{exaregconexp}
The simplest non-trivial example of a regular conditional expectation is perhaps the following:

\medskip

\begin{exam}[Functions of i.i.d. sequences]
\label{exaregconexpexa}
Refer to the setting explained along Example \ref{exalinpro} on page \pageref{exalinpro} and consider the following observation: if $(\Omega^{-},\mathcal{F}^{-})$ and $(\Omega^{+},\mathcal{F}^{+})$ denote respectively the space of complex-valued sequences of the form $(a_{k})_{k\leq 0}$ and $(a_{k})_{k>0}$ ($k\in\mathbb{Z}$) with the sigma algebras $\mathcal{F}^{-}$ and $\mathcal{F}^{+}$ generated by the respective  finite dimensional cylinders,  then 
  $(\mathbb{C}^{\mathbb{Z}},\mathcal{F})=(\Omega^{-}\times \Omega^{+},\mathcal{F}^{-}\otimes \mathcal{F}^{+})$ and the projections (defined in an obvious way) $\pi^{-}: \Omega \to \Omega^{-}$ and $\pi^{+}: \Omega \to \Omega^{+}$ are measurable with respect to the respective sigma-algebras. Note also that $\mathcal{F}_{0}=(\pi^{-})^{-1}{
  \mathcal{F}^{-}}$.
  
  For every $\omega\in\Omega$, let $\omega^{+}:=\pi^{+}(\omega)$ and $\omega^{-}:=\pi^{-}(\omega)$, 
  consider the function $\delta_{\omega}:\Omega\to\Omega$ given by
  $$\delta_{\omega}(z)=(\omega^{-},z^{+})$$
  and define the measure of ``partial integration with respect to the future'' $\mathbb{P}_{\omega}:=\mathbb{P}\delta_{\omega}^{-1}$.  
   We claim that if the sigma algebras $\sigma(\xi_{k})_{k\leq 0}$ and $\sigma(\xi_{k})_{k> 0}$ are independent (with respect to $\mathbb{P}'$) or, equivalently, if $(\xi_{k})_{k\in\mathbb{Z}}$ is i.i.d. (consider the hypothesis  of stationarity)  then $\{\mathbb{P}_{\omega}\}_{\omega\in \Omega}$ is a decomposition of $E_{0}$. 
   
   Let us prove this: first note that, by the hypothesis of independence, $\mathbb{P}=\mathbb{P}^{-}\times\mathbb{P}^{+}$, where $\mathbb{P}^{-}$ (respectively $\mathbb{P}^{+}$) is the measure in $(\Omega^{-},\mathcal{F}^{-})$ (respectively $(\Omega^{+},\mathcal{F}^{+})$) induced by $(\xi_{k})_{k\leq 0}$ (respectively $(\xi_{k})_{k>0}$) by the procedure explained in Example \ref{exalinpro}.
   
Now  fix $A\in\mathcal{F}$, and let us give an explicit formula for $\mathbb{P}_{\omega}({A})$:
   \begin{equation}
   \label{pomeequfuniid}
   \mathbb{P}_{\omega}A=\mathbb{P}[\delta_{\omega}\in A]=\mathbb{P}\{z\in\Omega:(\omega_{-},z_{+})\in A\}=\mathbb{P}_{+}\{y \in \Omega_{+}: (\omega_{-},y)\in A\}
   \end{equation}
where we used Fubini's Theorem (see \cite{bilpromea} Theorem 18.3, see also Theorems 18.1 and 18.2 there\footnote{These are theorems related to real-valued functions, but this poses no serious restriction. The reader may as well replace ``$\mathbb{C}$'' by ``$\mathbb{R}$'' in this example and refer to Example \ref{examarpro} to cover the complex-valued case.}) to guarantee the validity of (\ref{pomeequfuniid}). By Fubini's theorem again, the function $u:(\Omega^{-},\mathcal{F}^{-})\to [0,\infty)$ given by 
\begin{equation}
\label{equparintseccom}
u(x)=\mathbb{P}^{+}\{y \in \Omega^{+}: (x,y)\in A\}
\end{equation}
is $\mathcal{F}^{-}-$measurable. Since $\omega\mapsto \mathbb{P}_{\omega}(A)$ is the same as $\omega\mapsto u\circ \pi^{-}(\omega)$, it follows that $\omega\mapsto \mathbb{P}_{\omega}(A)$ is $\mathcal{F}_{0}-$measurable.

Now, $\mathcal{F}_{0}=\{B\times \Omega^{+}:B\in\mathcal{F}^{-}\}$ (to see this use, for instance, the $\pi-\lambda$ theorem), and a further application of Fubini's theorem shows that for every $B\in \mathcal{F}^{-}$
$$\int_{B\times\Omega^{+}}\mathbb{P}_{\omega}(A)\,d\mathbb{P}(\omega)= \int_{B\times\Omega^{+}}\mathbb{P}^{+}\{y \in \Omega^{+}: (\omega^{-},y)\in A\}\,d\mathbb{P}(\omega)=$$
\begin{equation}
\label{equverconexp}
\int_{B} \mathbb{P}^{+}\{y \in \Omega^{+}: (x,y)\in A\}\,d\mathbb{P}^{-}(x)=\mathbb{P}(A\cap (B\times\Omega^{+})).
\end{equation}

These facts show that for every $A\in\mathcal{F}$, $\omega\mapsto \mathbb{P}_{\omega}(A)$ defines a version of $\mathbb{P}[A|\mathcal{F}_{0}]$, and an application of Proposition \ref{proredregcondis} shows that, indeed, $\{\mathbb{P}_{\omega}\}_{\omega\in \Omega}$ is a decomposition of $E_{0}$.
\end{exam}

\subsection{Functions of Stationary Markov Chains}
\label{examarpro}

To begin with, let us recall the notion of a {transition probability matrix:}
 
 \medskip
 
 \begin{dfn}[Transition Probability Matrix]
 \label{traprodef}
 Given two measurable spaces $(\Omega,\mathcal{F})$ and $(\mathcal{K},\mathcal{G})$, a {\bfseries transition probability matrix} between $(\Omega,\mathcal{F})$ and $(\mathcal{K},\mathcal{G})$ is a function
 \begin{equation}
 \label{traprodefequ}
 P:\Omega\times \mathcal{G}\to [0,1]
 \end{equation}
 with the property that for every $\omega\in \Omega$, $P(\omega,\cdot)$ is a probability measure in $\mathcal{G}$ and for every $G\in\mathcal{G}$, $P(\cdot,G)$ is $\mathcal{F}-$measurable. If $(\Omega,\mathcal{F})=(\mathcal{K},\mathcal{G})$, we say that $P$ is a transition probability matrix {\bfseries in} $(\Omega,\mathcal{F})$.
 \end{dfn}

We also introduce the following terminology.

\medskip

\begin{dfn}[Markov Chains]
\label{marchadef}
Assume that for every $k\in\mathbb{Z}$, a measurable space $(S_{k},\mathcal{S}_{k})$ is given, and let $(\xi_{k})_{k\in\mathbb{Z}}$ be a sequence of random elements $\xi_{k}:\Omega\to S_{k}$   defined on the  same probability space $(\Omega,\mathcal{F},\mathbb{P})$ (Definition \ref{raneledef}). We say that $(\xi_{k})_{k\in\mathbb{Z}}$ is a {\bfseries Markov chain}  if for every $k\in\mathbb{Z}$ there exists a transition probability matrix (Definition \ref{traprodef}) $P_{k}:S_{k}\times\mathcal{S}_{k+1}\to [0,1]$  such that for every $A\in\mathcal{S}_{k+1}$,
\begin{equation}
\label{marchapro}
\omega\mapsto {P}_{k}(\xi_{k}(\omega),A)
\end{equation}
defines a version of $\mathbb{P}[\xi_{k+1}\in A|\sigma(\xi_{j})_{j\leq k}]:= E[I_{A}\circ \xi_{k+1}| \sigma(\xi_{j})_{j\leq k}]$ (the conditional expectation is taken with respect to $\mathbb{P}$). The Markov chain has a {\bfseries fixed state space} if $(S_{k},\mathcal{S}_{k})=(S_{0},\mathcal{S}_{0})$ for every $k\in\mathbb{Z}$. If the state space is fixed, the Markov chain  is {\bfseries stationary} if $(\xi_{k})_{k\in\mathbb{Z}}$ is stationary (the law of $(\xi_{n},\cdots,\xi_{n+k})$ is the same for every $n\in\mathbb{Z}$ if $k\in\mathbb{Z}$ is fixed), and it is {\bfseries homogeneous} if (we can choose) $P_{k}=P_{0}$  for all $k\in\mathbb{Z}$.  
\end{dfn}

\medskip

\begin{remark}
\label{staimphom}
Every stationary Markov chain is homogeneous: in the context of Definition \ref{marchadef}, given $k\in\mathbb{Z}$ and $A,B\in\mathcal{S}:=\mathcal{S}_{0}$,
$$\mathbb{P}([\xi_{k+1}\in A]\cap [\xi_{k}\in B])=\mathbb{P}([\xi_{1}\in A]\cap [\xi_{0}\in B])= \int_{\Omega}P_{0}(\xi_{0}(\omega),A)I_{B}(\xi_{0}(\omega))d\mathbb{P}(\omega)=$$
$$\int_{S}P_{0}(x,A)I_{B}(x)\,d\mathbb{P}\xi_{0}^{-1}(x)=\int_{S}P_{0}(x,A)I_{B}(x)\,d\mathbb{P}\xi_{k}^{-1}(x)=$$
$$\int_{\Omega}P_{0}(\xi_{k}(\omega),A)I_{B}(\xi_{k}(\omega))\,d\mathbb{P}(\omega),$$
so that, necessarily, $P_{k}(\xi_{k}(\omega),A)=P_{0}(\xi_{k}(\omega),A)$, $\mathbb{P}-$a.s. And we can replace $P_{k}=P_{0}$.
\end{remark}
Note also that condition (\ref{marchapro}) implies, in particular, that for every $A\in\mathcal{S}$ and every $k\in\mathbb{Z}$
\begin{equation}
\label{genmarchaequ} 
\mathbb{P}[\xi_{k+1}\in A|\sigma(\xi_{j})_{j\leq k}]=\mathbb{P}[\xi_{k+1}\in A|\sigma(\xi_{k})].
\end{equation}

\medskip

\begin{dfn}[Generalized Markov Chain]
\label{genmarcha}
In the context of Definition \ref{marchadef}, if we can verify (\ref{genmarchaequ}) (regardless of whether the family of transitions matrices $(P_{k})_{k\in\mathbb{Z}}$ satisfying (\ref{marchapro}) exists), we call $(\xi_{k})_{k\in\mathbb{Z}}$ is a {\bfseries generalized Markov chain}. The generalized Markov chain has a fixed state space if for every $k\in\mathbb{Z}$, $(S_{k},\mathcal{S}_{k})=(S_{0},\mathcal{S}_{0})$.
\end{dfn}

Every generalized Markov chain whose state space is a complete and separable metric space (with is Borel sigma algebra) is a Markov chain:

\medskip

\begin{prop}[Existence of Markov Kernels]
\label{proeximarker}
If $(S,d)$ is a complete and separable metric space with Borel sigma algebra $\mathcal{S}$ and $\xi_{1},\xi_{2}$ are random elements on $S$ defined on the same probability space $(\Omega,\mathcal{F},\mathbb{P})$,  there exists a transition probability matrix $P$ in $(S,\mathcal{S})$ (Definition \ref{traprodef}) such that the map $\Omega\to [0,1]$ given by
$$\omega\mapsto P(\xi_{1}(\omega),A)$$
defines a version of $\mathbb{P}[\xi_{2}\in A|\sigma(\xi_{1})]$. 
\end{prop}

{\bf Proof:} This follows from Exercise 1 in \cite{bau}, Section 44. \qed


\medskip
Our first example in this section is the following:

\medskip

\begin{exam}[Functions of Stationary Markov Chains]
\label{examarproexa}
 
To motivate the construction that follows, start by  noticing that Example  \ref{exaregconexpexa} can be extended to the context in which $\sigma(\xi_{k})_{k\leq 0}$ and $\sigma(\xi_{k})_{k\geq 0}$ are not necessarily independent, provided that if we replace $\mathbb{P}^{+}$ by a measure $\mathbb{P}^{x}$ in (\ref{equparintseccom}) the function $u(x)$ is (still) $\mathcal{F}^{-}-$measurable, and that we can (still) verify the equalities in (\ref{equverconexp}) (dropping the third term) with $\mathbb{P}^{+}$ replaced by $\mathbb{P}^{x}$. 

This is the case for instance in the context of {\it functions of stationary
 , homogeneous Markov Chains on 
a complete and separable metric space $S$}. In what follows $(S,d)$ denotes a 
complete and separable metric space with metric $d$ and Borel sigma algebra $\mathcal{S}$ and, for every $k\in\mathbb{N}^{*}$, $\mathcal{S}^{k}$ denotes the product sigma algebra in $S^{k}$. The rest of the notation copies that in Example \ref{exalinpro} and Section \ref{exaregconexp}: we will work (again) over the space $(\Omega,\mathcal{F})$ of $S-$valued sequences $(a_{k})_{k\in\mathbb{Z}}$ with the sigma algebra generated by the finite-dimensional cylinders, and we will use the decomposition $(\Omega,\mathcal{F})=(\Omega^{-}\times\Omega^{+},\mathcal{F}^{-}\otimes\mathcal{F}^{+})$ as in Section \ref{exaregconexp} (with $\mathbb{C}$ replaced by $S$). Again, $\mathcal{F}_{0}=\sigma(x_{k})_{k\leq0}$, where $x_{j}:\Omega\to S$ is the projection in the $j-$th coordinate.

Let us start by explaining the construction of the processes under consideration. 
\subsubsection*{Construction of stationary, homogeneous Markov Chains}

Let $P:S\times\mathcal{S}\to [0,1]$ be a transition probability matrix in $(S,\mathcal{S})$ (Definition \ref{traprodef}). 
We will also assume that $P(\cdot,\cdot)$ admits an invariant probability measure $\mathbb{P}$ (see item 6. below). Our goal is to construct a probability measure $\mathbb{P}_{\mathbb{Z}}$ on $(S^{\mathbb{Z}},\mathcal{S}^{\mathbb{Z}})$ such that the coordinate functions $x_{k}:S^{\mathbb{Z}}\to S$ define a stationary Markov Chain (with $\mathbb{P}_{\mathbb{Z}}x_{0}^{-1}=\mathbb{P}$) on $(S^{\mathbb{Z}},\mathcal{S}^{\mathbb{Z}},\mathbb{P}_{\mathbb{Z}})$ with state space $S$ and transition probability $P$:  for every $k\in\mathbb{Z}$, 
$$\mathbb{P}_{\mathbb{Z}}(x_{k+1}\in A|x_{k})=P(x_{k},A).$$
The construction can be summarized as follows:
\begin{enumerate}
\item Given $A_{0},A_{1},\cdots, A_{k}\in\mathcal{S}$, $x\in S$ and $n\in\mathbb{Z}$, and denoting  the integral of a measurable function $f:S\to\mathbb{C}$ with respect to the measure $P(x.\cdot)$ by
$$\int_{S}f(y)P(x,\,dy)$$
define $P_{0}^{k}(x, A_{0}\times\cdots\times A_{n})$ by
$$P_{0}^{k}(x, A_{0}\times\cdots\times A_{n}):=I_{A_{0}}(x)\int_{A_{1}}\cdots\int_{A_{k-1}} P(y_{k-1},A_{k})P(y_{k-2},dy_{k-1})\cdots P(x,dy_{1})$$
\begin{equation}
\label{equmarchalawrec}
=\int_{S}\cdots\int_{S}I_{A_{0}\times A_{1}\times\cdots\times A_{k}}(x,y_{1},\cdots ,y_{k})P(y_{k-1},dy_{k})\cdots P(x,dy_{1}).
\end{equation}
The $\pi-\lambda$ theorem and Kolmogorov's existence theorem\footnote{It is important to point out that Kolmogorov's existence theorem is not guaranteed without special assumptions on the structure of the underlying measurable space (see \cite{andandjes} for counterexamples). The validity of Kolmogorov's existence theorem for the case of complete and separable metric spaces is, on the other side, a well established fact.} guarantee that $P_{0}^{k}(x,\cdot)$, thus defined over rectangular sets, extends in a unique way to a measure $P_{0}^{\infty}(x,\cdot)$ on $\mathcal{S}^{\mathbb{N}}$ (more explicitly:  Theorem 3.1 in \cite{bilpromea} allows us to see that (\ref{equmarchalawrec}) defines a unique probability measure  in $\mathcal{S}^{k+1}$, and the $\pi-\lambda$ theorem guarantees that if $A\in \mathcal{S}^{k}$, $P_{0}^{k}(x,A\times S)=P_{0}^{k-1}(x,A)$. An application of Proposition III-3-3 in \cite{nev} implies that $P_{0}^{\infty}(x,\cdot)$ exists and is unique). Even more  (see  Proposition V-2-1 in \cite{nev}), for every $A\in\mathcal{S}^{\mathbb{N}}$, the function $S\to[0,1]$ given by
$$x\mapsto P_{0}^{\infty}(x,A)$$
is $\mathcal{S}-$measurable. Thus $P_{0}^{\infty}:S\times \mathcal{F}^{+}\to [0,1]$ is a transition probability between $(S,\mathcal{S})$ and $(\Omega^{+},\mathcal{F}^{+})$. 

\item We can extend $P_{0}^{\infty}$ to a transition probability between $(\Omega^{-},\mathcal{F}^{-})$ and $(\Omega^{+},\mathcal{F}^{+})$ in the following way: given $\omega^{-}\in \Omega^{-}$ and $A\in\mathcal{F}^{+}$
\begin{equation}
\label{traproomeminomeplu}
\mathbb{P}_{0}^{\infty}(\omega^{-},A):= P_{0}^{\infty}(x_{0}(\omega^{-}),A).
\end{equation}
$\mathbb{P}_{0}^{\infty}$, thus extended, is clearly $\mathcal{F}^{-}$ measurable for every fixed $A$, showing that it is (indeed) a transition probability matrix $\Omega^{-}\times \mathcal{F}^{+}\to[0,1]$.

\item Notice that for every $k\in\mathbb{N}^{*}$, $P_{0}^{\infty}$ restricts to a transition probability $P_{0}^{k}:S\times\mathcal{S}^{k}\to[0,1]$ in the obvious way: if $\pi_{k}:S^{\mathbb{N}}\to S^{k}$ is the natural projection, $P_{0}^{k}(x,A):=P_{0}^{\infty}(x,\pi_{k}^{-1}(A))$. 

More explicitly, note that $P_{0}^{1}(x,A)=P(x,A)$ and for general $k$, $P_{0}^{k}(x,A)$ is given by the last line of (\ref{equmarchalawrec}) with $I_{A_{0}\times\cdots\times A_{k}}$ replaced by $I_{A}$ (apply the $\pi-\lambda$ theorem to the $\lambda-$system of sets in $\mathcal{S}^{k+1}$ for which this  holds).

\item Given a probability measure $\mu$ on $\mathcal{S}$, $n\in\mathbb{Z}$ and $k\in\mathbb{N}$, define the probability measure $\mu_{n}^{n+k}$ in $\mathcal{S}^{k}$ in the following way: for every $A\in \mathcal{S}^{k}$
$$\mu_{n}^{n+k}A:=\int_{S}P_{0}^{k}(x,A)\,d\mu(x).$$ 

\item If we assume that $\mathbb{P}$, a probability measure on $\mathcal{S}$, is a {\it stationary probability measure} for $P$. This is, that $\mathbb{P}_{0}^{1}=\mathbb{P}$, then it is easy to see that for every $n\in\mathbb{Z}$, $k\in\mathbb{N}$ and $A\in\mathcal{S}$, 
\begin{equation}
\label{comequ}
\mathbb{P}_{n}^{n+k}(S\times\cdots\times S\times A)=\mathbb{P}(A).
\end{equation}

It follows from this that for every simple function $f(y)=\sum_{j=1}^{r}a_{j}I_{A_{j}}$ ($A_{j}\in\mathcal{S}$):
$$\int_{S}\cdots \int_{S}f(y)P(y_{n-1},y)\cdots P(x,y_{n-l+1})d\mathbb{P}(x)=\int_{S}f(x)\,d\mathbb{P}(x)$$
and by  an approximation argument analogous to the one leading to Proposition \ref{proredregcondis} the same holds for every $\mathbb{P}$-integrable function $f$.

\item In particular the following holds: for every $A_{n},\cdots, A_{n+k}\in\mathcal{S}$, if we denote 
$$f(y):=\int_{S}\cdots \int_{S}I_{A_{n}\times\cdots \times A_{n+k}}(y,y_{n+1},\cdots ,y_{n+k})P(y_{n+k-1},dy_{n+k})\cdots P(y, dy_{n-1})$$
then
$$\mathbb{P}_{n-l}^{n+k}(S^{l}\times A_{n}\times\cdots \times A_{n+k})=\int_{S}\cdots \int_{S}f(y_{n})P(y_{n-1},y_{n})\cdots P(x,y_{n-l+1})d\mathbb{P}(x)=$$
$$\int_{S}f(x)\,d\mathbb{P}(x)=\int_{S}P_{0}^{k}(x,A_{n}\times\cdots \times A_{n+k}) \,d\mathbb{P}(x)=:\mathbb{P}_{0}^{k}(A_{n}\times\cdots \times A_{n+k}),$$
and it follows by a further application of the $\pi-\lambda$ theorem and Kolmogorov's existence theorem that there exists a unique probability measure $\mathbb{P}_{\mathbb{Z}}$ on $(\Omega,\mathcal{F})$ such that for every $k\in\mathbb{N}$ and every set of the form
\begin{equation}
\label{genfindimcyl}
H_{k,A}:=[(x_{-k},\cdots ,x_{k})\in A]
\end{equation}
where $A\in\mathcal{S}^{2k+1}$, 
$$\mathbb{P}_{\mathbb{Z}}(H_{k,A})=\mathbb{P}_{0}^{2k+1}(A)$$ 
(note that the sets of the form (\ref{genfindimcyl}) indeed generate $\mathcal{F}$).
\end{enumerate}

The coordinate functions $(x_{k})_{k\in\mathbb{Z}}$ give, in this setting,  a stationary Markov chain defined on $(\Omega,\mathcal{F},\mathbb{P}_{\mathbb{Z}})$ with transition probability $P$ and law $\mathbb{P}$ (see \cite{nev}, V-2 for more details on this). It is not hard to see in particular that for every $n\in\mathbb{Z}$, $k\in\mathbb{N}$ and  $f\in L^{1}_{\mathbb{P}}$, 
$$E[f(x_{n+k})|\sigma(x_{n})](\omega)=\int_{S} f(y)P_{0}^{k}(x_{n}(\omega),dy),$$
where $P_{0}^{k}(x,dy)$ denotes (in this case) the marginal distribution 
$$P_{0}^{k}(x,A)=P_{0}^{k}(x,S\times\cdots\times S\times A),$$ 
for $A\in\mathcal{S}$. More generally, given any function $f:S^{k}\to\mathbb{C}$ such that $f\circ (x_{1},\dots,x_{k})$ is $\mathbb{P}_{\mathbb{Z}}-$integrable:
\begin{equation}
\label{maropegenequ}
E[f\circ (x_{n+1},\dots,x_{n+k})|\sigma(x_{n})](\omega)=\int_{S}\cdots\int_{S}f(z_{1},\cdots,z_{k})P(z_{k-1},dz_{k})\cdots P(x_{n}(\omega),dz_{1}).
\end{equation}

For future reference, we will introduce the notation
\begin{equation}
\label{maropefunequ}
(P_{0}^{k}f)(z_{0}):=\int_{S}\cdots\int_{S}f(z_{1},\cdots,z_{k})P(z_{k-1},dz_{k}),\cdots,P(z_{0},dz_{1}).
\end{equation}

where $f:S^{k}\to\mathbb{C}$ is an appropriate function (in particular $E[f\circ (x_{n+1},\dots,x_{n+k})|\sigma(x_{n})]=P_{0}^{k}f\circ x_{n}$).

\subsubsection{Regularity}
We will show now that, again, $E_{0}=E[\,\cdot\,|\mathcal{F}_{0}]$ is regular. 
 
\begin{enumerate}
\item To do so we proceed as follows: given $\omega\in\Omega$, let $\mathbb{P}_{\omega}$ be the probability measure  on $\mathcal{F}=\mathcal{F}^{-}\otimes\mathcal{F}^{+}$ given in the following way: for $A\in \mathcal{F}$,
\begin{equation}
\label{defregconexpmarchaequ}
\mathbb{P}_{\omega}(A)=\mathbb{P}_{0}^{\infty}(\omega^{-},\{y\in\Omega^{+}:(\omega^{-},y)\in A\})
\end{equation}
where $\mathbb{P}_{0}^{\infty}$ is given by (\ref{traproomeminomeplu}). We proceed now to verify that for every $A\in\mathcal{F}$, $\omega\mapsto \mathbb{P}_{\omega}(A)$ is a version of $\mathbb{P}_{\mathbb{Z}}[A|\mathcal{F}_{0}]$ which (again) is sufficient to prove the regularity of $E[\,\cdot\,|\mathcal{F}_{0}]$ in virtue of Proposition \ref{proredregcondis}.

\item{\it $\mathcal{F}_{0}-$ measurability.} To see that $\omega\mapsto \mathbb{P}_{\omega}(A)$ is $\mathcal{F}_{0}-$measurable note that, by the $\pi-\lambda$ theorem applied to the $\lambda-$system of sets $A\in\mathcal{F}$ such that $\omega\mapsto \mathbb{P}_{\omega}(A)$ is $\mathcal{F}_{0}-$measurable, it suffices to see that this is the case under the assumption that $A=A^{-}\times A^{+}\in\mathcal{F}^{-}\times\mathcal{F}^{+}$. But it is easy to see that, in this case
$$\mathbb{P}_{\omega}(A)=I_{A^{-}}(\omega^{-})\mathbb{P}_{0}^{\infty}(\omega^{-},A^{+})$$
which defines an $\mathcal{F}_{0}-$measurable function of $\omega$ because the function $f_{A}:(\Omega^{-},\mathcal{F}^{-})\to [0,1]$ given by 
$$f_{A}(u)= I_{A^{-}}(u)\mathbb{P}_{0}^{\infty}(u,A^{+})$$
is $\mathcal{F}^{-}-$measurable and 
$$\mathbb{P}_{\omega}(A)=f_{A}\circ \pi^{-}(\omega).$$


\item{\it Integral equation.} To check that for every $A\in\mathcal{F}_{0}$ and $B\in\mathcal{F}$ 
\begin{equation}
\label{intequ}
\int_{\Omega}\mathbb{P}_{\omega}(B)I_{A}(\omega)\,d\mathbb{P}_{\mathbb{Z}}(\omega)=\mathbb{P}_{\mathbb{Z}}(A\cap B)
\end{equation}
we start by noticing the following: if  we can check (\ref{intequ}) for 
\begin{equation}
\label{fincylgenf0}
A'=[(x_{-k},\dots, x_{0})\in A_{-k}'\times\cdots \times A_{0}']
\end{equation} 
fixed and every set $B$ of the form 
\begin{equation}
\label{findimcylgenf}
B=[(x_{-l},\dots,x_{l})\in A_{-l}\times\cdots\times A_{0}\times B_{1}\times\cdots \times B_{l}]
\end{equation}
then, by the $\pi-\lambda$ theorem, (\ref{intequ}) holds for every $B\in\mathcal{F}$ whenever $A'$ is a finite dimensional cylinder of the form (\ref{fincylgenf0}). Then, since for fixed $B\in\mathcal{F}$, (\ref{intequ}) holds for every finite dimensional cylinder $A'$ of the form (\ref{fincylgenf0}) and these generate $\mathcal{F}_{0}$, a new application of the $\pi-\lambda$ theorem gives (\ref{intequ}) for every $A\in \mathcal{F}_{0}$.

\item Thus it suffices to check (\ref{intequ}) for $A'$, $B$ as in (\ref{fincylgenf0}) and (\ref{findimcylgenf}). Note that without loss of generality we can assume that $k=l$, and that in this case, taking $C_{j}:=A_{j}\cap A_{j}'$ ($j=-k,\cdots,0$),
$$\mathbb{P}_{\omega}(B)I_{A'}(\omega)=I_{C_{-k}\times\cdots \times C_{0}}(x_{-k}(\omega),\cdots,x_{0}(\omega)){P}_{0}^{k}(x_{0}(\omega),B_{1}\times\cdots\times B_{k})=\mathbb{P}_{\omega}(C),$$
where
$$C=[(x_{-k},\cdots,x_{k})\in C_{-k}\times \cdots\times C_{0}\times B_{1}\times\cdots \times{B}_{k}].$$
In conclusion, it suffices to see that if $B$ is any cylinder of the form (\ref{findimcylgenf}):
$$\int_{\Omega}\mathbb{P}_{\omega}(B)\,d\mathbb{P}_{\mathbb{Z}}(\omega)=\mathbb{P}_{\mathbb{Z}}(B).$$
Let us do this for the case $k=1$ (the general case is analogous):
$$\int_{\Omega}\mathbb{P}_{\omega}(B)d\mathbb{P}_{\mathbb{Z}}(\omega)=\int_{\Omega} I_{A_{-1}\times A_{0}}(x_{-1}(\omega),x_{0}(\omega))P_{0}^{1}(x_{0}(\omega), B_{1})d\mathbb{P}_{\mathbb{Z}}(\omega)=$$
$$\int_{\Omega}\left(\int_{S}I_{A_{-1}\times A_{0}\times B_{1}}(x_{-1}(\omega),x_{0}(\omega),z_{1})P(x_{0}(\omega),dz_{1})\right)\,d\,\mathbb{P}_{\mathbb{Z}}(\omega)=$$
$$\int_{S}\int_{S}\int_{S}I_{A_{-1}\times A_{0}\times B_{1}}(z_{-1},z_{0},z_{1})P(z_{0},dz_{1}))P(z_{-1},dz_{0})d\mathbb{P}(z_{-1})=\mathbb{P}(B)$$
as desired.




\end{enumerate}
The result of this construction can be summarized in the following way: 

\medskip

\begin{prop}[Functions of Stationary Markov Chains and Regular Conditional Expectations]
\label{profunstamarcha}
{\it If $(\Omega,\mathcal{F},\mathbb{P}_{\mathbb{Z}})$ is the probability space  constructed above, $T:\Omega\to\Omega$ is the left shift (specified again by $x_{k}\circ T=x_{k+1}$) and for some $p\geq 1$, $f:\Omega^{-}\to \mathbb{C}$ belongs to $L^{p}_{\mathbb{P}_{\mathbb{Z}}}$ (where $f$ is extended to $\Omega$ in the obvious way: $\tilde{f}(\omega)=f(\omega^{-})$), then
\begin{enumerate}
\item[(a).] If  $\mathcal{F}_{0}:=\sigma(x_{k})_{k\leq 0}$, then $(\mathcal{F}_{k})_{k\in\mathbb{Z}}:=(T^{k}\mathcal{F}_{0})_{k\in\mathbb{Z}}=(\sigma(x_{j})_{j\leq k})_{k\in\mathbb{Z}}$ is a $T-$filtration.
\item[(b).] If for every $k\in\mathbb{Z}$, $X_{k}:=T^{k}f:=f\circ T^{k}$, then the stationary sequence $(X_{k})_{k\in\mathbb{Z}}$ is $(\mathcal{F}_{k})_{k\in\mathbb{Z}}-$adapted.
\item[(c).] The conditional expectation $E_{0}=E[\,\cdot\,|\mathcal{F}_{0}]$ is regular, and for every $X\in L^{1}_{\mathbb{P}_{\mathbb{Z}}}$, a version of $E_{0}X$ is given by
$$E_{0}[X](\omega)=\int_{\Omega}X(z)\,d\mathbb{P}_{\omega}(z)$$
where, for every $\omega\in\Omega$, $\mathbb{P}_{\omega}$ is given by (\ref{defregconexpmarchaequ}).
\end{enumerate}}
\end{prop}


%
\end{exam}

Our last example shows how to represent a stationary sequence of random functions (on a complete and separable metric space) as a function of a Markov chain, a construction that allows us to see that  a stationary process admits regular conditional expectations with respect to ``the past''.

\medskip

\begin{exam}[Stationary Sequences as Functions of Markov Chains]
\label{staseqasfunmarcha}
Under the setting introduced in Example \ref{examarproexa},  consider now  the $\mathcal{F}/\mathcal{F}^{-}$ measurable function $\xi_{0}=\pi^{-}$ and the $\mathcal{F}^{-}/\mathcal{S}$ measurable function $x_{0}^{-}$: the restriction of $x_{0}$ to $\Omega^{-}$. Note that, if for every $k\in\mathbb{Z}$, $\xi_{k}:=T^{k}\xi_{0}$ then, since $\sigma(\xi_{0})=\sigma(x_{j})_{j\leq 0}=\mathcal{F}_{0}$,  we have that for every $k\in\mathbb{Z}$, $\sigma(\xi_{k})=:\mathcal{F}_{k}$ and therefore, since $(\mathcal{F}_{k})_{k\in\mathbb{Z}}$ is increasing,  $\sigma((\xi_{j})_{j\leq k})=\sigma(\xi_{k})$. 

In particular, if $\mathbb{P}$ is {any} probability measure in $(\Omega,\mathcal{F})$, then for any $k\in\mathbb{Z}$ and any $\mathcal{F}^{-}/\mathcal{C}$ measurable function $f:\Omega^{-}\to \mathbb{C}$
\begin{equation}
\label{chimarcha}
\mathbb{P}[f(\xi_{k+1})|\sigma(\xi_{j})_{j\leq k}]=\mathbb{P}[f(\xi_{k+1})|\sigma(\xi_{k})],
\end{equation} 
provided that $f(\xi_{k+1})\in L^{1}_{\mathbb{P}}(\mathcal{F})$. By taking $f=I_{A}$  for any given $A\in \mathcal{F}^{-}$ we see, by an application of Proposition \ref{proeximarker}, that (\ref{chimarcha}) implies that $(\xi_{k})_{k\in\mathbb{Z}}$ is a Markov chain\footnote{The state space $(\Omega^{-},\mathcal{F}^{-})$ is generated by a complete and separable metric space by the standard fact that the countable product of such spaces can be metrized in such a way that it has those two properties.} under (any) $\mathbb{P}$. If $(\xi_{k})_{k\in\mathbb{Z}}$ is stationary (under $\mathbb{P}$) it is also homogeneous (Remark \ref{staimphom}). 

Let now $(X'_{k})_{k\in\mathbb{Z}}$ be a sequence of random elements in $S$ defined on a probability space $(\Omega',\mathcal{F}',\mathbb{P}')$. Assume that $(X'_{k})_{k\in\mathbb{Z}}$ is stationary, so that the probability measure on $\mathcal{S}^{\mathbb{Z}}$ specified by
$$\mathbb{P}_{\mathbb{Z}}((x_{n},\cdots,x_{n+k})\in A):=\mathbb{P}'((X'_{n},\cdots,X'_{n+k})\in A)$$
for every $A\in \mathcal{S}^{k+1}$ makes $(x_{k})_{k\in\mathbb{Z}}$ a copy (in distribution) of $(X_{k})_{k\in\mathbb{Z}}$. Under $\mathbb{P}_{\mathbb{Z}}$, the Markov chain $(\xi_{k})_{k\in\mathbb{Z}}$ is stationary. 

If we apply the previous observations to $\mathbb{P}_{\mathbb{Z}}$, and consider $f:=x_{0}^{-}$ we get that, for every $k$, $f(T^{k}\xi_{0})=f(\xi_{k})=x_{k}$, and therefore $x_{k}=f(\xi_{k})$ is a {\it function} of the (stationary) Markov chain $(\xi_{k})_{k\in\mathbb{Z}}$. Since the finite dimensional distributions of $(x_{k})_{k\in\mathbb{Z}}$ (under $\mathbb{P}_{\mathbb{Z}}$) are the same as those of $X_{k}'$ (under $\mathbb{P}'$), we see that {\it every stationary process in a complete and separable metric space is equivalent (in distribution) to a function of a Markov chain}. Under this equivalence, the ``past'' sigma algebra is regular: there exists a family of probability measures $\{\mathbb{P}_{\omega}\}_{\omega\in\Omega}$  such that, under $\mathbb{P}_{\mathbb{Z}}$
$$E[x_{k}|\sigma(x_{j})_{j\leq 0}](\omega)=E[f(\xi_{k})|\sigma(\xi_{0})](\omega)=\int_{\Omega} f(z)d\mathbb{P}_{\omega}(z).$$
As a matter of fact, the analysis in Example \ref{examarproexa} shows that 
$$E[f(\xi_{k})|\sigma(\xi_{0})](\omega)=P_{0}^{k}f (\xi_{0}(\omega)),$$
where $P_{0}^{k}f$ is given by (\ref{maropefunequ}) (considering $f$ as constant in $(z_{1},\cdots,z_{k-1})$) via the transition probability matrix guaranteed for $(\xi_{k})_{k\in\mathbb{Z}}$ by Proposition \ref{proeximarker} and Remark \ref{staimphom}.
\end{exam}

\medskip

\begin{remark}[Nonstationary Case]
\label{uniregnat}
An analysis of the arguments  in Example \ref{examarproexa} shows that the assumption of stationarity is superfluous in the following sense: the role of $\mathbb{P}$ (the invariant measure of $P$) in the construction carried along that example is to guarantee that the coordinate functions {\it indeed} define a stationary process {\it and} that the compatibility conditions of Kolmogorov's existence theorem hold. We can drop the requirement of stationarity and still carry on the given construction if we start from ``marginal'' probability measures $\{\mathbb{P}_{n}\}_{n\in\mathbb{Z}}$ on $\mathcal{S}$ (representing the distribution of $\xi_{n}$), a family of transition measures $\{P_{n}\}_{n\in\mathbb{Z}}$ on $(S,\mathcal{S})$ (representing the transitions $\mathbb{P}(\xi_{n+1}\in A|\sigma(\xi_{n}))$) and if, following the arguments  in Example \ref{examarproexa}, (\ref{comequ}) holds for every $(n,k)\in\mathbb{Z}\times\mathbb{N}$ if we replace $\mathbb{P}$ by $\mathbb{P}_{n+k}$. 

This is the case if $(X'_{k})_{k\in\mathbb{Z}}$ is {\it any} sequence of random  elements in a complete and separable metric space $S$ and $\mathbb{P}_{n}$ is the law of $\xi_{n}=(\cdots,X'_{n-1}, X'_{n})$ in $(\Omega^{-},\mathcal{F}^{-})$. By following the construction along Example \ref{staseqasfunmarcha}, this gives a representation of any sequence of random elements on $(S,\mathcal{S})$ as a sequence of functions of a (not necessarily stationary) Markov chain. 
\end{remark}
\medskip

 
\section{Regularity and Quenched Convergence}
\label{regandquecon} 
The natural question at this point is the following: suppose, in the context of Definition \ref{regconexpdef}, that $E_{0}$ is regular, and assume that $X_{n}$ converges in the quenched sense to $X$ as $n\to\infty$. Can we say anything about the convergence of $(X_{n})_{n\in\mathbb{N}}$ with respect to the measures in the decomposition of $E_{0}$? The following proposition provides an answer sufficiently good for our purposes. 
\medskip

\begin{prop}[Regularity and Quenched Convergence]
\label{prounicon}
In the context of Definition \ref{defquecon}, assume that $(S,d)$ is  separable. If $E_{0}$ is regular and $Y_{n}$ converges to $Y$ in the quenched sense,  there exists a set $\Omega_{0}\subset \Omega$ with $\mathbb{P}\Omega_{0}=1$ such that for all $f:S\to\mathbb{R}$ continuous and bounded and all $\omega\in \Omega_{0}$
\begin{equation}
\label{defunicon}
\int_{\Omega}f\circ Y_{n}(z)d\mathbb{P}_{\omega}(z)\to_{n} \int_{\Omega'}f\circ Y(z)d\mathbb{P}'(z).
\end{equation} 
In particular, $Y_{n}$ converges to  $Y$ in the quenched sense if and only if for $\mathbb{P}-$a.e $\omega$, $Y_{n}\Rightarrow Y$ with respect to $\mathbb{P}_{\omega}$.
\end{prop}
 
Notice that if $\{\mathbb{P}_{\omega}\}_{\omega\in\Omega}$ is a decomposition of $E_{0}$ then, by the definition of quenched convergence, and denoting again by $E^{\omega}$ the integration with respect to $\mathbb{P}_{\omega}$,
$$E^{\omega}f(Y_{n})\to Ef(Y)$$
as $n\to\infty$ for every $\omega\in\Omega_{f}$, where $\mathbb{P}\Omega_{f}=1$. Proposition \ref{prounicon} states that if $(S,d)$ is separable, $\Omega_{f}$ can be chosen {\it independent} of $f$, namely $\Omega_{f}:=\Omega_{0}$ for all $f\in \mathbf{C}^{b}(S)$. The set $\Omega_{0}$ depends, nonetheless, on $(Y_{n})_{n}$.
 
 {\bf Proof of Proposition \ref{prounicon}:} Consider functions $U_{k,\epsilon}$ as in the statement 2. of Theorem \ref{refporthesta}. As remarked in the paragraph above there exists, for all $k\in \mathbb{N}$ and $\epsilon>0$ ($\epsilon\in \mathbb{Q}$), a set $\Omega_{k,\epsilon}\subset \Omega$ with $\mathbb{P}\Omega_{k,\epsilon}=1$ such that for all $\omega\in \Omega_{k,\epsilon}$, 
 $E^{\omega}U_{k,\epsilon}(Y_{n})\to EU_{k,\epsilon}(Y)$ as ${n\to\infty}$. Now take $\Omega_{0}:=\bigcap_{k,\epsilon} \Omega_{k,\epsilon}$ and use Theorem \ref{refporthesta}.\qed
 

\medskip


\section{Product Spaces and Regularity}
\label{queconprospa}

We finish our discussion about regular conditional expectations with the following result, showing that the notion of regular conditional expectation behaves well under the product of probability spaces.

\medskip

\begin{prop}[Product Spaces and Regularity]
\label{proregpro}
Let $(\Theta,\mathcal{B}, \lambda)$ and $(\Omega,\mathcal{F},\mathbb{P})$ be probability spaces, and let $\mathcal{B}_{0}\subset \mathcal{B}$, $\mathcal{F}_{0}\subset \mathcal{F}$ be sub-sigma algebras such that $E[\,\cdot\,|\mathcal{B}_{0}]$ and $E[\,\cdot\,|\mathcal{F}_{0}]$ are regular (Definition \ref{regconexpdef}). If  $\{\lambda_{\theta}\}_{\theta\in\Theta}$ and $\{\mathbb{P}_{\omega}\}_{\omega\in \Omega}$ are, respectively, decompositions of $E[\,\cdot\,|\mathcal{B}_{0}]$ and $E[\,\cdot\,|\mathcal{F}_{0}]$, then the conditional expectation $E[\,\cdot\,|\mathcal{B}_{0}\otimes\mathcal{F}_{0}]$ with respect to $\mathcal{B}_{0}\otimes\mathcal{F}_{0}$ and $\lambda\times\mathbb{P}$ admits the decomposition  $\{\lambda_{\theta}\times\mathbb{P}_{\omega}\}_{(\theta,\omega)\in \Theta\times \Omega}$.
\end{prop}

{\bf Proof:} We will proceed in two steps.

{\itshape\bfseries Step 1. }{\it  Assume that $\mathcal{B}_{0}=\mathcal{B}$}. In this case we will prove that for any $E\in\mathcal{B}\otimes\mathcal{F}$ the function
\begin{equation}
\label{regconexpprocas1}
\tilde{I}_{E}(\theta,\omega)= \int_{\Omega}I_{E}(\theta,z)\,d\mathbb{P}_{\omega}(z)
\end{equation} 
(which is well defined for every $\theta$ by Theorem 18.1 in \cite{bilpromea}) defines a version of $\mathbb{P}[{E}|\mathcal{B}\otimes\mathcal{F}_{0}]$. Note that by Proposition \ref{proredregcondis} this proves also the desired conclusion for any $f\in L^{1}_{\lambda\times\mathbb{P}}$, and that in the special case in which for every $\theta\in\Theta$, $\{\theta\}\in\mathcal{B}$ (and therefore  $\lambda_{\theta}=\delta_{\theta}$, the Dirac measure at $\theta$, defines a decomposition for $E[\,\cdot\,|\mathcal{B}]=Id$, the identity map on $L^{1}_{\lambda}$) there is consistency with the given conclusion.


To prove that (\ref{regconexpprocas1}) defines a $\mathcal{B}\otimes\mathcal{F}_{0}-$measurable function, note first that if  $E=A\times B$ is a rectangular set, then (\ref{regconexpprocas1}) is equal to the function
$$(\theta,\omega)\mapsto I_{A}(\theta)\mathbb{P}_{\omega}(B),$$
which is clearly $\mathcal{B}\otimes\mathcal{F}_{0}$ measurable. 

Now consider the family ${\mathcal{G}}$ of sets $E \in \mathcal{B}\otimes\mathcal{F}$ such that (\ref{regconexpprocas1}) is $\mathcal{B}\otimes\mathcal{F}_{0}$-measurable. Since for any family  $\{E_{n}\}_{n}{\subset}{\mathcal{G}}$ of mutually disjoint sets the choice $E=\cup_{n}E_{n}$ gives that 

$$\tilde{I}_{E}(\theta,\omega)=\sum_{n}\tilde{I}_{E_{n}}(\theta,\omega)$$
(apply the monotone convergence theorem) and ${\mathcal{G}}$ includes the set $\Theta\times\Omega$, ${\mathcal{G}}$ is a $\lambda-$system. Since ${\mathcal{G}}$ includes the finite unions of disjoint rectangles it follows, by the $\pi-\lambda$ theorem, that ${\mathcal{G}}=\mathcal{B}\otimes\mathcal{F}$. This proves the $\mathcal{B}\otimes\mathcal{F}_{0}-$measurability of (\ref{regconexpprocas1}) for every $E\in \mathcal{B}\otimes\mathcal{F}$. 

Now, by Fubini's theorem and the definition of $\{\mathbb{P}_{\omega}\}_{\omega\in\Omega}$, given any rectangular set $E'=A'\times B'\in\mathcal{B}\otimes{\mathcal{F}}_{0}$
$$\int_{\Theta\times\Omega}\tilde{I}_{E}(\theta,\omega)I_{E'}(\theta,\omega)d(\lambda\times \mathbb{P})(\theta,\omega)=\int_{A'}\int_{B'}\tilde{I}_{E}(\theta,\omega)d\mathbb{P}(\omega)d\lambda(\theta)=$$
$$\int_{A'}\int_{B'}{I}_{E}(\theta,\omega)d\mathbb{P}(\omega)d\lambda(\theta)=(\lambda\times\mathbb{P})(E\cap E'),$$
and a further application of the $\pi-\lambda$ theorem shows that the equality between the extremes holds for any $E'\in\mathcal{B}\otimes\mathcal{F}_{0}$, which proves that (\ref{regconexpprocas1}) is indeed a version of $\mathbb{P}[E|\mathcal{B}\otimes\mathcal{F}_{0}]$.

{\bfseries\itshape Step 2.} {\it General Case.} For the general case note first the following: by the case treated in the previous step and Proposition \ref{proredregcondis},  given any $f\in L^{1}_{\lambda\times\mathbb{P}}$ the function 
$$\tilde{f}(\theta,\omega)=\int_{\Omega}f(\theta,z)\,d\mathbb{P}_{\omega}(z)$$
is a version of $E[f|\mathcal{B}\otimes\mathcal{F}_{0}]$. Also
$$E[f\,|\mathcal{B}_{0}\otimes\mathcal{F}_{0}]=E[E[f\,|\mathcal{B}\otimes\mathcal{F}_{0}]|\mathcal{B}_{0}\otimes\mathcal{F}_{0}]=E[\tilde{f}\,|\mathcal{B}_{0}\otimes\mathcal{F}_{0}]$$
 
 $(\lambda\times\mathbb{P})-$a.s. It follows by a second application of the Step 1 and Proposition \ref{proredregcondis} with $\mathcal{B}\otimes\mathcal{F}_{0}$ in the role of $\mathcal{B}\otimes\mathcal{F}$ and with $\mathcal{B}_{0}$ in the role of $\mathcal{F}_{0}$, that
$$(\theta,\omega)\mapsto \int_{\Theta}\int_{\Omega}f(x,z) d\mathbb{P}_{\omega}(z)\,d\lambda_{\theta}(x)$$
defines a version of $E[f|\mathcal{B}_{0}\otimes\mathcal{F}_{0}]$.\qed

\chapter{Quenched Asymptotics of Normalized Fourier Averages}
\label{resandcom}
In this chapter we will introduce the results on asymptotic distributions to be proved along this monograph. The main results are theorems \ref{quecltfoutra}, \ref{nonqueconthe}, \ref{invpriavefre} and \ref{barinvprihancon} (theorems \ref{queconznundratdec} and \ref{queconznmaxwoo} can be seen as versions of the previous ones refined by the introduction of additional structure).

For some of the results, including the main ones, we will limit our discussion to the presentation of the statements and to the comments necessary to clarify their meaning. We will nonetheless provide proofs of some of the corollaries and ``secondary'' results whenever they can be reached in a straightforward manner from the discussions already made.

This chapter is organized as follows: in Section \ref{secquecltfoutra} we present the {\it Central Limit Theorem for Fourier transforms at (a.e-)fixed frequencies} (Theorem \ref{quecltfoutra}), which extends to the quenched setting Theorem \ref{cltpelwu} and opens the door to several questions regarding the validity of this quenched convergence in stronger forms.

Then, in Section \ref{necrancen}, we address the first issue in the direction of these questions: the necessity of the ``random'' centering of the normalized ergodic averages in order to guarantee the conclusion of Theorem \ref{quecltfoutra}. We will state a general result (Theorem \ref{nonqueconthe}) showing that this is indeed a necessary condition, but we will still address, in Section \ref{secnonrancennec}, particular cases in which this normalization is irrelevant.

Finally, in Section \ref{secquefclt}, we will address the problem of extending the quenched central limit theorems under consideration to corresponding quenched invariance principles. We will state a result (Theorem \ref{invpriavefre}) showing that this is indeed possible in the sense of {\it averaged frequencies}. The (stronger) version for fixed frequencies remains open, but a special case (Theorem \ref{barinvprihancon}), and some of its consequences, are discussed in Section \ref{queinvpriavefre}.

Some sections have a part dedicated to ``general comments''. The purpose of these discussions is to clarify the meaning of the results previously given, to describe some of the relations between them, and to motivate the discussions that follow both in the corresponding as in further sections.

\section{The quenched CLT for Fourier Transforms} 
\label{secquecltfoutra}

The purpose of this section is to present the most general version of the quenched central limit theorem for Fourier Transforms availabe in this monograph. The result is the following.  

\medskip

\begin{thm}[The Quenched Central Limit Theorem for Fourier Transforms]
\label{quecltfoutra}
Let $(X_{k})_{k\in\mathbb{Z}}=(T^{k}X_{0})_{k\in\mathbb{Z}}$ be a square-integrable ergodic process (Definition \ref{staproL2}) adapted to an increasing $T-$filtration $(\mathcal{F}_{k})_{k\in\mathbb{Z}}$ (Definition \ref{defadafil}). Assume that $\mathcal{F}_{\infty}$ (Definition \ref{taisigalgdef}) is countably generated (Definition \ref{deffcougen}), that $E_{0}=E[\,\cdot\,|\mathcal{F}_{0}]$ is regular (Definition \ref{regconexpdef}), denote by $S_{n}(\theta)$ the $n-$th discrete Fourier Transform of $(X_{k})_{k\in\mathbb{Z}}$ (Definition \ref{defdisfoutra}) and let 
\begin{equation}
\label{norrancenavedisfoutraequ}
Y_{n}(\theta):=\frac{1}{\sqrt{n}}(S_{n}(\theta)-E_{0}S_{n}(\theta)).
\end{equation}
Then there exist $I\subset[0,2\pi)$ with $\lambda(I)=1$ such that the following holds:
\begin{enumerate}
\item For every $\theta\in I$, there exists a nonnegative number $\sigma(\theta)$ such that 
\begin{equation}
\label{sigtheequ}
\sigma^{2}(\theta)=\lim_{n}E_{0}|Y_{n}(\theta)|^{2},\, \mbox{\,\,\,\,\,\,\,\it $\mathbb{P}-$a.s. and in $L^{1}_{\mathbb{P}}$.}
\end{equation}

\item If $N_{1}, N_{2}$ denote independent standard normal random variables and $i:=\sqrt{-1}$, then for every $\theta\in I$, the process $Y_{n}(\theta)$ converges in the quenched sense (Definition \ref{defquecon}) with respect to $\mathcal{F}_{0}$ to
\begin{equation}
\label{limdisquecltequ}
Y(\theta)=\frac{\sigma(\theta)}{\sqrt{2}}(N_{1}+iN_{2}),
\end{equation}
(or, what is the same, $Y_{n}(\theta)$ convergence in the quenched sense to a bivariate normal, centered variable with covariance matrix (\ref{covmatasynor})).
\end{enumerate}

In addition, $\theta\mapsto \sigma^{2}(\theta)$ is the spectral density (Definition \ref{defspeden}) of the process $(X_{k}-E_{-\infty}X_{k})_{k\in\mathbb{Z}}$, where $E_{-\infty}$ denotes the conditional expectation $E[\,\cdot\,|\mathcal{F}_{-\infty}]$ with respect to  $\mathcal{F}_{-\infty}$. 
\end{thm}

\medskip

Before moving on to further comments, let us state the following Corollary, whose proof is given in full detail to facilitate further discussions.

\medskip

\begin{cor}
\label{corquecltprospa}
In the context of Theorem \ref{quecltfoutra}, and denoting by  $Y_{n}:[0,2\pi)\times\Omega\to \mathbb{C}$ and  
$Y:[0,2\pi)\times\Omega'\to \mathbb{C}$ the functions defined respectively by $Y_{n}(\theta,\omega)=Y_{n}(\theta)(\omega)$ and $Y(\theta,\omega')=Y(\theta)(\omega')$, there exists a set $\Omega_{0}\subset \Omega$ with $\mathbb{P}\Omega_{0}=1$ such that for every $\omega\in\Omega_{0}$, $Y_{n}\Rightarrow Y$ under $\lambda\times\mathbb{P}_{\omega}$.
\end{cor}

{\bf Proof:} First note that if $\mathcal{B}_{0}=\{\emptyset,[0,2\pi)\}$ is the trivial sigma-algebra and $\lambda_{\theta}:=\lambda$ for all $\theta\in [0,2\pi)$ then, by Proposition \ref{proregpro}, the family of measures
$$\{\lambda_{\theta}\times\mathbb{P}_{\omega}\}_{\omega\in\Omega}$$
is a decomposition of $E[\,\cdot\,|\mathcal{B}_{0}\otimes\mathcal{F}_{0}]$.

Now, the functions $Y_{n}, Y$ given in the statement of Corollary \ref{corquecltprospa} are clearly measurable with respect to the respective product sigma algebras and therefore, in virtue of Proposition \ref{proregpro}, we can read the statement of Theorem \ref{quecltfoutra} in the following way: {\it for any continuous and bounded function $f:\mathbb{C}\to\mathbb{R}$}
\begin{equation}
\label{queconasconprospa}
E[f\circ Y_{n}|\mathcal{B}\otimes\mathcal{F}_{0}]\to_{n} E[f\circ Y(\theta,\cdot)]\mbox{\,\,\,\,\,\,\,\,\, \it $\lambda\times\mathbb{P}-$a.s.}
\end{equation}
Let us explain this in detail: note that by an application of Proposition \ref{proregpro}, the function at the right-hand side in (\ref{queconasconprospa}) is a version of the conditional expectation of $f\circ Y$ with respect to the sigma field $\mathcal{B}\otimes \{\emptyset,\mathcal{F}'\}$. Therefore this function is $\mathcal{B}\otimes\mathcal{F}'$-measurable and, since it is constant over $\Omega'$ for $\theta$ fixed, it is $\mathcal{B}-$measurable. By regarding it as constant on $\Omega$ for $\theta$ fixed, it can be considered $\mathcal{B}\otimes\mathcal{F}-$measurable. This shows that the set where the convergence in (\ref{queconasconprospa}) occurs belongs to $\mathcal{B}\otimes\mathcal{F}$. 

Now,  a further application of Proposition \ref{proregpro} shows that  
$$(\theta,\omega)\mapsto \int_{\Omega}f\circ Y_{n}(\theta,z)\,d\mathbb{P}_{\omega}(z)$$ defines a version of $E[f\circ Y_{n}|\mathcal{B}\otimes\mathcal{F}_{0}]$, and since for $\lambda-$a.e. fixed $\theta$,
$$\int_{\Omega}f\circ Y_{n}(\theta,z)\,d\mathbb{P}_{\omega}(z)\to E[f\circ Y(\theta,\cdot)],$$
$\mathbb{P}-$a.s., we deduce that the set where the convergence in (\ref{queconasconprospa}) occurs has, indeed, product measure one.

It follows from Proposition \ref{inhquecon} that $Y_{n}$ converges to  $Y$ in the quenched sense with respect to $\mathcal{B}_{0}\otimes\mathcal{F}_{0}$. The conclusion follows at once from the observation at the beginning of this proof and Proposition \ref{prounicon}.\qed

\subsection*{General Comments}

Note that the convergence in (\ref{queconasconprospa}) resembles the convergence that follows from Theorem \ref{pelwufclt} by evaluating the corresponding random functions at $t=1$. As we shall see,   Corollary \ref{corquecltprospa} can indeed be extended to a quenched invariance principle without imposing any further hypothesis to the processes under consideration (see Theorem \ref{invpriavefre} below). At the moment of writing this monograph this is not the case  for Theorem \ref{quecltfoutra}, whose extension to an invariance principle will be possible for us only at the expense of further assumptions.

With regards to the statement of Theorem \ref{quecltfoutra}, 
the following comments are worth at this point.

\begin{enumerate}
\item First, note that Theorem \ref{quecltfoutra} is apparently a re-statement of Theorem \ref{cltpelwu}: it basically emerges  from that result by replacing  ``$X_{k}$'' by ``$X_{k}-E_{0}X_{k}$'' and  ``convergence in distribution'' by ``quenched convergence''. Note nevertheless that the process $(X_{k}-E_{0}X_{k})_{k\in\mathbb{Z}}$ is generally non-stationary (all its entries are zero for $k\leq 0$), and therefore that substitution brings us outside of the hypotheses of Theorem \ref{cltpelwu}. 

\item Another one of the hypotheses of Theorem \ref{cltpelwu} is missing from the statement of Theorem \ref{quecltfoutra}: the regularity of $(X_{k})_{k\in\mathbb{Z}}$ (Definition \ref{regcondef}), but this actually ``can be obtained''  from the theory already developed via a simple substitution, as we proceed now to explain.

First, the process $$(X_{-\infty,k})_{k\in\mathbb{Z}}:=(X_{k}-E_{-\infty}X_{k})_{k\in\mathbb{Z}}$$ 
is stationary and regular. Indeed: $(X_{-\infty,k})_{k\in\mathbb{Z}}$ is stationary by (\ref{equinvprotai}), and an application of Proposition \ref{regcon} (which is actually implicit in the statement of Theorem \ref{quecltfoutra}) shows that it is regular.

Now, since $E_{0}E_{-\infty}=E_{-\infty}$, a simple computation shows that (see the notation in Definition \ref{defdisfoutra}) 
\begin{equation}
\label{equregnotneed}
S_{n}((X_{k})_{k},\theta,\cdot)-E_{0}S_{n}((X_{k})_{k},\theta,\cdot)=S_{n}((X_{-\infty,k})_{k},\theta,\cdot)-E_{0}S_{n}((X_{-\infty,k})_{k},\theta,\cdot)
\end{equation}
and therefore we can study the asymptotics of $Y_{n}(\theta)$ assuming, via the substitution of $X_{k}$ by $X_{k}-E_{-\infty}X_{k}$ for all $k\in\mathbb{Z}$, that $(X_{k})_{k\in\mathbb{Z}}$ is stationary, centered, and regular.
\item Now consider the following observation: in the context of Theorem \ref{quecltfoutra}, the process
$$Z_{n}(\theta):=\frac{1}{\sqrt{n}}S_{n}(\theta)$$
satisfies 
\begin{equation}
\label{equrelznyn}
Z_{n}(\theta)=Y_{n}(\theta)+\frac{E_{0}S_{n}(\theta)}{\sqrt{n}},
\end{equation}
and since $Y_{n}(\theta)$ converges in the quenched sense, and therefore in distribution to (\ref{limdisquecltequ}), we have the following corollary.
\end{enumerate}

\medskip

\begin{cor}
\label{corobsconpro}
Under the hypotheses of Theorem \ref{quecltfoutra}, the conclusion of Theorem \ref{cltpelwu} remains true (without necessarily assuming the regularity of $(X_{k})_{k\in\mathbb{Z}}$) if for $\lambda-$a.e $\theta$, 
\begin{equation}
\label{concontozerproequ}
\frac{E_{0}S_{n}(\theta)}{\sqrt{n}}\Rightarrow_{n} 0,
\end{equation} 
in which case $\theta\mapsto \sigma^{2}(\theta)$ is the spectral density of the (regular) process $(X_{k}-E_{-\infty}X_{k})_{k\in\mathbb{Z}}$ (where ``$E_{-\infty}$'' is as in the last statement of Theorem \ref{quecltfoutra}). In particular,  the conclusion of Theorem \ref{cltpelwu} follows under the hypotheses of Theorem \ref{quecltfoutra} if $(X_{k})_{k\in\mathbb{Z}}$ is regular.
\end{cor}

{\bf Proof:} Only the last statement requires a proof. To do so we will prove that, under the hypothesis of Theorem \ref{quecltfoutra}, the hypothesis of regularity in Theorem \ref{cltpelwu} imply the fulfillment of (\ref{concontozerproequ}).

Using the notation introduced in theorems \ref{cltpelwu} and \ref{quecltfoutra} we have by orthogonality that, for every $\theta\in I$ 
$$E|Z_{n}(\theta)|^{2}=E|Z_{n}(\theta)-E_{0}Z_{n}(\theta)|^{2}+E|E_{0}Z_{n}(\theta)|^{2}=E[E_{0}|Y_{n}(\theta)|^{2}]+E|E_{0}Z_{n}(\theta)|^{2}$$
and it follows from Theorem \ref{spedenasalim}, Theorem \ref{quecltfoutra}, and Fatou's lemma 
that for $\lambda-$a.e $\theta$
$$\sigma^{2}(\theta)=\limsup_{n}E|Z_{n}(\theta)|^{2}\geq  
\liminf_{n}E[E_{0}|Y_{n}(\theta)|^{2}]+\limsup_{n}E|E_{0}Z_{n}(\theta)|^{2}\geq $$
$$ E[\liminf_{n}E_{0}|Y_{n}(\theta)|^{2}]+\limsup_{n}E|E_{0}Z_{n}(\theta)|^{2}=
\sigma^{2}(\theta)+\limsup_{n}E|E_{0}Z_{n}(\theta)|^{2},$$
which implies that $E_{0}Z_{n}(\theta)$ converges to zero in $L^{2}_{\mathbb{P}}$. This clearly implies (\ref{concontozerproequ}).\qed

\medskip

\section{The Random Centering}
\label{necrancen}

This leaves us with a question about the ``missing'' element  on the statement in Theorem \ref{cltpelwu}: the random centering ``$-E_{0}S_{n}(\theta)$'' in the definition of $Y_{n}(\theta)$. 

More precisely, consider the following observations: every process satisfying the hypotheses of Theorem \ref{cltpelwu} satisfies the hypotheses of Theorem \ref{quecltfoutra}, and by the arguments following the statement of Theorem \ref{quecltfoutra}, the processes involved in the statement of Theorem \ref{quecltfoutra} can be assumed to satisfy the hypotheses of Theorem \ref{cltpelwu}. 

Even more, in Corollary \ref{corobsconpro} we obtained the convergence in distribution of the normalized discrete Fourier transforms
\begin{equation}
\label{noravedisfoutraequ}
Z_{n}(\theta):= \frac{1}{\sqrt{n}}S_{n}(\theta)
\end{equation}
by using the convergence in distribution of $Y_{n}(\theta)$ and the ``{ad hoc}'' hypothesis for the remainder, but the following question is still to be addressed.

{\bf Question:} {\it can we actually prove that $(Z_{n}(\theta))_{n\in\mathbb{N}}$ converges in the quenched sense under the hypotheses of Theorem \ref{cltpelwu}?}

\subsection{Necessity of the Random Centering}
\label{secnecrancen}
To begin the discussion regarding the question above note that by (\ref{equrelznyn}), and since $Y_{n}(\theta)$ admits the same quenched limit as the limit (in distribution) of $Z_{n}(\theta)$,  the ``perturbation'' to quenched convergence, if any, is due to the behavior of ${E_{0}S_{n}(\theta)}/{\sqrt{n}}$ under $\mathbb{P}_{\omega}$. 

This can actually be described in a very precise way, as stated by the following theorem. 

\medskip

\begin{thm}[Possible Quenched Limits for the Non-centered Normalized Averages]
\label{quecondes}
In the context of Theorem \ref{quecltfoutra}, given $\theta\in I$ and denoting by $Z_{n}(\theta)$ a (fixed) version of the random variable in (\ref{noravedisfoutraequ}) ($n\in\mathbb{N}$) and by $E_{0}Z_{n}(\theta)$ a (fixed)  version of $E[Z_{n}(\theta)|\mathcal{F}_{0}]$, there exists $\Omega_{\theta}\subset \Omega$ with $\mathbb{P}\Omega_{\theta}=1$ such that, for $\omega\in \Omega_{\theta}$ the following are equivalent
\begin{enumerate}
\item $Z_{n}(\theta)$ is convergent in distribution under $\mathbb{P}_{\omega}$.
\item There exists
\begin{equation}
\label{limdetconznpome}
L_{\theta}(\omega)=\lim_{n}{E_{0}[Z_{n}(\theta)]}(\omega),
\end{equation}
and $Z_{n}(\theta)\Rightarrow Y(\theta)+L_{\theta}(\omega)$ under $\mathbb{P}_{\omega}$.
\end{enumerate}
\end{thm}

The proof of this theorem is deferred to Section \ref{secprothequecondes}, but we will use it at this point to prove the following corollary.

\medskip

\begin{cor}
\label{sufnecconqueconzn}
In the context of Theorem \ref{quecltfoutra}, denoting by $Z_{n}(\theta)$ the random variable (\ref{noravedisfoutraequ}), and assuming that $(X_{k})_{k\in\mathbb{Z}}$ is regular (Definition \ref{regcondef}), the following are equivalent for $\theta\in I$.
\begin{enumerate}
\item $Z_{n}(\theta)$ converges in the quenched sense as $n\to\infty$.
\item $E_{0}Z_{n}(\theta)\to_{n} 0$, $\mathbb{P}-$a.s.,
\end{enumerate}
in which case the (quenched) limit of $Z_{n}(\theta)$ is $Y(\theta)$.
\end{cor}

{\bf Proof:}  Fix $\theta\in I$. Since, by Theorem \ref{cltpelwu}, $Z_{n}(\theta)\Rightarrow Y(\theta)$ (under $\mathbb{P}$), the only possible quenched limit of $Z_{n}(\theta)$ is certainly $Y(\theta)$ (see the paragraph following Remark \ref{remdefquecon}). 

Now, by Proposition \ref{prounicon}, the quenched convergence of $Z_{n}(\theta)$ to $Y(\theta)$ is equivalent to the following: there exists a set $\Omega_{\theta,1}\subset \Omega$ with $\mathbb{P}\Omega_{\theta,1}=1$ such that for every $\omega\in \Omega_{\theta,1}$
$$Z_{n}(\theta)\Rightarrow Y(\theta)$$
under $\mathbb{P}_{\omega}$ as $n\to\infty$.

Now note that
\begin{equation}
\label{decznthe}
Z_{n}(\theta)=Y_{n}(\theta)+E_{0}Z_{n}(\theta).
\end{equation}

By Proposition \ref{prounicon} and Theorem \ref{quecltfoutra}, there exists $\Omega_{\theta,2}$ with $\mathbb{P}\Omega_{\theta,2}=1$ such that for every $\omega\in \Omega_{\theta,2}$
$$Y_{n}(\theta)\Rightarrow Y(\theta) $$
under $\mathbb{P}_{\omega}$ as $n\to\infty$. The conclusion follows considering $$\omega\in \bigcap_{k=0}^{2}\Omega_{\theta,k}$$
where $\Omega_{\theta,0}$ is the set specified in Theorem \ref{quecondes} and applying  Proposition \ref{procontypcom} (use the complex version of Proposition \ref{proquecon} (respectively, Proposition \ref{degcascon}) when $\sigma(\theta)>0$ (respectively, when $\sigma(\theta)=0$)).\qed




We return to the question above, that about the quenched convergence (in general) of $Z_{n}(\theta)$ for $\theta\in I$. The actual answer is {\it no}, as our next main result  shows.

\medskip

\begin{thm}[An Example of non-Quenched Convergence]
\label{nonqueconthe}
There exist $\mathcal{F}$, $\mathcal{F}_{0}$, $T$, and $(X_{k})_{{k}\in\mathbb{Z}}$ as in the hypotheses of Theorem \ref{cltpelwu} such that $E_{0}:=E[\,\cdot\,|\mathcal{F}_{0}]$ is regular and for any  decomposition $\{\mathbb{P}_{\omega}\}_{\omega\in \Omega}$ of $E_{0}$ (Definition \ref{regconexpdef}) 
$$Z_{n}(\theta)= \frac{1}{\sqrt{n}}S_{n}(\theta)$$
admits no limit in distribution under $\mathbb{P}_{\omega}$ for every $\theta\in[0,2\pi)$ and $\mathbb{P}-$a.e  $\omega$.
\end{thm}

\subsection*{General Comments}
With regards to the results in this section it is important to observe the following: for  the process $(X_{n})_{n\in\mathbb{Z}}$ to be constructed along the proof of Theorem \ref{nonqueconthe}, if $Y_{n}(\theta)$ is given by (\ref{norrancenavedisfoutraequ}), $Y(\theta)$ is given by (\ref{limdisquecltequ}) and $Z_{n}(\theta)$ is given by (\ref{noravedisfoutraequ}), then  certainly
$$Z_{n}(\theta)\Rightarrow{Y}(\theta)$$
as $n\to \infty$ for $\lambda-$almost every $\theta$. Theorem \ref{nonqueconthe} not only states that this convergence is {\it not} quenched, but it states that $Z_{n}(\theta)$ cannot converge when started at $\mathbb{P}-$a.e $\omega$, this is,  $Z_{n}(\theta)$ does not admit a limit (in distribution) under $\mathbb{P}_{\omega}$ for $\mathbb{P}-$a.e $\omega$. As a matter of fact, we will see that for this process
\begin{equation}
\label{plimsupinfone}
\mathbb{P}[\limsup_{n}|E_{0}Z_{n}(\theta)|=\infty]=1,
\end{equation}
which makes impossible the convergence under $\mathbb{P}_{\omega}$ for $\mathbb{P}-$almost every $\omega$ in virtue of Theorem \ref{quecondes}.

This enforces the intuitive idea that $\mathcal{F}_{0}$ represents the ``deterministic part'' of the processes in question. Note again that, even if we can prove that (\ref{limdetconznpome}) exists for $\mathbb{P}-$a.e $\omega$, we {\it cannot} a priori conclude that $Z_{n}(\theta)$ converges in the quenched sense, because according to our definition of quenched convergence and Proposition \ref{prounicon}, the asymptotic distribution of $Z_{n}(\theta)$ under $\mathbb{P}_{\omega}$ must be independent of $\omega$. Of course, this is more a limitation of our definition of quenched convergence (Definition \ref{defquecon}) than an inherent pathology of the behavior of a ($\mathbb{P}-$convergent) process under the measures $\mathbb{P}_{\omega}$.

\subsection{Cases of Quenched Convergence without Random Centering}
\label{secnonrancennec}

Now consider the following observation: by the proof of Corollary \ref{corobsconpro}  and (\ref{equregnotneed}), for every $\theta\in I$, $E_{0}(Z_{n}(\theta)-E_{-\infty}Z_{n}(\theta))\to 0$ in $L^{1}_{\mathbb{P}}$. It follows from Fatou's lemma (see the proof of Theorem 16.4 on \cite{bilpromea})  that, if we assume the condition
\begin{equation}
\label{weaconqueconzn}
\sup_{n}|E_{0}(Z_{n}(\theta)-E_{-\infty}Z_{n}(\theta))|\in L^{1}_{\mathbb{P}},
\end{equation}
then
$$E[\limsup_{n}|E_{0}(Z_{n}(\theta)-E_{-\infty}Z_{n}(\theta))|]\leq \limsup_{n}E[|E_{0}(Z_{n}(\theta)-E_{-\infty}Z_{n}(\theta))|]=0,$$
which is possible if and only if $E_{0}(Z_{n}(\theta)-E_{-\infty}Z_{n}(\theta))\to 0$, $\mathbb{P}-$a.s. Thus the following result follows from Corollary \ref{sufnecconqueconzn}.

\medskip

\begin{cor}
In the context of Theorem \ref{quecltfoutra}, denote by $Z_{n}(\theta)$ the random variable given in (\ref{noravedisfoutraequ}). Then the validity of condition (\ref{weaconqueconzn}) for $\theta\in I$ implies that $(Z_{k}(\theta)-E_{-\infty}Z_{k}(\theta))_{k\in\mathbb{N}}$ converges to $Y(\theta)$ in the quenched sense. In particular, the condition 
\begin{equation}
\label{weaconqueconznreg}
\sup_{n}|E_{0}Z_{n}(\theta)|\in L^{1}_{\mathbb{P}}
\end{equation}
for $\theta\in I$ implies that
$(Z_{k}(\theta))_{k\in\mathbb{N}}$ converges to $Y(\theta)$ in the quenched sense if $(X_{k})_{k\in\mathbb{Z}}$ is regular. 
\end{cor}

Let us give two more results regarding the quenched convergence of $(Z_{n}(\theta))_{n\in\mathbb{N}}$ for $\theta\in I$ in terms of decay of correlations, whose proof will be given in Section \ref{prothenonrancen}.

\medskip

\begin{thm}
\label{queconznundratdec}
In the context of Theorem \ref{quecltfoutra}, denote by $Z_{n}(\theta)$ the random variable given in (\ref{noravedisfoutraequ}). If the condition
\begin{equation}
\label{conratdecdif}
\sum_{k\in \mathbb{N}^{*}}\frac{|E_{0}[X_{k}-X_{k-1}]|^{2}}{k}<\infty, \mbox{\,\,\,\,\,\,\,$\mathbb{P}$-{\it a.s.}}
\end{equation}
holds,  there exists $J\subset I$ with $\lambda(J)=1$ such that, for every $\theta\in J$, $Z_{n}(\theta)\Rightarrow_{n} Y(\theta)$ in the quenched sense.
\end{thm}





Our last theorem in this direction is related to the {\it Maxwell and Woodroofe condition}, and its proof is essentially an application of results found by Cuny and Merlev\'{e}de in \cite{cunmer}.  The statement is the following.

\medskip

\begin{thm}[Quenched Convergence under the Maxwell-Woodroofe Condition]
\label{queconznmaxwoo}
In the context of Theorem \ref{quecltfoutra}, and given $\theta\in I$, denote by $Z_{n}(\theta)$ the random variable given in (\ref{noravedisfoutraequ}), then the {\upshape Maxwell and Woodroofe condition}
\begin{equation}
\label{maxwoocon}
\sum_{k\in \mathbb{N}^{*}}\frac{||E_{0}S_{k}(\theta)||_{_{\mathbb{P},2}}}{k^{3/2}}<\infty
\end{equation}
implies the quenched convergence of $Z_{n}(\theta)$ to $Y(\theta)$.
\end{thm}

\medskip

\begin{remark}
It is possible to relax the assumption ``$\theta\in I$'' to ``$e^{2i\theta}\notin Spec_{p}(T)$'' in the hypotheses of Theorem \ref{queconznmaxwoo} by using a direct martingale approximation also presented in \cite{cunmer}. See the proof of Theorem 6 in \cite{barpel} for details.
\end{remark}

\section{Quenched Functional Central Limit Theorem}
\label{secquefclt}
Finally, let us address the question of the validity of the quenched Central Limit Theorem in its functional form. 

To begin with, let us recall the definition of the space $(S,d)$ of complex valued cadlag functions on $[0,\infty)$: Definition \ref{skospacomvalcadfun}), and that a {random element} of  $S$ is (by definition) a measurable function $W:\Omega'\to S$ where $(\Omega',\mathcal{F}',\mathbb{P}')$ is a probability space and $S$ is endowed with its Borel sigma algebra $\mathcal{S}$. By an adaptation of the theory for $D[[0,\infty)]$ (see for instance Theorem 16.6 in \cite{bilconpromea}), $\mathcal{S}$ is also the sigma algebra generated by the finite dimensional cylinders
\begin{equation}
\label{deffindimcyl}
H_{t_{1}\dots t_{k},A}:=[\pi_{t_{1}\dots t_{k}}\in A],
\end{equation}
where $A$ is a Borel set in $\mathbb{C}^{k}$, $0\leq t_{1}\leq\cdots \leq t_{k}$, and $\pi_{t_{1}\dots t_{k}}$ is given by (\ref{proope}). 

It follows (see the argument in \cite{bilconpromea}, p.84) that if $(\Omega',\mathcal{F}',\mathbb{P}')$ is a probability space, $W:\Omega'\to S$ is a random element of $S$ if and only if for every $t\geq 0$, $\pi_{t}\circ W$ (i.e., the function $\omega'\mapsto W(\omega')(t)$) is a random variable in $(\Omega',\mathcal{F}',\mathbb{P}')$.

\subsection*{The Question}

Here the problem is the following: consider the setting in the hypothesis of Theorem \ref{quecltfoutra}, and for $(\theta,\omega)\in[0,2\pi)\times\Omega$, consider the function $W_{n}:[0,\infty)\times\Omega\to\mathbb{C}$ given by
\begin{equation}
\label{funqueinvpri}
W_{n}(\theta,\omega)(t):=\frac{S_{\left\lfloor nt\right\rfloor}(\theta,\omega)-E_{0}[S_{\left\lfloor nt\right\rfloor}(\theta,\cdot)](\omega)}{\sqrt{n}}.
\end{equation}
This is: for fixed $\theta,\omega$ and $n$, $W_{n}(\theta,\omega)$ takes the value 
$$\frac{S_{k}(\theta,\omega)-E_{0}[S_{k}(\theta,\cdot)](\omega)}{\sqrt{n}}$$
whenever $t\in [k/n,(k+1)/n)$. 

Note that for fixed $(\theta,\omega)$, $W_{n}(\theta,\omega)$ is an element of $S$, and that there are two ways in which we can regard $W_{n}$ as a random element of $S$:

\begin{enumerate}
\item {\it Fixed frequency approach.} For fixed $\theta\in [0,2\pi)$, consider the function $W_{n}(\theta):\Omega\to S$
\begin{equation}
\label{funqueinvprifixfre}
W_{n}(\theta)(\omega)=W_{n}(\theta,\omega).
\end{equation}
Then $W_{n}(\theta)$ is a random element of $S$.
\item {\it Averaged frequency approach.} Consider the product space $([0,2\pi)\times\Omega, \mathcal{B}\otimes\mathcal{F},\lambda\times\mathbb{P})$. Then the function $W_{n}: [0,2\pi)\times\Omega\to S$ is a random element of $S$.
\end{enumerate}

Our goal is to give results on the quenched convergence of $W_{n}$ from both the fixed frequency and the averaged frequency points of view. Note that, by the discussion in Section \ref{raneleprospa} (see the discussion following Theorem \ref{invpriavefre}), results for $\lambda-$almost every fixed frequency imply results for averaged frequencies. 

\medskip

\subsection{The Invariance Principle for Averaged Frequencies}
 
Our first result concerns the validity of the quenched Invariance Principle under the averaged frequency approach.  It is the following:

\medskip

\begin{thm}[The Quenched Invariance Principle for Averaged Frequencies]
\label{invpriavefre}
In the setting of Theorem \ref{quecltfoutra}, let $B_{1}$, $B_{2}$ be independent standard Brownian motions on $[0,\infty)$ defined on some probability space $(\Omega',\mathcal{F}',\mathbb{P}')$.  Consider the trivial sigma-algebra $\mathcal{B}_{0}:=\{\emptyset,[0,2\pi)\}\subset \mathcal{B}$, and let $S$ be the space of cadlag complex valued functions with the Skorohod distance (Definition \ref{skospacomvalcadfun}). Then the sequence $(W_{n})_{n\in\mathbb{N}^{*}}$ of random elements  of $S$ specified by (\ref{funqueinvpri}) converges in the quenched sense with respect to $\mathcal{B}_{0}\otimes\mathcal{F}_{0}$ to the random function $B:[0,2\pi)\times\Omega'\to S$ specified by
\begin{equation}
\label{defasyfunequ}
B(\theta,\omega')= \frac{\sigma(\theta)}{\sqrt{2}}(B_{1}(\omega')+iB_{2}(\omega')).
\end{equation}
Equivalently, for any decomposition $\{\mathbb{P}_{\omega}\}_{\omega\in\Omega}$ of $E_{0}$ (Definition \ref{regconexpdef}), there exists $\Omega_{0}\subset \Omega$ with $\mathbb{P}\Omega_{0}=1$ such that for every $\omega\in\Omega_{0}$
\begin{equation}
\label{queconpomeequ}
W_{n}\Rightarrow B \mbox{\,\,\,\,\,\,\,\,\it under $\lambda\times\mathbb{P}_{\omega}$.}
\end{equation}
\end{thm}

This theorem should be compared with Theorem \ref{pelwufclt}: it plays a role with respect to this theorem similar to that of Theorem \ref{quecltfoutra} with respect to Theorem \ref{cltpelwu}.

\subsection{Invariance Principles for Almost Every Fixed Frequencies}
\label{queinvpriavefre}
Of course, we would like to give an extension of Theorem \ref{quecltfoutra} in the direction of an invariance principle valid for $\lambda-$a.e fixed  frequency, which in particular would imply the convergence stated in Theorem \ref{invpriavefre}.

To be more precise, note that if we are able to prove that for $\lambda-$almost every fixed $\theta$ the sequence $(W_{n}(\theta))_{n\geq 0}$ of random elements of $S$ (defined on $(\Omega,\mathcal{F},\mathbb{P)}$) converges in the quenched sense to $B(\theta,\cdot)$ with respect to $\mathcal{F}_{0}$ then, by an argument similar to that in the proof of Corollary \ref{corquecltprospa}, the quenched convergence of $W_{n}$ with respect to $\mathcal{B}_{0}\otimes\mathcal{F}_{0}$ follows at once.

The validity of the quenched invariance principle for $\lambda-$almost every $\theta$ is a problem under current research.\footnote{At the moment of writing this monograph, the author ignores whether this stronger form of the invariance principle can be proved without assumptions additional to those in Theorem \ref{quecltfoutra}.} In this work, we will give a result in the direction of Hannan-like conditions guaranteeing its fulfillment.\footnote{But other approaches are possible. For instance via the results  in \cite{cunmer} (see the proof of Theorem \ref{queconznmaxwoo} for an illustration of the use of these results).} 

\subsection*{Motivation}
{In what follows, the notations and the assumptions are those given in Theorem \ref{quecltfoutra}}.

To illustrate our last results we start by considering the {\it Hannan} condition: recall the definition (\ref{defpkequ}) of the projection operators $\mathcal{P}_{k}$ ($k\in\mathbb{Z}$). We say that $(X_{k})_{k\in\mathbb{Z}}$ satisfy the {\itshape Hannan Condition} if
\begin{equation}
\label{hancon}
\sum_{n\in\mathbb{N}}||\mathcal{P}_{0}X_{n}||_{2}<\infty.
\end{equation} 

Cuny and Voln\'{y} showed, in \cite{cunvol}, that in the context of Theorem \ref{invpriavefre}, condition (\ref{hancon}) guarantees that $W_{n}(0)$ converges to 
$$B'(0)=\sigma(0) B_{1}.$$
where 
$$\sigma^{2}(0)=\lim_{n}E_{0}[|Y_{n}(0)|^{2}]$$
$\mathbb{P}-$a.s. (see the notation in Theorem \ref{quecltfoutra}).

In spite of the fact that this is a quenched result for (only) one frequency, and that the quenched asymptotic distribution of $W_{n}(\theta,\cdot)$ does not correspond to a two-dimensional Brownian motion (but to a one-dimensional one), we will see that this condition is actually strong enough to guarantee the quenched convergence of $W_{n}(\theta,\cdot)$ at {every} $\theta\neq 0$ provided that $e^{2i\theta}\notin Spec_{p}(T)$.

Our main result in this direction depends on  the following condition 
\begin{equation}
\label{imhancon}
\sum_{n\geq 0}||\mathcal{P}_{0}(X_{n+1}-X_{n})||_{2}<\infty,
\end{equation}
which is clearly a ``weak'' version of the {Hannan} condition\footnote{To see that this condition is strictly weaker than the Hannan condition consider the process 
$$X_{k}:=\sum_{j\geq 1}\frac{1}{j}x_{k-j}$$
where $(x_{j})_{j\in\mathbb{Z}}$ are the coordinate functions in $\mathbb{R}^{\mathbb{Z}}$, seen as an i.i.d sequence in $L^{2}$, $T$ is the left shift, and $\mathcal{F}_{0}=\sigma(x_{k})_{k\leq 0}$ (see Example \ref{exalinpro} in page \pageref{exalinpro}).} (\ref{hancon}). The result is the following.

\medskip
 
\begin{thm}[A Quenched Invariance Principle for Fixed Frequencies]
\label{barinvprihancon}
With the notation and assumptions of Theorem \ref{invpriavefre}, and assuming (\ref{imhancon}), if $e^{2i\theta}\notin Spec_{p}(T)$,  then $W_{n}(\theta,\cdot)$ converges in the quenched sense to 
\begin{equation}
\label{quelimfixfreequ}
\omega'\mapsto \frac{\sigma(\theta)}{\sqrt{2}}(B_{1}(\omega')+iB_{2}(\omega')).
\end{equation}
where $\sigma(\theta)$ is given as in Theorem \ref{quecltfoutra} (see (\ref{sigtheequ})).
\end{thm}

\medskip

It is worth to further specify a case in which the set of frequencies where the asymptotic distribution is as in (\ref{quelimfixfreequ}) can be easily described. To motivate the following Theorem recall that $T$ is weakly mixing if and only if $Spec_{p}(T)=\{1\}$ (see \cite{qua}, Section 8 for a review of this and other related facts).

Now, as a subgroup of $\mathbb{T}$, $Spec_{p}(T)$ is finite (actually: closed) if and only if there exists $m\in\mathbb{N}^{*}$ such that
\begin{equation}
\label{eigangfinspe}
Spec_{p}(T):=\{e^{2\pi k i/m}\}_{k=0}^{m-1}.
\end{equation}
In other words $Spec_{p}(T)$ is finite if and only it it consists of the points in the unit circle given by the {rational rotations} by an angle of $2\pi/m$ or, what is the same, by the $m-$th roots of unity.

Our last result is the following.

\medskip

\begin{cor}
\label{weamixcas}
Assume that $Spec_{p}(T)$ is finite and its elements are the $m-$th roots of unity.
Under the hypothesis and the notation in Theorem \ref{barinvprihancon}, $W_{n}(\theta)$ converges in the quenched sense to (\ref{quelimfixfreequ}) for all $\theta\in [0,2\pi)$ such that $e^{2im\theta}\neq 1$. If $T$ is in particular weakly mixing, (\ref{quelimfixfreequ}) describes the asymptotic quenched limit of $W_{n}(\theta)$ for all $\theta\neq 0,\pi$.
\end{cor}

{\bf Proof:} Immediate from (\ref{eigangfinspe}) and Theorem \ref{barinvprihancon}.\qed

\part{Proofs}
\section*{General Setting}
\label{gensetpro}
In addition to the notation introduced at the beginning, the following setting will be fixed throughout this part of the monograph: $(\Omega,\mathcal{F},\mathbb{P})$ will be a fixed probability space. $T:\Omega\to \Omega$ will be a fixed invertible, bimeasurable measure-preserving transformation (Definition \ref{defmeapretra}). As before, $T$ will (also) denote the Koopman operator associated to the map $T$ (Definition \ref{defkooope}), and $Spec_{p}(T)$ will denote its point spectrum (Definition \ref{poispe}).  $(\mathcal{F}_{k})_{k\in\mathbb{Z}}$ will be a fixed $T-$filtration (Definition \ref{defadafil}) where $\mathcal{F}_{0}$ is countably generated (Definition \ref{deffcougen}), and given $k\in\mathbb{Z}\cup\{-\infty,\infty\}$, we will denote by $E_{k}$ the conditional expectation with respect to $\mathcal{F}_{k}$, where $\mathcal{F}_{\pm\infty}$ are given via Definition \ref{taisigalgdef}. 

We will assume that  $E_{0}:=E[\,\cdot\,|\mathcal{F}_{0}]$ is regular,  $\{\mathbb{P}_{\omega}\}_{\omega\in\Omega}$ will be a fixed decomposition of $E_{0}$ (Definition \ref{regconexpdef}), and for a given $\omega\in \Omega$, $E^{\omega}$ will denote integration with respect to $\mathbb{P}_{\omega}$.  Given $p>0$, we will also denote by $Id:L^{p}_{\mathbb{P}}\to L^{p}_{\mathbb{P}}$ the identity function (the domain of $Id$ will be clear from the context). When needed, we will use {explicitely} the version of $E_{0}$ given by integration with respect to $\mathbb{P}^{\omega}$: $E_{0}X(\omega):=E^{\omega}X$ for every $X\in L^{1}_{\mathbb{P}}$. Such restriction will not be assumed without explicit indication. 

Finally, $B_{1}, B_{2}$ will denote independent standard Brownian motions defined on some probability space $(\Omega',\mathcal{F}',\mathbb{P}')$, and $N_{j}=B_{j}(1)$ ($j=1,2$) denote independent standard normal random variables on $(\Omega',\mathcal{F}',\mathbb{P}')$.

\subsection*{Dot Product}
In what follows, we will use the notation $a\cdot b$ to denote the dot product between vectors in $\mathbb{R}^{n}$ ($n\in\mathbb{N}^{*}$). Thus if $a=(a_{1},\cdots,a_{n})$ and $b=(b_{1},\dots,b_{n})$ are elements of $\mathbb{R}^{n}$
$$a\cdot b:=a_{1}b_{1}+\cdots + a_{n}b_{n}.$$
In particular, if $z=z_{1}+iz_{2}$ and $w=w_{1}+iw_{2}$ are complex numbers (with respective real and imaginary parts $z_{1},w_{1}$ and $z_{2},w_{2}$)
$$z\cdot w:=z_{1}w_{1}+z_{2}w_{2}.$$

\newpage

\section*{Structure of the Arguments}

Our goal in this part of the monograph is to prove the results stated, but not proved, in Chapter \ref{resandcom}. The general structure of the forthcoming arguments is the following.

\begin{enumerate}
\item{\it Martingale case.} We will start by addressing the martingale case. More precisely, we will prove that if $\theta\in[0,2\pi)$ is such that $e^{2i\theta}\notin Spec_{p}(T)$ (see Definition \ref{poispe}) and $D_{0}(\theta)\in L^{2}_{\mathbb{P}}(\mathcal{F}_{0})\ominus L^{2}_{\mathbb{P}}(\mathcal{F}_{-1})$, then the conclusion no. 1. in Theorem \ref{barinvprihancon} holds replacing $X_{0}$ by $D_{0}(\theta)$ in the statement of this theorem. In this case $\sigma^{2}(\theta)=E|D_{0}(\theta)|^{2}$.
\item {\it Martingale approximations, proof of Theorem \ref{quecltfoutra}.} The next step is the following: given a  stationary process $(X_{k})_{k}=(T^{k}X_{0})_{k}$ with $X_{0}\in L^{2}_{\mathbb{P}}$, we will construct a random element $D_{0}:[0,2\pi)\times\Omega\to L^{2}_{\mathbb{P}}$ with the property that for $\lambda-$a.e $\theta$, $D_{0}(\theta,\cdot)\in L^{2}_{\mathbb{P}}(\mathcal{F}_{0})\ominus L^{2}_{\mathbb{P}}(\mathcal{F}_{-1})$. We will then prove Theorem \ref{quecltfoutra} by showing that for $\lambda-$a.e $\theta$ and $\mathbb{P}-$a.e $\omega$
\begin{equation}
\label{traineclt}
||(Id-E_{0})S_{n}((X_{k})_{k},\theta)-S_{n}((T^{k}D_{0}(\theta,\cdot))_{k},\theta)||_{\mathbb{P}_{\omega},2}=o(\sqrt{n})
\end{equation}
(see Definition \ref{defdisfoutra} for the notation) and then applying Theorem  \ref{limlimlem} together with the martingale case and the discussions made before. 

\item{\it Proof of theorems \ref{invpriavefre} and \ref{barinvprihancon}}. To achieve the proof of these two theorems we will first show that for $\mathbb{P}-$a.e $\omega$
\begin{equation}
\label{trainefcltavefre}
||\max_{1\leq k\leq n}|(Id-E_{0})S_{n}((X_{k}(z))_{k},\theta)-S_{n}((D_{k}(\theta,z))_{k},\theta)|||_{\lambda\otimes\mathbb{P}_{\omega},2}=o(\sqrt{n})
\end{equation}
where $D_{k}(\theta,z):=D_{0}(\theta,T^{k}z)$ (this map will be $\mathcal{B}\otimes\mathcal{F}_{\infty}-$measurable). This will give the proof of Theorem \ref{invpriavefre} by an approximation argument again and the martingale results in Section \ref{marcassec}.

Then we will see that, under the conditions  in the hypothesis of Theorem \ref{barinvprihancon}, (\ref{trainefcltavefre}) holds for $\lambda-$a.e $\theta\in (0,2\pi)$ fixed (actually, for every $\theta$ with $e^{2i\theta}\notin Spec_{p}(T)$) replacing $\lambda\otimes\mathbb{P}_{\omega}$ by $\mathbb{P}_{\omega}$, which again implies the functional form of Theorem \ref{quecltfoutra} by the martingale version previously proved.

\item {\it Proof of Theorem \ref{nonqueconthe}. }We will then prove Theorem \ref{nonqueconthe} by specializing our study to the case explained in Example \ref{exalinpro}: we will see that there exist a sequence $(a_{k})_{k\in\mathbb{Z}}\in l^{2}(\mathbb{N})$ generating a linear process with the property announced in Theorem \ref{nonqueconthe}.

\item{\it Proof of theorems \ref{queconznundratdec} and \ref{queconznmaxwoo}. }The proofs of these results end the content of this monograph. We will achieve them by using the characterization of quenched convergence without random centering given in Corollary \ref{sufnecconqueconzn} (which is proved in previous sections), together with suitable interpretations of results present in the existing literature applied to the processes under our consideration.

\end{enumerate}

\chapter{Martingale Case}
\label{marcas}
This chapter is devoted to present the martingale theorems (Theorem \ref{queinvprimar} and Corollary \ref{marcasprospa}) which will be used  to prove the results on quenched asymptotics presented in Chapter \ref{resandcom} via suitable martingale  approximations and transport theorems.

In Section \ref{preresmar}, we introduce some results from the existing literature which will allow us to carry out the proof of Theorem \ref{queinvprimar} by specializing to the case under our consideration. Section \ref{marcassec} presents the aforementioned proofs of the martingale case.

\section{Preliminary Results}
\label{preresmar}
In this short section we present some preliminary facts needed to prove Theorem \ref{queinvprimar} below, from which all the  proofs of the (positive) results announced in Chapter \ref{resandcom} will follow via suitable martingale approximations. With the exception of Lemma \ref{rel16cunmerpel} (proved first by Cuny et.al in \cite{cunmerpel}), the results presented here pertain to the classical literature, but we decided to include their statements due to their very specific role among the proofs of our main theorems. The setting is that explained in page \pageref{gensetpro}.

Our first result is a lemma that will allows us, among other things, to characterize the asymptotic finite-dimensional distributions of the normalized discrete Fourier transforms of a martingale at a frequency  not associated to an element of $Spec_{p}(T)$ (in the sense just to be stated). 

\medskip

\begin{lemma}
\label{rel16cunmerpel}
Let $\theta\in [0,2\pi)$ be such that $e^{-2i\theta}\notin Spec_{p}(T)$, let $p\geq 1$ and let $Y\in L^{p}_{\mathbb{P}}$. Then for every $z\in\mathbb{C}$ 
\begin{equation}
\label{cunmerpelrel16equ}
\lim_{n}\frac{1}{n}\sum_{k=0}^{n-1} E_{k-1}(z\cdot  (T^{k}Ye^{ik\theta}))^{2}= \frac{|z|^2}{2\,}{E|Y|^{2}} \mbox{\,\,\,\,\,\it $\mathbb{P}-$a.s. and in $L^{p}_{\mathbb{P}}$,}
\end{equation} 
where the (probability one) set $\Omega_{\theta}$ of pointwise convergence does not depend on $z$. 
\end{lemma}

{\bf Proof:}\footnote{For an alternative explanation of this proof see the proof of relation (16) in \cite{cunmerpel}.} Let $z=z_{1}+iz_{2}$, and note first that
\begin{equation}
\label{simmatpro}
E_{k-1}({z}\cdot  (T^{k}Ye^{ik\theta}))^{2}=T^{k}E_{-1}(z\cdot (Ye^{ik\theta}))^{2}.
\end{equation}
Now, using Euler's formula and the double-angle identities, it is an elementary (though somewhat tedious) exercise in trigonometry to prove that, if $z=z_{1}+iz_{2}$ and $Y=Y_{1}+iY_{2}$ (where $z_{j}, Y_{j}$, $j=1,2$ are real-valued) and $k\in\mathbb{N}$,
$$(z\cdot (Ye^{ik\theta}))^{2}=\left(\frac{z_{1}^{2}+z_{2}^{2}}{2}(Y_{1}^{2}+Y_{2}^{2})\right)+\left((z_{2}^{2}-z_{1}^{2})Y_{1}Y_{2}-(Y_{2}^{2}-Y_{1}^{2})z_{1}z_{2}\right)\sin(2k\theta)+$$
\begin{equation}
\label{gencunmerpelrel16equ}
\left((z_{1}Y_{1}+z_{2}Y_{2})^{2}-(Y_{2}z_{1}-Y_{1}z_{2})^{2}\right)\frac{\cos(2k\theta)}{2}\,,
\end{equation}
thus there exist real constants (depending on $z$) $a_{j}, b_{j}$ ($j=1,2,3$) such that
$$(z\cdot (Ye^{ik\theta}))^{2}=\frac{|z|^{2}}{2\,}|Y|^{2}+$$
\begin{equation}
\label{gencunmerpelrel16equsim}
\left(a_{1}Y_{1}^{2}+a_{2}Y_{2}^{2}+a_{3}Y_{1}Y_{2}\right)\cos (2k\theta)+\left(b_{1}Y_{1}^{2}+b_{2}Y_{2}^{2}+b_{3}Y_{1}Y_{2}\right)\sin (2k\theta).
\end{equation}
The conclusion follows at once from (\ref{simmatpro}), (\ref{gencunmerpelrel16equsim}), Theorem \ref{ergthedisfoutra} and Corollary \ref{corcontozer}, by taking $\Omega_{\theta}$ as the set of probability one where, according to the notation on Theorem \ref{ergthedisfoutra}
$$S_{n}(E_{-1}|Y|^2,0)/n\to E|Y|^{2},\,\,\,S_{n}(E_{-1}[Y_{1}Y_{2}],2\,\theta)/n\to 0, \,\,\,\mbox{\it and \,\,\,\,}S_{n}(E_{-1}|Y|^{2},2\,\theta)/n\to 0 $$
as $n\to\infty$.\qed

The next two theorems are very classical. We will use them to prove our martingale limit theorems in the setting of discrete Fourier transforms in the quenched sense.

\bigskip

\begin{thm}[The Lindeberg-L\'{e}vy Theorem for Martingales]
\label{bilmarclt}
For each $n\in\mathbb{N}^{*}$, let $\Delta_{n1},\dots,\Delta_{nk},\dots$ be a sequence of real-valued martingale differences with respect to some increasing filtration $\mathcal{F}_{0}^{n}\subset \cdots\subset\mathcal{F}_{k}^{n}\subset\dots$. Define, for $1\leq k\leq n$, $\sigma_{nk}:=E[\Delta_{nk}^{2}||\mathcal{F}_{n}^{k-1}]$. If for some $\sigma\geq 0$ the following two conditions hold
\begin{enumerate}
\item $\sum_{k\geq 0}\sigma_{nk}^{2}\Rightarrow \sigma^{2}$ {\it as $n\to \infty$,}
\item $\sum_{k\geq 0}E[\Delta_{nk}^{2}I_{[\Delta_{nk}\geq \epsilon]}]\to 0$ {\it as $n\to \infty$,}
\end{enumerate}
then $Z_{n}:=\sum_{k\geq 0}\Delta_{nk}\Rightarrow \sigma N$ where $N$ is a standard normal random variable.
\end{thm}

{\bf Proof:} \cite{bilpromea}, p.476. \qed

\bigskip

\begin{thm}[The Functional form of Theorem \ref{bilmarclt}]
\label{bilinvpri}
For each $n\in\mathbb{N}^{*}$, let $\Delta_{n1},\dots,\Delta_{nk},\dots$ be a sequence of real-valued martingale differences with respect to some increasing filtration $\mathcal{F}_{0}^{n}\subset \cdots\subset\mathcal{F}_{k}^{n}\subset\dots$ and defines, for $1\leq k\leq n$, $\sigma_{nk}^{2}:=E[\Delta_{nk}^{2}||\mathcal{F}_{n}^{k-1}]$. If for some $\sigma\geq 0$ the following two conditions hold for every $t\geq 0,\epsilon>0$
\begin{enumerate}
\item $\sum_{k\leq nt}\sigma_{nk}^{2}\Rightarrow_{n} \sigma^2 t,$
\item $\sum_{k\leq nt}E[\Delta_{nk}^{2}I_{[\Delta_{nk}\geq \epsilon]}]\to_{n} 0,$
\end{enumerate}
then the random functions $X_{n}(t):=\sum_{k\leq nt}\Delta_{nk}$ converge in distribution to $\sigma W$ in the sense of $D[[0,\infty)]$, where $W$ is a standard Brownian motion.
\end{thm}

{\bf Proof:} This is a slight reformulation of Theorem 18.2 in \cite{bilconpromea}, (pp. 194-195): the case $\sigma>0$ follows by a simple renormalization, and to cover the case $\sigma=0$, note that the convergence  (18.6) in \cite{bilconpromea} becomes a simple consequence of the definition given there of $\zeta_{nk}$ and the hypothesis (corresponding to $\sigma=0$)
$$\sum_{k\leq nt}\sigma_{nk}^{2}\Rightarrow 0$$
for every $t\geq 0$.\qed

\section{Martingale Case}
\label{marcassec}

As already mentioned, all of the positive results in Section \ref{resandcom} follow from the following theorem via suitable martingale approximations. 

\medskip

\begin{thm}[The Quenched Invariance Principle for the Discrete Fourier Transforms of a Martingale]
\label{queinvprimar}
Under the setting introduced in page \pageref{gensetpro}, and given $\theta\in [0,2\pi)$ such that $e^{-2i\theta}\notin Spec_{p}(T)$ (Definition \ref{poispe}), assume that $D_{0}(\theta)\in L^{2}_{\mathbb{P}}(\mathcal{F}_{0})\ominus L^{2}_{\mathbb{P}}(\mathcal{F}_{-1})$ is given, and define the $(\mathcal{F}_{k-1})_{k\in \mathbb{N}^{*}}-$adapted martingale $(M_{k}(\theta))_{k\in\mathbb{N}}$ by
\begin{equation}
\label{nmardef}
M_{n}(\theta):=\sum_{k=0}^{n-1}T^{k}D_{0}(\theta)e^{ik\theta}
\end{equation}
for all $n\in \mathbb{N}$. Then the sequence $(V_{k}(\theta))_{k\in\mathbb{N}^{*}}$ of random elements of  $D[[0,\infty),\mathbb{C}]$ defined by
\begin{equation}
\label{defofvnthe}
V_{n}(\theta)(t):= M_{\left\lfloor nt\right \rfloor}(\theta)/\sqrt{n} 
\end{equation}
for every $n\in\mathbb{N}^{*}$, converges in the quenched sense with respect to $\mathcal{F}_{0}$ to the random function $B(\theta):\Omega'\to D[[0,\infty),\mathbb{C}]$ given by
\begin{equation}
\label{asyvnthemarcas}
B(\theta)(\omega')= [E|D_{0}(\theta)|^2/2]^{1/2}(B_{1}(\omega')+iB_{2}(\omega')).
\end{equation}
\end{thm}

\medskip

\begin{remark}
\label{remequmarcas}
Before proceeding to the proof it is worth noticing the following: the conclusion of  Theorem \ref{quecltfoutra}, specialized to this case, is a statement about the asymptotic distribution of the random variables $V_{n}(\theta)(1)$. Now, by Corollary \ref{dedmerpel} and the orthogonality  under $E_{0}$ of $(T^{k}D_{0}(\theta))_{k\in\mathbb{N}}$, \footnote{Note that if $(k,r)\in\mathbb{N}\times\mathbb{N}^{*}$ is given then, since $T^{r}{D_{0}(\theta)}\in L^{2}_{\mathbb{P}}(\mathcal{F}_{r})\ominus L^{2}_{\mathbb{P}}(\mathcal{F}_{r-1})$,
$$E_{0}[T^{k}D_{0}(\theta)T^{k+r}\overline{D_{0}(\theta)}]=E_{0}[T^{k}D_{0}(\theta)E_{k}T^{k+r}\overline{D_{0}(\theta)}]=E_{0}T^{k}[D_{0}(\theta)E_{0}T^{r}\overline{D_{0}(\theta)}]=0.$$} 
$$E[|D_{0}(\theta)|^{2}]=\lim_{n}\frac{1}{n}\sum_{k=1}^{n-1}E_{0}T^{k}|D_{0}(\theta)|^{2}=\lim_{n}\frac{1}{n}E_{0}|M_{n}(\theta)-E_{0}M_{n}(\theta)|^{2}$$
so that the equality (\ref{sigtheequ}) is certainly verified in this case.
\end{remark}

{\bf Proof of Theorem \ref{queinvprimar}:} Let us start by sketching the argument of the proof: we will see that there exists $\Omega_{\theta}\subset \Omega$ with $\mathbb{P}\Omega_{\theta}=1$ such that for every $\omega\in \Omega_{\theta}$ the following holds:

\begin{enumerate}
\item [\bf a.]{\it The sequence of random functions $(V_{n}(\theta))_{n}$ in $D[[0,\infty),\mathbb{C}]$ is tight with respect to $\mathbb{P}_{\omega}$.} To prove this, we will actually prove the convergence in distribution of both the real and imaginary parts of $(V_{n}(\theta))_{n}$  to a Brownian motion via Theorem \ref{bilinvpri} (see the ``Criteria for Tightness'' in section \ref{condoinfc}).
\item[\bf b.] {\it The finite dimensional asymptotic distributions under $\mathbb{P}_{\omega}$ of  $(V_{n}(\theta))_{n}$ converge to those of two independent Brownian motions with the scaling $E[(D_{0}(\theta))^2]^{1/2}/\sqrt{2}$.} under $\mathbb{P}_{\omega}$. For this we will proceed via the Cramer-Wold theorem, using some of the results already presented.
\end{enumerate}

We go now to the details: first, we will assume, making it explicit only when necessary, that $E_{0}$ is the version of $E[\,\cdot\,|\mathcal{F}_{0}]$ given by integration with respect to the decomposing probability measures $\{\mathbb{P}_{\omega}\}_{\omega\in \Omega}$ (see Definition \ref{regconexpdef}).

Now denote, for every $k\in \mathbb{N}$
\begin{equation}
\label{defdkmarcas}
D_{k}(\theta):=T^{k}D_{0}(\theta).
\end{equation}

Let $\Omega_{\theta,1}'$ be the set of probability one guaranteed by Lemma \ref{rel16cunmerpel} for the case $Y=D_{0}(\theta)$. By Remark \ref{remproone}, there exists a set $\Omega_{\theta,1}$ with $\mathbb{P}\Omega_{\theta,1}=1$ such that for every $\omega\in \Omega_{\theta,1}$
$$\lim_{n}\frac{1}{n}\sum_{k=0}^{n-1} E_{k-1}(z\cdot  (D_{k}(\theta)e^{ik\theta}))^{2}= \frac{|z|^{2}}{2\,}E|D_{0}(\theta)|^{2} $$
$\mathbb{P}_{\omega}$-a.s. for all $z\in\mathbb{C}$.

For such $\omega$'s the first hypothesis of Theorem \ref{bilinvpri} is verified  by the triangular arrays $(Re(M_{k}(\theta)/\sqrt{n}))_{1\leq k \leq n}$  and $(Im(M_{k}(\theta)/\sqrt{n}))_{1\leq k \leq n}$ ($n\in\mathbb{N}^{*}$) with respect to $\mathbb{P}_{\omega}$, because they arise from the particular choices $z=1$ and $z=i$ respectively.

To verify the second hypothesis in Theorem \ref{bilinvpri} we start from the $\mathbb{P}-$a.s. inequality
$$E_{0}\left[\frac{1}{n}\sum_{k=0}^{n-1} \left((Re(D_{k}(\theta)e^{ik\theta}))^2I_{[|Re(D_{k}(\theta)e^{ik\theta})|\geq \epsilon \sqrt{n}]]}+(Im(D_{k}(\theta)e^{ik\theta}))^2I_{[|Im(D_{k}(\theta)e^{ik\theta})|\geq \epsilon\sqrt{n}]}\right)\right]\leq $$
\begin{equation}
\label{inetigreaima}
E_{0}\left[\frac{1}{n}\sum_{k=0}^{n-1} |D_{k}(\theta)|^2I_{[|D_{k}(\theta)|\geq \epsilon\sqrt{n}]}\right].
\end{equation}
Now, given $\eta>0$ there exists $N\geq 0$ such that $\mu_{N}:= {E}[|D_{0}(\theta)|^{2}I_{[|D_{0}(\theta)|^{2}\geq \epsilon^{2} N]}]<\eta$, and therefore
$$\limsup_{n}\frac{1}{n}\sum_{k= 0}^{n-1}E_{0}T^{k}[|D_{0}(\theta)|^{2}I_{[|D_{0}(\theta)|^{2}\geq \epsilon^{2} n]}]\leq$$
\begin{equation}
\label{inelimequzer}
\limsup_{n}\frac{1}{n}\sum_{k=0}^{ n-1}E_{0}T^{k}[|D_{0}(\theta)|^{2}I_{[|D_{0}(\theta)|^{2}\geq \epsilon^{2} N]}]=\mu_{N}\leq \eta
\end{equation}
over a set $\Omega_{\theta,\epsilon,\eta}$ with $\mathbb{P}\Omega_{\theta,\epsilon,\eta}=1$, where we made use of Corollary \ref{dedmerpel}. Without loss of generality, (\ref{inetigreaima}) holds for all $\omega\in \Omega_{\theta,\epsilon,\eta}$.

Denote by $Z_{n}^{\epsilon}$ the random variable at the left-hand side of the inequality (\ref{inetigreaima}) and note that, if we define 
\begin{equation}
\label{defome2pri}
\Omega_{\theta,2}=\bigcap_{\epsilon>0,\eta>0}\Omega_{\theta,\epsilon,\eta}
\end{equation}
where the intersection runs over rational $\epsilon$, $\eta$, then $\mathbb{P}\Omega_{\theta,2}=1$, and for every $\epsilon>0$ and every $\omega\in \Omega_{\theta,2}$
$$\lim_{n}Z_{n}^{\epsilon}(\omega)=0.$$
or, what is the same,  for all $\omega\in \Omega_{\theta,2}$
$$\frac{1}{n}\sum_{k=0}^{n-1} \left((Re(D_{k}(\theta)e^{ik\theta}))^2I_{[|Re(D_{k}(\theta)e^{ik\theta})|\geq \epsilon \sqrt{n}]]}+(Im(D_{k}(\theta)e^{ik\theta}))^2I_{[|Im(D_{k}(\theta)e^{ik\theta})|\geq \epsilon\sqrt{n}]}\right)$$ 
goes to $0$ in $L^{1}_{\mathbb{P}_{\omega}}$ as $n\to \infty$. 

Thus, if $\Omega_{\theta,3}$ is a set of probability one such that $(Re(M_{k}(\theta))_{k\in \mathbb{N}^{*}}$ and $(Im(M_{k}(\theta)))_{k\in \mathbb{N}^{*}}$ is a $(\mathcal{F}_{k-1})_{k\in\mathbb{N}^{*}}-$ adapted martingale in $L^{2}_{\mathbb{P}_{\omega}}$ for all $\omega\in \Omega_{\theta,3}$ (Corollary \ref{cormarundpome}), the hypotheses 1. and 2. in Theorem \ref{bilinvpri} are verified for all $\omega$ in the set $\Omega_{\theta}$ defined by
\begin{equation}
\label{defometheequ}
\Omega_{\theta}:=\bigcap_{k=1}^{3}\Omega_{\theta,k}.
\end{equation} 
Since $\mathbb{P}\Omega_{\theta}=1$ this finishes the proof of {\bf a.}

To prove {\bf b.} {we will show that for any given $n\in \mathbb{N}$, any $\omega\in \Omega_{\theta}$, and any $0\leq t_{1}\leq \cdots\leq t_{n}$,  the $\mathbb{C}^{n}=\mathbb{R}^{2n}-$valued process
$$(V_{n}(\theta)(t_{1}), V_{n}(\theta)(t_{2})-V_{n}(\theta)(t_{1}),\cdots, V_{n}(\theta)(t_{n})-V_{n}(\theta)(t_{n-1}))$$
has the same asymptotic distribution as
$$\mathbf{B}^{\theta}(t_{1},\cdots, t_{n}):=$$
$$[E|D_{0}(\theta)|^2/2]^{1/2}(B_{1}(t_{1}),B_{2}(t_{1}), B_{1}(t_{2})-B_{1}(t_{1}), B_{2}(t_{2})-B_{2}(t_{1}),\cdots, B_{2}(t_{n})-B_{2}(t_{n-1}))$$
under $\mathbb{P}_{\omega}$ and therefore, by the mapping theorem (\cite{bilconpromea}, Theorem 2.7), the finite dimensional asymptotic distributions of $V_{n}(\theta)$ under $\mathbb{P}_{\omega}$ and those of (\ref{asyvnthemarcas}) under $\mathbb{P}'$ are the same.}

For simplicity we will assume $n=2$. The argument generalizes easily to an arbitrary $n\in\mathbb{N}$. 

Our goal is thus to prove that for all $\omega\in \Omega_{\theta}$ and all $0\leq s\leq t$ the asymptotic distribution of 
\begin{equation}
\label{tesvec}
\mathbf{V}_{n}^{\theta}(s,t):=(V_{n}(\theta)(s),V_{n}(\theta)(t)-V_{n}(\theta)(s)) 
\end{equation}
(a $\mathbb{C}^{2}=\mathbb{R}^{4}$-valued process) is the same under $\mathbb{P}_{\omega}$ as that of  

\begin{equation}
\label{asydistesvec}
\mathbf{B}^{\theta}(s,t):=[E|D_{0}(\theta)|^2/2]^{1/2}(B_{1}(s),B_{2}(s), B_{1}(t)-B_{1}(s), B_{2}(t)-B_{2}(s))
\end{equation}

under $\mathbb{P}'$.

To prove the convergence in distribution of (\ref{tesvec}) to (\ref{asydistesvec}) we will use the Cramer-Wold theorem. This is, we will see that for any $\omega\in\Omega_{\theta}$, any $0\leq s\leq t$, and any 
\begin{equation}
\label{defu}
\mathbf{u}=(a_{1},a_{2},b_{1},b_{2})\in \mathbb{R}^{4}
\end{equation}
the asymptotic distribution under $\mathbb{P}_{\omega}$ of the stochastic process $(U_{n})_{n\in\mathbb{N}^{*}}$ defined by
\begin{equation}
\label{varcrawoldev}
U_{n}:=\mathbf{u}\cdot\mathbf{V}^{\theta}_{n}(s,t)
\end{equation}
is that of a normal random variable with variance
\begin{equation}
\label{asyvarcrawol}
\sigma_{\mathbf{u},s,t}^{2}(\theta):=\frac{E[|D_{0}(\theta)|^{2}]}{2}((a_{1}^{2}+a_{2}^{2})s+(b_{1}^{2}+b_{2}^{2})(t-s)).
\end{equation}

To do so we will verify the hypotheses of Theorem \ref{bilmarclt}. Fix $\mathbf{u}$ as above and note that
$$U_{n}=\sum_{k=0}^{\left\lfloor ns\right\rfloor}\eta_{nk}(a_{1},a_{2})+\sum_{k=\left\lfloor ns\right\rfloor+1}^{\left\lfloor nt\right\rfloor}\eta_{nk}(b_{1},b_{2})$$ 
where 
\begin{equation}
\label{defetank}
\eta_{nk}(x_{1},x_{2})=\frac{1}{\sqrt{n}}(x_{1},x_{2})\cdot e^{ik\theta}T^{k}D_{0}(\theta).
\end{equation} 

By the construction of $\Omega_{\theta}$, for every $0\leq r$, every $x_{1}, x_{2}$ and every $\omega\in \Omega_{\theta}$,  $(\eta_{nk}(x_{1},x_{2}))_{0\leq k\leq \left\lfloor nr\right\rfloor}$ is a triangular array of $(\mathcal{F}_{k})_{k}-$ adapted (real-valued)  martingale differences under $\mathbb{P}_{\omega}$, and by Lemma \ref{rel16cunmerpel} combined with Remark \ref{remproone} we can assume that
\begin{equation}
\label{firhyptover}
\sum_{k\leq ns}E_{k-1}[\eta_{nk}^{2}(a_{1},a_{2})]+\sum_{ns< k\leq nt}E_{k-1}[\eta_{nk}^{2}(b_{1},b_{2})]\to_{n} \sigma_{\mathbf{u},s,t}^{2}(\theta) 
\end{equation}

$\mathbb{P}_{\omega}-$a.s.\footnote{More precisely: redefine $\Omega_{\theta}$ above by intersecting it with the set $\Omega_{\theta}'$ of elements $\omega$ for which the convergence in Lemma \ref{rel16cunmerpel}  happens $\mathbb{P}_{\omega}-$a.s.} This verifies the first hypothesis in Theorem \ref{bilmarclt} under $\mathbb{P}_{\omega}$ for all $\omega\in \Omega_{\theta}$ for the triangular array defining $U_{n}$. 

It remains to prove that if $\omega\in \Omega_{\theta}$ then
\begin{equation}
\label{sechyptover}
\sum_{k\leq ns}E_{0}[\eta_{nk}^{2}(a_{1},a_{2})I_{[|\eta_{nk}(a_{1},a_{2})|>\epsilon]}](\omega)\to 0.
\end{equation} 
 
This is, that for all $\omega\in \Omega_{\theta}$
$$\sum_{k\leq ns}\eta_{nk}^{2}(a_{1},a_{2})I_{[|\eta_{nk}(a_{1},a_{2})|>\epsilon]}\to 0$$
in $L^{1}_{\mathbb{P}_{\omega}}$.

To do so we depart from the Cauchy-Schwartz inequality to get that
$$\eta_{nk}^{2}(x_{1},x_{2})\leq \frac{1}{n}(x_{1}^{2}+x_{2}^{2})T^{k}|D_{0}(\theta)|^{2},$$ 
so that the sum in (\ref{sechyptover}) is bounded by 
$$\frac{1}{n}\sum_{k\leq ns}E_{0}T^{k}[(a_{1}^{2}+a_{2}^{2})|D_{0}(\theta)|^{2}I_{[(a_{1}^{2}+a_{2}^{2})|D_{0}(\theta)|^{2}\geq \epsilon^{2} n]}].$$

This obviously goes to zero when $a_{1}=a_{2}=0$. Otherwise it is the same as 
$$(a_{1}^{2}+a_{2}^{2})\frac{1}{n}\sum_{k\leq ns}E_{0}T^{k}[|D_{0}(\theta)|^{2}I_{[|D_{0}(\theta)|^{2}\geq {\epsilon^{2} n}/{(a_{1}^{2}+a_{2}^{2})}]}],$$

which, again, goes to zero as $n\to \infty$ for every $\omega\in \Omega_{\theta}$. \qed

\medskip

\begin{remark}
\label{remtwopar}
When necessary, {\it specially when discussing quenched convergence in the product space} $([0,2\pi)\times\Omega,\mathcal{B}\otimes\mathcal{F})$, we will specify the dependence on $\omega\in \Omega$ of a given family $\{Y(\theta)\}_{\theta\in\Theta}$ of functions $Y(\theta):\Omega\to S$ parametrized by $\theta$ by seeing them as  sections of  functions depending on two parameters. So if, for instance, $D_{0}(\theta)$ is the function introduced in Theorem \ref{queinvprimar}, we will write
$$D_{0}(\theta,\omega):=D_{0}(\theta)(\omega)$$
and so on.
\end{remark}

The following result basically follows from Theorem \ref{queinvprimar} via Theorem \ref{domcontheconexp}. We state it in a language that will be convenient for our forthcoming proofs.

\medskip

\begin{cor}[The Averaged-frequency Quenched Invariance Principle for Martingales]
\label{marcasprospa}
Assume that $D_{0}(\theta)\in L^2_{\mathbb{P}}(\mathcal{F}_{0})\ominus L^{2}_{\mathbb{P}}(\mathcal{F}_{-1})$ is given for every $\theta\in [0,2\pi)$, and that the function $(\theta,\omega)\mapsto D_{0}(\theta,\omega)$  is $\mathcal{B}\otimes\mathcal{F}-$measurable (see Remark \ref{remtwopar}). Then, with the notation in Theorem \ref{queinvprimar}, and assuming that $\mathcal{F}$ is countably generated, there exists $\Omega_{0}\subset \Omega$ with $\mathbb{P}\Omega_{0}=1$ such that for all $\omega_{0}\in \Omega_{0}$, the distribution of $(\theta,\omega)\mapsto V_{n}(\theta,\omega)$ under $\lambda\times\mathbb{P}_{\omega_{0}}$ converges to that of $(\theta,\omega')\mapsto B(\theta,\omega')$ under $\lambda\times\mathbb{P}'$.
\end{cor}

{\bf Proof:}  First, the $\mathcal{B}\otimes \mathcal{F}-$(resp. $\mathcal{B}\otimes \mathcal{F}'-$)measurability of $(\theta,\omega)\mapsto V_{n}(\theta,\omega)$(resp. $(\theta,\omega')\mapsto B(\theta,\omega')$) follows at once from Remark \ref{remcrimead0infc} (page \pageref{remcrimead0infc}).

We claim that given any continuous and bounded function $f:D[[0,\infty),\mathbb{C}]\to \mathbb{R}$
\begin{equation}
\label{} 
\lim_{n}E[f\circ V_{n}|\mathcal{B}\otimes\mathcal{F}_{0}](\theta,\omega)=Ef(B(\theta))
\end{equation}
$\lambda\times\mathbb{P}-$a.s., where the expectation at the left-hand side (resp. right-hand side) denotes integration with respect to $\mathbb{P}$ (resp. $\mathbb{P}'$). 

Before proceeding to the proof of (\ref{}), let us explain why this implies the desired conclusion:

\begin{enumerate}
\item First, note that (\ref{}) can be considered an equality of $\mathcal{B}\otimes\mathcal{F}_{0}$ measurable functions, the $\mathcal{B}-$measurable function at the right being considered as constant in $\Omega$ for fixed $\theta$. 

\item It follows by an application of Theorem \ref{domcontheconexp} that, for any given $\mathcal{B}_{0}\subset \mathcal{B}$
\begin{equation}
\label{semqueconequ}
\lim_{n}E[f\circ V_{n}|\mathcal{B}_{0}\otimes\mathcal{F}_{0}]=E[Ef(B(\theta))|\mathcal{B}_{0}\otimes \mathcal{F}_{0}]
\end{equation}
$\lambda\times\mathbb{P}-$a.s. 

\item If $\mathcal{B}_{0}=\{\emptyset,[0,2\pi)\}$ is the trivial sigma algebra then (see Example \ref{triexareg} in page \pageref{triexareg}) if we define $\lambda_{\theta}:=\lambda$ for all $\theta\in [0,2\pi)$, $\{\lambda_{\theta}\}_{\theta\in [0,2\pi)}$ is a decomposition of $E[\,\cdot\,|\mathcal{B}_{0}]$ and it follows, from Proposition \ref{proregpro}, that (\ref{semqueconequ}) is nothing but the statement of convergence $V_{n}\Rightarrow B$ under $\lambda\times \mathbb{P}_{\omega}$ for $\mathbb{P}-$a.e $\omega$: this is the desired conclusion.

\end{enumerate} 

{\it Proof of (\ref{}).} To prove (\ref{}) we proceed as follows: first, the set
$$\{(\theta,\omega):\lim_{n}(E[f\circ V_{n}|\mathcal{B}\otimes\mathcal{F}_{0}](\theta,\omega)-Ef(B(\theta)))=0\}$$
is $\mathcal{B}\otimes\mathcal{F}$ measurable, and to see that it has product measure one it suffices to see that for $\lambda-$a.e fixed $\theta$
\begin{equation}
\label{redfin}
\mathbb{P}[\lim_{n}(E[f\circ V_{n}|\mathcal{B}\otimes\mathcal{F}_{0}](\theta,\cdot)-Ef(B(\theta)))]=1.
\end{equation}

Let $I$ be the set
$$I:=\{\theta\in [0,2\pi): e^{2i\theta}\notin Spec_{p}(T)\},$$
which satisfies $\lambda(I)=1$ according to Proposition \ref{prospecou} ($\mathcal{F}$ is countably generated). We claim that (\ref{redfin}) holds for every $\theta\in I$.

To see why this claim is true, note that by Proposition \ref{proregpro} and Example \ref{triexareg} again, if $\delta_{\theta}$ denotes the Dirac measure at $\theta$, then
$$\{\delta_{\theta}\times \mathbb{P}_{\omega}\}_{(\theta,\omega)\in [0,2\pi)\times\Omega}$$ is a decomposition of $E[\,\cdot\,|\mathcal{B}\times\mathcal{F}_{0}]$, and Theorem \ref{queinvprimar} gives that for every $\theta\in I$ there exists $\Omega_{\theta}$ with $\mathbb{P}\Omega_{\theta}=1$ such that for every $\omega\in \Omega_{\theta}$
$$
\lim_{n}E[f\circ V_{n}|\mathcal{B}\otimes\mathcal{F}_{0}](\theta,\omega)=\lim_{n}\int_{[0,2\pi)\times \Omega}f\circ V_{n}(\alpha,z) \,d(\delta_{\theta}\times \mathbb{P}_{\omega})(\alpha,z)=$$
$$\lim_{n}\int_{\Omega}f(V_{n}(\theta,z))\,d\mathbb{P}_{\omega}(z)=\lim_{n}E[f(V_{n}(\theta))|\mathcal{F}_{0}](\omega)= Ef(B(\theta))$$
as desired.\qed



\chapter{Proofs of Theorems \ref{quecltfoutra}, \ref{invpriavefre} and \ref{barinvprihancon}}

The exposition is divided as follows: Section \ref{marapp} presents the martingale  approximation results leading to the proof of the theorems stated in the title of this chapter. This section is divided into two parts: ``Approximation Lemmas'' (Section \ref{applem}), giving a presentation  of the abstract martingale approximation results that will be used to construct the proofs of the corresponding theorems, and ``The Approximating Martingales'' (Section \ref{appmar}), in which we present the actual martingales to be used along the rest of the chapter.

Section \ref{proquecltfoutra} presents the proof of Theorem \ref{quecltfoutra} which, in analogy with the forthcoming proofs, consists of verifying the hypothesis of the corresponding lemma from Section \ref{applem} via the martingales introduced in Section \ref{appmar}. The key step is a further, ``concrete'' approximation lemma (Lemma \ref{genmarapplem}), whose proof at some point makes use of a technique analogous to that used to prove Theorem \ref{ergthedisfoutra}. With such lemma and the previous results at hand, the proof of the aforementioned theorem is reduced to a few, almost obvious, lines.

Section \ref{protheinvpri} is devoted to the proofs of theorems \ref{invpriavefre} and \ref{barinvprihancon}. The reason to present these proofs in the same section lies in the fact that, as the reader will see, the corresponding arguments can be considered ``branches'' of the same decomposition of the difference between the process and the approximating martingales (Lemma \ref{firapp}), and in particular to stress the ``smoothing'' role of Hunt and Young's inequality (Theorem \ref{hunyou}) in the proofs involving ``averaged'' (as opposed to ``fixed'') frequencies.

The chapter finishes with a note (see page \pageref{notlimlimlem}) pointing out that the use of Theorem \ref{limlimlem} along these proofs is not essential.

\section{Martingale Approximations}
\label{marapp}

In this section we will give a series of approximation lemmas whose verification will imply the results stated as theorems \ref{quecltfoutra}, \ref{invpriavefre} and \ref{barinvprihancon}. For the sake of clarity, we will limit our discussion in this section to state and prove the aforementioned lemmas and in particular to explain {\it why} these imply the corresponding results stated in Chapter \ref{resandcom}. We will also present, without further analysis, the martingales used along the proofs. The actual verification of the hypotheses in these lemmas under the hypotheses of the corresponding theorems via the given martingales is deferred to later sections.  

\subsection{Approximation Lemmas}
\label{applem}
Our first approximation lemma is the following.

\medskip

\begin{lemma}[Approximation Lemma for Theorem \ref{quecltfoutra}]
\label{applemquecltfoutra}
Under the hypotheses and notation in Theorem \ref{quecltfoutra}, assume that there exists $I'\subset [0,2\pi)$ with $\lambda(I')=1$ satisfying the following: for every $\theta\in I'$, there exists $D_{0}(\theta)\in L^{2}(\mathcal{F}_{0})\ominus L^{2}(\mathcal{F}_{-1})$ with the property that, if we denote $M_{n}(\theta):=\sum_{k=0}^{n-1}T^{k}D_{0}(\theta)e^{ik\theta}$ ($n\in\mathbb{N}^{*}$),
\begin{equation}
\label{applemquecltfoutraequ}
\lim_{n}\frac{1}{n}E_{0}|S_{n}(\theta)-E_{0}S_{n}(\theta)-M_{n}(\theta)|^{2}=0 
\end{equation}
$\mathbb{P}-$a.s. and in $L^{1}_{\mathbb{P}}$. Then the conclusion of Theorem \ref{quecltfoutra} holds with $I=I'\setminus\{\theta:e^{2i\theta}\in Spec_{p}(T)\}$ and
\begin{equation}
\label{equsigthequecltfoutra}
\sigma^{2}(\theta)=E|D_{0}(\theta)|^2.
\end{equation}
\end{lemma}

Before proving this lemma let us point out the following interesting fact: assume that, for $\theta\in [0,2\pi)$, $D_{0}(\theta)$ and $D'_{0}(\theta)$ are given as in Lemma \ref{applemquecltfoutra}, and let $(M_{n}(\theta))_{n\in\mathbb{N}^{*}}$ and $(M_{n}'(\theta))_{n\in\mathbb{N}^{*}}$ be the corresponding $(\mathcal{F}_{n-1})_{n\in\mathbb{N}^{*}}-$adapted martingales. Then, according to Corollary \ref{dedmerpel} and the footnote in Remark \ref{remequmarcas}
$$E|D_{0}(\theta)-D'_{0}(\theta)|^{2}=\lim_{n}E_{0}\frac{1}{n}|M_{n}(\theta)-M'_{n}(\theta)|^{2}\leq$$
\begin{equation}
\label{unimardifequ} 2\limsup_{n}\frac{1}{n}(E_{0}|S_{n}(\theta)-E_{0}S_{n}(\theta)-M_{n}(\theta)|^{2}+E_{0}|S_{n}(\theta)-E_{0}S_{n}(\theta)-M'_{n}(\theta)|^{2}).
\end{equation}

In particular, we have the following uniqueness result.

\medskip

\begin{prop}[Uniqueness of $D_{0}(\theta)$]
\label{uniappmar}
In the context of Lemma \ref{applemquecltfoutra}, and given $\theta\in [0,2\pi)$ (not necessarily in $I'$), there exists {\it at most} one function $D_{0}(\theta)\in L^{2}_{\mathbb{P}}(\mathcal{F}_{0})\ominus L^{2}_{\mathbb{P}}(\mathcal{F}_{-1})$ satisfying (\ref{applemquecltfoutraequ}).
\end{prop}

{\bf Proof:} Combine (\ref{applemquecltfoutraequ}) with (\ref{unimardifequ}).\qed

We proceed now to the proof of Lemma \ref{applemquecltfoutra}.

{\bf Proof of Lemma \ref{applemquecltfoutra}:} First, note that $\lambda(I)=1$ by Proposition \ref{prospecou}.

Let now  $\{\mathbb{P}_{\omega}\}_{\omega\in\Omega}$ be a decomposition of $E_{0}$ (Definition \ref{regconexpdef}). According to (\ref{applemquecltfoutraequ}), there exists, for $\theta\in I$, $\Omega_{\theta}\subset \Omega$ with $\mathbb{P}\Omega_{\theta}=1$ such that, for all $\omega\in \Omega_{\theta}$
\begin{equation}
\label{applemquecltfoutraequreg}
\lim_{n}\frac{1}{n}||S_{n}(\theta)-E_{0}S_{n}(\theta)-M_{n}(\theta)||_{_{\mathbb{P}_{\omega},2}}^{2}=0.
\end{equation}
and  the quenched convergence stated in Theorem \ref{quecltfoutra} follows at once from Theorem \ref{queinvprimar} (taking $t=1$), Proposition \ref{prounicon}, and Corollary \ref{corlimlimlem} (replacing $X_{r,n}:=M_{n}(\theta)/\sqrt{n}$ and $X_{n}:=(S_{n}(\theta)-E_{0}S_{n}(\theta))/\sqrt{n}$).

Now, by orthogonality under $E_{0}$ (see the footnote in Remark \ref{remequmarcas}) and Corollary \ref{dedmerpel}, 
$$E[|D_{0}(\theta)|^{2}]=\lim_{n}\frac{1}{n}\sum_{k=0}^{n-1}E_{0}T^{k}|D_{0}(\theta)|^{2}=\lim_{n}\frac{1}{n}E_{0}|M_{n}(\theta)|^{2}$$
in the $\mathbb{P}-$a.s and $L^{1}_{\mathbb{P}}$ senses, which implies  by (\ref{applemquecltfoutraequreg}) and the Minkowski inequality that  
\begin{equation}
\label{staonequecltfoutramar}
\lim_{n}\frac{1}{n}E_{0}|S_{n}(\theta)-E_{0}S_{n}(\theta)|^{2}=E[|D_{0}(\theta)|^{2}]
\end{equation}
$\mathbb{P}-$a.s. and in $L^{1}_{\mathbb{P}}$: this is the statement {\it 1.} in Theorem \ref{quecltfoutra}. 

Finally, to see that $\theta\mapsto \sigma^{2}(\theta)$ necessarily defines a version of the spectral density of $(X_{k}-E_{-\infty}X_{k})_{k\in\mathbb{Z}}$ we proceed as follows: integrating (\ref{staonequecltfoutramar}) and using the $L^{1}_{\mathbb{P}}$ convergence we get that, for $\theta\in I$
$$E[|D_{0}(\theta)|^{2}]=\lim_{n}\frac{1}{n}E|S_{n}(\theta)-E_{0}S_{n}(\theta)|^2=$$
$$\lim_{n}\frac{1}{n}E|S_{n}(\theta)-E_{-\infty} S_{n}(\theta)-E_{0}(S_{n}(\theta)-E_{-\infty}S_{n}(\theta))|^2=$$
\begin{equation}
\label{equvarmarspeden}
\lim_{n}\frac{1}{n}E|(S_{n}(\theta)-E_{-\infty} S_{n}(\theta))|^{2},
\end{equation}
where for the last equality we used the fact that
$$\lim_{n}\frac{1}{{n}}E|E_{0}(S_{n}(\theta)-E_{-\infty} S_{n}(\theta))|^{2}=0$$
(see the proof of Corollary \ref{corobsconpro}). The conclusion follows from (\ref{equvarmarspeden}), Theorem \ref{spedenasalim} and the fact that $\lambda(I)=1$.
\qed

Our next two approximation lemmas make use of an additional parameter, ``$r$'', whose presence will allow us in particular to carry on the proofs of Theorems \ref{invpriavefre} and \ref{barinvprihancon} without restricting ourselves explicitly to the set $I$ in Lemma \ref{applemquecltfoutra}.\footnote{Nevertheless, we will work under this restriction when carrying on the actual proofs.}
 


\medskip

\begin{lemma}[Approximation Lemma for Theorem \ref{barinvprihancon}]
\label{appleminvprifixfre}
With the notation and conventions introduced on page \pageref{gensetpro}, and  with the additional notation (\ref{funqueinvpri}) and (\ref{funqueinvprifixfre}), let $\theta\in[0,2\pi)$ be such that $e^{2i\theta}\notin Spec_{p}(T)$. Assume given, for every $r\in \mathbb{N}$, a function $D_{r,0}(\theta)\in L^{2}_{\mathbb{P}}(\mathcal{F}_{0})\ominus L^{2}_{\mathbb{P}}(\mathcal{F}_{-1})$, and given $n\in\mathbb{N}^{*}$ denote by $M_{r,n}(\theta)$ the function
$$M_{r,n}(\theta):=\sum_{k=0}^{n-1}T^{k}D_{r,0}(\theta)e^{ik\theta}.$$
Then the hypothesis
\begin{equation}
\label{appleminvprifixfreequ}
\lim_{r}\limsup_{n}E_{0}[\frac{1}{n}\max_{1\leq k\leq n}|S_{k}(\theta)-E_{0}S_{k}(\theta)-M_{r,k}(\theta)|^2]=0,\mbox{\,\,\,\,\,\,\,\,$\mathbb{P}-$a.s.},
\end{equation}
implies the existence of 
\begin{equation}
\label{limdorsqu}
\sigma^{2}(\theta):=\lim_{r}{E[|D_{r,0}^{2}(\theta)|]},
\end{equation}
and if we denote 
\begin{equation}
\label{asylimapplembarinvpriequ}
B(\theta)(\omega'):=(\sigma^{2}(\theta)/2)^{1/2}(B_{1}(\omega')+iB_{2}(\omega')),
\end{equation} 
then $W_{n}(\theta)$ converges in the quenched sense (with respect to $\mathcal{F}_{0}$)  to $B(\theta)$ as $n\to\infty$.
 \end{lemma} 

Before proceeding to the proof of Lemma \ref{appleminvprifixfre}, let us point out the following.

\medskip
 
\begin{remark}[Consistency of the Notation (\ref{limdorsqu})]
\label{remconnot}
Notice that, in the context of Lemma \ref{applemquecltfoutra}, if for $\theta\in I$  the hypotheses of Lemma \ref{appleminvprifixfre} are verified, then necessarily 
$$\lim_{r}E[|D_{r,0}(\theta)|^2]=E|D_{0}(\theta)|^{2},$$
where $D_{0}(\theta)$ is chosen according to Lemma \ref{applemquecltfoutra}.
 
To see this just note that, for such $\theta$, the conclusion of Lemma \ref{applemquecltfoutra} follows from Lemma \ref{appleminvprifixfre} by evaluating (\ref{asylimapplembarinvpriequ}) at $t=1$, and compare the corresponding random variables thus obtained.\footnote{If $X,Y$ are nonzero random variables, $X=Y$ in distribution, and $a,b\geq0$ are constants with $aX=bY$ in distribution, then $a=b$: for every $M>0$, $0=bE[|Y|I_{[|bY|\leq M]}]-aE[|X|I_{[|aX|\leq M]}]=(b-a)E[|X|I_{[|aX|\leq M]}]$.}
\end{remark}
 
 {\bf Proof of Lemma \ref{appleminvprifixfre}:} Start by recalling the notation and criteria introduced in Section \ref{condoinfc}, specially in the numeral {\it \ref{resdom}.}, and define $V_{r,n}$ as in (\ref{defofvnthe}) with $D_{r,0}$ in place of $D_{0}$ for every $(r,n)\in\mathbb{N}\times\mathbb{N}^{*}$. 
 
For $m\geq 1$, the Skorohod metric $d_{m}$ on $D[[0,m],\mathbb{C}]$, is dominated by the uniform (product) metric. Thus for every $(m,n)\in\mathbb{N}^{*}\times \mathbb{N}^{*}$
$$d_{m}(r_{m}W_{n}(\theta,\omega),r_{m}V_{r,n}(\theta,\omega))\leq \frac{\sqrt{m}}{\sqrt{n'}}\max_{1\leq k\leq n' }|S_{k}(\theta,\omega)-E_{0}S_{k}(\theta)(\omega)-M_{r,k}(\theta,\omega)|.$$

where $n'=mn$. It follows from (\ref{appleminvprifixfreequ}) that there exists $\Omega_{0,1}\subset \Omega$ with $\mathbb{P}\Omega_{0,1}=1$ such that if $\omega\in\Omega_{0,1}$
\begin{equation}
\label{applembarinvpriequ}
\lim_{r}\limsup_{n}||d_{m}(r_{m}W_{n}(\theta),r_{m}V_{r,n}(\theta))||_{\mathbb{P}_{\omega},2}=0.
\end{equation}

Now, according to Theorem \ref{queinvprimar} and Proposition \ref{prounicon}, there exists $\Omega_{0,2}\subset \Omega$ with $\mathbb{P}\Omega_{0,2}=1$ with the following property: for every $\omega\in \Omega_{0,2}$,
$$V_{r,n}(\theta)\Rightarrow_{n} B_{r}(\theta)$$
under $\mathbb{P}_{\omega}$ where $B_{r}(\theta)$ is the random element with domain in $(\Omega',\mathcal{F}',\mathbb{P}')$ defined by 
\begin{equation}
\label{defbrequ}
B_{r}(\theta)(\omega'):=(E[|D_{r,0}(\theta)|^{2}]/2)^{1/2}(B_{1}(\omega')+iB_{2}(\omega')).
\end{equation}
Since for every fixed $m\geq 0$, $B_{r}$ is $\mathbb{P}'-$a.e continuous at $m$,  the observations in Section \ref{condoinfc} (numeral {\it \ref{resdom}.} again) imply that for every $m\in\mathbb{N}$ and every $\omega\in \Omega_{0,2}$ 
\begin{equation}
\label{conresvmfixfre}
r_{m}V_{r,n}(\theta)\Rightarrow_{n} r_{m}B_{r}(\theta)
\end{equation}
under $\mathbb{P}_{\omega}$ as $n\to\infty$.

Let $\Omega_{0}:=\Omega_{0,1}\cap\Omega_{0,2}$. According to Theorem \ref{queinvprimar} and Corollary \ref{corlimlimlem}, (\ref{applembarinvpriequ}) together with (\ref{conresvmfixfre}) imply the following: given $\omega\in\Omega_{0}$ and  $m>0$, there exists  a random element $\hat{B}^{m}(\theta)$ of $D[[0,m],\mathbb{C}]$ such that 
\begin{equation}
\label{asylimauxb}
r_{m}W_{n}(\theta)\Rightarrow_{n} \hat{B}^{m}
\end{equation}
under $\mathbb{P}_{\omega}$, and $r_{m}B_{r}(\theta)\Rightarrow_{r} \hat{B}^{m}(\theta)$ under $\mathbb{P}'$. 

We claim that, actually, there exists
$$\sigma^{2}(\theta):=\lim_{r}{E[|D_{r,0}^{2}(\theta)|]}$$
from where it follows easily that, if $B(\theta)$ is given by (\ref{asylimapplembarinvpriequ}), the distribution of $\hat{B}^{m}(\theta)$ is the same as that of $r_{m}B(\theta)$, and the conclusion will follow at once from {\it \ref{resdom}.} in Section \ref{condoinfc}, 
 (\ref{asylimauxb}) and Proposition \ref{prounicon}, because $\mathbb{P}\Omega_{0}=1$.

{\it Proof of the existence of (\ref{limdorsqu})}. To prove the existence of the limit (\ref{limdorsqu}) notice first that, by Theorem \ref{theconindoinfty} there exists, for every $m>0$,  a number $0<t<m$ such that $r_{m}B_{r}(\theta)(t)\Rightarrow_{r} \hat{B}^m(\theta)({t})$. For any of such $t$ we get the existence of a random variable $N(\theta,t)$ such that, if $N_{1}$ and $N_{2}$ are i.i.d standard normal variables
$$\left(\frac{tE|D_{0,r}(\theta)|^{2}}{2}\right)^{1/2}(N_{1}+iN_{2})\Rightarrow_{r} N({\theta},t)$$
and the existence of the limit in (\ref{limdorsqu}) follows at once from Proposition \ref{procontypcom} in page \pageref{procontypcom}.

Finally, note  that $\sigma(\theta)$ is indeed given by (\ref{sigtheequ}) in accordance to Remark \ref{remconnot} and the statement of Lemma \ref{applemquecltfoutra}.\qed

\medskip

Before proceeding to the next approximation lemma let us anticipate the fact that, under (\ref{imhancon}), the hypotheses of Lemma \ref{appleminvprifixfre} will be verified for every $\theta\in [0,2\pi)$ provided that $e^{2i\theta}\notin Spec_{p}(T)$. On proving this, we will encounter some ``intermediate'' approximations that will lead us to verify the hypotheses of Lemma \ref{appleminvpriavefre} below assuming {\it only} the hypotheses of Theorem \ref{invpriavefre}. 

Our next approximation lemma is the ``two-parameters'' version of the previous one.

\medskip
 
 \begin{lemma}[Approximation Lemma for the Proof of Theorem \ref{invpriavefre}] 
\label{appleminvpriavefre}
Under the setting in page \pageref{gensetpro}, denote by $\mathcal{D}_{\infty,\mathbb{C}}$ the Borel sigma-algebra of $D[[0,\infty),\mathbb{C}]$. Assume that for every $(r,\theta)\in\mathbb{N}\times [0,2\pi)$, $D_{r,0}(\theta)$ is given as in the statement of Lemma \ref{appleminvprifixfre} and that the function
$(\theta,\omega)\mapsto D_{r,0}(\theta)(\omega)$  is $\mathcal{B}\otimes\mathcal{F}_{\infty}/\mathcal{D}_{\infty,\mathbb{C}}$ measurable, and denote by $E_{0}$ the version of $E[\,\cdot\,|\mathcal{F}_{0}]$ given by integration with respect to $\{\mathbb{P}_{\omega}\}_{\omega\in\Omega}$ (Definition \ref{regconexpdef}): $E_{0}X(\omega)=E^{\omega}X$ for every $X\in L^{1}_{\mathbb{P}}$. Then the hypotheses
\begin{enumerate}
\item There exists $I'\subset [0,2\pi)$ with $\lambda(I')=1$ such that, for every $\theta\in I'$
$$\sigma^{2}(\theta):=\lim_{r} E|D_{r,0}(\theta)|^{2}$$
is well defined.
\item The equality 
\begin{equation}
\label{appleminvpriavefreequ}
\lim_{r}\limsup_{n}\int_{0}^{2\pi}E^{\omega}[\frac{1}{n}\max_{1\leq k\leq n}|S_{k}(\theta,\cdot)-E_{0}[S_{k}(\theta,\cdot)]-M_{r,k}(\theta,\cdot)|^2] d\lambda(\theta)=0
\end{equation}
(see also Remark \ref{remtwopar}) holds for $\mathbb{P}-$a.e $\omega$.
\end{enumerate}
imply  (together) the conclusion of Theorem \ref{invpriavefre}.

\end{lemma}

{\bf Proof:} First, the assumption that $E_{0}[Z](\omega)=E^{{\omega}}Z$ where $E^{\omega}$ denotes integration with respect to $\mathbb{P}_{\omega}$ guarantees the $\mathcal{B}\otimes\mathcal{F}_{0}-$measurability of the integrand (see {\it Step 2.} in the proof of Proposition \ref{proregpro}). In particular, by \cite{bilpromea}, Theorem 18.1-(ii), the given integral makes sense for every $\omega\in\Omega$.

Now note that, since  for every $(k,r,\theta)\in\mathbb{Z}\times\mathbb{N}\times[0,2\pi)$ the random variables $T^{k}D_{r,0}(\theta,\cdot)$ and $T^{k}X_{0}$ are $\mathcal{F}_{\infty}-$measurable, we can assume that $\mathcal{F}=\mathcal{F}_{\infty}$. Since $\mathcal{F}_{\infty}$ is countably generated ($\mathcal{F}_{0}$ is), Corollary \ref{marcasprospa} (page \pageref{marcasprospa}) guarantees that there exist $\Omega_{0,1}$ with $\mathbb{P}\Omega_{0,1}=1$ such that for every $\omega\in\Omega_{0,1}$, $(\theta,\omega)\mapsto V_{r,n}(\theta,\omega)$ converges  to (\ref{defbrequ}) under $\lambda\times\mathbb{P}_{\omega}$.

The same arguments as in the proof of Lemma \ref{appleminvpriavefre} 
guarantee that there exists $\Omega_{0}\subset \Omega$ with $\mathbb{P}\Omega_{0}=1$ such that, for every fixed $r\in\mathbb{N}$, the sequence of $\mathcal{B}\otimes\mathcal{F}/\mathcal{D}_{\infty,\mathbb{C}}-$measurable functions  $(\theta,\omega)\mapsto V_{r,n}(\theta,\omega)$ satisfy $V_{r,n}\Rightarrow_{n} B_{r}$ under $\lambda\times\mathbb{P}_{\omega}$ (where $(\theta,\omega')\mapsto B_{r}(\theta,\omega')$ is given by (\ref{defbrequ})), that there exists a random function $\hat{B}\in D[[0,\infty),\mathbb{C}]$ (defined on some unspecified probability space) such that  $B_{r}\Rightarrow_{r} \hat{B}$ (under $\lambda\times\mathbb{P}'$) and that, for every $\omega\in\Omega_{0}$, $W_{n}\Rightarrow_{n} \hat{B}$ under $\lambda\times\mathbb{P}_{\omega}$. 

To prove that we can take $\hat{B}=B$, where $B$ is as in the statement of Theorem \ref{invpriavefre},   note  that the $\lambda-$a.e well definition  of (\ref{limdorsqu}) guarantees that $B_{r}$ converges to 
$$B(\theta,\omega')=(\sigma^{2}(\theta)/2)^{1/2}(B_{1}(\omega')+B_{2}(\omega')),$$
$\lambda\times \mathbb{P}'-$a.s. and that, according to Remark \ref{remconnot} and Lemma \ref{applemquecltfoutra}, $\theta\mapsto\sigma^{2}(\theta)$ is certainly a version of the spectral density of $(X_{k}-E_{-\infty}X_{k})_{k\in\mathbb{Z}}$.\qed

\subsection{The Approximating Martingales}
\label{appmar}

We finish this section introducing the martingales used along the proofs of the results established in this chapter. We will defer any discussion about the martingales themselves to later sections.

For every $(r,n,\theta)\in \mathbb{N}\times\mathbb{N}^{*}\times[0,2\pi)$, denote 
 
\begin{equation}
\label{appmardefequ}
\begin{array}{lcr}
D_{r,0}(\theta):=\sum_{k=0}^{r}\mathcal{P}_{0}X_{k}e^{ik\theta}, &  \,\,\,\,\,\,\,\,& M_{r,n}(\theta):=\sum_{k=0}^{n-1}T^{k}D_{r,0}(\theta)e^{ik\theta}\\
 &  &\\
D_{0}(\theta):=\lim_{r}D_{r,0}(\theta), & \,\,\,\,\,\,\,\, & M_{n}(\theta):=\sum_{k=0}^{n-1}T^{k}D_{0}(\theta)e^{ik\theta}.
\end{array}
\end{equation}

{\it When necessary}, we will indicate the dependence on $X_{0}$, $T$, and $\mathcal{P}_{0}$ by denoting 
\begin{equation}\label{impnotequ}
D_{r,0}(\theta)=D_{r,0}(X_{0},T,\mathcal{P}_{0},\theta),
\end{equation}
and so on.

where $D_{0}(\theta)$ is defined as a limit  in the $L^{2}_{\mathbb{P}}$ sense, provided that such limit exists. 

In consonance with the notation introduced in Definition \ref{defdisfoutra}, we treat the case $\theta=0$ denoting, for every $r\in\mathbb{N}$, 
\begin{equation}
\label{appmardef0equ}
D_{r,0}:=D_{r,0}(0), \,\, D_{0}:=D_{0}(0){\,\,\,\,\,\mbox{and}
\,\,\,\,\,\,}M_{r,n}:=M_{r,n}(0), \, M_{n}:=M_{n}(0).
\end{equation}

\section{Proof of Theorem \ref{quecltfoutra}}
\label{proquecltfoutra}
The following lemma will be of fundamental importance to prove the validity of the hypotheses of Lemma \ref{applemquecltfoutra}.
\medskip

\begin{lemma}[Almost Surely Approximation Lemma]
\label{genmarapplem} In the context of Theorem \ref{quecltfoutra}, and with the notation (\ref{appmardefequ}) and (\ref{appmardef0equ}), fix  $\,\theta\in [0,2\pi)$ and assume that $D_{r,0}(\theta)$
converges $\mathbb{P}-$a.s as $r\to \infty$ and that $\sup_{r\in\mathbb{N}}|D_{r,0}(\theta)|\in L^{2}_{\mathbb{P}}$. Then $D_{0}(\theta)$ is well defined and
\begin{equation}
\label{genmarapplemequ}
\lim_{n}\frac{1}{n}E_{0}|S_{n}(\theta)-E_{0}S_{n}(\theta)-M_{n}(\theta)|^{2}=0.
\end{equation}
$\mathbb{P}-$a.s. and in $L^{1}_{\mathbb{P}}$.
\end{lemma}

{\bf Proof}. We will proceed in two steps.

{\it Step 1. Assume $\theta=0$.} We will prove that, if $\sup_{r\in\mathbb{N}}|D_{r,0}|\in L^{2}_{\mathbb{P}}$ and $D_{r,0}$ converges $\mathbb{P}-$a.s. as $r\to\infty$, then 
$$\lim_{n}\frac{1}{n}E_{0}|S_{n}-E_{0}S_{n}-M_{n}|^{2}=0,$$
$\mathbb{P}-$a.s.

Let $D_{0}:=\lim_{r}D_{r,0}$ (in the  $\mathbb{P}-$a.s sense). To see that $\lim_{r}D_{r,0}=D_{0}$ in $L^{2}_{\mathbb{P}}$ note that   
$|D_{r,0}-D_{0}|\leq 2\sup_{r\in\mathbb{N}}|D_{r,0}|$ and therefore, by the dominated convergence theorem and the hypotheses on $(D_{r,0})_{r\in\mathbb{N}}$, 
$$\lim_{r}E|D_{0}-D_{r,0}|^{2}=0.$$
as desired. Note also that, by a similar argument 
\begin{equation}
\label{contozertaiser}
\lim_{N}E[\sup_{j\geq N}|D_{0}-D_{j,0}|^2]=0.
\end{equation}

Notice now that  $D_{0}\in L^{2}_{\mathbb{P}}(\mathcal{F}_{0})\ominus L^{2}_{\mathbb{P}}(\mathcal{F}_{-1})$, because this is a closed subspace of $L^{2}_{\mathbb{P}}$. To prove  (\ref{genmarapplemequ}) let us start in the following way: given $n\in\mathbb{N}^{*}$, we have
$$S_{n}-E_{0}S_{n}-M_{n}=\sum_{k=0}^{n-1}(E_{k}S_{n}-E_{k-1}S_{n}-T^{k}D_{0})=\sum_{k=0}^{n-1}\mathcal{P}_{k}T^{k}(S_{n-k}-D_{0}).$$
The term at the right-hand side is a decomposition of the term  at the left-hand side as a sum of orthogonal functions with respect to $E_{0}$ (see the footnote on Remark \ref{remequmarcas}), and therefore
 $$\frac{1}{n}E_{0}[|S_{n}-E_{0}S_{n}-M_{n}|^2]=\frac{1}{n}\sum_{k=1}^n E_{0}[|\mathcal{P}_{k}T^{k}(S_{n-k}-D_{0})|^2]=\frac{1}{n}\sum_{k=1}^{n}E_{0}T^{k}|D_{0}-D_{n-k,0}|^{2}.$$ 
 Now fix $N\in\mathbb{N}^{*}$. For every $n\geq N$ we can decompose
 $$\frac{1}{n}\sum_{k=1}^{n}E_{0}T^{k}|D_{0}-D_{n-k,0}|^{2}= \frac{1}{n}\sum_{k=1}^{n-N}E_{0}T^{k}|D_{0}-D_{n-k,0}|+\sum_{k=n-N+1}^{n}E_{0}T^{k}|D_{0}-D_{n-k,0}|\leq $$
$$\frac{1}{n}\sum_{k=1}^{n-N}E_{0}T^{k}\sup_{j\geq N}|D_{0}-D_{k,0}|^{2}+ \frac{2}{n}\,\sum_{k=n-N+1}^{n}E_{0}T^{k}\sup_{j\geq 0}|D_{j,0}|^{2}.$$
Note that, since the last summand contains (only) the last $N$ elements of the $n+1-$th ergodic average  in the statement of Theorem \ref{ergtheadafil} corresponding  to the random variable $\sup_{j\geq 0}|D_{j,0}|^{2}$ we have, according to such result combined with the estimates above,  that
$$\limsup_{n\in\mathbb{N}}\frac{1}{n}E_{0}|S_{n}-E_{0}S_{n}-M_{n}|^{2}\leq E_{0}P_{0}[\sup_{j\geq N}|D_{0}-D_{k,0}|^{2}],$$
both in the $\mathbb{P}-$a.s. and $L^{1}_{\mathbb{P}}$ senses, where $P_{0}$ is the orthogonal projection over the subspace of $T-$invariant functions (see Remark \ref{ptkooopecas}). The conclusion follows via (\ref{contozertaiser}) and the continuity in $L^{1}_{\mathbb{P}}$ of $E_{0}P_{0}$ by letting $N\to\infty$.

\medskip

{\it Step 2: general case.} The general case follows from the previous one via the following argument: assume that the hypotheses in Lemma \ref{genmarapplem} hold for a given $\theta\in [0,2\pi)$, and denote by $\tilde{E}$ the integration with respect to $\lambda\times\mathbb{P}$.  Let $\tilde{X}_{0}$ be the extension to the product space specified by  Definition \ref{extprospa}, let  $\tilde{T}_{\theta}$ be the extension map in (\ref{exttprospathe}) and, for every $k\in\mathbb{Z}$, let $\tilde{\mathcal{F}}_{k}=\mathcal{B}\otimes\mathcal{F}_{k}$ and  $\tilde{E}_{k}:=E[\,\cdot\,|\tilde{\mathcal{F}}_{k}]$, so that for  $Y\in L^{1}_{\mathbb{P}}$, $$\tilde{E}_{k}\tilde{Y}(u,\omega):=E[\tilde{Y}|\tilde{\mathcal{F}}_{k}](u,\omega)=e^{iu}E_{k}Y(\omega)=\widetilde{E_{k}Y}(u,\omega).$$
Note also that, since $|\tilde{Y}(u,\omega)|^{2}=|Y(\omega)|^{2}$, then $\tilde{E}_{k}|\tilde{Y}|^{2}(u,\omega)=E_{k}|Y|^{2}(\omega)$.

It is not hard to see that $(\tilde{\mathcal{F}}_{k})_{k\in\mathbb{Z}}$ is a $\tilde{T}_{{\theta}}-$filtration (Definition \ref{defadafil}),  that $\tilde{X}_{0}\in L^{2}_{\lambda\times\mathbb{P}}(\tilde{\mathcal{F}}_{0})$, and that if we follow the definitions in (\ref{appmardefequ}) and (\ref{appmardef0equ}) (see also (\ref{impnotequ})) with $\tilde{X}_{0}$, $\tilde{T}_{{\theta}}$, and  $\tilde{E}_{k}$ in place of $X_{0}$, $T$ and $E_{k}$ then we get that 
$${D}_{r,0}(\tilde{X}_{0},\tilde{T}_{\theta},\tilde{\mathcal{P}}_{0},0)=\tilde{D}_{r,0}(X_{0},T,\mathcal{P}_{0},\theta),$$
and similarly for ${M}_{r,n}(\tilde{X}_{0},\tilde{T}_{\theta},\tilde{P}_{0},0)$ and ${S}_{n}(\tilde{X}_{0},\tilde{T}_{\theta},0)$.

In particular, $\tilde{E}[\sup_{r\in\mathbb{N}}|\tilde{D}_{r,0}(\theta)|^{2}]= E[\sup_{r\in\mathbb{N}}|{D}_{r,0}(\theta)|^{2}]<\infty$ and $$\tilde{D}_{r,0}(\theta)(u,\omega)=e^{iu}D_{r,0}(\theta)(\omega)$$
converges $\lambda\times\mathbb{P}-$a.s.

Finally, by the case already studied $(\theta=0)$, we have that for $\lambda\times\mathbb{P}-$a.e $(u,\omega)$, 
\begin{equation}
\label{frozertothe}
0=
\lim_{n}\frac{1}{n}\tilde{E}_{0}|\tilde{S}_{n}(\theta)-\tilde{E}_{0}\tilde{S}_{n}(\theta)-\tilde{M}_{n}(\theta)|^{2}(u,\omega)=\lim_{n}\frac{1}{n}{E}_{0}|{S}_{n}(\theta)-{E}_{0}{S}_{n}(\theta)-{M}_{n}(\theta)|^{2}(\omega),
\end{equation}
which implies the desired conclusion by fixing $u$ in such a way that the first equality holds $\mathbb{P}-$a.s.\qed
Lemma \ref{genmarapplem} completes the set of tools needed to reach the proof of Theorem \ref{quecltfoutra}.

{\bf Proof of Theorem \ref{quecltfoutra}:} Note that, according to (\ref{relconexpkooopeequ}),
$$T^{k}\mathcal{P}_{-k}X_{0}=\mathcal{P}_{0}X_{k}$$
for every $k\in\mathbb{Z}$ and therefore, since $T$ is measure preserving
\begin{equation}
\label{boul2sumpkxo}
||X_{0}||_{_{\mathbb{P},2}}^{2}=\sum_{k\geq 0}||\mathcal{P}_{-k}X_{0}||_{_{\mathbb{P},2}}^{2}=\sum_{k\geq 0}||\mathcal{P}_{0}X_{k}||_{_{\mathbb{P},2}}^2.
\end{equation}
An application of Proposition \ref{proconsupL2} (page \pageref{proconsupL2}) combined with lemmas \ref{applemquecltfoutra} and \ref{genmarapplem} gives the conclusion in Theorem \ref{quecltfoutra}.\qed

\section{Proof of Theorems \ref{invpriavefre} and \ref{barinvprihancon}}
\label{protheinvpri}

We move on now to the construction of proofs for theorems \ref{invpriavefre} and \ref{barinvprihancon}. 

Following the explanations in Section \ref{applem}, our goal is to prove the approximations (\ref{appleminvprifixfreequ}) and (\ref{appleminvpriavefreequ}) in lemmas \ref{appleminvprifixfre} and \ref{appleminvpriavefre}. Our first step towards this goal is to prove the following decomposition.

\medskip
 
\begin{lemma}
\label{firapp}
In the setting on page \pageref{gensetpro}, for all $(n,r,\theta)\in \mathbb{N}\times\mathbb{N}^{*}\times [0,2\pi)$, $X_{0}\in L^{2}_{\mathbb{P}}(\mathcal{F}_{0})$, and with the notation (\ref{appmardefequ}), the following equality holds :
\begin{equation}
\label{firappequ}
\begin{array}{lcl}
S_{n}(\theta)-E_{0}S_{n}(\theta)-M_{r,n}(\theta)&=  &-e^{i(n-1)\theta}\left(\sum_{k=1}^{r}  (T^{n-1}E_{0}X_{k}-E_{0}T^{n-1}E_{0}X_{k}) e^{ik\theta}\right)\\
 & & \\
 & +&e^{ir\theta}\sum_{k=2}^{n-1}(T^{k}E_{-1}X_{r}-E_{0}T^{k}E_{-1}X_{r})e^{ik\theta}\\
 & & \\
  & -&D_{r,0}(\theta).
\end{array}
\end{equation}
\end{lemma}

{\bf Proof:} Fix $(n,r,\theta)\in \mathbb{N}\times\mathbb{N}^{*}\times [0,2\pi)$. We depart from the following decomposition of $X_{0}$ (the array is intended to make visible the rearrangements):
\begin{equation}
\label{prolem11}
\begin{array}{rll}
X_{0}=E_{0}X_{0} =&(E_{0}-E_{-1})X_{0} &+ \,\,\,\,\,\,\,\,\,\,\,\,\,\,\,\,\,\,E_{-1}X_{0}\\
 +& (E_{0}-E_{-1})X_{1}e^{i\theta} &- (E_{0}-E_{-1})X_{1}e^{i\theta}+ \\
 +& (E_{0}-E_{-1})X_{2}e^{i2\theta} &- (E_{0}-E_{-1})X_{2}e^{i2\theta}+ \\
  & & \vdots \\
  +& (E_{0}-E_{-1})X_{r}e^{ir\theta} &- (E_{0}-E_{-1})X_{r}e^{ir\theta}\\
  & & \\
=&\sum_{k=0}^{r}(\mathcal{P}_{0}X_{k})e^{ik\theta}&-\sum_{k=1}^{r}(E_{0}X_{k}e^{ik\theta}-E_{-1}X_{k-1}e^{i(k-1)\theta})\\
& & \\
 & & +E_{-1}X_{r}e^{ir\theta}.
\end{array}
\end{equation}
  
Now, using the equality 
$$\sum_{j=0}^{n-1}e^{ij\theta}T^{j}\sum_{k=1}^{r}(E_{0}X_{k}e^{ik\theta}-E_{-1}X_{k-1}e^{i(k-1)\theta})=e^{i(n-1)\theta}T^{n-1}\sum_{k=1}^{r}E_{0}X_{k}e^{ik\theta}-\sum_{k=0}^{r-1}E_{-1}X_{k}e^{ik\theta}$$
we get, from (\ref{prolem11}), that
\begin{equation}
\label{prolem12}
\begin{array}{rll}
S_{n}(\theta)= & M_{r,n}(\theta)-( e^{i(n-1)\theta}T^{n-1}\sum_{k=1}^{r}E_{0}X_{k}e^{ik\theta}-\sum_{k=0}^{r-1}E_{-1}X_{k}e^{ik\theta})&
\\

 & & \\
+& \sum_{j=0}^{n-1}e^{ij\theta} T^{j}E_{-1}X_{r}e^{ir\theta}& 
\end{array}
\end{equation}
and that 
\begin{equation}
\label{prolem13}
\begin{array}{rll}
E_{0}S_{n}(\theta)= & D_{r,0}(\theta)-( E_{0}e^{i(n-1)\theta}T^{n-1}\sum_{k=1}^{r}E_{0}X_{k}e^{ik\theta}-\sum_{k=0}^{r-1}E_{-1}X_{k}e^{ik\theta})& \\
 & & \\
+& \sum_{j=0}^{n-1}e^{ij\theta}E_{0}T^{j}E_{-1}X_{r}e^{ir\theta} &  
\end{array}.
\end{equation}

(\ref{firappequ}) follows from (\ref{prolem12}) and (\ref{prolem13}) (see also Proposition \ref{relconexpkooope}).\qed

\medskip

The next step towards (\ref{appleminvpriavefreequ}) lies in the use of appropriate upper bounds for the terms at the right-hand side in (\ref{firappequ}): for a given $(n,r,\theta)\in \mathbb{N}\times\mathbb{N}^{*}\times [0,2\pi)$ let us denote by $A_{r,n}:[0,2\pi)\times\Omega\to \mathbb{C}$ and $B_{r,n}:[0,2\pi)\times\Omega\to\mathbb{C}$ the $\mathcal{B}\otimes{\mathcal{F}_{\infty}}-$measurable functions
\begin{equation}
\label{B}
A_{r,n}(\theta,\omega):=\sum_{k=1}^{r}  (T^{n-1}E_{0}X_{k}(\omega)-E_{0}T^{n-1}E_{0}X_{k}(\omega)) e^{ik\theta},
\end{equation}
\begin{equation}
\label{C}
B_{r,n}(\theta,\omega):=\sum_{k=0}^{n-1}(T^{k}E_{-1}X_{r}(\omega)-E_{0}T^{k}E_{-1}X_{r}(\omega))e^{ik\theta}.
\end{equation}
Then we have the following lemma.

\medskip

\begin{lemma}
\label{secappequ}
In the context of Lemma \ref{firapp}, and with the notation (\ref{B}) and (\ref{C}), there exists a constant $C>0$ such that, if $E_{0}$ is given by the regular version $E_{0}X(\omega)=E^{\omega}X$ ($X\in L^{1}_{\mathbb{P}}$) then
\begin{enumerate}
\item  For all  $(n,r,\theta)\in \mathbb{N}\times\mathbb{N}^{*}\times [0,2\pi)$,  $\alpha \in \mathbb{R}$, and $\omega\in \Omega$
\begin{equation}
\label{firappb}
E_{0}\left[\max_{k\leq n}\left|A_{r,k}(\theta,\cdot)\right|^2\right](\omega)\leq 4\alpha^2+ 4\sum_{j=0}^{n-1 }(T^{j}+E_{0}T^{j})|(E_{0}S_{r}(\theta))I_{[|E_{0}S_{r}(\theta)|>\alpha]}|^2(\omega).
\end{equation}

\item 
For all $\omega\in \Omega$
\begin{equation}
\label{firappc}
\int_{0}^{2\pi}{E}_{0}\left[\max_{k\leq n}\left|B_{r,k}(\theta,\cdot)\right|^2\right](\omega)\,d\lambda(\theta)\leq C\sum_{j=2}^{n-1}E_{{0}}|E_{j-1}X_{j+r}-E_{0}X_{j+r}|^{2}(\omega).
\end{equation}
\end{enumerate}
\end{lemma}

\medskip


{\bf Proof of Lemma \ref{secappequ}:} 
We will prove (\ref{firappb}) using a truncation argument: let $U_{\alpha}$ be the (non-linear) operator given by $U_{\alpha}Y:= Y I_{|Y|\geq \alpha}$, and fix the version of $E_{0}$ given by $E_{0}X(\omega)=E^{\omega}X$ ($X\in L^{1}_{\mathbb{P}}$), then for all $\omega\in \Omega$
$$\max_{k \leq n}| A_{r,k}(\theta,\cdot)|(\omega)=\max_{ k \leq n} |(Id-E_{0})(T^{k-1}E_{0}S_{r}(\theta))|^2(\omega)\leq$$
$$4\alpha^{2}+ 2\max_{k \leq n} |(Id-E_{0})T^{k-1}U_{\alpha}(E_{0}S_{r}(\theta))|^2(\omega)\leq$$ $$4(\alpha^2+\sum_{j=0}^{n-1}T^{j}|U_{\alpha}(E_{0}S_{r}(\theta))|^2(\omega)+\sum_{j=0}^{N-1}E_{0}T^{j}|U_{\alpha}(E_{0}S_{r}(\theta))|^2)(\omega),$$ 
where we used Jensen's inequality. This clearly implies (\ref{firappb}).

Let us now prove (\ref{firappc}): by Theorem \ref{hunyou} there exists a constant $C$ such that
$$\int_{0}^{2\pi}\max_{k\leq n}|B_{r,k}(\theta,z)|^{2}\,d\lambda(\theta)\leq C\int_{0}^{2\pi}|\sum_{j=2}^{n-1}(T^{j}E_{-1}X_{r}(z)-E_{0}T^{j}E_{-1}X_{r}(z))e^{ij\theta}|^{2}d\lambda(\theta)=$$
$$C\sum_{j=2}^{n-1}|E_{j-1}X_{j+r}(z)-E_{0}X_{j+r}(z)|^{2}.$$
  
The conclusion follows at once by integrating with respect to $E^{\omega}$ over these inequalities and using Tonelli's theorem.\qed

\subsection{Proof of Theorem \ref{invpriavefre}}

Under the hypothesis of Theorem \ref{invpriavefre}, if we can prove that there exists $\Omega_{0}\subset \Omega$ with $\mathbb{P}\Omega_{0}=1$ such that for all $\omega\in \Omega_{0}$, (\ref{appleminvpriavefreequ}) holds, then, combining this with the proof of Theorem \ref{quecltfoutra} (see Section \ref{proquecltfoutra}) and Lemma \ref{appleminvpriavefre}, the conclusion given in Theorem \ref{invpriavefre} will hold as well.
  
Let us do so: by Lemma \ref{firapp}, it is sufficient to prove that there exists $\Omega_{0}$ with $\mathbb{P}\Omega_{0}=1$ such that  if for $(k,r,\theta)\in \mathbb{N}\times\mathbb{N}^{*}\times [0,2\pi)$ we replace $Z_{r,k}(\theta,\omega):=A_{r,k}(\theta,\omega)$ or $Z_{r,k}(\theta,\omega):=B_{r,k}(\theta,\omega)$, then
\begin{equation}
\label{apphan}
\lim_{r} \limsup_{n} \int_{0}^{2\pi}{E}_{0}\left[\frac{1}{n}\max_{1\leq k \leq n} |Z_{r,k}(\theta,\cdot)|^2\right](\omega)\,d\lambda(\theta)=0.
\end{equation}
for all $\omega\in \Omega_{0}$.
\medskip

{\it Proof of (\ref{apphan}) with $Z_{r,k}(\theta,\omega):=A_{r,k}(\theta,\omega)$}: if we fix the version of $E_{0}$ given by $E_{0}X(\omega)=E^{\omega}X$ ($X\in L^{1}_\mathbb{P}$) then it is clear that  for any $\omega\in \Omega$

\begin{equation}
\label{majrhsb}
|E_{0}S_{r}(\theta)I_{[|E_{0}S_{r}(\theta)|>\alpha]}|(\omega)\leq \left|\left(\sum_{j=0}^{r-1}E_{{0}}|X_{j}|\right)I_{[\sum_{j=0}^{r-1}E_{0}|X_{j}|>\alpha]}\right|(\omega),
\end{equation}

and it follows by an application of Theorem \ref{ergthedisfoutra} (ergodic case), combined with Corollary \ref{dedmerpel} and (\ref{firappb}) (fixing first $\alpha>0$ so that the expectation of the random variable at the right  in (\ref{majrhsb}) is less than any fixed $\eta>0$), that
\begin{equation}
\label{btozer}
\lim_{n }{E}_{0}\left[\frac{1}{n}\max_{1\leq k \leq n}{|A_{r,k}(\theta,\cdot)|^{2}}\right]=0 \mbox{\,\,\,\,\,\,\,\, \it $\mathbb{P}-$ a.s.}
\end{equation}
Note that here the (probability one) set $\Omega_{0,1}$ of convergence does not depend on $\theta$ and, even more, the convergence is uniform in $\theta$ for any fixed $\omega\in \Omega_{0,1}$. It follows that for every $\omega\in \Omega_{0,1}$ 
$$\limsup_{n}\int_{0}^{2\pi}{E}_{0}\left[\frac{1}{n}\max_{1\leq k \leq n}{|A_{r,k}(\theta,\cdot)|^{2}}\right](\omega)\,d\lambda(\theta)\leq$$
$$ \int_{0}^{2\pi}\limsup_{n}{E}_{0}\left[\frac{1}{n}\max_{1\leq k \leq n}{|A_{r,k}(\theta,\cdot)|^{2}}\right](\omega)d\lambda(\theta)=0$$
as desired.

\medskip

{\it Proof of (\ref{apphan}) with $Z_{r,n}(\theta,\cdot):=B_{r,n}(\theta,\cdot)$}: again, fix the version of $E_{0}$ given by $E_{0}X(\omega)=E^{\omega}X$. We depart from (\ref{firappc}) and note that, if for every $j\in\mathbb{Z}$, $X_{-\infty,j}:=X_{j}-E_{-\infty}X_{j}$ then, by (\ref{equinvprotai})
$$\sum_{k=2}^{n-1}E_{0}|(E_{k-1}-E_{0})X_{k+r}|^{2}=\sum_{k=2}^{n-1}E_{0}|(E_{k-1}-E_{0})X_{-\infty,k+r}|^{2}=$$
$$\sum_{k=2}^{n-1}E_{0}T^{k-1}|(E_{0}-E_{-k+1})X_{-\infty,r+1}|^2=\sum_{k=1}^{n-2}(E_{0}T^{k}|E_{0}X_{-\infty,r+1}|^2-|E_{0}X_{-\infty,k+r+1}|^2)\leq$$
$$\sum_{k=1}^{n-2}E_{0}T^{k}|E_{0}X_{-\infty,r+1}|^2 $$
$\mathbb{P}-$a.s. It follows from (\ref{firappc}) and 
Corollary \ref{dedmerpel} that
  \begin{equation}
  \label{ctozer}
  \limsup_{n\to \infty }\int_{0}^{2\pi}{E}_{0}[\frac{1}{n}\max_{1\leq k \leq n}|B_{r,k}(\theta,\cdot)|^2](\omega)\,d\lambda(\theta)\leq C||E_{0}X_{-\infty,r+1}||_{_{\mathbb{P},2}}^2 =C||E_{-(r+1)}X_{-\infty,0}||_{_{\mathbb{P},2}}^2
  \end{equation}
  
$\mathbb{P}-$a.s. over a set $\Omega_{0,2,r}$ independent of $\theta$ and therefore, by the regularity condition (\ref{reg}), (see also (\ref{regconequ}))
  $$\lim_{r}\limsup_{n}\frac{1}{n}\int_{0}^{2\pi}{E}_{0}[\max_{1\leq k \leq n}|B_{r,k}(\theta,\cdot)|^2]\,d\lambda(\theta)=0$$
for all $\omega\in\Omega_{0,2}:=\cap_{r\in\mathbb{N}}\Omega_{0,2,r}$. 

To conclude, take $\Omega_{0}:=\Omega_{0,1}\cap\Omega_{0,2}$. \qed

\subsection{Proof of Theorem \ref{barinvprihancon}}

Let us start by recalling the following (Doob's) maximal inequality (\cite{revyor}, p.53): {\it if $p>1$ is given and $(M_{k})_{k\in \mathbb{N}^{*}}$ is a positive submartingale in $L^{p}_{\mu}$ then}
\begin{equation}
\label{doo}
||M_{n}||_{p,\mu}\leq||\max_{0\leq k\leq n}M_{k}||_{p,\mu}\leq\frac{p}{p-1}||M_{n}||_{p,\mu}.
\end{equation}
A combination of Doob's maximal inequality (\ref{doo}) with Corollary \ref{cormarundpome} gives the following result.

\medskip

\begin{lemma}
\label{doomaxl2que}
With the notation and conventions in page \pageref{gensetpro}, if $(M_{k})_{k\in\mathbb{N}^{*}}$ is a $(\mathcal{F}_{k-1})_{k\in\mathbb{N}^{*}}-$adapted martingale in $L^{2}_{\mathbb{P}}$ then
\begin{equation}
\label{condoo}
E_{0}[\max_{0\leq k\leq n}|M_{k}|]^2 \leq 4E_{0}|M_{n}|^2, \mbox{\,\,\,\,\,$\mathbb{P}$-\it a.s.}
\end{equation}
\end{lemma}

To prove Theorem \ref{barinvprihancon} we will need some additional estimates which will allow us to exploit the structure brought by (\ref{imhancon}).
\medskip

\begin{lemma}
\label{leminvprihancon}
In the context of Theorem \ref{barinvprihancon}, and under the notation on page \pageref{gensetpro}, consider the random variables ${B}_{n,r}(\theta,\cdot)$ given by (\ref{C}).
Then for all $(r,n,\theta)\in \mathbb{N}\times\mathbb{N}^{*}\times[0,2\pi)$,  
$$|1-e^{i\theta}|\left(E_{0}\left[\max_{k\leq n}\left|B_{r,k}(\theta,\cdot)\right|^2\right]\right)^{\frac{1}{2}} \leq $$
\begin{equation}
\label{secappc}
2\sum_{k=1}^{n-4}(\sum_{j=1}^{n-2}E_{0}T^{j}|\mathcal{P}_{0}(X_{k+r+1}-X_{k+r})|^{2})^{\frac{1}{2}} + (E_{0}|Y(n,r,\theta)|^{2})^{1/2}
\end{equation}
$\mathbb{P}-$a.s., where the residual $Y(n,r,\theta)$  is such that, under (\ref{reg}): 
\begin{equation}
\label{asynegerrc}
\lim_{r}\limsup_{n}\frac{1}{n}E_{0}|Y(n,r,\theta)|^{2}=0,\mbox{\,\,\,\,\,\it $\mathbb{P}$-a.s. 
}
\end{equation}

\end{lemma}

{\bf Proof:} We start by computing 
\begin{equation}
\label{brdec}
\begin{array}{rll}
(1-e^{i\theta})B_{n,r}(\theta)&=&(TE_{-1}X_{r}-E_{0}TE_{-1}X_{r})e^{i\theta} \\
& & \\
-(TE_{-1}X_{r}-E_{0}TE_{-1}X_{r})e^{2i\theta}&+& (T^{2}E_{-1}X_{r}-E_{0}T^{2}E_{-1}X_{r})e^{2i\theta}\\
 
  &\vdots &\\
-(T^{n-2}E_{-1}X_{r}-E_{0}T^{n-2}E_{-1}X_{r})e^{(n-1)i\theta}& +&(T^{n-1}E_{-1}X_{r}-E_{0}T^{n-1}E_{-1}X_{r})e^{(n-1)i\theta}\\
 &&\\
-(T^{n-1}E_{-1}X_{r}-E_{0}T^{n-1}E_{-1}X_{r})e^{in\theta}& =&\\
 \end{array}$$
$$-e^{i\theta}\sum_{k=0}^{n-2}(T^{k}(E_{-1}X_{r}-TE_{-1}X_{r})-E_0T^{k}(E_{-1}X_{r}-TE_{-1}X_{r}))e^{ik\theta}$$
$$-(T^{n-1}E_{-1}X_{r}-E_{0}T^{n-1}E_{-1}X_{r})e^{in\theta}.
\end{equation}

Let us stop now to make the following digression: assume that $Y_{0}\in L^2_{\mathbb{P}}$ is $\mathcal{F}_{0}-$measurable and let $Y_{j}:=T^{j}Y_{0}$ ($j\in\mathbb{Z}$) and  
$$S(Y_{0},n,\theta):=\sum_{k=0}^{n-1}Y_{k}e^{ik\theta}$$
the $n-$th discrete Fourier transform of $(Y_{j})_{j\in\mathbb{Z}}$.

Such $Y_{0}$ admits the decomposition (see (\ref{decl2f0}) and (\ref{relconexpkooopeequ}))
$$Y_{0}=\sum_{l=0}^{\infty}\mathcal{P}_{-l}Y_{0}+E_{-\infty}Y_{0}=\sum_{l=0}^{\infty}T^{-l}\mathcal{P}_{0}Y_{l}+E_{-\infty}Y_{0}.$$
Since these series are convergent in the $L^{2}_{\mathbb{P}}-$sense,  it follows that
$$E_{0}[S(Y_{0},n,\theta)]=\sum_{k=0}^{n-1}(\sum_{l=0}^{\infty}E_{0}T^{k}\mathcal{P}_{-l}Y_{0}+E_{-\infty}T^{k}Y_{0})e^{ik\theta}=$$
$$\sum_{k=0}^{n-1}(\sum_{l=k}^{\infty}E_{0}T^{k}\mathcal{P}_{-l}Y_{0}+E_{-\infty}T^{k}Y_{0})e^{ik\theta}=
\sum_{k=0}^{n-1}(\sum_{l=k}^{\infty}T^{k}\mathcal{P}_{-l}Y_{0}+E_{-\infty}T^{k}Y_{0})e^{ik\theta},$$
and it follows that
$$(Id-E_{0})S(Y_{0},n,\theta)=\sum_{k=0}^{n-1}T^{k}(\sum_{l=0}^{k-1}\mathcal{P}_{-l}Y_{0})e^{ik\theta}=\sum_{k=0}^{n-1}\sum_{l=0}^{k-1}T^{k-l}\mathcal{P}_{0}Y_{l}e^{ik\theta}= $$
\begin{equation}
\label{equstanotboo}
\sum_{k=0}^{n-1}\sum_{j=1}^{k}T^{j}\mathcal{P}_{0}Y_{k-j}e^{ik\theta}=\sum_{k=0}^{n-2}\sum_{j=1}^{n-k-1}T^{j}\mathcal{P}_{0}Y_{k}e^{i(k+j)\theta}.
\end{equation}

To continue towards the proof of (\ref{secappc}), apply (\ref{equstanotboo}) with $Y_{0}=(Id-T)E_{-1}X_{r}$ (so that $\mathcal{P}_{0}Y_{0}=-\mathcal{P}_{0}X_{r+1}$ and $\mathcal{P}_{0}Y_{k}=-\mathcal{P}_{0}(X_{r+k+1}-X_{r+k})$ for $k\geq 1$) to arrive at the identity
$$
(1-e^{i\theta})B_{r,n}(\theta)=-(T^{n-1}E_{-1}X_{r}-E_{0}T^{n-1}E_{-1}X_{r})e^{in\theta}+e^{i\theta}\sum_{j=1}^{n-1}T^{j}\mathcal{P}_{0}X_{r+1}e^{ij\theta}$$ $$+\,\,e^{i\theta}\sum_{k=1}^{n-4}\sum_{j=1}^{n-k-1}T^{j}\mathcal{P}_{0}(X_{r+k+1}-X_{r+k})e^{i(k+j)\theta}$$

so that, for a fixed $n\geq 4$
\begin{equation}
\label{almsecappc}
\max_{0\leq k\leq n}|(1-e^{i\theta})B_{k,r}(\theta)|\leq \sum_{j=1}^{n-4}\max_{1\leq k \leq n-2}|\sum_{l=1}^{n}T^{l}\mathcal{P}_{0}(X_{r+j+1}-X_{r+j})e^{il\theta}|+Y(n,r,\theta)
\end{equation}
where 
$$Y(n,r,\theta):=\max_{1\leq k \leq n}|T^{k-1}E_{-1}X_{r}-E_{0}T^{k-1}E_{-1}X_{r}|+\max_{1\leq k \leq n}|\sum_{j=1}^{k-1}T^{j}\mathcal{P}_{0}X_{r+1}e^{ij\theta}|$$
(to see that the ``max'' can be taken over $n\geq 1$ note that $B_{r,0}(\theta,\cdot)=0$).

We will prove (\ref{asynegerrc}) now. To do so we notice that by orthogonality under $E_{0}$ 
, Jensen's inequality, and Doob's maximal inequality (\ref{condoo}),
$$E_{0}(Y(n,r,\theta))^{2}\leq {8}\sum_{j=0}^{n-1}E_{0}T^{j}|E_{-1}X_{r}|^{2}+{8}E_{0}[|\sum_{j=1}^{n-1}T^{j}\mathcal{P}_{0}X_{r+1}e^{ij\theta}|^{2}]$$
$$={8}\sum_{j=0}^{n-1}E_{0}T^{j}|E_{-1}X_{r}|^{2}+{8}\sum_{j=1}^{n-1}E_{0}T^{j}|\mathcal{P}_{0}X_{r+1}|^2$$
$\mathbb{P}-$a.s, 
and we use Corollary \ref{dedmerpel} to conclude that
$$\limsup_{n}\frac{1}{n}E_{0}(Y(n,r,\theta))^{2}\leq 8||E_{-(r+1)}X_{0}||_{\mathbb{P},2}^{2}+8||\mathcal{P}_{-(r+1)}X_{0}||_{\mathbb{P},2}^{2}\leq 16||E_{-(r+1)}X_{0}||_{\mathbb{P},2}^{2}$$
{$\mathbb{P}$-a.s}. 
Using this (\ref{asynegerrc}) follows clearly from (\ref{reg}).

To finish the proof of (\ref{secappc}) we appeal to (\ref{almsecappc}) and we note that, by the conditional Minkowski's inequality, orthogonality and Doob's maximal inequality (\ref{condoo}):
$$(E_{0}[\sum_{j=1}^{n-4}\max_{1\leq k \leq n-2}|\sum_{l=1}^{k}T^{l}\mathcal{P}_{0}(X_{r+j+1}-X_{r+j})e^{il\theta}|]^2)^{1/2}\leq$$ 
$$\sum_{j=1}^{n-4}(E_{0}[\max_{1\leq k \leq n-2}|\sum_{l=1}^{k}T^{l}\mathcal{P}_{0}(X_{r+j+1}-X_{r+j})e^{il\theta}|]^2)^{1/2}\leq$$
$$2\sum_{j=1}^{n-4}(E_{0}[ |\sum_{l=1}^{n-2}T^{l}\mathcal{P}_{0}(X_{r+j+1}-X_{r+j})e^{il\theta}|]^2)^{1/2}=$$
$$2\sum_{j=1}^{n-4}(\sum_{l=1}^{n-2}E_{0}T^{l}|\mathcal{P}_{0}(X_{r+j+1}-X_{r+j})|^2)^{1/2}$$
$\mathbb{P}-$a.s. 
\qed

The following  lemma completes the box of tools needed to complete our proof of Theorem \ref{barinvprihancon}.

\medskip

\begin{lemma}
\label{conashan}
Under the condition (\ref{imhancon}), the series
\begin{equation}
\label{sersup}
\sum_{j=1}^{\infty}(\sup_{n}\frac{1}{n}\sum_{l=1}^{n}E_{0}T^{l}|\mathcal{P}_{0}(X_{j+r+1}-X_{j+r})|^2)^{\frac{1}{2}}
\end{equation}
converges $\mathbb{P}$-a.s.
\end{lemma}

{\bf Proof:} Remember that $L^{2,\infty}_{\mu}$ is a Banach space with the topology given by the norm $|||\cdot|||_{\mu,2}$ (see Section \ref{dunschweal2}). Thus by the inequality (\ref{equl2wl2}) the desired conclusion follows at once taking $Q=E_{0}T$ (see also (\ref{eque0tk})), $Y_{j}=\mathcal{P}_{0}(X_{j+r+1}-X_{j+r})$, 
and using condition (\ref{imhancon})
.\qed

\medskip

{\bf Proof of Theorem \ref{barinvprihancon}:} Our goal is to verify that, under the hypotheses of Theorem \ref{barinvprihancon}, the approximation (\ref{appleminvprifixfreequ}) in Lemma \ref{appleminvprifixfre} holds. Note that, by (\ref{firappequ}) and (\ref{btozer}), it suffices to prove that, under (\ref{imhancon}),
\begin{equation}
\label{apphanfixfre}
\lim_{r} \limsup_{n} {E}_{0}\frac{1}{n}\max_{1\leq k \leq n} |B_{r,k}(\theta,\cdot)|^2=0 \mbox{\,\,\,\,\,\, {\it $\mathbb{P}$-a.s.}}
\end{equation}
for every $\theta\neq 0$.    
    
{\bf Reduction:} We start from the following observation: by (\ref{relconexpkooopeequ}) and (\ref{equinvprotai}), the definition of $B_{r,n}(\theta,\cdot)$ remains unchanged if we replace $X_{r}$ by $X_{-\infty,r}:=X_{r}-E_{-\infty}X_{r}$. {\it Thus we can assume, without loss of generality, that $(X_{k})_{k\in\mathbb{Z}}$ is regular} (Definition \ref{regcondef}, see also (\ref{regconequ})).

{\it Proof of (\ref{apphanfixfre}) under the assumption of regularity (see the reduction above)}: it follows from (\ref{secappc}) that
$$\left(\frac{1}{n}E_{0}\left[\max_{k\leq n}\left|B_{r,k}(\theta,\cdot)\right|^2\right]\right)^{\frac{1}{2}}\leq$$
$$ |1-e^{i\theta}|^{-1}(2\sum_{j=1}^{n-4}(\frac{1}{n}\sum_{l=1}^{n}E_{0}T^{l}|\mathcal{P}_{0}(X_{j+r+1}-X_{j+r})|^{2})^{\frac{1}{2}}+(\frac{1}{n}E_{0}(Y(n,r,\theta))^2)^{1/2}),$$
$\mathbb{P}$-a.s. 

So from (\ref{asynegerrc}), Lemma \ref{conashan}, the dominated convergence theorem, and Corollary \ref{dedmerpel}  we get that
\begin{equation}
\label{ctozerhan1}
\limsup_{n\to \infty} \left(\frac{1}{n}E_{0}\left[\max_{k\leq n}\left|B_{r,k}(\theta,\cdot)\right|^2\right]\right)^{\frac{1}{2}}\leq 2|1-e^{i\theta}|^{-1}\sum_{j=r+1}^{\infty}||\mathcal{P}_{0}(X_{j+1}-X_{j})||_{_{\mathbb{P},2}}+o_{r}(1)
\end{equation}
{$\mathbb{P}$-a.s.} and therefore, by condition (\ref{imhancon}) again\footnote{Note that (\ref{imhancon}) was used already when applying Lemma \ref{conashan}.}
\begin{equation}
\label{ctozer2}
\lim_{r}\limsup_{n} \left(\frac{1}{n}E_{0}\left[\max_{k\leq n}\left|B_{r,k}(\theta,\cdot)\right|^2\right]\right)^{\frac{1}{2}}=0
\end{equation}
$\mathbb{P}$- a.s., as desired.\qed

\medskip

\begin{remark}
\label{ermconbfixfre}
The set $\Omega_{0}$ of convergence in the last statement can be chosen independent of $\theta$ (this requires some care, but the general strategy is to use the representation $\omega\mapsto E^{\omega}X$ of $E[X|\mathcal{F}_{0}]$ along all of the arguments). It is not clear, on the other side, whether the convergence is uniform in $\theta$ for a fixed $\omega\in \Omega_{0} $ (due to the factor $|1-e^{i\theta}|^{-1}$). Contrast this with the statement following (\ref{btozer}). 
 \end{remark}

\subsection{A Note on Theorem \ref{limlimlem}}
\label{notlimlimlem}

The proofs presented along this chapter can be carried out using the following (more restrictive) classical version of Theorem \ref{limlimlem} (see Theorem 3.2 in \cite{bilconpromea} for a proof).

\medskip

\begin{thm}[Transport Theorem]
\label{limlimlemcla}
{\it Let $(S,{d})$ be a complete and separable metric space. Assume that for all natural numbers $r,n$, $X_{r,n}$  and $X_{n}$ are random elements of $S$ defined on the same probability space $(\Omega,\mathcal{F},\mathbb{P})$, that $X_{r,n}\Rightarrow_{n} Z_{r}$, and that $Z_{r}\Rightarrow _{r} X$. Then the hypothesis
\begin{equation}
\label{liminf}
\lim_{r}\limsup_{n}\mathbb{P}[d(X_{r,n},X_{n})\geq \epsilon]=0 \mbox{\,\,\,\,\,\,\,\it for all $\epsilon>0$,}
\end{equation}
implies that $X_{n}\Rightarrow_{n}X$.}
\end{thm}

The adaptation of the arguments above to the further restriction imposed by this theorem poses  no serious difficulty: it suffices to see that the martingales $D_{r,0}$ converge in the appropriate way, as $r\to\infty$, for each one of the theorems proved along this section, and in particular to use the results presented in Chapter \ref{marcas} to deduce the asymptotic distributions associated to the processes generated via the martingale differences $D_{0}=\lim_{r}D_{r,0}$. The details are left to the reader.

\chapter{Proofs Related to the Random Centering}
\label{chapronecrancen}

In this chapter we present proofs of the results stated, but not proved, in Section \ref{necrancen}. 

Section \ref{secprothequecondes}, devoted to the proof of Theorem \ref{quecondes}, addresses the problem of the meaning of convergence under $\mathbb{P}_{\omega}$, for a fixed (and appropriately chosen) $\omega$, of a stochastic process $(Y_{n})_{n}$ for which $Y_{n}-E_{0}Y_{n}$ converges in the quenched sense. The main novelty is Proposition \ref{queconrev}, a relaxed form of Theorem \ref{quecondes} from which this result follows easily.

Section \ref{thefoutraofalin} can be considered as pertaining to an exposition on the general theory of quenched convergence, but we present it in this part of the monograph because, first, it can be considered as a relatively straightforward specialization of the results presented in Chapter \ref{resandcom}, and second, the exposition is written with the purpose of giving a precise interpretation of the processes there considered (adapted linear processes) under the light of the hypotheses in theorems \ref{cltpelwu} and \ref{quecltfoutra}, paving the way to the arguments present in Section \ref{prothenonquecon} towards the proof of Theorem \ref{nonqueconthe}. The main result in this section is Proposition \ref{asyvart}.

Section \ref{prothenonquecon} presents the proof of Theorem \ref{nonqueconthe}. We start with two approximation lemmas (lemmas \ref{lemxij} and \ref{lemest}) and then proceed to give an instance of the process announced by Theorem \ref{nonqueconthe} by following a construction due to Voln\'{y} and Woodroofe.

Finally, in Section \ref{prothenonrancen}, we present the proofs of theorems \ref{queconznundratdec} and \ref{queconznmaxwoo}. Besides presenting these proofs, this section aims to illustrate how the techniques involved in the proofs of some previous results in this monograph can  (and should) be used to expand the family of theorems on quenched asymptotics for the Discrete Fourier Transforms of dependent sequences by combining the estimates present in the existing literature for non-rotated partial sums with the martingale limit theorems developed along this work. 

The notation is that introduced at the beginning and on page \pageref{gensetpro}.    

\section{Proof of Theorem \ref{quecondes}}
\label{secprothequecondes}

To begin with, suppose that we know that the integrable process $(Y_{n})_{n\in\mathbb{N}}$ is such that that $Y_{n}-E_{0}Y_{n}\Rightarrow Y$ in the quenched sense. Thus (Proposition \ref{prounicon}) there exists $\Omega_{0}\subset \Omega$ with $\mathbb{P}\Omega_{0}=1$ such that for every $\omega\in \Omega_{0}$, $Y_{n}-E_{0}Y_{n}\Rightarrow Y$ under $\mathbb{P}_{\omega}$. 

{\bf Question:} {\it What are the possible limit laws for $(Y_{n})_{n\in\mathbb{N}}$ under $\mathbb{P}_{\omega}$ for a fixed $\omega$?}

To answer this question we depart from the following auxiliary result. Remember that, for a metric space $S$, $\mathbf{C}^{b}(S)$ denotes the space of functions $h:S\to\mathbb{R}$ that are continuous and bounded.
\medskip

\begin{lemma}
\label{barpellem18}
Under the setting introduced in page \pageref{gensetpro}, if $X$ is $\mathcal{F}-$measurable and $Z$ is $\mathcal{F}_{0}-$measurable\footnote{These are {\it not} $\mathbb{P}-$equivalence classes, but actual ``versions'' of $X$ and $Z$.}, there exists $\Omega(X,Z)\subset \Omega$ with $\mathbb{P}(\Omega(X,Z))=1$ such that, for every $h\in\mathbf{C}^{b}({\mathbb{C}^{2}})$ 
\begin{equation}
\label{lemconexp}
E^{\omega}[h(X,Z)]=E^{\omega}[h(X,Z(\omega))]
\end{equation}
for all $\omega\in \Omega(X,Z)$.
\end{lemma}

{\bf Proof:} Let $A\in\mathcal{F}_{0}$ be given, and assume first that $Z=I_{A}$. Then, clearly
$$h(X,Z)=h(X,1)I_{A}+h(X,0)I_{\Omega\setminus A}$$
and therefore
$$E_{0}h(X,Z)=E_{0}[h(X,1)]I_{A}+E_{0}[h(X,0)]I_{\Omega\setminus A}$$
$\mathbb{P}-$a.s., which implies, via the representation $E_{0}Y(\omega)=E^{\omega}Y$ ($Y\in L^{1}_{\mathbb{P}}$), that there exists $\Omega(h,X,Z)$ with $\mathbb{P}\Omega(h,X,Z)=1$ such that (\ref{lemconexp}) holds for every $\omega\in\Omega(h,X,Z)$. This argument can be easily extended to the case of simple functions $Z=\sum_{k=1}^{n}a_{k}I_{A_{k}}$ with $A_{k}\in\mathcal{F}_{0}$ ($k=1,\dots, n$). 

Now assume that $Z$ is an arbitrary $\mathcal{F}_{0}-$measurable function,  let $(Z_{n})_{n\in\mathbb{N}}$ be a sequence of simple functions with $\lim_{n}Z_{n}(\omega)=Z(\omega)$ (\cite{bilpromea}, p.254), and let $$\Omega(h,X,Z):=\cap_{n\in\mathbb{N}}\Omega(h,X,Z_{n}).$$

Clearly $\mathbb{P}\Omega(h,X,Z)=1$, and the dominated convergence theorem, together with the definition of $\Omega(h,X,Z)$ and the continuity of the bounded function $h$ imply that for every $\omega\in \Omega(h,X,Z)$
$$E^{\omega}(h(X,Z))=\lim_{n}E^{\omega}h(X,Z_{n})=\lim_{n}E^{\omega}h(X,Z_{n}(\omega))=E^{\omega}h(X,Z(\omega)).$$

Finally, let $(h_{n})_{n\in\mathbb{Z}}$ be a family of Urysohn functions as in the statement $2.$ of Theorem \ref{refporthesta}, and let $$\Omega(X,Z):=\cap_{n\in\mathbb{N}}\Omega(h_{n},X,Z).$$ 
Again, $\mathbb{P}\Omega(X,Z)=1$, and by Theorem \ref{refporthesta} (replacing  $X_{n}$ by $(X,Z)$ for all $n\in\mathbb{N}$ and $X$ by $(X,Z(\omega))$) and the definition of $\Omega(X,Z)$, for  every $\omega\in \Omega(X,Z)$ and every $h\in \mathbf{C}^{b}(\mathbb{C}^{2})$, $E^{\omega}h(X,Z)=E^{\omega}h(X,Z(\omega))$. \qed

This lemma, in combination with Proposition \ref{procontypcom} gives the following result.

\medskip

\begin{prop}[Possible Limit Laws for a Fixed Starting Point]
\label{queconrev}
With the notation of Lemma \ref{barpellem18}, assume that $(Z_{n})_{n}$ is a sequence of functions in $L^{1}_{\mathbb{P}}$ such that $Z_{n}-E_{0}Z_{n} \Rightarrow_{n} Y$ under $\mathbb{P}_{\omega}$ for all $\omega\in \Omega_{0,1}$ ($\Omega_{0,1}\subset \Omega$ is {\it any} given set, not even assumed measurable), and let $\Omega_{0,2}:=\cap_{n}\Omega(Z_{n},E_{0}Z_{n})$, where $\Omega(Z_{n},E_{0}Z_{n})$ is the set specified in the conclusion of Lemma \ref{barpellem18} \footnote{Of course, here we are implicitely fixing versions of $Z_{n}$ and $E_{0}Z_{n}$. We will leave these details to the reader.}. Then, given $\omega\in \Omega_{0}:=\Omega_{0,1}\cap\Omega_{0,2}$, $Z_{n}\Rightarrow Z_{\omega}$ under $\mathbb{P}_{\omega}$ if and only if $L(\omega)=\lim_{n\to\infty}E_{0}Z_{n}(\omega)$ exists, in which case 
\begin{equation}
\label{quelimcenvsnocen}
Z_{\omega}=Y+ L(\omega)
\end{equation}
(in distribution).
\end{prop}

\medskip


{\bf Proof:} Given $\omega\in \Omega_{0}$ and any bounded and continuous function $h:\mathbb{C}\to\mathbb{R}$, the hypotheses imply that
$$E^{\omega}h(Z_{n}-E_{0}Z_{n}(\omega))=E^{\omega}h(Z_{n}-E_{0}Z_{n})\to Eh(Y),$$
as $n\to\infty$, so that  $Z_{n}-E_{0}Z_{n}(\omega)\Rightarrow Y$ under $\mathbb{P}_{\omega}$ (Portmanteau's Theorem). From $Z_{n}=Z_{n}-E_{0}Z_{n}(\omega)+E_{0}Z_{n}(\omega)$ the conclusion follows  via Proposition \ref{procontypcom} in page \pageref{procontypcom}.\qed

\medskip




We can proceed now to the proof of Theorem \ref{quecondes}
\medskip


{\bf Proof of Theorem \ref{quecondes}:} 
We  appeal to Proposition \ref{queconrev}, replacing $Z_{n}$ by $Z_{n}(\theta)$ and taking 
$$\Omega_{\theta,1}:=\{\omega\in \Omega: Y_{n}(\theta)=Z_{n}(\theta)-E_{0}Z_{n}(\theta)\Rightarrow Y({\theta}) \mbox{\, \,\,\it under  $\mathbb{P}_{\omega}$}\}$$
and $\Omega_{\theta,2}= \cap_{n\in\mathbb{N}}\Omega(Z_{n}(\theta),E_{0}Z_{n}(\theta))$. This gives that,  for  any $\omega\in \Omega_{\theta}:=\Omega_{\theta,1}\cap\Omega_{\theta,2}$ fixed, there exists $L_{\theta}(\omega):=\lim_{n}E_{0}Z_{n}(\theta,\omega)$ and therefore, if $Z_{n}(\theta)\Rightarrow_{n} Z_{\omega}(\theta)$ under $\mathbb{P}_{\omega}$, $Z_{\omega}({\theta})=Y({\theta})+L_{\theta}(\omega)$ (in distribution) under $\mathbb{P}_{\omega}$. 
With this observation, the conclusion follows easily from Proposition \ref{procontypcom} and the fact that $\mathbb{P}\Omega_{\theta}=1$.\qed

As explained at the end of Section \ref{secnecrancen} (see the ``General Comments''), it follows from Corollary \ref{sufnecconqueconzn} that if we can provide an example of a regular process $(X_{n})_{n}$ (Definition \ref{regcondef}) satisfying the hypothesis of Theorem \ref{quecltfoutra} for which
\begin{equation}
\label{protoprov}
\mathbb{P}\left(\limsup_{n\to \infty}\left|\frac{1}{\sqrt{n}}E_{0}S_{n}(\theta)\right|>0\right)>0\mbox{\,\,\,\, \it for   $\theta$ in a set $I'$ with $\lambda(I')>0$} 
\end{equation}
we will prove, in particular, the necessity of the random centering ``$-E_{0}Z_{n}(\theta)$'' for a nonnegligible subset of $I$ (namely $I\cap I'$).

In their paper \cite{volwoononque}, Voln\'{y} and Woodroofe provide an example of a sequence $(X_{n})_{n}$ for which a quenched CLT holds for $(Y_{n}(0))_{n}$ but not for $(Z_{n}(0))_{n}$. We will adapt their construction to give an example satisfying (\ref{plimsupinfone}) for every $\theta\in[0,2\pi)$ (this clearly implies (\ref{protoprov})). By the discussion at the end of Section \ref{secnecrancen} again, this will make the proof of Theorem \ref{nonqueconthe}.

The main novelty adapting the example in \cite{volwoononque}, which arises from a careful construction of a sequence $(a_{n})_{n\in\mathbb{N}}$ of nonnegative coefficients of a linear process is that, in order to guarantee the validity of the ``inductive step'' defining $a_{n+1}$ from $a_{1},\dots, a_{n}$, one needs to prove that a certain type of convergence is uniform in $\theta$ (see Lemma \ref{lemxij} below). While it would be sufficient to prove this uniform convergence for $\theta$ in (for instance) an open subinterval $I'$ of $[0,2\pi)$ in order to construct a valid example, a compactness argument allows us to do it for $I'=[0,2\pi)$.


\section{Theorem \ref{quecltfoutra} for Adapted Linear Processes}
\label{thefoutraofalin}

Let $(a_{k})_{k\in\mathbb{N}}\in l^{2}(\mathbb{N})$ be given, thus $a_{k}\in\mathbb{C}$ for every $k\in\mathbb{N}$, and $\sum_{k}{|a_{k}|^{2}}<\infty$. As explained in Section \ref{thefoutraofanint}, Carleson's theorem (Theorem \ref{carthe}) guarantees the convergence, for  $\lambda-$a.e $\theta$,  of the  series
$$\sum_{j\geq 0} a_{j}e^{ij\theta}$$
and therefore there exists a well defined ($\lambda-$equivalence class of) function(s) $f:[0,2\pi)\to\mathbb{C}$ given  by
\begin{equation}
\label{detpart}
\theta\mapsto f({\theta})=\lim_{n}\sum_{j=0}^{n-1} a_{j}e^{ij\theta}            
\end{equation}
(see Definition \ref{deflimfun}) and $f(\theta)$, thus defined, is a $2\pi-$periodic function, square integrable over $[0,2\pi)$, and satisfying $\hat{f}(n)=a_{n}$ for every $n\in\mathbb{N}$, where $\hat{f}$ denotes the Fourier transform (\ref{deffoutra}).

For every $k\in\mathbb{Z}$, denote by 
\begin{equation}
\label{deffkequ}
f_{k}(\theta):=\sum_{j=0}^{k-1}a_{j}e^{ij\theta}
\end{equation}
(thus $f_{k}=0$ if $k\leq 0$), and consider the setting explained along Example \ref{exaregconexpexa} in page \pageref{exaregconexpexa}. If we regard $(a_{k})_{k\in\mathbb{N}}$ as an element of $l^{2}(\mathbb{Z})$ with $a_{k}=0$ for $k<0$ then, as  explained in Example \ref{exalinpro},  the process $(X_{k})_{k\in\mathbb{Z}}$ given by (\ref{deflinpro}) is $(\mathcal{F}_{k})_{k\in\mathbb{Z}}-$adapted. Since in this setting $(x_{k})_{k\in\mathbb{Z}}$, the sequence of coordinate functions, is i.i.d.,  Kolmogorov's zero-one law (\cite{bilpromea}, Theorem 22.3) implies that for every $k\in\mathbb{Z}$
$$E[X_{k}|\mathcal{F}_{-\infty}]=E[X_{k}]=0$$
and it follows  that $(X_{k})_{k\in\mathbb{Z}}$ is regular (see \ref{regconequ}). 

In conclusion, these processes satisfy the hypotheses of Theorem \ref{quecltfoutra} and are regular. This will be the basis for the construction of the example stated by Theorem \ref{nonqueconthe}.  
   
To begin our discussions, start by noting that, in the context just explained, we have the following two expressions for the Fourier Transforms $S_{n}(\theta)$ (Definition \ref{defdisfoutra}) of $(X_{k})_{k\in\mathbb{Z}}$:

\begin{equation}
\label{snexp1}
S_{n}(\theta)=\sum_{j=-\infty}^{n-1}(f_{-j+n}-f_{-j})(\theta)x_{j}e^{ij\theta},
\end{equation}

\begin{equation}
\label{snexp2}
S_{n}(\theta)=\sum_{k=0}^{\infty}a_{k}\sum_{j=0}^{n-1}e^{ij\theta}x_{j-k}.
\end{equation}

in the $L^{2}_{\mathbb{P}}-$sense\footnote{The changes in the order of summation involved can be easily justified in this case. We skip this detail.}. Now let us denote, for all $k\in \mathbb{Z}$,  

\begin{equation}
\label{defzeta}
\zeta_{-k}(\theta):=\sum_{j=0}^{k}e^{-ij\theta}x_{-j} 
\end{equation}

(note that $\zeta_{-k}=0$ if $k<0 $). Then from  (\ref{snexp1}) and (\ref{snexp2}) the following two equalities follow respectively: 

\begin{equation}
\label{exprojsn1}
E_{0}S_{n}(\theta)=\sum_{j\in \mathbb{N}}x_{-j}(f_{j+n}-f_{j})(\theta)e^{-ij\theta},
\end{equation}

\begin{equation}
\label{exprojsn2}
E_{0}S_{n}(\theta)=\sum_{j\in \mathbb{N}} a_{j}(\zeta_{-j}-\zeta_{-j+n})(\theta)e^{ij\theta}.
\end{equation}

In particular
$$E_{0}|S_{n}(\theta)-E_{0}S_{n}(\theta)|^2=E_{0}|\sum_{j=1}^{n-1}e^{ij\theta}x_{j}f_{-j+n}(\theta)|^2=$$
$$||x_{0}||_{_{\mathbb{P},2}}^2\sum_{j=1}^{n-1}|f_{n-j}(\theta)|^2
=||x_{0}||_{_{\mathbb{P},2}}^2\sum_{j=1}^{n-1}|f_{j}(\theta)|^2,$$
so that, by Theorem \ref{quecltfoutra}, there exists $I\subset (0,2\pi)$ with $\lambda(I)=1$ such that for every $\theta\in I$
$$\lim_{n}\frac{E_{0}|S_{n}(\theta)-E_{0}S_{n}(\theta)|^2}{n}=\lim_{n\to \infty}\frac{1}{n}\sum_{j=1}^{n-1}||x_0||_{_{\mathbb{P},2}}^{2}|\, f_{j}(\theta)|^2={|\,||x_0||_{_2} f(\theta)|^2}.$$

Using this fact, we get the following version of Theorem \ref{quecltfoutra}:

\bigskip

\begin{prop}[Theorem \ref{quecltfoutra} for Adapted Linear Processes]
\label{asyvart}
With the notations and the setting in Example \ref{exalinpro} and Example \ref{exaregconexpexa}, and given a $(\mathcal{F}_{k})_{k\in\mathbb{Z}}-$adapted linear process (\ref{deflinpro}) (thus $a_{k}=0$ if $k<0$), there exists $I\subset (0,2\pi)$ such that for every $\theta\in I$, $(S_{n}(\theta)-E_{0}S_{n}(\theta))/\sqrt{n}$ is asymptotically normally distributed under $\mathbb{P}_{\omega}$ for $\mathbb{P}$-a.e $\omega$ $(\mathbb{P}_{\omega}$ is given by (\ref{pomeequfuniid})), and its asymptotic distribution  (under $\mathbb{P}_{\omega}$) corresponds  to a complex-valued normal random variable with independent real and imaginary parts, each with mean zero and variance
$$\sigma_{\theta}^2=\frac{|\,||x_0||_{_{\mathbb{P},2}} f(\theta)|^2}{2},$$
where $f$ is given by (\ref{detpart}).
\end{prop}

\section{Proof of Theorem \ref{nonqueconthe}}
\label{prothenonquecon}

We finally address here the construction leading to the proof of Theorem \ref{nonqueconthe}.  The notation along the following arguments is borrowed from the previous sections in this chapter. In particular, $\zeta_{n}$ is defined by (\ref{defzeta}) for every $n\in\mathbb{Z}$, and $T$ is the left shift in $\mathbb{R}^{\mathbb{Z}}$ which, under the setting of Example \ref{exaregconexpexa}, is weakly mixing (see \cite{qua}, p.13). 

By \cite{cunmerpel}, p.4075 (Section 4.1) applied to the sequence $(\delta_{1j})_{j\in \mathbb{Z}}$ ($\delta_{ij}$ denotes the Kronecker $\delta-$function) and the fact that $T$ is weakly mixing, the following law of the iterated logarithm holds\footnote{Note that, for the linear process $(x_{n})_{n\in \mathbb{Z}}$ (the coordinate functions), which corresponds to convolution with the sequence $(\delta_{1j})_{j\in\mathbb{Z}}\in l^{2}(\mathbb{Z})$, the spectral density with respect to Lebesgue measure is the constant function $||x_{0}||_{_{\mathbb{P},2}}^{2}/2\pi$.}: for every $t\in (0,2\pi)\setminus\{\pi\}$
\begin{equation}
\label{lilperlin}
\limsup_{n\to \infty}\frac{|\zeta_{-n}(\theta)|}{\sqrt{n \log \log n}}=||x_{0}||_{_{\mathbb{P},2}}.
\end{equation}

$\mathbb{P}-$almost surely.

If $\theta=0$ or $\theta=\pi$, and  $x_{0}$ is real-valued and symmetric ($\mathbb{P}[x_{0}\leq t]=\mathbb{P}[x_{0}\geq -t]$), the L.I.L. as stated above holds with $||x_{0}||_{_{\mathbb{P},2}}$ replaced by $\sqrt{2}\,||x_{0}||_{_{\mathbb{P},2}}$ (\cite{bilpromea}, Theorem 9.5). Assume this from now on.

The equality (\ref{lilperlin}) clearly implies that $$\limsup_{n}\frac{|\zeta_{-n}(\theta)|}{\sqrt{n}}= \infty$$
 {$\mathbb{P}$-a.s}. The following lemma states that the divergence occurs ``at comparable speeds'' for every $\theta$. 

\medskip

\begin{lemma}
\label{lemxij}
Consider  $\zeta_{-k}$ as defined by (\ref{defzeta}). Then for every $\lambda\in \mathbb{R}$ and every $0<\eta\leq 1$ there exists $N\in\mathbb{N}$ satisfying

$$\mathbb{P}\left(\max_{1\leq n\leq N}\frac{|\zeta_{-n}(\theta)|}{\sqrt{n}}> \lambda\right)\geq 1-\eta $$

for all $\theta\in [0,2\pi)$. In particular

$$\mathbb{P}\left(\max_{1\leq n\leq m}\frac{|\zeta_{-n}(\theta)|}{\sqrt{n}}\geq \lambda\right)\geq 1-\eta $$

for all $m\geq N$.
\end{lemma}

{\bf Proof:} Fix $\lambda\in \mathbb{R}$ and $0<\eta\leq 1$. Let\footnote{We work over the interval $[0,2\pi]$ (instead of $[0,2\pi)$). This has no effect for the validity of the conclusion and is assumed in order to take advantage of compactness, as will be clear along the proof.} $\theta\in [0,2\pi]$ and $\epsilon>0$ be given and define

$$E_{\epsilon,m}(\theta):=\left[\inf_{|\delta|<\epsilon}\left\{\max_{1\leq n\leq m}\frac{|\zeta_{-n}(\theta+\delta)|}{\sqrt{n}}\right\}>\lambda\right]$$ 
and

$$E_{m}(\theta):=\left[\max_{1\leq n\leq m}\frac{|\zeta_{-n}(\theta)|}{\sqrt{n}}>\lambda\right] .$$

Note that, for fixed  $m$, the sequence of sets $E_{\epsilon,m}(\theta)$ is decreasing with respect to $\epsilon$ ($\epsilon_{1}<\epsilon_{2}$ implies that $E_{\epsilon_{2},m}(\theta)\subset E_{\epsilon_{1},m}(\theta)$), and that the (random) function $\theta\mapsto \max_{1\leq n\leq m}{|\zeta_{-n}(\theta)|}/{\sqrt{n}}$ is continuous for all $m$. In particular

\begin{equation}
\label{decsets}
\bigcup_{\epsilon>0}E_{\epsilon,m}(\theta)=E_{m}(\theta),
\end{equation}

where the sets in the union increase as $\epsilon$ decreases to $0$.

Now, there exists a minimal $N(\theta)$ such that 
$ \mathbb{P}(E_{N({\theta})}(\theta))>1-\eta$. To see this note that the family $\{E_{k}(\theta)\}_{k\geq 0}$ is increasing with $k$, and its union contains the set 
$$[\limsup_{n}{|\zeta_{-n}(\theta)|}/{\sqrt{n}}>\lambda]$$
which has $\mathbb{P}-$measure $1$ by (\ref{lilperlin}).

It follows from (\ref{decsets}) that there exists an $\epsilon_{\theta}$ such that
\begin{equation}
\label{eqept}
\mathbb{P}(E_{\epsilon_{\theta},N(\theta)}(\theta))>1-\eta \,.
\end{equation}
Now, the family of sets $\{(\theta-\epsilon_{\theta}, \theta+\epsilon_{\theta})\}_{\theta\in [0,2\pi]}$ is an open cover of $[0,2\pi]$, and therefore it admits an open subcover $\{(\theta_{j}-\epsilon_{j}, \theta_{j}+\epsilon_{j})\}_{j=1}^{r}$ where $\epsilon_{j}:=\epsilon_{\theta_{j}}$. Let $N=\max\{N(\theta_{1}),\dots,N(\theta_{r})\}$. We claim that for every $\theta\in [0,2\pi]$
$$\mathbb{P}(E_{N}(\theta))> 1-\eta\,.$$
Indeed, given $\theta\in[0,2\pi]$, with $\theta_{j}-\epsilon_{j}<\theta<\theta_{j}+\epsilon_{j}$,

$$E_{N}(\theta)\supset E_{N(\theta_{j})}(\theta)= \left[\max_{1\leq n \leq N(\theta_{j})}\frac{|\zeta_{-n}(\theta_{j}+\theta-\theta_{j})|}{\sqrt{n}}>\lambda\right]\supset E_{\epsilon_{j},N(\theta_{j})}(\theta_{j})\,,$$

and the conclusion follows from (\ref{eqept}) and the definition of $E_{N}(\theta)$. \qed

\bigskip

 Let us now move to the following observation: if $(n_{k})_{k\in \mathbb{N}}$ is a strictly increasing sequence of natural numbers and if $(a_{k})_{k\in\mathbb{N}}$ is square summable and satisfies $a_{j}=0$ if $j\notin \{n_{k}\}_{k}$  then, using (\ref{exprojsn2}) we have, for every given $k\in\mathbb{N}$,
\begin{equation}
\label{parsnt}
E_{0}S_{n}(\theta)=\sum_{j=0}^{k}e^{in_{j}\theta}a_{n_{j}}(\zeta_{-n_{j}}-\zeta_{-n_{j}+n})(\theta)+ \sum_{j=k+1}^{\infty}e^{in_{j}\theta}a_{n_{j}}(\zeta_{-n_{j}}-\zeta_{-n_{j}+n})(\theta)=:$$
$$A_{k}(n,\theta)+B_{k}(n,\theta)
\end{equation}
so that 

$$\mathbb{P}\left(\max_{n_{k-1}<n\leq n_{k}}\frac{|E_{0} S_{n}(\theta)|}{\sqrt{n}}\geq 2^{k}\right)\geq$$
$$ \mathbb{P} \left(\max_{n_{k-1}<n\leq {n_{k}}}\frac{|A_{k}(n,\theta)|}{\sqrt{n}} \geq 2^{k+1}\right)- \mathbb{P}\left(\max_{n_{k-1}<n\leq {n_{k}}}\frac{|B_{k}(n,\theta)|}{\sqrt{n}}\geq 2^{k}\right)
\geq $$

\begin{equation}
\label{finestsn}
\mathbb{P} \left(\max_{n_{k-1}<n\leq {n_{k}}}\frac{|A_{k}(n,\theta)|}{\sqrt{n}} \geq 2^{k+1}\right)- \mathbb{P}\left(\max_{n_{k-1}<n\leq {n_{k}}}|B_{k}(n,\theta)|\geq 2^{k}\right).
\end{equation}

Now, if $n_{k-1}< n \leq n_{k}$ then, actually

$$A_{k}(n,\theta)= \sum_{j=0}^{k-1}e^{in_{j}\theta}a_{n_{j}}\zeta_{-n_{j}}(\theta)+e^{in_{k}}a_{n_{k}}(\zeta_{-n_{k}}-\zeta_{-n_{k}+n})(\theta).$$
The first summand at the right-hand side in this expression is bounded by 
$$\sum_{j=1}^{k-1}\sum_{r=0}^{n_j}|a_{n_{j}}||\xi_{-r}|$$

and therefore there exists $\gamma_{k}>0$ such that

\begin{equation}
\label{firestakn}
\mathbb{P}\left(\left|\sum_{j=0}^{k-1}e^{in_{j}\theta}a_{n_{j}}\zeta_{-n_{j}}(\theta)\right|>\gamma_{k}\right)\leq \left(\frac{1}{2}\right)^{k+2}
\end{equation}

for all $\theta\in [0, 2\pi]$.

All together (\ref{parsnt}), (\ref{finestsn}) and (\ref{firestakn}) give the following result.

\bigskip

\begin{lemma}
\label{lemest}
Let $(n_{k})_{k\in\mathbb{N}}$ be a strictly increasing sequence of natural numbers and  let $(a_{j})_{j\in\mathbb{N}}\in l^{2}(\mathbb{N})$ be a square summable sequence of nonnegative numbers with $a_{j}=0$ for $j\notin \{n_{k}\}_{k}$. Then for every sequence of real numbers $(\gamma_{k})_{k}$ satisfying (\ref{firestakn}) the following inequality holds

$$\mathbb{P}\left(\max_{n_{k-1}< n \leq n_{k}}\frac{|E_{0}S_{n}(\theta)|}{\sqrt{n}}\geq 2^{k}\right)\geq \mathbb{P}\left({a_{n_{k}}\max_{n_{k-1}< n \leq n_{k}}\frac{|(\zeta_{-n_{k}}-\zeta_{-n_{k}+n})(\theta)|}{\sqrt{n}}\geq {\gamma_{k}+2^{k+1}}} \right)$$
 \begin{equation}
\label{genestpe0}
- \mathbb{P}\left(\max_{n_{k-1}<n\leq {n_{k}}}|B_{k}(n,\theta)|\geq 2^{k}\right)-\left(\frac{1}{2}\right)^{k+2}
\end{equation}
for all $\theta\in [0,2\pi]$.
\end{lemma}

This completes the set of pieces needed to construct the example stated by Theorem \ref{nonqueconthe}.

\medskip

{\bf Proof of Theorem \ref{nonqueconthe}:} 
Following \cite{volwoononque}, assume that $||x_{0}||_{_{\mathbb{P},2}}=1$ and let $(n_{j})_{j\geq 0}$, $(a_{j})_{j\geq 0}$, and $(\gamma_{j})_{j\geq 0}$  be defined inductively as follows:  $n_{0}=1$, $\gamma_{0}=0$, $a_{0}=0$, $a_{1}=\frac{1}{2}$, and given $n_{0},\cdots, n_{k-1}$, $a_{0}, \dots, a_{n_{k-1}}$ and $\gamma_{0},\dots,\gamma_{k-1}$ , let $\gamma_{k}$ be such that
$$\mathbb{P}\left(\left|\sum_{j=1}^{k-1}a_{n_{j}}e^{in_{j}\theta}\zeta_{-n_{j}}(\theta)\right|> \gamma_{k}\right)\leq \left(\frac{1}{2}\right)^{k+2}, $$
 
(see(\ref{firestakn})) and let $n_{k}>n_{k-1}$ be such that

\begin{equation}
\label{almthe}
\mathbb{P}\left(\max_{n_{k-1}< n\leq n_{k}}\frac{|(\zeta_{-n_{k-1}}-\zeta_{-n_{k-1}+n})(\theta)|}{\sqrt{n}}\geq \frac{\gamma_{k}+2^{k+1}}{a_{n_{k-1}}}\right)\geq 1-\left(\frac{1}{2}\right)^{k+1}
\end{equation}

for all $\theta\in[0,2\pi]$. The choice of $n_{k}$ is possible according to Lemma \ref{lemxij} ($|(\zeta_{-n_{k-1}}-\zeta_{-n_{k-1}+n})(\theta)|$ and $|\zeta_{-n}(\theta)|$ have the same distribution). Then define $a_{n_{k}}=\frac{1}{2^{k}\sqrt{n_{k-1}}}$ and $a_{j}=0$ for $n_{k-1}<j<n_{k}-1$.

The sequences $(a_{j})_{j\geq 0}$ and $(\gamma_{k})_{k}$, thus defined, satisfy the hypotheses of Lemma \ref{lemest} and therefore, by the estimates (\ref{genestpe0}) and (\ref{almthe}), 
$$\mathbb{P}\left(\max_{n_{k-1}< n \leq n_{k}}\frac{|E_{0}S_{n}(\theta)|}{\sqrt{n}}\geq 2^{k}\right)\geq 
1-\left(\frac{1}{2}\right)^{k+2}-\mathbb{P}\left(\max_{n_{k-1}<n\leq {n_{k}}}|B_{k}(n,\theta)|\geq 2^{k}\right)$$
for all $\theta\in [0,2\pi]$.

We claim that, under the present conditions,
\begin{equation}
\label{estbk}
\mathbb{P}\left(\max_{n_{k-1}<n\leq {n_{k}}}|B_{k}(n,\theta)|\geq 2^{k}\right)\leq \left(\frac{1}{2}\right)^{k+2}
\end{equation}
for $k\geq 3$.

Fix $k\geq 3$ and note that, for fixed $\theta$, $(|\zeta_{-n}(\theta)|)_{n\in \mathbb{N}}$ is an $L^2_{\mathbb{P}}$ submartingale (with respect to $(\mathcal{G}_{n})_{n\in\mathbb{Z}}$, where $\mathcal{G}_{k}=\sigma((x_{-j})_{j\leq k})$) and therefore, by Doob's maximal inequality (\ref{doo}):

$$E\left(\max_{k\leq n}\left|\zeta_{-k}(\theta)\right|\right)\leq ||\max_{k\leq n}\left|\zeta_{-k}(\theta)\right|||_{_{\mathbb{P},2}}\leq 2\,||\zeta_{-n}(\theta)||_{_2}= 2\,\sqrt{n+1}.$$

This gives

$$E\left(\max_{n_{k-1}<n\leq n_{k}}\left|B_{k}(n,\theta)\right|\right)\leq \sum_{j=k+1}^\infty a_{n_{j}}E\left(\max_{k\leq n_{k}-n_{k-1}}\left|\zeta_{-k}(\theta)\right|\right)\leq$$

$$\sum_{j=k+1}^\infty \frac{1}{2^{j-1}}\sqrt{\frac{n_{k}-n_{k-1}+1}{n_{j-1}}}\leq \frac{1}{2^{k-1}}, $$
  and therefore, by Markov's inequality (\cite{bilpromea}, p.276, (21.12))
  
  $$\mathbb{P}\left(\max_{n_{k-1}<n\leq n_{k}}\left|B_{k}(n,\theta)\right|\geq 2^{k}\right)\leq \frac{1}{2^{^{2k-1}}}\leq \left(\frac{1}{2}\right)^{k+2}$$
  
 as claimed.
 
 To finish the proof we observe that a combination of (\ref{genestpe0}), (\ref{almthe}) and (\ref{estbk}) gives, under the present choices of $(a_{k})_{k}$ and $(n_{k})_{k}$, that
 
 $$\mathbb{P}\left(\max_{n_{k-1}< n \leq n_{k}}\frac{|E_{0}S_{n}(\theta)|}{\sqrt{n}}< 2^{k}\right)\leq \left(\frac{1}{2}\right)^{k+1}$$
 so that, by the first Borel-Cantelli Lemma (\cite{bilpromea}, Theorem 4.3)
  $$\max_{n_{k-1}< n \leq n_{k}}\frac{|E_{0}S_{n}(\theta)|}{\sqrt{n}}\geq 2^{k} \mbox{   \it except for finitely many $k$'s,}$$
$\mathbb{P}-$a.s. This clearly implies that $\limsup_{n}{|E_{0}S_{n}(\theta)|}/{\sqrt{n}}=\infty$ {$\mathbb{P}-$a.s.}  \qed

\section{Proof of Theorems  \ref{queconznundratdec} and \ref{queconznmaxwoo}}
\label{prothenonrancen}

In this section we address the proofs of theorems \ref{queconznundratdec} and  \ref{queconznmaxwoo}. As the reader may expect at this point, these are just consequences of the fact that the hypotheses in these theorems are sufficient to verify the validity of item {\it 2.} in Corollary \ref{sufnecconqueconzn} (page \pageref{sufnecconqueconzn}). 

It is important to point out that, for proving Theorem \ref{queconznmaxwoo}, we will use  again the extensions to the product space described in Definition \ref{extprospa} and (\ref{exttprospathe}). Together with the proofs of Theorem \ref{ergthedisfoutra}, Theorem \ref{ergtheadafil} and Lemma \ref{genmarapplem} (``{\it Step 2.}''), this will serve as a further illustration of how this method of {\it lifting to the product space} allows us to translate estimates on (non-rotated) partial sums to corresponding results for Discrete Fourier Transforms.

\subsection{Proof of Theorem \ref{queconznundratdec}}

We will use the criterion given in Corollary \ref{sufnecconqueconzn}: we will prove that there exists $I'\subset [0,2\pi)$ with $\lambda(I')=1$ such that $(1-e^{i\theta})E_{0}S_{n}(\theta)/\sqrt{n}\to_{n} 0$, $\mathbb{P}-$a.s for every $\theta\in I'$. The conclusion follows by taking $J=I\cap I'$, where $I$ is the set guaranteed by Theorem \ref{quecltfoutra} and using item {\it 2.} in Corollary \ref{sufnecconqueconzn} (note that, for $\theta\in I$, $e^{i\theta}\neq 1$).

To do so we use an argument similar in spirit to the one leading to the proof of Lemma \ref{leminvprihancon} (the decompositions are way simpler here). Thus note that, if $n\in \mathbb{N}^{*}$ and $\theta\in [0,2\pi)$ are given
$$(1-e^{i\theta}){E}_{0}\frac{S_{n}(\theta)}{\sqrt{n}}=\frac{E_{0}
(S_{n}(\theta))-e^{i\theta}E_{0}(S_{n}(\theta))}{\sqrt{n}}=$$
\begin{equation}
\label{tertozer}
\frac{1}{\sqrt{n}}e^{i\theta}X_{0}-e^{in\theta}\frac{1}{\sqrt{n}}E_{0}X_{n-1}+\frac{1}{\sqrt{n}}\sum_{k=1}^{n-1}E_{0}[X_{k}-X_{k-1}]e^{ik\theta}.
\end{equation}

We will analyze each term in the last sum separately: the first term in the above expression,
$e^{it}X_{0}/\sqrt{n}$,  is trivially
convergent to zero for every $\omega\in \Omega$. 

Now, the conditional Jensen's inequality gives that:
\begin{equation}
\label{jenineapp}
|e^{in\theta}\frac{1}{\sqrt{n}}E_{0}X_{n-1}|^{2}\leq\frac{1}
{n}E_{0}|X_{n-1}|^{2}
\end{equation}
$\mathbb{P}-$a.s., and if we write
$$\frac{1}{n}E_{0}|X_{n-1}|^{2}=\frac{1}{n}\sum_{j=0}^{n-1}E
_{0}T^{j}|X_{0}|^{2}-\frac{1}{n}\sum_{j=0}^{n-2}E_{0}T^{j}|X_{0}|^{2},$$
we see that $E_{0}|X_{n-1}|^{2}/n\to_{n} 0$, $\mathbb{P}-$a.s. by Theorem \ref{ergtheadafil}, and therefore $e^{in\theta}E_{0}X_{n-1}/\sqrt{n}\to 0$, $\mathbb{P}-$a.s. by (\ref{jenineapp}).

To prove the convergence of the third term note that, since we are under the assumption (\ref{conratdecdif}), an argument similar to the one leading to the proof of the $\mathbb{P}-$a.s convergence of  (\ref{ranfunfixfre})  for $\lambda-$a.e $\theta$ (page \pageref{ranfunfixfre}) implies that there exists $I'\subset [0,2\pi)$ with $\lambda(I')=1$ such that for every $\theta\in I'$
$$\sum_{k\in \mathbb{N}^{*}}\frac{E_{0}(X_{k}-X_{k-1})}{k^{1/2}}e^{ik\theta}$$
converges $\mathbb{P}-$a.s, and the  Kronecker lemma (\cite{dur}, Theorem 2.5.5) implies that, for $\theta\in I'$
$$\frac{1}{\sqrt{n}}\sum_{k=1}^{n-1}E_{0}(X_{k}
-X_{k-1})e^{ik\theta}\rightarrow_{n} 0$$

$\mathbb{P}-$a.s. The conclusion follows from these arguments, as explained at the beginning, via (\ref{tertozer}).
\qed

\subsection{Proof of Theorem \ref{queconznmaxwoo}} 
The proof of theorem \ref{queconznmaxwoo} is, as announced at the beginning of this section, an application of the results in \cite{cunmer}. The proof that we will present here depends on the following lemma:

\medskip

\begin{lemma}
\label{cunmerlem}
In the context on page \pageref{gensetpro}. If 
\begin{equation}
\label{maxwoocon}
\sum_{k\in\mathbb{N}^{*}}\frac{||E_{0}S_{k}||_{{\mathbb{P},2}}}{k^{3/2}}<\infty,
\end{equation}
then $E_{0}S_{k}/\sqrt{n}\to_{n} 0$, $\mathbb{P}-$a.s.
\end{lemma}

{\bf Proof:} This is a direct consequence of Theorem 4.7 in \cite{cunmer}: with the notation of that paper, take $X_{0}$ in place of $f$, $\psi=1$ (the constant function), $E_{0}T$ in place of $T$, and use (\ref{eque0tk}).\qed 

{\bf Proof of Theorem \ref{queconznmaxwoo}:} Given $\theta\in I$, let $\tilde{T}_{\theta}$ and $\tilde{X_{0}}$ be the extensions of $X_{0}$ and $T$ recalled along the ``{\it Step 2.}'' in the proof of Lemma \ref{genmarapplem} (page \pageref{genmarapplem}). Keeping the notation introduced there note that, by  an argument analogous to the one leading to the chain of equalities (\ref{frozertothe}),
\begin{equation}
\sum_{k\in\mathbb{N}^{*}}\frac{||\tilde{E}_{0}\tilde{S}_{k}(\theta)||_{_{\lambda\times\mathbb{P},2}}}{k^{3/2}}=\sum_{k\in\mathbb{N}^{*}}\frac{||{E}_{0}S_{k}(\theta)||_{_{\mathbb{P},2}}}{k^{3/2}},
\end{equation}
so that, if the last series is convergent, Lemma \ref{cunmerlem} implies that there exists $\tilde{\Omega}_{0}\subset [0,2\pi)\times\Omega$ with $(\lambda\times\mathbb{P})\tilde{\Omega}_{0}=1$ such that for every $(u,\omega)\in \tilde{\Omega}_{0}$
$$0=\lim_{n}\frac{\tilde{E}_{0}\tilde{S}_{k}(\theta)}{\sqrt{n}}(u,\omega)= \lim_{n\to\infty}e^{iu}\frac{{E}_{0}S_{k}(\theta)}{\sqrt{n}}(\omega),$$
and it follows  that $E_{0}S_{k}(\theta)/\sqrt{n}\to_{n} 0$, $\mathbb{P}-$a.s. The conclusion follows again via Corollary \ref{sufnecconqueconzn}.\qed

\medskip

\end{document}